\documentclass[12pt,twoside]{amsart}
\usepackage{amssymb}
\usepackage{amscd}
\usepackage[abbrev,alphabetic]{amsrefs}
\usepackage{hyperref}
\usepackage{comment}
\usepackage{array,multirow,tabularx,longtable}
\usepackage{tikz}
\usetikzlibrary{cd}
\usepackage{here}
\usepackage{multirow}
\usepackage{ascmac}
\usepackage[margin=1.25in]{geometry}

\usepackage{fancybox}
\usepackage{ascmac}

\title[Fano threefolds in positive characteristic II]
{Fano threefolds in positive characteristic II} 
\author{Hiromu Tanaka} 
\address{Graduate School of Mathematical Sciences, 
The University of Tokyo, 
3-8-1 Komaba, Meguro-ku, Tokyo 153-8914, JAPAN} 
\email{tanaka@ms.u-tokyo.ac.jp}

\newcommand{\Bl}[0]{{\operatorname{Bl}}}

\newcommand{\Br}[0]{{\operatorname{Br}}}
\newcommand{\et}[0]{{\operatorname{\acute{e}t}}}

\newcommand{\codim}[0]{{\operatorname{codim}}}

\newcommand{\Hilb}[0]{{\operatorname{Hilb}}}
\newcommand{\Univ}[0]{{\operatorname{Univ}}}
\newcommand{\NE}[0]{{\operatorname{NE}}}
\newcommand{\red}[0]{{\operatorname{red}}}

\newcommand{\Coker}[0]{{\operatorname{Coker}}}
\newcommand{\Ker}[0]{{\operatorname{Ker}}}

\renewcommand{\Im}[0]{{\operatorname{Im}}}

\newcommand{\Proj}[0]{{\operatorname{Proj}}}
\newcommand{\Spec}[0]{{\operatorname{Spec}}}

\newcommand{\Bs}[0]{{\operatorname{Bs}}}
\newcommand{\Supp}[0]{{\operatorname{Supp}}}
\newcommand{\Pic}[0]{{\operatorname{Pic}}}
\newcommand{\mult}[0]{{\operatorname{mult}}}
\newcommand{\Ex}[0]{{\operatorname{Ex}}}

\renewcommand{\line}[0]{{\operatorname{line}}}
\newcommand{\conic}[0]{{\operatorname{conic}}}

\renewcommand{\min}[0]{{\operatorname{min}}}
\newtheorem{thm}{Theorem}[section]
\newtheorem{lem}[thm]{Lemma}
\newtheorem{cor}[thm]{Corollary}
\newtheorem{prop}[thm]{Proposition}
\newtheorem{claim}[thm]{Claim}    
\newtheorem{clm}[thm]{Claim}    
\newtheorem*{claim*}{Claim}         
\newtheorem{step}{Step}

\theoremstyle{definition}

\newtheorem{dfn}[thm]{Definition}

\newtheorem{rem}[thm]{Remark}
      
\newtheorem{nota}[thm]{Notation}         
\newtheorem{nasi}[thm]{}

\makeatletter
  
  \@addtoreset{equation}{thm}
  \makeatother

\newcommand{\cred}{\color{black}}

\newcommand{\MO}{\mathcal{O}}

\newcommand{\R}{\mathbb{R}}
\newcommand{\Q}{\mathbb{Q}}
\newcommand{\Z}{\mathbb{Z}}
\newcommand{\F}{\mathbb{F}}
\newcommand{\A}{\mathbb{A}}
\newcommand{\G}{\mathbb{G}}
\renewcommand{\P}{\mathbb{P}}

\newcommand{\m}{\mathfrak{m}}

\newcommand{\wt}{\widetilde}

\begin{document}

\maketitle

\begin{abstract}
Let $X$ be a smooth Fano threefold  over an algebraically closed field of positive characteristic. 
Assume that $|-K_X|$ is very ample and each of  the index and the Picard number is equal to one. We prove that $3 \leq g \leq 12$ 
and $g \neq 11$ for the genus $g$ of $X$. Moreover, we show that there exists no smooth curve on $X$ along which the blowup is Fano. 
\end{abstract}

\tableofcontents

\section{Introduction}

This article is the second part of our series of papers. 
The purpose of the first paper \cite{TanI} and this paper is to classify Fano threefolds in positive characteristic with $\rho(X)=1$. 
For the index $r_X$ of $X$, 
the classification for the case $r_X \geq 2$ has been carried out by \cite{Meg98} (cf. \cite[Theorem 2.18, Theorem 2.23]{TanI}). 
Hence it is enough to classify the case when $\rho(X)=r_X=1$, i.e., $\Pic\,X$ is generated by $\MO_X(K_X)$. 
In the first part \cite{TanI}, we classified Fano threefolds $X$ such that
$\rho(X)=r_X=1$ and $|-K_X|$ is not very ample. 
In this article, we treat the case when $\rho(X)=r_X=1$ and $|-K_X|$ is very ample. 
More precisely, the main theorem is as follows.

\begin{thm}[Corollary \ref{c-pt-not-ample}, Theorem \ref{t-g-bound}, Theorem \ref{t-blowup-main}, Theorem \ref{t g neq 11}]\label{intro-main}
Let $k$ be an algebraically closed field of characteristic $p>0$. 
Let $X$ be a three-dimensional smooth projective variety over $k$ 
such that $|-K_X|$ is very ample. 
Assume that 
$\Pic\,X$ is generated by $\MO_X(-K_X)$ as an abelian group. 
We define $g$ by $(-K_X)^3 = 2g-2$. 
Then the following hold. 
\begin{enumerate}
    \item $3 \leq g \leq 12$. 
    \item $g \neq 11$. 
    \item Let $\Gamma$ be a point  or a smooth curve on $X$. 
    Then $-K_Y$ is not ample for  the blowup $Y \to X$ along $\Gamma$. 
\end{enumerate}
\end{thm}

\noindent 
Under the same notation as in Theorem \ref{intro-main}, 
the author does not know whether there exists a line on $X$. 


\subsection{Description of the proof}\label{ss-overview}

In what follows, we shall explain some of ideas of the proof of  Theorem \ref{intro-main}. 

\medskip

\textbf{Theorem \ref{intro-main}(1)}: 
We now overview the proof of Theorem \ref{intro-main}(1), 
although it is very similar to the one in characteristic zero. 
Let $X \subset \P^{g+1}$ be an anti-canonically embedded Fano threefold. 
Assume that there exists a conic $\Gamma$ on $X$, i.e., 
$\Gamma$ is a curve on $X$ satisfying ${\cred -K_X \cdot \Gamma=2}$, 
and hence $\Gamma \simeq \P^1$. 
Let $\sigma : Y \to X$ be the blowup along $\Gamma$. 
We focus on a typical case: 
\begin{itemize}
\item $|-K_Y|$ is base point free, $-K_Y$ is big, 
and 
\item for the morphism $\psi : Y \to Z$ obtained by the Stein factorisation of $\varphi_{|-K_Y|}$, 
$\psi$ is a birational morphism to a projective normal threefold $Z$ such that $\dim \Ex(\psi) = 1$. 
\end{itemize}
In particular, $\psi : Y \to Z$ is a flopping contraction. 
Let $\psi^+ : Y^+ \to Z$ be its flop. 
It is easy to see that $\rho(Y^+)=2$ and $-K_{Y^+}$ is nef and big. 
Since  $K_{Y^+}$ is $\psi^+$-numerically trivial, 
there exists a unique $K_{Y^+}$-negative extremal ray $R$ of $\overline{\NE}(Y^+)$. 
Let $\tau : Y^+ \to W$ be the  contraction of $R$. 
To summarise, we get the following diagram: 
\begin{equation}\label{intro-e1-2ray}
\begin{tikzcd}
Y \arrow[d, "\sigma"'] \arrow[rd, "\psi"]& & Y^+ \arrow[ld, "\psi^+"'] \arrow[d, "\tau"]\\
X & Z & W.
\end{tikzcd}
\end{equation}
The classification of 
extremal rays \cite{Kol91} 
enables us to apply a case study depending on the type of $R$. 
We then get the upper bound $g \leq 12$ by computing 
some intersection numbers on $Y^+$. 
In order to assure the existence of a conic $\Gamma$ used in the above argument, 
we shall apply the same argument as above by starting with a blowup at a general point $P \in X$ instead of a blowup along a conic $\Gamma$. 
For a more detailed exposition, we refer to Subsection \ref{ss-2ray-intro}. 

\medskip

\textbf{Theorem \ref{intro-main}(2)}: 
The idea of the proof of Theorem \ref{intro-main}(2) (i.e., $g \neq 11$) 
is to construct a $W(k)$-lift of the diagram (\ref{intro-e1-2ray}) 
in reverse order. 
Assume $g=11$. 
Then we shall see that $W = \P^3_k$ and $\tau : Y^+ \to W = \P^3_k$ is a blowup along a smooth rational curve $B$ (Theorem \ref{t-conic-flop}). 
By using the Euler sequence for $\P^3_k$, 
we can check $H^1(B, N_{B/W})=0$, 
which implies  that the closed embedding $B \subset W = \P^3_k$ lifts to $W(k)$.  
Then $Y^+$ has a $W(k)$-lift $\wt{Y}^+$. 
It is easy to find a $W(k)$-lift 
$\wt{\psi}^+ : \wt{Y}^+ \to \wt{Z}^+$
of  $\psi^+ : Y^+ \to Z$. 
The next step is to construct a $W(k)$-lift 
$\wt{\psi} : \wt{Y} \to \wt{Z}$ of $\psi: Y \to Z$. 
To this end, we shall prove the existence of the flop 
of the flopping contraction $\wt{\psi}^+ : \wt{Y}^+ \to \wt{Z}^+$ 
(Proposition \ref{p flop lift}). 
This flop $\wt{\psi} : \wt{Y} \to \wt{Z}$ will be shown to be 
a $W(k)$-lift of $\psi: Y \to Z$ by using the uniqueness of flops. 
Finally, it is easy to get a $W(k)$-lift $\wt{\sigma} : \wt{Y} \to \wt{X}$ of $\sigma : Y \to X$. 
This leads to a contradiction, 
because the characteristic-zero fibre $X_0$ of $\wt{X}$ is a Fano threefold 
with $\rho(X_0) = r_{X_0}=1$ and $(-K_{X_0})^3 = 20$. 
For more details, we refer to Section \ref{s-g11}.









\medskip 

\textbf{Theorem \ref{intro-main}(3)}: 
Let us overview the proof of Theorem \ref{intro-main}(3) for the case when $\Gamma$ is a smooth curve. 
Let $X \subset \P^{g+1}$ be an anti-canonically embedded Fano threefold. 
Let $\Gamma$ be a smooth curve on $X$ and take the blowup $\sigma : Y \to X$ along $\Gamma$. 
Suppose that $-K_Y$ is ample. 
It suffices to derive a contradiction. 
Since $-K_Y$ is ample, 
we have the extremal ray $R$ of $\NE(Y)$ not corresponding to $\sigma: Y \to X$. 
Let $\tau: Y \to W$ be the contraction of $R$. 
\begin{itemize}
    \item 
    If $\tau$ is of type C or D, then we define $D$ as the pullback of an ample generator of $\Pic\,W (\simeq \Z)$. 
    \item  
    If $\tau$ is of type E, then we set $D :=\Ex(\tau)$, 
    which is a prime divisor on $Y$. 
\end{itemize}
For $E := \Ex(\sigma)$, it is easy to see that we can write  
\[
D \sim -\alpha K_Y - \beta E
\]
for some  positive integers $\alpha, \beta$ (Lemma \ref{l-beta=r}). 

We now consider the case when $\tau$ is of type D, i.e., $W \simeq \P^1$. 
In this case, we have $D^3 = D^2 \cdot (-K_Y)=0$. 
For $h := h^1(\Gamma, \MO_{\Gamma})$ and $d := -K_X \cdot \Gamma$, 
we can check that the equalities $D^3 = D^2 \cdot (-K_Y)=0$ induce the following Diophantine equations 
(Lemma \ref{l-678-easy-bound}): 
\begin{equation}\label{e1-main}
(g-2-d+h) \alpha^2 -(d-2h+2) \alpha\beta +(h-1)\beta^2=0
\end{equation}
\begin{equation}\label{e2-main}
d (-\alpha+\beta) + (2h-2) (\alpha+\beta) =0. 
\end{equation}
Although we have five parameters $g, d, h, \alpha, \beta$, 
all of these parameters are bounded, e.g., $3 \leq g \leq 12, 0 \leq h \leq 10, 1 \leq \beta \leq 3$. 
A key part is to restrict the possibilities for $\alpha$. 
We only overview how to prove $\alpha \leq 3$. 
Note that (\ref{e1-main}) implies $(h-1)\beta^2 \in \alpha \Z$. 
Since $\alpha$ and $\beta$ are  coprime for the case of  type $D$, 
we get $h-1 \in \alpha \Z$. 
Then, by (\ref{e2-main}), $d(\alpha -\beta) =(2h-2)(\alpha+\beta) \in \alpha \Z$. 
Since $\alpha$ and $\beta$ are coprime, so are $\alpha$ and $\alpha -\beta$. 
Hence $d \in \alpha \Z$. 
By $h-1 \in \alpha \Z$ and $d \in \alpha \Z$, 
(\ref{e1-main}) implies $(h-1)\beta^2 \in \alpha^2 \Z$, and hence $h-1 \in \alpha^2 \Z$. 
Finally, by using $0 \leq h \leq 10$, we get $h=1$ or $\alpha \in \{ 2, 3\}$. 
If $h=1$, then (\ref{e2-main}) implies $\alpha =\beta \in \{1, 2, 3\}$. 
Similarly, we can settle the case when $\tau: Y \to W$ is not of type $E_1$ 
almost by using elementary number-theoretic methods.

The most complicated part is the case when $\tau : Y \to W$ is of type $E_1$. 
In this case, we have $\beta =r_W$ for the index $r_W$ of the Fano threefold $W$ and 
we can write $\alpha +1 =\beta \gamma$ for some $\gamma \in \Z_{>0}$. 
In particular, we have $-K_W \sim \beta H_W$ for some Cartier divisor $H_W$ on $W$. 
In what follows, we assume  $h \geq 2$ and $\gamma \geq 2$ for simplicity. 
We have the following Diophantine equations (Lemma \ref{l-678-E1-ABC}). 
\begin{enumerate}
    \item[(A)] $\frac{1}{2} \beta H_W^3 = (g-d +h -2)
\gamma^2 - (d -2h+2)\gamma   + (h-1)$. 
\item[(B)] $H_W^3 = (2g-2d +2h -4)\gamma^3 -3(d -2h+2)\gamma^2+ 3(2h-2)\gamma +(d+2h-2)$. 
\item[(C)] $\alpha H_W^3 = (\beta \gamma -1)H_W^3 = (d -2h+2)\gamma^2 -2(2h-2)\gamma  -(d+2h-2)$. 
\end{enumerate}
We shall treat the following two cases separately. 
\begin{enumerate}
\item[(I)]  $h \leq 4$. 
\item[(II)] $h\geq 5$. 
\end{enumerate}

(I) Assume $h \leq 4$. Note that there are only $17$ possibilities for $(\beta, H_W^3)$ 
(Remark \ref{r-E1-deg-genus}, Remark \ref{r-w-bound}). 
By $0 \leq h \leq 4$, we have at most $85$ possibilities for $(h, \beta, H_W^3)$. 
Actually, we shall apply case study by fixing this triple in many cases. 
Here we overview the proof only for the case when 
$h=2, \beta=2,$ and $H_W^3 =3$. 
By (C), we get 
\[
(2\gamma -1) \cdot 3 = (d-2) \gamma^2 -4\gamma -(d+2). 
\]
This, together with  $\gamma \geq 2$, implies 
\[
d = \frac{2\gamma^2+10\gamma -1}{\gamma^2-1}= 2 + \frac{10\gamma +1}{\gamma^2-1}. 
\]
We then obtain $10\gamma +1 \in (\gamma^2-1)\Z \subset (\gamma -1)\Z$. 
Taking modulo $\gamma -1$, we get $11 \equiv 10\gamma +1 \equiv 0 \mod \gamma-1$. 
Hence $\gamma =2$ or $\gamma =12$. 
If $\gamma =12$, then we would get a contradiction:  
$\Z \ni \frac{10\gamma +1}{\gamma^2-1} = \frac{121}{143} \not\in \Z$. 
Thus $\gamma =2$. 
It holds that  $d = 2 + \frac{10\gamma +1}{\gamma^2-1} = 9$. 
By using (A), we get 
\begin{eqnarray*}
3 &=& \frac{1}{2} \beta H_W^3\\
&=& (g-d +h -2)
\gamma^2 - (d -2h+2)\gamma   + (h-1)\\
&=& (g-9 +2 -2) \cdot 4 - (9 -4+2) \cdot 2   + (2-1)\\
&=& 4(g-9) -13. 
\end{eqnarray*}
Then it holds that  $g =13$, which contradicts $ 3 \leq g \leq 12$ (Theorem \ref{intro-main}(1)). 

\medskip

(II) Assume $h \geq 5$. 
As in (I), we may fix $(\beta, H_W^3)$. 
Then the essential part is to restrict the possibilities for 
$(\gamma, h)$. We only overview how to show $\gamma \leq 5$ and $h \leq 6$. 
A key observation is to consider the following quadratic function: 
\[
f(x) :=  (g-d +h -2)x^2 - (d -2h+2)x   + (h-1). 
\]
Note that (A) implies $f(\gamma) = \frac{1}{2} \beta H_W^3$. 
By $h-1 \geq 1$, $g-d +h -2 =\frac{1}{2}(-K_Y)^3 >0$,  and 
$d -2h+2 =(-K_Y)^2 \cdot E >0$, 
we have $f(-3) \geq 9 + 3+ 1 =13 >11$. 
For the symmetry axis $x=s$ of the quadratic function $f(x)$ 
(explicitly, $s:= \frac{1}{2} \cdot \frac{d-2h+2}{g-d +h -2} \in \R_{>0}
$), 
we obtain  
\[
f(2s+3) = f(s +(s+3)) = f(s-(s+3)) = f(-3) >11. 
\]
This, together with $1 \leq f(\gamma) = \frac{1}{2} \beta H_W^3 \leq 11$, 
implies $\gamma < 2s +3$. 
It is easy to check that  $d \leq h+g-3$ (Lemma \ref{l-678-easy-bound}). 
Moreover, it is not so hard to treat the case when $d = h+g -3$ 
(Lemma \ref{l-d=h+g-3}).  
In what follows, we consider the case when $d \leq  h+g-4$. 
By $g \leq 12$ and $h \geq 5$, we get 
\[
s =\frac{1}{2} \cdot \frac{d-2h+2}{g-d +h -2} \leq \frac{1}{2} \cdot \frac{(h+g-4)-2h+2}{g-(h+g-4) +h -2} 
\]
\[
= \frac{1}{4}(-h+g-2) 
\leq \frac{1}{4}(-5+12-1) = \frac{3}{2}. 
\]
Therefore, $\gamma < 2s +3 \leq 2 \cdot \frac{3}{2} +3 =6$. Thus $\gamma \leq 5$. 
By (C), we obtain 
\begin{eqnarray*}
0&<& \alpha H_W^3  \\
&=& (d -2h+2)\gamma^2 -2(2h-2)\gamma  -(d+2h-2)\\
&<& \gamma ( (d -2h+2)\gamma -2(2h-2)). 
\end{eqnarray*}
This, together with $d \leq h+g-3$, implies 
\[
2(2h-2) <  (d -2h+2)\gamma \leq (h+g-3 -2h+2)\gamma = (g-1-h)\gamma.
\]
Thus $4h-3 \leq (g-1-h) \gamma$, and hence 
\begin{equation}\label{e1-678-h-gamma}
h \leq \frac{(g-1)\gamma +3}{\gamma+4} \leq 
\frac{11\gamma +3}{\gamma +4} 
\overset{(*)}{\leq}  
\frac{11 \cdot 5 +3}{5 +4} =\frac{58}{9} <7, 
\end{equation}
 where the inequality $(*)$ follows from $\gamma \leq 5$ and 
$\frac{11\gamma +3}{\gamma +4} = 11 - \frac{41}{\gamma+4}$. 
Thus $h \leq 6$.

\begin{rem}
It is easy to give a computer-assisted proof which is much shorter than that of this paper. 
Actually, if $\tau$ is of type $E_1$, then 
all the cases except for $(g, d, h, \beta, \gamma, H_W^3)=(8, 4, 1, 2, 2, 4)$ are ruled out by using (A) and (C), together with suitable bounds on the parameters, e.g., $6 \leq g \leq 12, 0 \leq h \leq g-2 \leq 10$ (Remark \ref{r computer}). 
Only for the case $(g, d, h, \beta, \gamma, H_W^3)=(8, 4, 1, 2, 2, 4)$, 
we need to use results on algebraic geometry (cf. the proof of Lemma \ref{l-678-E1-h=1}). 
In this paper, we shall give a human-readable proof for the sake of logical completeness. 
\end{rem}

\begin{rem}
The contents outside Section \ref{s-nonFano} 
and Section \ref{s-g11} 
are quite similar to the argument in  characteristic zero \cite{Tak89}, \cite{IP99}. 
However, we include whole the proofs, 
as there are several obstructions appearing only in positive characteristic. 
\end{rem}





\vspace{5mm}

\textbf{Acknowledgements:} 
The author would like to thank 
Tatsuro Kawakami, Yuri Prokhorov, Teppei Takamatsu, 
Burt Totaro, Jakub Witaszek, Fuetaro Yobuko
for useful comments, constructive suggestions, and answering questions. 
The author was funded by JSPS KAKENHI Grant numbers JP22H01112 and JP23K03028. 

\section{Preliminaries}\label{s-prelim}

\subsection{Notation}\label{ss-notation}

In this subsection, we summarise notation used in this paper. 

\begin{enumerate}
\item We will freely use the notation and terminology in \cite{Har77} and \cite{KM98}. 
In particular, $D_1 \sim D_2$ means linear equivalence of Weil divisors. 
\item 
Throughout this paper, 
we work over an algebraically closed field $k$ 
of characteristic $p>0$ unless otherwise specified. 
\item For an integral scheme $X$, 
we define the {\em function field} $K(X)$ of $X$ 
as the local ring $\MO_{X, \xi}$ at the generic point $\xi$ of $X$. 
For an integral domain $A$, $K(A)$ denotes the function field of $\Spec\,A$. 

\item 
For a scheme $X$, its {\em reduced structure} $X_{\red}$ 
is the reduced closed subscheme of $X$ such that the induced closed immersion 
$X_{\red} \to X$ is surjective. 
\item 
For a field $\kappa$, we say that $X$ is a {\em variety} over $\kappa$ if 
$X$ is an integral scheme which is separated and of finite type over $\kappa$. 
We say that $X$ is a {\em curve} (resp. a {\em surface}, resp. a {\em threefold}) 
if $X$ is a variety over $\kappa$ of dimension one (resp. {\em two}, resp. {\em three}). 
Unless otherwise specified, a variety means a variety over $k$. 
\item 
 We say that $X$ is a {\em Fano threefold}  
if $X$ is a three-dimensional smooth projective variety over $k$ such that $-K_X$ is ample. 
The {\em index} $r_X$ of a Fano threefold $X$ is defined as the largest positive integer $r$ such that 
there exists a Cartier divisor $H$ on $X$ satisfying $-K_X \sim r H$. 
\item 
Our notation will not distinguish between invertible sheaves and
Cartier divisors. For example, given a Fano threefold $X$, 
the equality $\Pic\,X = \Z K_X$ means that $\Pic\,X$ is isomorphic to $\Z$ which is generated by $\MO_X(K_X)$ (i.e., $\rho(X)=r_X=1$). 
\item 
Given a closed subscheme $\Gamma$ of $\P^N_k$, 
$\langle \Gamma \rangle$ denotes the smallest linear subvariety of $\P^N_k$ 
containing $\Gamma$, which is nothing but the intersection of all the hyperplanes containing $\Gamma$. 
By abuse of notation, $\langle \Gamma \rangle =\P^r$ means that 
$\langle \Gamma \rangle$ is of dim $r$. 
\item 
Given a closed subscheme $Z$ of $\P^N_k$, 
we set $\MO_Z(r) := \MO_{\P^N}(r)|_Z$ for every $r \in \Z$. 
\item 
Given a closed subscheme $X$ of $\P^N_k$, 
we say that $X \subset \P^N_k$ is {\em an intersection of quadrics} if 
the scheme-theoretic equality 
\[
X = \bigcap_{X \subset Q} Q
\]
holds, where $Q$ runs over all the quadric hypersurfaces of $\P^N_k$ containing $X$. Here the scheme-theoretic equality $X = \bigcap_{X \subset Q}Q$ means the equality $I_X = \sum_{X \subset Q} I_Q$ of the corresponding coherent ideals. 
\item For $n \in \Z_{\geq 0},$ we set $\F_n := \P_{\P^1}( \MO_{\P^1} \oplus \MO_{\P^1}(n))$. 
\item Given a projective normal variety $X$ and a Cartier divisor $D$ on $X$, 
 $\kappa(X, D)$ denotes the {\em Iitaka dimension}, which is characterised by the following properties: 
\begin{itemize}
\item If $H^0(X,\MO_X(mD))=0$ for every $m \in \Z_{>0}$, then 
$\kappa(X, D) := -\infty$. 
\item Assume that $H^0(X, \MO_X(m_1D)) \neq 0$ for some $m_1 \in \Z_{>0}$. 
Then there exist $a, b \in \Q_{>0}$ and $m_2 \in \Z_{>0}$ such that 
the following holds for every $m \in \Z_{>0}$: 
\[
a m^{\kappa(X, D)} \leq h^0(X, \MO_X(mm_1m_2D)) \leq b m^{\kappa(X, D)}. 
\]
\end{itemize}
\item $W(k)$ denotes the ring of Witt vectors of $k$. 
It is well known that $W(k)$ is a discrete valuation ring such that $pW(k)$ is the maximal ideal and $W(k)/pW(k) \simeq k$. 
\item 
{\cred 
Let $X$ be a closed subscheme in $\P^N_k$ and 
take a closed point $P$ of $X$. 
Throughout this paper, the {\em tangent space} $T_PX$  of $X \subset \P^N_k$ at $P$ 
means the embedded tangent space in $\P^N_k$. 
Recall that $T_PX$ satisfies the following properties (a)-(c). 
\begin{enumerate}
\item $T_PX$ is a linear subvariety of $\P^N_k$ passing through $P$. 
\item $T_P(X \cap Y) = T_PX \cap T_PY$ for  closed subschemes $X$ and $Y$ of $\P^N_k$ containing $P$. 
\item 
If $P = [1:a_1: \cdots :a_N] \in \P^N_k$ and 
$F(X_0, ..., X_N) \in k[X_0, ..., X_N]$ is 
a homogeneous polynomial $F(X_0, ..., X_N) \in k[X_0, ..., X_N]$ 
satisfying $F(1, a_1, ..., a_N)=0$ (i.e., $P\in V(F)$), 
then 
\[
T_P(V(F)) \cap \A^N_k = 
\left\{ \frac{\partial f}{\partial x_1}(a) (x_1-a_1) + \cdots + 
 \frac{\partial f}{\partial x_N}(a)(x_1-a_1) =0\right\},  
\]
where $x_i := X_i/X_0$, $f(x_1, ..., x_N) := F(1, x_1, ..., x_N)$, $a:=(a_1, ..., a_n)$, and $\mathbb A^N_k := D_+(X_0)$. 
\end{enumerate}
Note that 
(a)-(c) give an alternative definition of $T_PX$ 
for the case when $P \in D_+(X_0)$.}
\end{enumerate}

\subsection{Brauer groups}

\begin{dfn}[\cite{CTS21}*{Definition 3.2.1}]\label{d-Brauer}
For a scheme $X$, 
the {\em Brauer group} $\Br(X)$ of $X$ is defined by 
\[
\Br(X) := H^2_{\et}(X, \G_m),
\]
where $H^2_{\et}(-)$ denotes the \'etale cohomology. 
For a ring $R$,  we set 
\[
\Br(R) := \Br(\Spec\,R) =H^2_{\et}(\Spec\,R, \G_m). 
\]
\end{dfn}

\begin{prop}\label{p-Br-vanish}
If $C$ is a curve over $k$, then $\Br(K(C))=0$. 
\end{prop}

\begin{proof}
The assertion follows from \cite[Theorem 1.2.14]{CTS21}, since $k$ is an algebraically closed field. 
\end{proof}

\begin{prop}\label{p-Br-Pic}
Let $\kappa$ be a field and let $\kappa^s$ be the separable closure of $\kappa$ inside the algebraic closure $\overline{\kappa}$. 
Let $X$ be a proper variety over $\kappa$ such that $H^0(X, \MO_X) = \kappa$. 
Set $X^s := X \times_{\kappa} \kappa^s$ and let $\alpha : X^s \to X$ be the projection. 
Then 
we have the following exact sequence: 
\[
0 \to \Pic X \xrightarrow{\alpha^*} (\Pic X^s)^{{\rm Gal}(\kappa^s/\kappa)} \to \Br(\kappa). 
\]
\end{prop}

\begin{proof}
In order to apply \cite[Proposition 5.4.2]{CTS21}, 
it suffices to check that $H^0(X^s, \MO_{X^s})^{\times} = (\kappa^s)^{\times}$. 
Recall that the equality $H^0(X, \MO_X) = \kappa$ means that the induced homomorphism 
$\MO_{\Spec\,\kappa} \to \theta_*\MO_X$ is an isomorphism for the canonical morphism $\theta : X \to \Spec\,\kappa$. 
This property is stable under flat base changes, and hence 
$H^0(X^s, \MO_{X^s}) = \kappa^s$. 
Therefore, we get $H^0(X^s, \MO_{X^s})^{\times} = (\kappa^s)^{\times}$. 
\end{proof}

\subsection{Fano threefolds}

In this subsection, we summarise some results on Fano threefolds which will be used in this paper.  

\begin{dfn}\label{d anti embedded}
We say that $X \subset \P^N$ is an {\em anti-canonically embedded Fano threefold} 
if \begin{enumerate}
\item $X$ is a Fano threefold, 
\item $X$ is a closed subscheme of $\P^N$, 
\item $\MO_X(-K_X) \simeq \MO_{\P^N}(1)|_X$, and 
\item the induced $k$-linear map 
\[
H^0(\P^N, \MO_{\P^N}(1)) \to H^0(X, \MO_{\P^N}(1)|_X)
\]
 is an isomorphism. 
\end{enumerate}
\end{dfn}

\begin{rem}
Let $X$ be a Fano threefold such that $|-K_X|$ is very ample. 
For $N := h^0(X, -K_X) -1$, 
let $j : X \hookrightarrow \P^N$ be the induced closed immersion. 
Then $j(X) \subset \P^N$ is an anti-canonically embedded Fano threefold. 
\end{rem}

\begin{prop}\label{p-RR}
Let $X \subset \P^N$ be an anti-canonically embedded Fano threefold with $\rho(X)=1$. 
Then the following hold. 
\begin{enumerate}
\item $\Pic\,X \simeq \Z$. 
\item 
$\bigoplus_{m=0}^{\infty} H^0(X, \MO_X(-mK_X))$ 
is generated by $H^0(X, \MO_X(-K_X))$ as a $k$-algebra. 
\item 
For all $i>0$ and $m \in \Z_{\geq 0}$, 
it holds that 
$H^i(X, -mK_X) = 0$ and 
\[
h^0(X, -mK_X)=  \frac{1}{12} m(m+1)(2m+1)(-K_X)^3 + 2m + 1. 
\]
In particular, ${\cred h}^0(X, -K_X) = \frac{1}{2} (-K_X)^3 +3 = g+2$ and $N=g+1$
for $g := \frac{1}{2} (-K_X)^3 +1$. Furthermore, $g \geq 3$. 
\end{enumerate}
\end{prop}

\begin{proof}
See \cite[Theorem 2.4, Corollary 2.6, Corollary 4.5, Theorem 6.2]{TanI}. 
\end{proof}

\begin{thm}\label{t-inter-quads}
Let $X \subset \P^{g+1}$ be an anti-canonically embedded Fano threefold. 
If $\rho(X)=1$ and $g \geq 5$, then $X$ is an intersection of quadrics. 
\end{thm}

\begin{proof}
See \cite[Theorem 1.2]{TanI}. 
\end{proof}

\begin{prop}\label{p-genus-345}
Let $X \subset \P^{g+1}$ be an anti-canonically embedded  Fano threefold with $\rho(X)=1$. 
Then the following hold. 
\begin{enumerate}
\item If $g=3$, then $X \subset \P^4$ is a smooth hypersurface of degree $4$. 
\item If $g=4$, then $X \subset \P^5$ is a complete intersection of 
a quadric hypersurface and a cubic hypersurface. 
\item If $g=5$, then $X \subset \P^6$ is a complete intersection of three quadric hypersurfaces. 
\end{enumerate}
\end{prop}

\begin{proof}
The assertion (1) follows from $\deg X = (-K_X)^3 = 2g -2 = 4$. 

Let us show (2). 
For $m \in \Z_{\geq 0}$, it holds by  Proposition \ref{p-RR} that 
\[
h^0(X, -mK_X) 
=  \frac{1}{12} m(m+1)(2m+1)(-K_X)^3 + 2m + 1, 
\]
and hence the following hold: 
\begin{itemize}
    \item $h^0(X, -K_X)  =  \frac{1}{2}(-K_X)^3 + 3 = g+2 = 6$. 
    \item $h^0(X, -2K_X)  =  \frac{5}{2}(-K_X)^3 + 5 = 5g = 20$. 
    \item $h^0(X, -3K_X)  =  7(-K_X)^3 + 7 = 14g-7 = 49$. 
\end{itemize}
Since $\bigoplus_{m=0}^{\infty} H^0(X, \MO_X(-mK_X))$ 
is generated by $H^0(X, -K_X)$ (Proposition \ref{p-RR}), 
the restriction map 
\[
\rho_m : H^0(\P^5, \MO_{\P^5}(m)) \to 
H^0(X, -mK_X)
\]
is surjective for every integer $m \geq 0$ \cite[Proposition 1.7]{Fuj75}. 
By 
\[
h^0(\P^5, \MO_{\P^5}(1)) = 6, \qquad 
h^0(\P^5, \MO_{\P^5}(2)) = 21, \qquad 
h^0(\P^5, \MO_{\P^5}(3)) = 56, 
\]
there exists  a quadric hypersurface $Y \subset \P^5$ containing $X$. 
Here $Y$ is an integral scheme, as otherwise $X$ is contained in a hyperplane, which is a contradiction. 
For the element $\zeta \in H^0(\P^5, \MO_{\P^5}(2))$ corresponding to $Y$ 
and 
the natural composite map 
\[
\theta : (k \zeta) \otimes H^0(\P^5, \MO_{\P^5}(1)) 
\hookrightarrow 
H^0(\P^5, \MO_{\P^5}(2)) \otimes H^0(\P^5, \MO_{\P^5}(1)) 
\to H^0(\P^5, \MO_{\P^5}(3)), 
\]
we have that 
\[
\dim \Im\,\theta \leq \dim H^0(\P^5, \MO_{\P^5}(1))=6.
\]
By $56-49 = 7 >6$, there exists an element 
\[
\xi \in \Ker (\rho_3 : H^0(\P^5, \MO_{\P^5}(3)) \to H^0(X, -3K_X))
\]
such that $\xi \not\in \Im\,\theta$. 
Let $Z$ be the cubic hypersurface corresponding to $\xi$. 
Let us show that $Z$ is an integral scheme. 
Otherwise, we get a decomposition $Z = W + W'$ as Weil divisors on $\P^5$, 
where $W$ and $W'$ are of degree $1$ and $2$, respectively. 
Then $W'$ should contain $X$, which implies $W' = Y$. 
However, this 
is a contradiction. 
Hence $Z$ is integral. 

Let us show that the scheme-theoretic equality
\[
X = Y \cap Z
\]
holds. 
The inclusion 
$X \subset Y \cap Z$ is clear. 
For a hyperplane $H$ on $\P^5$, we obtain 
\[ 
X \cdot H^3 =6 = Y \cdot Z \cdot H^3. 
\]
In particular, $(Y \cap Z)_{\red} = X$ and $Y \cap Z$ is generically reduced. 
As $Y \cap Z$ is Cohen--Macaulay, 
$Y \cap Z$ satisfies $R_0$ and $S_1$, and  hence $Y \cap Z$ is an integral scheme. 
Then the inclusion 
$X \subset Y \cap Z$ must be an equlity: $X = Y \cap Z$. 
Thus (2) holds. 

Let us show (3). We have 
$h^0(X, -2K_X) = 25$ and $h^0(\P^6, \MO(2)) = 28$. 
For a $k$-linear basis $\zeta_1, \zeta_2, \zeta_3$ of $H^0(\P^6, \MO(2) \otimes I_X)$, let $Q_1, Q_2, Q_3$ be the corresponding  quadric hypersurfaces. 
By Theorem \ref{t-inter-quads}, 
we have a scheme-theoretic equality $X = Q_1 \cap Q_2 \cap Q_3$. 
By counting dimension, this intersection must be a complete intersection. 
Thus (3) holds. 
\qedhere

\end{proof}

\section{Extremal rays}\label{s-ext-ray}

In this section, we shall establish some foundational results on contractions $\tau : V \to W$ of extremal rays on a smooth projective threefold $V$. 
For example, we prove that 
the following sequence is exact for a curve $\ell$ on $V$ contracted by $\tau$ 
(Proposition \ref{p-cont-ex} in Subsection \ref{ss-Pic-gp}): 
\[ 
0 \longrightarrow \Pic(W) \xrightarrow{\tau^*} \Pic(V) \xrightarrow{(-\cdot \ell)} \mathbb{Z}.  
\]
In Subsection \ref{ss-type-ext-ray}, we shall recall the definitions of types and lengths of extremal rays, e.g., 
\begin{itemize}
\item $\tau$ is {\em of type} C if $\dim W=2$. 
\item $\tau$ is {\em of type} D if $\dim W =1$. 
\item $\tau$ is {\em of type} E if $\dim W =3$. 
\end{itemize}
For the case when $V$ is a weak Fano threefold with $\rho(V)=2$ and $H^1(V, \MO_V)=0$, 
we shall compute intersection numbers $D^3, D^2 \cdot (-K_V), D \cdot (-K_V)^2$ and 
lengths of extremal rays (Subsection \ref{ss-weak-Fano}), where $D$ is defined as follows: 
\begin{itemize}
\item $D := \tau^*L$ for a line $L$ on $W \simeq \P^2$ if $\tau$ is of type C. 
\item $D := \tau^*P$ for a point $P$ on $W \simeq \P^1$  if $\tau$ is of type D. 
\item $D:= \Ex(\tau)$, which is a prime divisor,  if $\tau$ is of type E. 
\end{itemize}
If $\tau$ is of type E, then this was established in \cite{Kol91}. 
If $\tau$ is of type C (resp. type D), 
then the proof is essentially contained in Subsection \ref{ss-conic-bdl} (resp. Subsection \ref{ss-dP-I} and Subsection  \ref{ss-dP-II}).


\subsection{Types of extremal rays}\label{ss-type-ext-ray}

\begin{dfn}\label{d-cont}
Let $V$ be a smooth projective threefold. 
Let $R$ be a $K_V$-negative extremal ray of $\overline{\NE}(V)$. 
By \cite[(1.1.1)]{Kol91}, there exists a unique morphism $\tau \colon  V \to W$, called the {\em contraction (morphism)} of $R$, 
to a projective normal variety $W$ such that (1) and (2) hold. 
  \begin{enumerate}
    \item $\tau_*\MO_V = \MO_W$. 
    \item For every curve $C$ on $V$, $[C]\in R$ if and only if $\tau(C)$ is a point.
  \end{enumerate}
\end{dfn}

\begin{dfn}\label{d-length}
Let $V$ be a smooth projective threefold. 
Let $R$ be a $K_V$-negative extremal ray of $\overline{\NE}(V)$ 
and let $\tau: V \to W$ be the contraction of $R$. 
We set
\[ 
\mu_R := \min \left\{ -K_V\cdot C \mid 
C \mathrm{~is~a~rational~curve~on~} V \mathrm{~with~} [C]\in R \right\}, 
\]
which is called the {\em length} of an extremal ray $R$. 
We say that $\ell$ is an {\em extremal rational curve} (of $R$) 
if $\ell$ is a rational curve on ${\cred V}$ such that 
$[\ell] \in R$ and $-K_V \cdot \ell =\mu_R$. 
\end{dfn}

\begin{dfn}\label{d-type-ext-ray}
Let $V$ be a smooth projective threefold. 
Let $R$ be a $K_V$-negative extremal ray of $\overline{\NE}(V)$ 
and let $\tau : V \to W$ be the contraction of $R$. 
\begin{enumerate}
\item $R$ and $\tau$ are called of {\em type $C$} if $\dim W = 2$. 
\begin{itemize}
\item $R$ and $\tau$ are called of {\em type $C_1$} if 
$\tau$ is  not smooth. 
\item $R$ and $\tau$ are called of {\em type $C_2$} if 
$\tau$ is smooth. 
\end{itemize}
\item $R$ and $\tau$ are called of {\em type $D$} if $\dim W = 1$. 
Set $K :=K(W)$ and $V_K$ denotes the generic fibre of $\tau: V \to W$. 
\begin{itemize}
\item $R$ and $\tau$ are called of {\em type $D_1$} if $1 \leq K_{V_K}^2 \leq 7$. 
\item $R$ and $\tau$ are called of {\em type $D_2$} if $K_{V_K}^2 =8$. 
\item $R$ and $\tau$ are called of {\em type $D_3$} if $K_{V_K}^2 =9$. 
\end{itemize}
\item $R$ and $\tau$ are called of {\em type $E$} if $\dim W = 3$. 
Set $D:=\Ex(\tau)$. 
\begin{itemize}
\item $R$ and $\tau$ are called of {\em type $E_1$} if 
$\tau(D)$ is a curve. 
\item  $R$ and $\tau$ are called of {\em type $E_2$} if 
$\tau(D)$ is a point, $D \simeq \P^2$, and $\MO_V(D)|_D \simeq \MO_D(-1)$. 
\item $R$ and $\tau$ are called of {\em type $E_3$} if 
$\tau(D)$ is a point, $D$ is isomorphic to a smooth quadric surface $Q$ on $\P^3$, and $\MO_V(D)|_D \simeq \MO_{\P^3}(-1)|_Q$. 
\item $R$ and $\tau$ are called of {\em type $E_4$} if 
$\tau(D)$ is a point, $D$ is isomorphic to a singular quadric surface $Q$ on $\P^3$, and $\MO_V(D)|_D \simeq \MO_{\P^3}(-1)|_Q$. 
\item  $R$ and $\tau$ are called of {\em type $E_5$} if 
$\tau(D)$ is a point, $D \simeq \P^2$, and $\MO_V(D)|_D \simeq \MO_D(-2)$. 
\end{itemize}
\item $R$ and $\tau$ are called of {\em type $F$} if $\dim W = 0$. 
\end{enumerate}
\end{dfn}

\subsection{Del Pezzo fibrations I}\label{ss-dP-I}


The purpose of this subsection is to study extremal rays of type $D$. 
For flexibility, we allow the case when the base curve $W$ is not necessarily projective. 

\begin{dfn}
We say that $\tau : V \to W$ is a {\em del Pezzo fibration} 
if \begin{enumerate}
    \item $V$ is a smooth threefold, 
    \item $W$ is a smooth curve, 
    \item $\tau_*\MO_V = \MO_W$, $\rho(V/W)=1$, and 
    \item $-K_V$ is $f$-ample.
\end{enumerate}
\end{dfn}

\begin{prop}\label{p-dP-fibres}
Let $V$ be a smooth threefold and 
let $\tau : V \to W$ be a del Pezzo fibration.  
Set $K :=K(W)$ and $\overline K := \overline{K(W)}$, where $\overline{K(W)}$ denotes the algebraic closure of $K(W)$.  
Let $V_K$ and $V_{\overline K}$ be the generic fibre and the geometric generic fibre of $\tau$, respectively. 
Then the following hold. 
\begin{enumerate}
\item For a closed point $w \in W$, the effective Cartier divisor $\tau^*w$ is a prime divisor. 
In particular, any scheme-theoretic fibre of $\tau$ is geometrically integral. 
\item $V_{\overline K}$ is a canonical del Pezzo surface, i.e., 
$V_{\overline K}$ is a projective normal surface 
such that $V_{\overline K}$ has at worst canonical singularities and $-K_{V_{\overline K}}$ is ample. 
\item $1 \leq K_{V_K}^2 =K_{V_{\overline K}}^2  \leq 9$. 
\item $\Pic\,V_{\overline K}$ is a finitely generated free abelian group. 
Furthermore, $\Pic\,V_K \simeq \Z$. 
\end{enumerate}
\end{prop}

\begin{proof}
The assertion (1) follows from the same argument as in \cite[Lemma 5.8]{Kaw21}. 
The assertion (2) holds by \cite[Theorem 15.2]{FS20} and \cite[Theorem 3.3]{BT22}.

Let us show (3). 
Let $\mu : Y \to V_{\overline K}$ be the minimal resolution of $V_{\overline K}$. 
By (2), we have $K_Y = \mu^* K_{V_{\overline K}}$ and $K_Y^2 = K_{V_{\overline K}}^2$. 
Since $Y$ is a smooth rational surface such that $-K_Y$ is nef and big, it holds that $1 \leq K_Y^2 \leq 9$. 
Thus (3) holds.

Let us show (4). 
By $\rho(V/W)=1$, we get $\rho(V_K) =1$ \cite[Lemma 6.6(2)]{Tan18b}. 
We have the injective group homomorphisms 
\[
\Pic\,V_K \hookrightarrow \Pic \,V_{\overline{K}} \hookrightarrow \Pic \,Y
\]
induced by the pullbacks. 
Hence it suffices to show that $\Pic\,Y$ is a finitely generated free $\Z$-module, 
which is well known. 
Thus (4) holds. 
\qedhere





\end{proof}


\begin{prop}\label{p-dP-fibres2}
Let $V$ be a smooth threefold and  $\tau : V \to W$ be a del Pezzo fibration.  
Then general fibres are geoemtrically canonical del Pezzo surfaces, i.e., 
there exist a non-empty open subset $W' \subset W$ such that 
if $F$ is an algebraically closed field and $\Spec\,F \to W'$ is a morphism, 
then $V \times_{W'} \Spec\,F$ is a canonical del Pezzo surface over $F$. 
\end{prop}

\begin{proof}
Let $V_{\overline K}$ be the geometric generic fibre, 
where $K := K(W)$ and $\overline K$ denotes its algebraic closure. 
Let $\mu : Y \to V_{\overline K}$ be the minimal resolution of $V_{\overline K}$. 
Note that this is a log resolution of $V_{\overline K}$ \cite[Section 4]{KM98} (i.e., $\Ex(\mu)$ is simple normal crossing).

We consider the following commutative diagram 
\[
\begin{CD}
\widetilde Y\\
@VV\widetilde \mu V\\
V_1 @>\alpha >> V\\
@VV\tau_1 V @VV\tau V\\
W_1 @>\beta >> W 
\end{CD}
\]
induced by a suitable quasi-finite morphism $\beta : W_1 \to W$ from a smooth curve $W_1$, 
where the square is cartesian and 
the base change of $\widetilde \mu : \widetilde Y \to V_1$ by $(-) \times_{W_1} \Spec\,\overline K$ coincides with $\mu : Y \to V_{\overline K}$. 
Then, for a general point $w_1$ of $W_1$ and the geometric point $\overline{w}_1$ 
corresponding to its algebraic closure, 
$\mu_{\overline{w}_1} : (\widetilde Y)_{\overline{w}_1} \to (V_1)_{\overline{w}_1}$ is a log resolution 
satisfying $K_{(\widetilde Y)_{\overline{w}_1}} = \mu_{\overline{w}_1}^* K_{(V_1)_{\overline{w}_1}}$. 
Then $(V_1)_{\overline{w}_1}$ is a canonical del Pezzo surface. 
We are done, because the image of $W_1 \to W$ contains a non-empty open subset $W'$ of $W$. 
\end{proof}

\subsection{Conic bundles}\label{ss-conic-bdl}

In this subsection, we recall some terminologies and results in \cite{Tan-conic}.

\begin{dfn}\label{d-cb}
We say that $\tau : V \to S$ is a {\em conic bundle} if 
$\tau : V \to S$ is a flat projective morphism of noetherian schemes 
such that the scheme-theoretic fibre 
$X_s :=\tau^{-1}(s)$ is isomorphic to a conic on $\mathbb P^2_{\kappa(s)}$ for every point $s \in S$.  
\end{dfn}

\begin{prop}\label{conic-embedding}
Let $\tau :V \to S$ be a conic bundle, where $V$ and $S$ are smooth varieties. 
Then the following hold. 
\begin{enumerate}
\item $\tau_*\MO_V = \MO_S$. 
\item $\tau_*\omega_V^{-1}$ is a locally free sheaf of rank $3$. 
\item $\omega_V^{-1}$ is very ample over $S$, and hence it defines a closed immersion 
$\iota: V \hookrightarrow  \mathbb P( \tau_*\omega_V^{-1})$ over $S$. 
\end{enumerate}
\end{prop}

\begin{proof}
See \cite[Lemma 2.5 and Proposition 2.7]{Tan-conic}. 
\end{proof}

\begin{rem}\label{r-cb-base-sm}
If $V$ is a smooth projective threefold and $R$ is a $K_V$-negative extremal ray of type $C$, 
then its contraction $\tau: V \to S$ is a conic bundle to a smooth projective surface $S$ \cite[Theorem 1.1]{Kol91}. 
\end{rem}

\begin{nasi}[Discriminant divisors]\label{n-disc}
Let $\tau : V \to S$ be a conic bundle 
from a smooth threefold $V$ to a smooth surface $S$. 
\begin{enumerate}
\item 
If $\tau$ is generically smooth, then we set $\Delta_{\tau}$ to be the {\em discriminant divisor} \cite[Definition 3.3]{Tan-conic}, which is an effective Cartier divisor on $S$. 
For a point $s \in S$, 
$V_s$ is not smooth if and only if $s \in \Delta_{\tau}$  \cite[Theorem 3.11]{Tan-conic}. 
Note that $\Delta_\tau$ is reduced  and normal crossing if $p \neq 2$  \cite[Proposition 7.2]{Tan-conic}. 
On the other hand,  $\Delta_\tau$ is not necessarily reduced if $p=2$ \cite[Example 7.5]{Tan-conic}. 
\item 
If $\tau$ is not generically smooth, then we set $\Delta_{\tau}$ to be the {\em discriminant bundle} \cite[Definition 5.9]{Tan-conic}, which is an invertible sheaf on $S$ (note that we adopt a slightly different notation from  \cite{Tan-conic}). 
\end{enumerate}
Even if $\tau$ is generically smooth, 
we can define the discriminant bundle of $\tau$ \cite[Definition 5.9]{Tan-conic}, which is isomorphic to $\MO_S(\Delta_{\tau})$ 
 \cite[Remark 5.11]{Tan-conic}. 
 On the other hand, the discriminant divisor has more information than the discriminant bundle, 
 we have introduced the terminologies as above.  
\end{nasi}


\begin{prop}\label{p-cb-formula}
Let $\tau : V \to S$ be a conic bundle, where $V$ is a smooth projective threefold and $S$ is a smooth projective surface. 
Then the following hold. 
\begin{enumerate}
\item $\tau_*K_V \sim -2S$. 
\item 
$-\tau_*(K_{V/S}^2) \equiv \Delta_{\tau}$. 
\item 
\begin{enumerate}
\item If $\tau$ is generically smooth, then $\Delta_{\tau}$ is an effective Cartier divisor on $S$. \item If $\tau$ is not genericaly smooth, then $p=2$ and $\Delta_{\tau} \equiv -K_S$. 
\end{enumerate}
\item 
For every Cartier divisor $D$ on $S$, it holds that 
\[
K_V^2 \cdot \tau^*D = -4 K_S \cdot D -\Delta_{\tau} \cdot D. 
\]
\end{enumerate}
\end{prop}

\begin{proof}
Let us show (1). 
For the function field $K := K(S)$, 
let  $V_K := V \times_S \Spec\,K$ be the generic fibre 
of $\tau: V \to S$.   
By Serre daulity and the Riemann--Roch theorem, it holds that 
\[
h^0(V_K, -K_{V_{K}}) = \chi(V_{K}, -K_{V_{K}}) = 
\deg (-K_{V_K}) + \chi(V_K, \MO_{V_K}) 
\]
\[
= \deg (-K_{V_K}) + h^0(V_K, \MO_{V_K}) > 0. 
\]
  Hence there is an effective Cartier divisor $T_0$ on $V_K$ such that $-K_{V_K} \sim T_0$.
  Moreover, $\deg T_0 = \deg (-K_{V_K})=2$.
  Let $T$ be the closure of $T_0$ in $V$, 
  which is an effective divisor on $V$.
  Then 
   there exists a divisor $F$ on $V$ such that 
    $-K_V \sim T + F$ and  $\tau(\Supp F) \neq S$. 
By $\deg T_0 = \deg (-K_{V_K})=2$, $\tau|_T : T \hookrightarrow V \xrightarrow{\tau} S$ is a generically finite morphism of degree two.
  Hence $\tau_*(-K_V) \sim \tau_*T=2S$.
 Thus (1) holds.

The assertion (2) follows from \cite[Theorem 5.15]{Tan-conic}. 
We obtain (3) by (\ref{n-disc}) and \cite[Proposition 3]{MS03}. 
The assertion (4) holds by the following: 
\begin{eqnarray*}
-\Delta_{\tau} \cdot D 
&\overset{{\rm (i)}}{=}& \tau_*(K_{V/S}^2) \cdot D \\
&\overset{{\rm (ii)}}{=}& K_{V/S}^2 \cdot \tau^*D\\
&=& (K_V -\tau^*K_S)^2 \cdot \tau^*D\\
&\overset{{\rm (iii)}}{=}& (K_V^2 -2K_V \cdot \tau^*K_S) \cdot \tau^*D\\
&\overset{{\rm (iv)}}{=}& K_V^2 \cdot \tau^*D  + 4 K_S \cdot D,  
\end{eqnarray*}
Here (i) holds by (2), (ii) follows from the projection formula, 
we get (iii) by $\tau^*K_S \cdot \tau^*K_S \cdot \tau^*D =0$, 
and (iv) holds by (1). 
\qedhere

\end{proof}

\subsection{Picard groups}\label{ss-Pic-gp}

\begin{prop}\label{p-cont-ex}
Let $V$ be a smooth projective threefold. 
Let $R$ be a $K_V$-negative extremal ray $R$ of $\overline{\NE}(V)$ and 
let $\tau : V \to W$ be the contraction of $R$. 
Let $\ell$ be a curve on $V$ satisfying $[\ell ] \in R$. 
Then the following sequence 
\[ 
0 \longrightarrow \Pic(W) \xrightarrow{\tau^*} \Pic(V) \xrightarrow{(-\cdot \ell)} \mathbb{Z}  \]
is exact,     where $(- \cdot \ell)(D):= D\cdot \ell$ for $D\in \Pic(V)$.
\end{prop}

\begin{proof}

\setcounter{step}{0}

\begin{step}\label{s1-cont-ex}
In order to show the assertion of Proposition \ref{p-cont-ex}, it is enough to prove the following $(*)$. 
\begin{enumerate}
\item[$(*)$] If $D$ is a nef Cartier divisor on $V$ such that $\overline{\NE}(V) \cap D^{\perp} = R$, then 
one of the following holds. 
\begin{enumerate}
\item[$(**)$] 
There exists a Cartier divisor $D_W$ on $W$ such that $D \sim \tau^*D_W$. 
\item[$(***)$] 
There exists $m \in \Z_{>0}$ such that both $|mD|$ and $|(m+1)D|$ are base point free. 
\end{enumerate}
\end{enumerate}
\end{step}

\begin{proof}[Proof of Step \ref{s1-cont-ex}]
It follows from $\tau_*\MO_V = \MO_W$ that $\tau^* : \Pic\,W \to \Pic\,V$ is injective. 
It is obvious that 
the composite group homomorphism 
\[
\Pic(W) \xrightarrow{\tau^*} \Pic(V) \xrightarrow{(-\cdot \ell)} \mathbb{Z} 
\]
is zero. 

Let us show that $(**)$ implies the assertion of Proposition \ref{p-cont-ex}. 
Assume $(**)$. 
Fix a Cartier divisor $D$ with $D \cdot \ell =0$. 
Then it suffices to find a Cartier divisor $D_W$ on $W$ such that $D \sim \tau^*D_W$. 
Fix an ample Cartier divisor $A_W$ on $W$. 
By the cone theorem \cite[Theorem 1.24]{KM98}, 
we can find a sufficiently large positive integer $m \gg 0$ such that $E := D +m \tau^*A_W$ is nef and $\overline{\NE}(V) \cap E^{\perp} = R$. 
By $(**)$, there exists a Cartier divisor $E_W$ on $W$ such that $E \sim \tau^*E_W$. 
Then the assertion of Proposition \ref{p-cont-ex} holds, because $D =E -m \tau^*A_W \sim \tau^*(E_W - mA_W)$. 

It suffices to prove that $(***)$ implies $(**)$. 
Assume $(***)$. 
Let $\varphi_{|mD|} : V \to W_m$ and $\varphi_{|(m+1)D|} : V \to W_{m+1}$ be the induced morphisms. 
Then their Stein factorisations coincide with $\tau: V \to W$. 
Hence there exist Cartier divisors $A$ and $B$ on $W$ such that  
$mD \sim \tau^*A$ and $(m+1)D \sim \tau^*B$. 
We then obtain $D  = (m+1)D - mD \sim \tau^*B -\tau^*A = \tau^*(B-A)$. 
Thus $(**)$ {\cred holds}. 
This completes the proof of Step \ref{s1-cont-ex}. 
\end{proof}

\begin{step}\label{s2-cont-ex}
The assertion of Proposition \ref{p-cont-ex} holds if $\dim W=0$. 
\end{step}

\begin{proof}[Proof of Step \ref{s2-cont-ex}]
By $\dim W =0$, $V$ is a Fano threefold with $\rho(X)=1$. 
Pick a nef Cartier divisor $D$ on $V$ such that $\overline{\NE}(V) \cap D^{\perp} = R$. 
Then $D \equiv 0$. 
We obtain $D \sim 0$ by Proposition \ref{p-RR}. 
Thus ($**$) holds. 
This completes the proof of Step \ref{s2-cont-ex}. 
\end{proof}

\begin{step}\label{s3-cont-ex}
The assertion of Proposition \ref{p-cont-ex} holds if $\dim W=1$. 
\end{step}

\begin{proof}[Proof of Step \ref{s3-cont-ex}]
By $\dim W =1$, $\tau: V \to W $ is a del Pezzo fibration. 
Let $D$ be a Cartier divisor on $V$ such that $D\cdot \ell$=0. 
By $\rho(V/W)=1$ and $\rho(V_K)=1$ (Proposition \ref{p-dP-fibres}), 
we obtain $D|_{X_K} \equiv 0$. 
Since $V_{\overline{K}}$ is a canonical del Pezzo surface (Proposition \ref{p-dP-fibres}), 
we have $D|_{V_{\overline{K}}} \sim 0$, which implies $D_{V_K} \sim 0$ by the flat base change theorem. 
  Hence we can write $D|_{V_K}=\mathrm{div}(\varphi)$ for some $\varphi \in K(V_K)=K(V)$. 
The support of the divisor $D-\mathrm{div}(\varphi)$ on $V$ is $\tau$-vertical, i.e.,  
$\tau(\Supp(D-\mathrm{div}(\varphi))) \neq W$. 
Therefore, $D-\mathrm{div}(\varphi) = \sum_{i=1}^r E_i$ 
for some prime divisors $E_1, ..., E_r$ such that 
each $\tau(E_i)$ is a point, 
where $E_i =E_j$ might hold even if $i \neq j$. 
For $w_i := \tau(E_i)$, we have $E_i =\tau^*w_i$ (Proposition \ref{p-dP-fibres}(1)). 
Therefore, 
\[
D \sim 
\sum_{i=1}^r E_i =
\sum_{i=1}^r \tau^*w_i = \tau^*\left( \sum_{i=1}^r w_i\right)
=\tau^*D_Y
\]
for $D_Y := \sum_{i=1}^r w_i$. 
Thus ($**$) holds. 
This completes the proof of Step \ref{s3-cont-ex}. 
\end{proof}

\begin{step}\label{s4-cont-ex}
The assertion of Proposition \ref{p-cont-ex} holds if $\dim W=2$. 
\end{step}

\begin{proof}[Proof of Step \ref{s4-cont-ex}]
Let $D$ be a Cartier divisor on $V$ such that $D\cdot \ell=0$. 
Recall that $\tau : V \to W$ is flat and $W$ is smooth \cite[Theorem 1.1]{Kol91}. 
For a closed point $w \in W$, the fibre $V_w$ is a (possibly reducible or non-reduced) conic on $\P^2$. 
Since $\MO_V(D)|_{V_w}$ is numerically trivial, 
we obtain $\MO_V(D)|_{V_w} \simeq \MO_{V_w}$, 
and hence $h^0(V_w, \MO_V(D)|_{V_w}) =1$ 
(if $Z:=V_w$ is a non-reduced conic, then 
we have $0 \to  I_{Z_{\red}} \to \MO_Z \to \MO_{Z_{\red}} \to 0$
\cite[Ch. III, Exercise 4.6]{Har77}).  
By the Grauert theorem \cite[Ch. III, Corollary 12.9]{Har77}, 
$\tau_*\MO_V(D)$ is an invertible sheaf and $\MO_V(D)$ is $\tau$-free, i.e., 
$\tau^*\tau_*\MO_V(D) \to \MO_V(D)$ is surjective. 
As both $\tau^*\tau_*\MO_V(D)$ and $\MO_V(D)$ are invertible sheaves, 
we get $\tau^*\tau_*\MO_V(D) \xrightarrow{\simeq} \MO_V(D)$. 
For a Cartier divisor $D_W$ satisfying $\tau_*\MO_V(D) \simeq \MO_W(D_W)$, 
we obtain $D \sim \tau^*D_W$.  
Thus ($**$) holds.  
This completes the proof of Step \ref{s4-cont-ex}. 
\qedhere


\end{proof}

\begin{step}\label{s5-cont-ex}
The assertion of Proposition \ref{p-cont-ex} holds if $\dim W=3$. 
\end{step}

\begin{proof}[Proof of Step \ref{s5-cont-ex}]
Fix a nef Cartier divisor $D$ on $V$ such that $\overline{\NE}(V) \cap D^{\perp} = R$. 
It suffices to show $(***)$. 
By $\dim W =3$, $\tau: V \to W$ is a birational morphism. 
Set $E := \Ex(\tau)$, which is a prime divisor such that $E \cdot \ell <0$ \cite[(1.1.2)]{Kol91}. 
Fix $n \in \Z_{>0}$ such that 
$A:=nD-E$ is an ample divisor on $V$. 
By the exact sequence 
  \[0 \longrightarrow \MO_V(mD-E) \longrightarrow \MO_V(mD) \longrightarrow \MO_E(mD) \longrightarrow 0,\]
  we obtain the exact sequence 
  \[H^0(V,\MO_V(mD))\longrightarrow H^0(E,\MO_E(mD)) \longrightarrow H^1(V,\MO_V(mD-E)).\]

  \begin{clm}\label{cl-cont-ex}
  The following hold for every sufficiently large integer $m \gg 0$. 
  \begin{enumerate}
      \item The set-theoretic inclusion $\mathrm{Bs}\,|mD|\subset E$ holds.
      \item $mD -E-K_V$ is ample.
      \item $|mD|_E|$ is base point free. 
      \item There exists $m_0 \in \Z_{>0}$ such that $m_0D \sim \tau^*D_W$ for some ample Cartier divisor $D_W$ on $W$
  \end{enumerate}
  \end{clm}

  \begin{proof}[Proof of Claim \ref{cl-cont-ex}]
Let us show (1). 
    Let $q,r\in \mathbb{Z}_{\geq 0}$ be the non-negative integers that  satisfy 
$m=qn+r$ and $0\le r < n$.
    Then 
    \begin{align*}
      mD = qnD + rD = qA + rD + qE.
    \end{align*}
    Since $n$ is fixed and there are only finitely possibilities for $r$, 
    $|qA+rD|$ is very ample for every $q \gg 0$. 
    In particular, $|qA+rD|$ is base point free for $q \gg 0$. 
    Hence we get $\mathrm{Bs}\,|mD|\subset E$ for every  $m \gg 0$. 
    Thus (1) holds. 

Since $K_V \cdot \ell <0$, $E \cdot \ell<0$, and 
$\overline{\NE}(V) \cap D^{\perp} =\R_{\geq 0}[\ell]$, 
the assertion (2) follows from the cone theorem. 
Then (2) implies (3) by 
\cite[Theorem 0.4 and Corollary 3.6]{Tan15} (note that $E$ is a projective canonical surface).

Let us show (4). 
It suffices to show that $D$ is semi-ample, i.e., 
$|m_0D|$ is base point free for some $m_0 \in \Z_{>0}$. 
By the cone theorem again, we see that $m_1D -E$ is ample for some $m_1\in \Z_{>0}$. 
By Keel's result \cite[Proposition 1.6]{Kee99}, 
it is enough to show that $|m_1D|_E|$ is semi-ample, which follows from (3). 
This completes the proof of Claim \ref{cl-cont-ex}. 
  \end{proof}

By Claim \ref{cl-cont-ex}(1)(3), it is enough 
to show $H^1(V,\MO_V(mD-E))=0$ for every $m \gg 0$. 
We first treat the case when $f(E) =:C$ is a smooth curve. 
In this case, $\tau$ is the blowup along $C$. 
In particular, $\dim \tau^{-1}(w) \leq 1$ for any $w \in W$. 
Since $mD-E -K_V$ is ample (Claim \ref{cl-cont-ex}(2)), 
it follows from \cite[Theorem 0.5]{Tan15} that $R^1\tau_*\MO_X(mD-E)=0$ 
for every $m \gg 0$. 
Recall that we have the following exact sequence induced by 
the corresponding Leray spectral sequence: 
\[
0 \to H^1(W, \tau_*\MO_V(mD-E)) \to H^1(V,\MO_V(mD-E)) \to 
H^0(W, R^1\tau_*\MO_V(mD-E))=0. 
\]
We can uniquely write $m = m_0q + r$ for the non-negative integers $q, r$ satisfying $0 \leq r <m_0$. 
In particular, if $m \gg 0$, then $q \gg 0$. 
It holds that 
\begin{eqnarray*}
H^1(W, \tau_*\MO_V(mD-E)) &\simeq& 
H^1(W, \tau_*\MO_V(q\tau^*D_W + rD -E)) \\
 &\simeq& 
H^1(W, \MO_W(qD_W) \otimes \tau_*\MO_V(rD -E)) \\
&=& 0, 
\end{eqnarray*}
where the last equality holds by the Serre vanishing theorem. 
This completes the proof for the case when $f(E)$ is a curve. 


Assume that $f(E)$ is a point. 
There is $m_1\in \mathbb{Z}_{>0}$ such that $mD-E$ is ample for every integer $m\geq m_1$.
By the Fujita vanishing theorem \cite[Theorem 3.1]{Fuj17}, 
we can find $m_2\in \mathbb{Z}_{>0}$ such that $H^1(V,\MO_V(mD-m_2E))=0$ for every  $m\geq m_1m_2$.
  By the exact sequence
  \[0\longrightarrow \MO_V(-E)\longrightarrow \MO_V\longrightarrow \MO_E\longrightarrow 0,\]
  we obtain the following exact sequences
  \[0\longrightarrow \MO_V(mD-2E)\longrightarrow \MO_V(mD-E) \longrightarrow \MO_E(mD-E) \longrightarrow 0,\]
  \[0\longrightarrow \MO_V(mD-3E)\longrightarrow \MO_V(mD-2E) \longrightarrow \MO_E(mD-2E) \longrightarrow 0,\]
  \begin{center}
    $\vdots$
  \end{center}
  where $\MO_E(mD-sE) := \MO_V(mD-sE)|_E$. 
In order to show $H^1(V,\MO_V(mD-E))=0$ for $m \gg 0$, it is enough to show that 
\[
H^1(E,\MO_E(mD-sE))=0
\]
 for all $1 \leq s\leq m_2$ and $m \gg 0$. As $m_2$ is fixed, 
 we may assume that $mD-sE$ is ample. 
Since  $\tau(E)$ is a point, i.e., $\tau$ is of type $E_2-E_5$, 
$E$ is a projective normal toric surface (Definition \ref{d-type-ext-ray}). 
As $\MO_E(mD-sE)$ is ample, we get $H^1(E,\MO_E(mD-sE))=0$ by \cite[Corollary 1.3]{Fuj07}. 
This completes the proof of Step \ref{s5-cont-ex}. 
\end{proof}

This completes the proof of Proposition \ref{p-cont-ex}. 
\end{proof}

\subsection{Del Pezzo fibrations II}\label{ss-dP-II}

\begin{lem}\label{l-del_pezzo-3}
Let $V$ be a smooth projective 
threefold. 
Let $R$ be an extremal ray of $\overline{\NE}(X)$ of type $D_3$ and 
let $\tau : V \to W$ be the contraction of $R$. 
For $K := K(W)$, set $V_K$ to be the generic fibre of $\tau$. 
Then the following hold. 
\begin{enumerate}
\item $V_K \simeq \mathbb P^2_K$. In particular, general fibres of $\tau$ are $\mathbb P^2$. 
\item There exist Cartier divisors $D$ on $V$ and $E$ on $W$ such that $K_{V} \sim 3D+f^*E$. 
\end{enumerate}
\end{lem}

\begin{proof}
Let us show (1). 
We have $K_{V_{\overline K}} \simeq \mathbb P^2_{\overline K}$. 
In particular, $V_K$ is a Severi-Brauer surface over $K$. 
Since $V_K$ has a $K$-rational point \cite[Ch. IV, Theorem 6.8]{Kol96},  
it follows from \cite[Theorem 5.1.3]{GS17} that $X_K \simeq \mathbb P^2_K$. 
Thus (1) holds.

Let us show (2).  
Fix a general closed point $0 \in W$. 
By (1), we can find a Cartier divisor $D_0$ on $V_0 := \tau^{-1}(0)$ such that $\omega_{V_0} \simeq \MO_{V_0}(3D_0)$. 
By \cite[Corollary 8.5.6(a)]{FGI05}, there exists a Cartier divisor $D'$ on $V' := V \times_W \Spec\,\MO_{W, 0}$ 
such that $D'|_{V_0} \sim D_0$. 
Taking the closure, we obtain a Cartier divisor $D$ on $V$ such that $D|_{V_0} \sim D_0$. 
It holds that $K_V -3D \sim f^*E$ for some Cartier divisor $E$ on $W$ (Proposition \ref{p-cont-ex}). 
Thus (2) holds. 
\end{proof}

\begin{lem}\label{l-del_pezzo-2}
Let $V$ be a smooth projective 
threefold. 
Let $R$ be an extremal ray of $\overline{\NE}(X)$ of type $D_2$ and 
let $\tau : V \to W$ be the contraction of $R$. 
For $K := K(W)$, set $V_K$ to be the generic fibre of $\tau$. 
Then the following hold. 
\begin{enumerate}
\item If $V_K$ is smooth over $K$ and $w$ is a general closed point of $W$,  then 
each of $V_{\overline K}$ and $V_w$ is isomorphic to $\mathbb P^1 \times \mathbb P^1$. 
\item If $V_K$ is not smooth over $K$ and $w$ is a general closed point of $W$, 
then $p =2$ and each of $V_{\overline K}$ and $V_w$ is isomorphic to $\mathbb P(1, 1, 2)$, which is isomorphic to a singular quadric surface on $\P^3$. 
\item There exist Cartier divisors $D$ on $V$ and $E$ on $W$ such that $K_{V} \sim 2D+\tau^*E$. 
\end{enumerate}
\end{lem}

\begin{proof}
Let us show (1). 
Assume that $V_K$ is smooth over $K$. 
Since $V_{\overline K}$ is a smooth del Pezzo surface over $\overline K$ with $K_{V_{\overline K}}^2 =8$, 
it holds that $V_{\overline K}$ is isomorphic to 
either $\mathbb P^1 \times \mathbb P^1$ or 
$\F_1 := \P_{\P^1}(\MO_{\P^1} \oplus \MO_{\P^1}(1))$. 
Note that if $V_{\overline K}$ is isomorphic to $\mathbb P^1 \times \mathbb P^1$, then so is $V_w$, because $-K_{V/W}$ will be divisible by two after taking suitable base change.


Suppose  $V_{\overline K} \simeq \F_1$. 
It suffices to derive a contradiction. 
Let $K^s$ be the separable closure of $K$ satisfying $K \subset K^s \subset \overline K$. 
For the unique $(-1)$-curve $\overline{\Gamma}$ on $V_{\overline K}$, 
we set $\Gamma^s \subset V_{K^s}$  to be the image of $\overline{\Gamma}$, which is a prime divisor on $V_{K^s}$. 
Since $V_{\overline K} \to V_{K^s}$ is a universal homeomorphism, 
$\Gamma^s$ is a unique curve on $V_{K_s}$ 
that satisfies $(\Gamma^s)^2 <0$ \cite[Lemma 2.3]{Tan18b}.

Since we have $\rho(V_{K^s})= \rho(V_{\overline K}) = 2$ \cite[Proposition 2.4]{Tan18b} and $\Gamma^s$ is a unique curve on $V_{K^s}$ satisfying $(\Gamma^s)^2 <0$, 
we can check that $\MO_{V_{K^s}}(\Gamma^s)$ is a unique primitive line bundle on $V_{K^s}$ such that 
$\kappa(V_{K^s},  \MO_{V_{K^s}}(\Gamma^s))=0$. 
In particular, $\MO_{V_{K^s}}(\Gamma^s)$ is a $G$-invariant element of $\Pic\,(V_{K^s})$ for $G := {\rm Gal}(K^s/K)$. 
It follows from Proposition \ref{p-Br-vanish}
that $\mathrm{Br}(K)=0$. 
Then we obtain $\Pic(V_K) \simeq \Pic(V_{K^s})^G$ (Proposition \ref{p-Br-Pic}). 
Hence we can write $\MO_{V_{K^s}}(\Gamma^s)\simeq \beta^*\mathcal{L}$ for some $\mathcal{L}\in \Pic(V_K)$, 
where $\beta \colon V_{K^s}\to V_K$ denotes the induced morphism. 
  By the flat base change theorem, we get 
  \[
  K^s \simeq H^0(V_{K^s},\MO_{V_{K^s}}(\Gamma^s)) 
  \simeq H^0(V_{K^s},\beta^*\mathcal{L}) \simeq H^0(V_K,\mathcal{L})\otimes _K K^s, 
  \]
  which implies  $H^0(V_K,\mathcal{L})\simeq K$. 
  Then there exists an effective divisor $\Gamma$ on $V_K$ such that $\mathcal{L}\simeq \MO_{V_K}(\Gamma)$.
  Then $\MO_{V_{K^s}}(\Gamma^s) \simeq \beta^*\MO_{V_{K^s}}(\Gamma) \simeq \MO_{V_{K^s}}(\beta^*\Gamma)$ and $\Gamma^s \sim \beta^*\Gamma$. 
By $(\Gamma^s)^2<0$, we get $\Gamma^2<0$. 
Then 
we can find a prime divisor $C$ on $V_K$ 
satisfying $C \subset \Supp\,\Gamma$ and $C^2 <0$. 
However, this contradicts $\rho(V_K) = 1$ (Proposition \ref{p-dP-fibres}). 
Thus (1) holds. 

Let us show (2). 
Assume that $V_K$ is not smooth over $K$. 
Since $V_{\overline K}$ has at worst canonical singularities 
(Proposition \ref{p-dP-fibres}), 
it follows from  \cite[Theorem 6.1]{Sch08}, 
that the possibilities for each singularity of $V_{\overline{K}}$ are 
$A_n(n=p^e-1,e\in \mathbb{N}_{>0}), D_n(n\ge 4)$, and $E_n(n=6,7,8)$.
Since $K_W^2=K_{V_{\overline{K}}}^2=8$ and $K_W^2+\rho(W)=10$, 
only the possibility is the case when $p=2$ and $V_{\overline K} \simeq \P(1, 1, 2)$. 
Then $V_{\overline K}$ is a singular quadric surface in $\P^3_{\overline K}$. 
By the same argument as in (1), also $V_w$ is a quadric surface in $\P^3_k$. 
Thus (2) holds. 


Let us show (3).  
Fix a general closed point $0 \in W$. 
Then there exists a Cartier divisor $D_0$ on $V_0$ such that $\omega_{V_0} \simeq \MO_{V_0}(2D_0)$. 
By \cite[Corollary 8.5.6(a)]{FGI05}, there exists a Cartier divisor $D'$ on $V' := V \times_W \Spec\,\MO_{W, 0}$ 
such that $D'|_{V_0} \sim D_0$. 
Taking the closure, we obtain a Cartier divisor $D$ on $V$ such that $D|_{V_0} \sim D_0$. 
It holds that $K_V -2D \sim f^*E$ for some Cartier divisor $E$ on $W$ (Proposition \ref{p-cont-ex}). 
\qedhere 

\end{proof}

\begin{rem}
By \cite[Proposition 14.7]{FS20}, there actually exists a del Pezzo fibration $\tau: V \to W$ 
from a smooth threefold $V$ such that the generic fibre is not smooth. 
\end{rem}

\begin{prop}\label{p-D-length}
Let $V$ be a smooth projective threefold and 
let $R$ be a $K_V$-negative extremal ray of $\overline{\NE}(V)$. 
  \begin{enumerate}
    \item If $R$ is of type $D_1$, then $\mu_R=1$.
    \item If $R$ is of type $D_2$, then $\mu_R=2$.
    \item If $R$ is of type $D_3$, then $\mu_R=3$. 
  \end{enumerate}
\end{prop}

\begin{proof}
Let $\tau : V \to W$ be the contraction of $R$. 
Set $V_K$ to be the generic fibre of $\tau$. 

Let us show (1). 
Assume that $\tau$ is of type $D_1$. 
By Definition \ref{d-type-ext-ray}, 
we have $1\le K_{V_K}^2 \le 7$.
Fix a general closed point $w \in W$. 
Then $V_w$ is a canonical del Pezzo surface with $1 \leq K_{V_w}^2 \le 7$ 
(Proposition \ref{p-dP-fibres2}).

  Let $\mu\colon Y\to V_w$ be the minimal resolution of $V_w$. 
Pick a $(-1)$-curve $E$ on $Y$, 
whose existence is guaranteed by $\rho(Y)  = 10 - K_{Y}^2 = 10-K_{V_w}^2 \geq 3$. 
By $K_Y=\mu^*K_{V_w}$, 
we obtain
  \begin{align*}
  -1 = K_Y \cdot E  = K_{V_w}\cdot \mu(E) =  K_V \cdot \mu(E).
  \end{align*}
  Hence ${\cred \mu_R} =1$. Thus (1) holds.

Let us show (2) and (3). 
Fix $r \in \{2, 3\}$ and assume that $\tau$ is of type $D_r$. 
The goal is to prove $\mu_R = r$. 
Fix an extremal rational curve $\ell_R$, which satisfies $-K_V \cdot \ell_R = \mu_R>0$. 
We can write 
\[
-K_V \sim r D + \tau^*E. 
\]
for some Cartier divisors $D$ on $V$ and $E$ on $W$ 
(Lemma \ref{l-del_pezzo-3}, Lemma \ref{l-del_pezzo-2}). 
In particular, $\mu_R = -K_V \cdot \ell_R = rD \cdot \ell_R \geq r$. 
Pick a general closed point $w \in W$. 
If $r=3$ (resp. $r=2$), then $V_w \simeq \P^2$ (resp. $V_w$ is a (possibly singular) quadric surface in $\P^3$). 
In any case, let $\ell$ be a line on $V_w$. 
Then we obtain $D \cdot \ell =1$. 
Hence $\mu_R \leq -K_V \cdot \ell = r D \cdot \ell =r$. 
Thus (2) and (3) hold. 
\qedhere

\end{proof}

\subsection{Extremal rays on weak Fano threefolds}\label{ss-weak-Fano}


\begin{dfn}
We say that $V$ is a {\em weak Fano threefold} if 
$V$ is a smooth projective threefold such that $-K_V$ is nef and big. 
\end{dfn}

In our application, we shall use 
$K_V$-negative extremal rays on weak Fano threefolds such that 
$\rho(V) =2$ and $H^1(V, \MO_V) =0$. 
In the rest of this section, we summarise properties on extremal rays under these assumptions. 


\begin{prop}[Type C]\label{p-typeC-intersec}
Let $V$ be a weak Fano threefold such that $\rho(V) =2$. 
Let $\tau : V \to W$ be a contraction of a $K_V$-negative extremal ray $R$ of type $C$. 
Then the following hold. 
\begin{enumerate}
    \item 
    $W \simeq \P^2$. 
    \item 
    For  a line $L$ on $W \simeq \P^2$ and $D :=\tau^*L$, the following hold: 
    \[
    D^3 =0, \qquad D^2 \cdot (-K_V) =2, \qquad D \cdot (-K_V)^2 = 12 - \deg \Delta_{\tau}. 
    \]
    \item \begin{itemize}
    \item 
    $\tau$ is of type $C_1$ $\Leftrightarrow$ $\mu_R=1$ $\Leftrightarrow$ $\Delta_{\tau}$ is ample. 
    \item 
    $\tau$ is of type $C_2$ $\Leftrightarrow$ $\mu_R=2$ $\Leftrightarrow$ $\deg \Delta_{\tau} =0$. 
    \end{itemize}
\end{enumerate}
\end{prop}

\begin{proof}
Let us show (1). 
Recall that $W$ is a smooth projective surface with $\rho(W)=1$. 
Fix an ample prime divisor $C$ on $W$. 
It is enough to show that $-K_W \cdot C >0$. 
By Proposition \ref{p-cb-formula}(4), we obtain 
\[
-4K_W \cdot C = (-K_V)^2 \cdot \tau^*C +\Delta_{\tau} \cdot C > \Delta_{\tau} \cdot C. 
\]
If $\tau$ is generically smooth, then $\Delta_{\tau}$ is an effective Cartier divisor on ${\cred W}$ (Proposition \ref{p-cb-formula}), 
and hence we get $-4K_W \cdot C  > \Delta_{\tau} \cdot C \geq 0$ by the ampleness of $C$. 
If $\tau$ is not generically smooth, then we have 
$-K_W \equiv \Delta_{\tau}$ (Proposition \ref{p-cb-formula}), 
which implies $-4K_W \cdot C >\Delta_{\tau} \cdot C = -K_W \cdot C$, and hence ${\cred -K_W} \cdot C >0$. 
Thus (1) holds.

Let us show (2). 
It is clear that $D^3 =0$. 
By $D^2 = \tau^*L^2 \equiv V_w$ for a fibre $V_w$ of $\tau: V \to W$, 
we get $D^2 \cdot (-K_V) = (-K_V) \cdot V_w= -\deg \omega_{V_w}=2$. 
The remaining one holds by 
\[
D \cdot (-K_V)^2 = \tau^*L \cdot K_V^2 = -4K_W \cdot L - \Delta_{\tau} \cdot L = 12 - \deg \Delta_{\tau}.
\]
where the 
{\cred second} equality holds by Proposition \ref{p-cb-formula}. 
Thus (2) holds.


Let us show (3). It is clear that $\mu_R=1$ if and only if there exists a non-smooth fibre of $\tau$. 
If $\tau : V \to W$ is generically smooth, then 
$\Delta_{\tau}$ is an effective divisor whose support parametrises the non-smooth fibres. 
If $\tau : V \to W$ is not generically smooth, then 
$\Delta_{\tau} \equiv -K_W$ (Proposition \ref{p-cb-formula}), and hence $\Delta_{\tau}$ is ample by (1). Thus (3) holds. 
\end{proof}

\begin{prop}[Type D]\label{p-typeD-intersec}
Let $V$ be a weak Fano threefold such that $\rho(V) =2$ and $H^1(V, \MO_V) =0$. 
Let $\tau : V \to W$ be a contraction of a $K_V$-negative extremal ray $R$ of type $D$. 
Then the following hold. 
\begin{enumerate}
    \item 
    $W \simeq \P^1$. 
    \item 
    For  a closed point $P$ on $W$ and $D :=\tau^*P$, the following hold: 
    \[
    D^3 = D^2 \cdot (-K_V) =0, \qquad 1 \leq D \cdot (-K_V)^2 \leq 9. 
    \]
    \item \begin{itemize}
    \item $\tau$ is of type $D_1$ $\Leftrightarrow$ $\mu_R=1$ $\Leftrightarrow$ $1 \leq D \cdot (-K_V)^2 \leq 7$.  
    \item $\tau$ is of type $D_2$ $\Leftrightarrow$ $\mu_R=2$ $\Leftrightarrow$ $D \cdot (-K_V)^2 = 8$.  
    \item $\tau$ is of type $D_2$ $\Leftrightarrow$ $\mu_R=3$ $\Leftrightarrow$ $D \cdot (-K_V)^2 = 9$.  
    \end{itemize}
\end{enumerate}
\end{prop}

\begin{proof}
By $H^1(W, \MO_W) \hookrightarrow H^1(V, \MO_V)=0$, (1) holds. 
Let us show (2). 
It is obvious that $D^3 = D^2 \cdot (-K_V) =0$. 
After replacing $P$ by a general closed point, 
we may assume that $D$ is a canonical del Pezzo surface (Proposition \ref{p-dP-fibres2}). 
Then we obtain $D \cdot (-K_V)^2 =(-K_D)^2$ and $1 \leq (-K_D)^2 \leq 9$. Thus (2) holds. 
The assertion (3) follows from (2) and Proposition \ref{p-D-length}. 
\end{proof}

Since the extremal ray of type $E_1$ induces a blowup along a smooth curve, 
some intersection numbers are computed by the following lemma.

\begin{lem}\label{l-blowup-formula}
Let $X$ be a smooth projective threefold. 
Set $g := \frac{(-K_X)^3}{2}+1 \in \Q$. 
\begin{enumerate}
\item 
Let $P$ be a point on $X$ and let $\sigma : Y \to X$ be the blowup at $P$. 
Set $E := \Ex(\sigma)$. 
Then $K_Y = \sigma^*K_X + 2E$ and the following hold. 
\begin{enumerate}
\item $(-K_Y)^3 = (-K_X)^3 -8 = 2g -10$.  
\item $(-K_Y)^2 \cdot E =4$.
\item $(-K_Y) \cdot E^2 =-2$.
\item $E^3 = 1$. 
\end{enumerate}
\item Let $\Gamma$ be a smooth curve on $X$ of genus $g(\Gamma)$  
and let $\sigma : Y \to X$ be the blowup along $\Gamma$. 
Set $E := \Ex(\sigma)$ and $\deg \Gamma := -K_X \cdot \Gamma$. 
Then  $K_Y = \sigma^*K_X + E$, $\deg N_{\Gamma/X} = 2g(\Gamma) -2 + \deg \Gamma$, and the following hold. 
\begin{enumerate}
\item $(-K_Y)^3  =  2g-2\deg \Gamma +2g(\Gamma)-4$. 
\item $(-K_Y)^2 \cdot E =\deg \Gamma -2g(\Gamma)+2$.
\item $(-K_Y) \cdot E^2 =2g(\Gamma)-2$.
\item $E^3 = -\deg \Gamma + 2-2g(\Gamma)$. 
\end{enumerate}
\end{enumerate}
\end{lem}

Although the assertion is well known, 
we include a sketch of a proof for the reader's convenience (cf. \cite[Lemma 2.2.14 and Lemma 4.1.2]{IP99}). 

\begin{proof}
The assertion (1) follows from $E \simeq \P^2$, $E|_E =\MO_E(-1)$ \cite[Ch. II, Theorem 8.24]{Har77}, and 
$K_Y = \sigma^*K_X+ 2E$. 
Let us show (2). 
We have $K_Y = \sigma^*K_X + E$, $\deg N_{\Gamma/X} = 2g(\Gamma) -2 + \deg \Gamma$, and $E^3 = -\deg N_{\Gamma/X}$  
(cf. \cite[Lemma 2.2.14]{IP99}). Then the following hold: 
\renewcommand{\labelenumi}{(\roman{enumi})}
\begin{itemize}
    \item $\sigma^*(-K_X)^3 = (-K_X)^3 = 2g -2$. 
    \item $\sigma^*(-K_X)^2 \cdot E = (\sigma^*(-K_X)|_E)^2 = 0$. 
    \item $\sigma^*(-K_X) \cdot E^2 = (\sigma^*(-K_X)|_E) \cdot (E|_E) = -\deg \Gamma$. 
    \item $E^3 = -\deg \Gamma + 2-2g(\Gamma)$. 
\end{itemize}
Thus (d) holds. 
The remaining assertions (a)--(c) hold by the following computation: 
\[
(-K_Y)^3 = ( -\sigma^*K_X -E)^3 = (2g-2) + 0 -3 \deg \Gamma  
-(-\deg \Gamma +2 -2g(\Gamma))
\]
\[
(-K_Y)^2 \cdot E =( -\sigma^*K_X -E)^2 \cdot E = 0 - 2(-\sigma^*K_X) \cdot E^2+E^3 = 
-2 (-\deg \Gamma) + (-\deg \Gamma +2 -2g(\Gamma))
\]
\[
(-K_Y) \cdot E^2 =( -\sigma^*K_X -E) \cdot E^2 = (-\sigma^*K_X) \cdot E^2 - E^3 
= -\deg \Gamma - (-\deg \Gamma + 2 -2g(\Gamma)). 
\]


\end{proof}

As the case of type $E_1$ is essentially contained in Lemma \ref{l-blowup-formula}, 
we here summarise the remaining cases of type $E$. 

\begin{prop}[Type $E_2-E_5$]\label{p-typeE-intersec}
Let $V$ be a smooth projective threefold. 
Let $\tau : V \to W$ be a contraction of a $K_V$-negative extremal ray $R$ of type $E$. 
Set $D := \Ex(\tau)$. 
If $E$ is not of type $E_1$, then $\tau(D)$ is a point and 
it satisfies Table \ref{table-typeE}. 
\begin{table}[h]
\caption{Type $E_2-E_5$}\label{table-typeE}
     \centering
{\renewcommand{\arraystretch}{1.35}%
      \begin{tabular}{|c|c|c|c|c|c|c|c|}
      \hline
type of $\tau$  & $D^3$  & $D^2 \cdot (-K_V)$ & $D \cdot (-K_V)^2$ & $K_V -\tau^*K_W$ & $\mu_R$ 
& Description of $D$\\      \hline
$E_2$  & $1$ & $-2$ & $4$ & $2D$ & $2$ & $\P^2$ \\           \hline
$E_3$  & $2$ & $-2$ & $2$ & $D$ & $1$ & a smooth quadric surface in $\P^3$\\           \hline
$E_4$  & $2$ & $-2$ & $2$ & $D$ & $1$ & a singular quadric surface in $\P^3$\\           \hline
$E_5$  & $4$ & $-2$ & $1$ & $\frac{1}{2}D$ & $1$ & $\P^2$\\           \hline
      \end{tabular}}
    \end{table}
\end{prop}

\begin{proof}
The assertion follows from \cite[Theorem 1.1]{Kol91} by using the adjunction formula. 
\end{proof}

\section{Two-ray game: foundations}

\subsection{Introduction}\label{ss-2ray-intro}

Let $X \subset \P^{g+1}$ be an anti-canonically embedded Fano threefold with $\Pic\,X =\Z K_X$. 
In this section, we summarise some foundational results on so-called two-ray games for $X$. 
The results in this section will be applied in Section \ref{s-line-I}-\ref{s-line-II}. 
The purpose of this subsection is to give an overview of its construction and some arguments. 
If the reader is familiar with the proof in characteristic zero, it is advised to skip this subsection. 


We start by the construction of  two-ray game. 
Let $\Gamma$ be either a point or a smooth rational curve on $X$. 
Take the blowup $\sigma : Y \to X$ along $\Gamma$. 
Under the situations appearing in our applications, we can show that $-K_Y$ is semi-ample and big. 
Then $|-mK_Y|$ is base point free for some $m \in \Z_{>0}$. 
Let $\psi : Y \to Z$ be the Stein factorisation of the induced morphism $\varphi_{|-mK_Y|} : Y \to \P^{h^0(Y, -mK_Y)-1}$. 
Since $-K_Y$ is big, $\psi : Y \to Z$ is a birational morphism to a projective normal threefold $Z$. 
There are the following three cases: 
\begin{itemize}
\item[(i)] $\dim \Ex(\psi)=2$ (divisorial contraction). 
\item[(ii)] $\dim \Ex(\psi)=1$ (flopping contraction). 
\item[(iii)] $\psi$ is an isomorphism, i.e., $-K_Y$ is ample. 
\end{itemize}
For simplicity, we now explain how to proceed for the case (ii). 
Since $\rho(Y)=2$, $\NE(Y)$ has exactly two extremal rays, 
which correspond to $\sigma$ and $\psi$. 
We take the flop $\psi^+:Y^+ \to Z$ of $\psi : Y \to Z$. 
It is known that $Y^+$ is smooth. 
By $-K_{Y} = \psi^*(-K_Z)$ and $-K_{Y^+} = (\psi^+)^*(-K_Z)$, $-K_{Y^+}$ is nef and big. 
By $\rho(Y^+)=2$, we have the $K_{Y^+}$-negative extremal ray $R$ not corresponding to $\psi^+: Y^+ \to Z$. 
Let $\tau : Y^+ \to W$ be the contraction of $R$. 
\[
\begin{tikzcd}
Y \arrow[d, "\sigma"'] \arrow[rd, "\psi"]& & Y^+ \arrow[ld, "\psi^+"'] \arrow[d, "\tau"]\\
X & Z & W.
\end{tikzcd}
\]
For the case (iii), the situation is similar. 
For $Z :=Y$, $Y^+ :=Y$, $\psi :={\rm id}$, and $\psi^+:={\rm id}$, 
let $\tau$ be the contraction of the extremal ray not corresponding to $\sigma$. 
We then get the above diagram. 

By using the above construction, we shall  prove the following results: 
\begin{enumerate}
\item $g \leq 12$ (Corollary \ref{c-g-bound}). 
\item $-K_Y$ is not ample, i.e., (iii) does not occur 
if $\Gamma$ is a point, a line, or a conic (Proposition \ref{p-pt-defect}, Proposition \ref{p-conic-defect}, Proposition \ref{p-line-defect}). 
\item If $g \geq 8$, then there exists a conic on $X$ (Theorem \ref{t-pt-to-conic}). 
\end{enumerate}
We now overview how to show (1) and (2).

(1) This inequality $g \leq 12$ will be proven by using the case when $\Gamma$ is a conic. 
Note that we may assume that $g \geq 8$, and hence there exists a conic by (3). 
We here only treat the case when (ii) holds. 
The possibilities for $g$ are restricted by solving suitable Diophantine equation. 
For example, if $\tau$ is of type D, then we obtain $D^3 =D^2 \cdot (-K_{Y^+})=0$ 
for $D := \tau^*({\rm point})$. 
It is easy to see that we can write $D \sim -\alpha K_{Y^+} -\beta E^+$ for some $\alpha, \beta \in \Z_{>0}$ (Lemma \ref{l-beta=r}), where $E^+$ denotes the proper transform $E^+$ of $E := \Ex(\sigma)$. 
For example, the equality $D^2 \cdot (-K_{Y^+}) =0$ induces 
an equation $(g-4) \alpha^2 -4 \alpha\beta -\beta^2=0$ by using the following relations: 
\begin{equation}\label{e1-Y-vs-Y^+}
(-K_{Y^+})^3 = (-K_Y)^3,\,\, (-K_{Y^+})^2 \cdot E^+ =  (-K_{Y})^2 \cdot E,\,\, 
(-K_{Y^+}) \cdot (E^+)^2 =  (-K_{Y}) \cdot E^2. 
\end{equation}
Indeed, these relations imply 
\begin{eqnarray*}
D^2 \cdot (-K_{Y^+}) 
&=& (-\alpha K_{Y^+} -\beta E^+)^2 \cdot (-K_{Y^+})\\
&=&(-\alpha K_{Y} -\beta E)^2 \cdot (-K_{Y})\\
&=& (-K_Y)^3 \alpha^2 -2(-K_Y)^2 \cdot E \alpha \beta+ (-K_Y) \cdot E^2\beta^2\\
&=& (2g-8)\alpha^2 -8 \alpha \beta-2\beta^2,
\end{eqnarray*}
where the last equality follows from Proposition \ref{p-conic-basic}(1).

(2) 
Opposed to (\ref{e1-Y-vs-Y^+}), $E^3$ and $(E^+)^3$ might not coincide. 
Actually, the proof of (2) is carried out by showing $E^3 \neq (E^+)^3$. 
Note that if $-K_Y$ is ample, then $\psi$ and $\psi^+$ are isomorphisms, and hence $E=E^+$, 
which implies $E^3 = (E^+)^3$. 
Therefore, if $E^3 \neq (E^+)^3$, then $-K_Y$ is not ample. 
We now overview how to check $E^3 \neq (E^+)^3$. 
For simplicity, assume that $\Gamma$ is a point. 
By $E^3=1$ (Lemma \ref{l-blowup-formula2}), it is enough to show $(E^+)^3 \neq 1$.  
In our situation, we have 
\[
D \sim -\alpha K_Y -\beta E^+
\]
for some $\alpha, \beta \in \Z_{>0}$. 
By (\ref{e1-Y-vs-Y^+}), we get 
\[
D^3 = \alpha^3 (-K_Y)^3 - 3 \alpha^2 \beta (-K_Y)^2 \cdot E + 3 \alpha \beta^2 (-K_Y) \cdot E^2 +\beta^3 (E^+)^3. 
\]
Note that  we have $(-K_Y)^3 =2g-10, (-K_Y)^2 \cdot E=4, (-K_Y) \cdot E^2=-2$ (Lemma \ref{l-blowup-formula2}).  
Since the triple $(g, \alpha, \beta)$ is explicitly given (Table \ref{table-2ray-pt} in Theorem \ref{t-pt-flop}), 
it is enough to compute $D^3$, in order to compute $(E^+)^3$. 
For all the cases except for the $E_1$ case, 
the value $D^3$ has been computed in Subsection \ref{ss-weak-Fano}. 
The $E_1$ case will be treated in Subsection \ref{ss-flopping}.

\medskip

In Subsection \ref{ss-div-cont} (resp. Subsection \ref{ss-flop}), we shall summarise some results which are used to establish (1) and (3) for the case (i) (resp. (ii)). 
For the case (i), 
the strategy is quite similar to the case (ii) explained above. 
For example, the starting point for the case (i) is to observe $D \sim -\alpha K_Y -\beta E$ for $D:=\Ex(\psi)$.

\subsection{Divisorial contractions}\label{ss-div-cont}

\begin{nota}\label{n-div-cont}
Let $X \subset \P^{g+1}$ be an anti-canonically embedded Fano threefold with 
$\Pic\,X = \Z K_X$. 
Let $\Gamma$ be either a point or a smooth rational curve on $X$. 
Let $\sigma : Y \to X$ be the blowup along $\Gamma$. 
Assume that 
\begin{enumerate}
\item[$(*)$] $-K_Y$ is semi-ample and big. 
\end{enumerate}
Let $\psi : Y \to Z$ be the birational morphism to a projective normal threefold $Z$ 
that is the Stein factorisation of $\varphi_{|-mK_Y|}$ 
for some (every) positive integer $m$ such that $|-mK_Y|$ is base point free. 
Assume that $\dim \Ex(\psi)=2$. Set $E := \Ex(\sigma)$, $D := \Ex(\psi)$, and $B := \psi(D) = \psi(\Ex (\psi))$. 
\begin{equation}
\begin{tikzcd}
& Y \arrow[ld, "\sigma"'] \arrow[rd, "\psi"]\\
X & & Z
\end{tikzcd}
\end{equation}
\end{nota}

\begin{prop}\label{p-div-cont1}
We use Notation \ref{n-div-cont}. 
Let $\zeta$ be a curve on $Y$ such that $-K_Y \cdot \zeta =0$. 
Then the following hold. 
\begin{enumerate}
\item $D$ is a prime divisor. 
\item $D \cdot \zeta <0$ and $-D$ is $\psi$-ample. 
\item $E \cdot \zeta >0$, $E$ is $\psi$-ample, and $D \neq E$. 
\item We can uniquely write 
\[
D \sim -\alpha K_Y -\beta E
\]
for some $\alpha, \beta \in \Z_{>0}$. 
\item $B$ is a curve. 
Furthermore, $(-K_Y)^2 \cdot D = 0$, $(-K_Y) \cdot D^2 <0$, and 
$\alpha (-K_Y)^3 = \beta (-K_Y)^2 \cdot E$. 
\item 
$(\psi|_D)_* \MO_D = \MO_B$, 
where {\cred $\psi|_D : D  \to B$} denotes the induced morphism. 
\item 
If $\xi$ is a curve on $D$ contained in a general fibre of $\psi|_D : D \to B$, 
then $K_Y \cdot \xi =0$ and $-D \cdot \xi \in \{1, 2\}$.   
\item $\beta \in \{1, 2\}$. 
\end{enumerate}
\end{prop}

\begin{proof}
Let us show (1) and (2). 
Since $-K_Y$ is big, we can write $-nK_Y = A + F$ for some $n \in \Z_{>0}$, 
where $A$ is an ample divisor and $F$ is an effective divisor. 
It holds that
\[
0 = (-nK_Y) \cdot \zeta = (A+F) \cdot \zeta > F \cdot \zeta. 
\]
Therefore, we can find a prime divisor $D'$ on $Y$, contained in $\Supp\,F$, 
such that $D' \cdot \zeta <0$. 
In particular, we obtain $D = \Ex(\psi) \subset D'$. 
Since $D'$ is a prime divisor, 
it follows from $\dim D=2$ that $D=D'$. 
It follows from $\rho(Y/Z)=1$ and $D \cdot \zeta <0$ that $-D$ is $\psi$-ample. 
Thus (1) and (2) hold.

Let us show (3). 
We have $K_Y = \sigma^*K_X + aE$ with $a>0$. 
By $K_Y \cdot \zeta =0$ and $\sigma^*K_X \cdot \zeta <0$, 
we obtain $E \cdot \zeta >0$. 
In particular, $D \neq E$. 
It follows from $\rho(Y/Z)=1$ and $E \cdot \zeta >0$ that $E$ is $\psi$-ample. 
Thus (3) holds. 

Let us show (4). 
By $\Pic\,Y = \Z K_Y \oplus \Z E$, 
there uniquely exist $\alpha, \beta \in \Z$ such that 
\[
D \sim -\alpha K_Y - \beta E. 
\]
It suffices to show $\alpha >0$ and $\beta >0$. 
If $\zeta$ is a general $\psi$-contracted curve, 
then we obtain 
\[
0 \overset{{\rm (2)}}{>} D \cdot \zeta = 
(-\alpha K_Y - \beta E) \cdot \zeta = -\beta E \cdot \zeta. 
\]
By (3), we get $\beta >0$. 
Since $D$ is effective and $-K_Y$ is big, we get $\alpha >0$. 
Thus (4) holds.


Let us show (5). 
We first prove that $B$ is a curve. 
Otherwise, the closed subset $\psi(D) = \psi(\Ex\,\psi)=B$ is a point. 
Then $(-K_Y)^2 \cdot D = (-K_Y) \cdot D^2 =0$. 
We obtain 
\[
0 = (-K_Y)^2 \cdot D = (-\alpha K_Y -\beta E) \cdot (-K_Y)^2 
= \alpha (-K_Y)^3 - \beta (-K_Y)^2 \cdot E
\]
and the following contradiction: 
\begin{eqnarray*}
0 &=& (-K_Y) \cdot D^2\\
&=& (-K_Y) \cdot (-\alpha K_Y -\beta E)^2\\
&=& \alpha^2 (-K_Y)^3 -2\alpha\beta (-K_Y)^2 \cdot E + \beta^2 (-K_Y) \cdot E^2\\
&=& \alpha \cdot \beta (-K_Y)^2 \cdot E -2\alpha\beta (-K_Y)^2 \cdot E + \beta^2 (-K_Y) \cdot E^2\\
&=&  -\alpha \beta (-K_Y)^2 \cdot E + \beta^2 (-K_Y) \cdot E^2\\
&<&0, 
\end{eqnarray*}
where the last inequality follows from $\alpha >0, \beta>0$, 
$(-K_Y)^2 \cdot E >0$, and $(-K_Y) \cdot E^2 <0$ (Lemma \ref{l-blowup-formula}). 
Hence $B$ is a curve. 

Since $-\ell K_Y|_D$ is linearly equivalent to a sum of fibres of $\psi|_D : D \to B$ for some $\ell \in \Z_{>0}$, 
we have $(- \ell K_Y)^2 \cdot D = (-\ell K_Y|_D)^2 = 0$. 
Since $-D$ is $\psi$-ample, we obtain $(-\ell K_Y) \cdot D^2 = (-\ell K_Y|_D) \cdot (D|_D)  <0$. 
By $D \sim -\alpha K_Y -\beta E$ and $(-K_Y)^2 \cdot D=0$, 
we get $\alpha (-K_Y)^3 = (D+\beta E) \cdot (-K_Y)^2 = \beta (-K_Y)^2 \cdot E$. 
Thus (5) holds.

Let us show (6). 
By the exact sequence 
\[
0 \to \MO_Y(-D) \to \MO_Y \to \MO_D \to 0, 
\]
it suffices to show $R^1\psi_*\MO_Y(-D)=0$. 
By (4), we have 
\[
-D -K_Y \sim (\alpha -1) K_Y+ \beta E.
\]
Since $K_Y$ is $\psi$-nef and $E$ is $\psi$-ample, 
$(\alpha -1) K_Y+ \beta E$ is $\psi$-ample, which implies 
$R^1\psi_*\MO_Y(-D)=0$ \cite[Theorem 0.5]{Tan15}. 
Thus (6) holds.

Let us show (7). 
Let $\nu : \widetilde D \to D$ be the normalisation of $D$. 
We have $\MO_{\widetilde D}(K_{\widetilde D} + C)  \simeq \nu^*\omega_D$ for the conductor $C$, 
where $C$ is an  effective $\Z$-divisor on $\widetilde D$ \cite[Proposition 2.3]{Rei94}. 
Let $\widetilde D \to \widetilde B$ be the Stein factorisation of 
the composite morphism $\widetilde D \xrightarrow{\nu} D \xrightarrow{\psi|_D} B$: 
\[
\begin{tikzcd}
\widetilde D \arrow[r, "\nu"]  \arrow[d] & D\arrow[d, "\psi|_D"]\\
\widetilde B \arrow[r] & B. 
\end{tikzcd}
\]
Pick a general fibre $\widetilde{\xi} \subset \widetilde{D}$ of $\widetilde D \to \widetilde B$ 
and we set $\xi := \nu(\widetilde{\xi})$. 
Note that 
$\widetilde{\xi}$ is contained in the smooth locus of $\widetilde D$ and 
$\widetilde{\xi} \to \xi$ is birational. 
We then obtain 
\[
-2 \leq -2 + C \cdot \widetilde{\xi} \leq  (K_{\widetilde D} + C) \cdot \widetilde{\xi} 
= \nu^*\omega_D \cdot \widetilde{\xi} = \omega_D \cdot \xi = (K_Y+D) \cdot \xi= D \cdot \xi <0. 
\]
In particular, we obtain 
$-D \cdot \xi \in \{ 1, 2\}$. Thus (7) holds. 
The assertion (8) immediately follows from (4) and (7). 
\end{proof}

Given $K_Y = \sigma^*K_X + aE$ as in Notation \ref{n-div-cont}, 
our goal is to find a curve $\zeta_X$ on $X$ satisfying $-K_X \cdot \zeta_X = a$ (e.g., if $\sigma :Y \to X$ is a blowup at a point, 
then $a=2$ and $\zeta_X$ is a conic on $X$). 
To this end, we consider the $\psi$-contracted curves $\zeta$. 
By $K_Y \cdot \zeta =0$, it is enough to find a curve $\zeta$ on $Y$ 
with $E \cdot \zeta =1$. 
The next proposition treats the other case.

\begin{prop}\label{p-div-cont2}
We use Notation \ref{n-div-cont}. 
We define $a \in \Z_{>0}$ by 
\[
K_Y = \sigma^*K_X  +a E. 
\]
Let $\nu : \widetilde D \to D$ be the normalisation of $D$ and let $\widetilde B$ 
be the Stein factorisation of $\nu \circ (\psi|_D)$: 
\[
\begin{CD}
\widetilde D @>\nu >> D\\
@VV\varphi V @VV\psi|_D V\\
\widetilde B @>\nu_B>> B. 
\end{CD}
\]
Let $C$ be the conductor of $\nu$, 
which is the effective $\Z$-divisor on $\widetilde D$ satisfying 
$\Ex(\nu) =\Supp\,C$ and $\MO_{\widetilde D}(K_{\widetilde D} + C) \simeq \nu^*\omega_D$.  
Assume that $E \cdot \zeta \neq 1$ 
for every curve $\zeta$ on $Y$ such that $K_Y \cdot \zeta =0$. 
Then the following hold. 
\begin{enumerate}
\item $\beta =1$. 
Furthermore, $-D \cdot \xi =2$ and $E \cdot \xi =2$ 
for every general fibre $\xi$ of $\psi|_D : D \to B$.  
\item 
$C$ is vertical over $\widetilde B$, i.e., $\varphi(\Supp\,C) \neq \widetilde B$. 
\item 
The induced morphism $\nu_B : \widetilde B \to B$ 
is the normalisation of $B$. 
\item $\varphi : \widetilde D \to \widetilde B$ is a $\P^1$-bundle, i.e., every scheme-theoretic fibre is isomorphic to $\P^1$. In particular, $\widetilde D$ is a smooth projective surface. 
\item $\omega_D^2 \in 4\Z$. 
\end{enumerate}
\end{prop}

\begin{proof}
Let us show (1)' below. 
Note that (1) follows from (1)', (2), and (3). 
\begin{enumerate}
\item[(1)'] $\beta =1$. 
Furthermore, $-D \cdot \xi =2$ and $E \cdot \xi =2$ 
if $\xi$ is a curve on $D$ contained in a general fibre of $\psi|_D : D \to B$. 
\end{enumerate}
By $D \sim -\alpha K_Y -\beta E$, $K_Y \cdot \xi =0$, and $-D \cdot \xi \in \{1, 2\}$ (Proposition \ref{p-div-cont1}(4)(7)), we get 
\[
\{1, 2\} \ni -D \cdot \xi = \alpha K_Y \cdot \xi +\beta E \cdot \xi  = \beta E \cdot \xi. 
\]
It follows from $E \cdot \xi \neq 1$ that 
$E \cdot \xi =2$, $\beta =1$, and $-D \cdot \xi =2$. 
Thus (1)' holds.


Let us show (2). 
For a general fibre $\widetilde{\xi} \subset \widetilde D$ of $\varphi : \widetilde D \to \widetilde B$ 
and its image $\xi = \nu(\widetilde{\xi})$,  
we obtain 
\[
-2 \leq -2 + C \cdot \widetilde{\xi} \leq (K_{\widetilde D} + C) \cdot \widetilde{\xi} 
= \nu^*\omega_D \cdot \widetilde{\xi} = \omega_D \cdot \xi = (K_Y+D) \cdot \xi= D \cdot \xi 
\overset{(1)'}{=}-2. 
\]
Therefore, $C \cdot \widetilde{\xi} =0$, i.e., 
every irreducible component of $C$ is contained in some fibre of $\varphi : \widetilde D \to \widetilde B$. 
Thus (2) holds. 
The assertion (3) follows from (2) and $(\psi|_D)_*\MO_D = \MO_B$  (Proposition \ref{p-div-cont1}(6)). 
Furthermore, we have $K_{\widetilde D} \cdot \widetilde{\xi}=-2$, and hence $\widetilde \xi \simeq \P^1$. 


Let us show (4). 
Pick a scheme-theoretic fibre $\widetilde{\zeta}$ of $\varphi : \widetilde D \to \widetilde B$, 
which is an effective Cartier divisor on $\widetilde D$. 
If $\widetilde{\zeta}$ is a prime divisor, then $\widetilde{\zeta} \simeq \P^1$, because a general fibre $\widetilde \xi$ is $\P^1$. 
By $\nu^*\omega_D^{-1} \cdot \widetilde{\zeta} =2$ and the $\varphi$-ampleness of  $\nu^*\omega^{-1}_D$, 
we may assume that 
$\widetilde{\zeta} =\widetilde{\zeta}_1 + \widetilde{\zeta}_2$, 
where 
each $\widetilde{\zeta}_i$ is a prime divisor 
(possibly $\widetilde{\zeta}_1 = \widetilde{\zeta}_2$). 
Although $\widetilde{\zeta}$ is a Cartier divisor, each $\widetilde{\zeta}_i$ might not be a Cartier divisor. 
We have $-\nu^*\omega_D \cdot \widetilde{\zeta}_i = 1$ for each $i \in \{1, 2\}$. 
Set $\zeta_i := \nu(\widetilde{\zeta}_i)_{\red}$, which is a 
curve on $D$ such that $\psi(\zeta_i)$ is a point. 
We have that 
\[
1 = -\nu^*(\omega_D) \cdot \widetilde{\zeta}_i = -\omega_D \cdot \zeta_i \times \deg(\widetilde{\zeta}_i \to \zeta_i). 
\]
In particular, 
$-D \cdot \zeta_i = {\cred -(K_Y+D) \cdot \zeta_i=} -\omega_D \cdot \zeta_i =1$. 
We then get 
\[
1= -D \cdot \zeta_i = \alpha K_Y \cdot \zeta_i  +\beta E \cdot \zeta_i  = \beta E \cdot \zeta_i, 
\]
which implies $E \cdot \zeta_i =1$. This contradicts our assumption. 
Thus (4) holds. 

Let us show (5). Fix a fibre $\widetilde \zeta$ of $\varphi :\widetilde{D} \to \widetilde{B}$. 
By (2) and (4), we have $C \equiv n \widetilde \xi$ for some $n \in \Z_{\geq 0}$. 
We then get 
\[
\omega_D^2 =  (K_{\widetilde D}+C)^2 
=(K_{\widetilde D}+n \widetilde \zeta)^2 
= K_{\widetilde D}^2 + 2n K_{\widetilde D} \cdot \widetilde \zeta 
= 8(1-g(\widetilde B)) -4n \in 4\Z. 
\]
Thus (5) holds. 
\qedhere


\end{proof}

\subsection{Computations for flops}\label{ss-flop}

\begin{nota}\label{n-flop}
Let $X \subset \P^{g+1}$ be an anti-canonically embedded Fano threefold with 
$\Pic\,X = \Z K_X$. 
Let $\Gamma$ be either a point or a smooth 
curve on $X$. 
Let $\sigma : Y \to X$ be the blowup along $\Gamma$. 
Assume that 
\begin{enumerate}
\item[$(*)$] $-K_Y$ is semi-ample and big. 
\end{enumerate}
Let $\psi : Y \to Z$ be the birational morphism to a projective normal threefold $Z$ 
that is the Stein factorisation of $\varphi_{|-mK_Y|}$ 
for some (every) positive integer $m$ such that $|-mK_Y|$ is base point free. 
Assume that either 
\begin{enumerate}
    \item[(I)] $\dim \Ex(\psi)=1$, or 
    \item[(II)] $\psi$ is an isomorphism, i.e., $-K_Y$ is ample. 
\end{enumerate}

(I) Assume $\dim \Ex(\psi)=1$. 
Note that $Z$ is Gorenstein, the flop $\psi^+ : Y^+ \to Z$ of $\psi$ exists, and $Y^+$ is smooth \cite[Theorem 1.1]{Tan-flop}. 
In particular, we have $K_Y = \psi^*K_Z$ and $K_{Y^+} =(\psi^+)^*K_Z$. 
Since $-K_{Y^+}$ is not ample but nef and big, 
there exists a unique $K_{Y^+}$-negative extremal ray $R$ of $\NE(Y^+)$. 
We have the contraction $\tau : Y^+ \to W$ of $R$. 
\begin{equation}
\begin{tikzcd}
Y \arrow[d, "\sigma"'] \arrow[rd, "\psi"]& & Y^+ \arrow[ld, "\psi^+"'] \arrow[d, "\tau"]\\
X & Z & W
\end{tikzcd}
\end{equation}
Set $E^+$ to be the proper transform of $E$ on $Y^+$. 

(II) Assume that $\psi$ is an isomorphism. 
In this case, we set $Z:=Y$, $Y^+ :=Y$, $\psi := {\rm id}$, $(\psi)^+ := {\rm id}$, 
and $E^+ := E$ for consistency. 
Let $R$ be the $K_Y$-negative extremal ray not corresponding to $\sigma$. 
Let $\tau: Y \to W$ be the contraction of $R$. 
\end{nota}



\begin{lem}\label{l-blowup-formula2}
We use Notation \ref{n-flop}. 
\begin{enumerate}
\item 
If the blowup centre $\Gamma$ is a point, then the following hold. 
\begin{enumerate}
\item $(-K_Y)^3 = (-K_{Y^+})^3 = 2g -10$. 
\item $(-K_Y)^2 \cdot E = (-K_{Y^+})^2 \cdot E^+ =4$.
\item $(-K_Y) \cdot E^2 = (-K_{Y^+}) \cdot (E^+)^2=-2$.
\item $E^3 = 1$. 
\end{enumerate}
\item 
If the blowup centre $\Gamma$ is a smooth curve of genus $g(\Gamma)$, then the following hold for $\deg \Gamma := -K_X \cdot \Gamma$. 
\begin{enumerate}
\item $(-K_Y)^3 = (-K_{Y^+})^3 =  2g-2\deg \Gamma +2g(\Gamma)-4$. 
\item $(-K_Y)^2 \cdot E = (-K_{Y^+})^2 \cdot E^+ =\deg \Gamma -2g(\Gamma)+2$.
\item $(-K_Y) \cdot E^2 = (-K_{Y^+}) \cdot (E^+)^2=2g(\Gamma)-2$.
\item $E^3 = -\deg \Gamma + 2-2g(\Gamma)$. 
\end{enumerate}
\end{enumerate}

\end{lem}


\begin{proof}
Since all the proofs are very similar, we only show (c) of (2).  
By Lemma \ref{l-blowup-formula}, we have $(-K_Y) \cdot E^2 =2g(\Gamma)-2$. 
Hence it is enough to prove $(-K_Y) \cdot E^2 = (-K_{Y^+}) \cdot (E^+)^2$. 
Pick $m \in \Z_{>0}$ such that $|-mK_Z|$ is very ample. 
Then we can find a general member $S \in |-mK_Z|$ which avoids $\psi(\Ex(\psi))$ and $\psi^+(\Ex(\psi^+))$, 
because these are zero-dimensional. 
For its pullbacks $T := \psi^*S \in |-mK_Y|$ and $T^+:=(\psi^+)^*S \in |-mK_{Y^+}|$, 
we have $T \simeq S \simeq T^+$ and hence $(E|_T)^2 = (E_Z|_S)^2 = (E^+|_{T^+})^2$ 
for $E_Z := \psi_*(E)$, 
which implies $(-K_Y) \cdot E^2 = (-K_{Y^+}) \cdot (E^+)^2$. 
\end{proof}

\begin{lem}\label{l-beta=r}
We use Notation \ref{n-flop}. 
\begin{enumerate}
\item 
$-E$ is $\sigma$-ample and $E$ is $\psi$-ample. 
$-E^+$ is $\psi^+$-ample and $E^+$ is $\tau$-ample. 
\item 
Assume that $\tau$ is of type $C$ or $D$. 
Let $\MO_W(1)$ be  the invertible sheaf  on $W$ that is an ample generator of $\Pic\,W$ 
(note that $W =\P^2$ (resp. $W=\P^1$) if $\tau$ is of type $C$ (resp. type $D$)). 
Let $D$ be a Cartier divisor on $Y^+$ satisfying $\MO_{Y^+}(D) \simeq \tau^*\MO_W(1)$.  
Then 
\begin{equation}\label{e1-beta=r}
D \sim -\alpha K_{Y^+} -\beta E^+\qquad 
\text{for \quad some}\quad \alpha, \beta \in \Z_{>0}.
\end{equation}
Moreover, $\alpha$ and $\beta$ are corpime, i.e., $\alpha \Z  +\beta \Z  = \Z$. 
\begin{enumerate}
\item[(a)] If $\tau$ is of type $C$, then $\beta \in \{1, 2\}$. 
\item[(b)] If $\tau$ is of type $D$, then $\beta \in \{1, 2, 3\}$. 
\end{enumerate}
Furthermore, $\beta$ is a divisor of the length $\mu$ of the extremal ray $R$. 
\item 
Assume that $\tau$ is of type $E$. 
Let $D$ be the $\tau$-exceptional prime divisor. 
Then 
\begin{equation}\label{e2-beta=r}
D \sim -\alpha K_{Y^+} -\beta E^+\qquad 
\text{for \quad some}\quad \alpha, \beta \in \Z_{>0}.
\end{equation}
\begin{enumerate}
\item[(c)] 
If $\tau$ is of type $E_1$, then 
$\alpha +1 \in \beta \Z$ and $\beta = r_W$ for the index $r_W$ of the Fano threefold $W$. 
\item[(d)] 
If $\tau$ is of type $E_2$, then 
$2\alpha +1 \in \beta  \Z$ and $\beta = r_W$ for the index $r_W$ of the Fano threefold $W$. 
\end{enumerate}
\end{enumerate}
\end{lem}


\begin{proof}
Let us show (1). 
We only treat the case when $\dim \Ex(\psi)=1$ (i.e., (I) in Notation \ref{n-flop}), as the other case 
(i.e., (II) in Notation \ref{n-flop}) is simpler. 
Since $-E$ is $\sigma$-ample, $-E$ is not $\psi$-nef, as otherwise $-E$ would be nef by Kleiman's ampleness criterion, which is absurd. 
By $\rho(Y/Z)=1$, $E$ is $\psi$-ample. 
Then $-E^+$ is $\psi^+$-ample, 
because the flop $Y \xrightarrow{\psi} Z \xleftarrow{\psi^+} Y^+$ 
switches $\psi$-ampleness and $\psi^+$-anti-ampleness. 
Again by Kleiman's ampleness criterion, $E^+$ is not $\tau$-nef, i.e., $E^+$ is $\tau$-ample.  
Thus (1) holds.

We now  prove (\ref{e1-beta=r}) and (\ref{e2-beta=r}) simultaneously. 
For both cases, $-D$ is  $\tau$-nef and $D$ is effective. 
It follows from (1) that $D \neq E^+$. 
 By $\Pic\,Y^+ = \Z K_{Y^+} \oplus \Z E^+$, we can write 
\[
D \sim -\alpha K_{Y^+} - \beta E^+
\]
for some $\alpha, \beta \in \Z$. 
If $\alpha <0$, then we would get 
\[
0 \geq \kappa(Y, -\beta E) =\kappa(Y^+, - \beta E^+) = \kappa(Y^+, \alpha K_{Y^+} +D) \geq \kappa(Y^+, -\alpha (-K_{Y^+}))  \geq 3.
\]
Hence $\alpha \geq 0$. 
Suppose $\alpha =0$. Then $D \sim -\beta E^+$. 
Since $D$ and $E^+$ are nonzero effective divisors, we  get $-\beta >0$. 
For a general curve $\zeta \subset D$ which is contracted by $\tau$, 
we get 
$0 > D \cdot \zeta = -\beta E^+ \cdot \zeta \geq 0$ by 
$D \neq E^+$, which is absurd. 
Thus we get $\alpha >0$. 
If $\beta \leq 0$, then we would obtain 
\[
2 \geq \kappa(Y^+, D) = \kappa(Y^+, -\alpha K_{Y^+} - \beta E^+) \geq  \kappa(Y^+, -\alpha K_{Y^+}) =3,
\]
which is a contradiction. 
Therefore, $\beta >0$. 
This completes the proof of (\ref{e1-beta=r}) and (\ref{e2-beta=r})

Let us show (2).  
Since $D$ is primitive in $\Pic\,Y$ (Proposition \ref{p-cont-ex}), 
$\alpha$ and $\beta$ are coprime. 
Let $\ell$ be an extremal rational curve of $R$, i.e., 
$\ell$ is a rational curve on $Y^+$ such that $D \cdot \ell =0$ and $-K_{Y^+} \cdot \ell = \mu$. 
We then obtain 
\[
\alpha \mu = \alpha (-K_{Y^+}) \cdot \ell = \beta E^+ \cdot \ell. 
\]
As $\alpha$ and $\beta$ are coprime, $\beta$ is a divisor of $\mu$. 
If $\tau$ is of type $C$ (resp. type $D$), 
then we have $\mu \in \{1, 2\}$ (resp. $\mu \in \{1, 2, 3\}$) 
by Proposition \ref{p-typeC-intersec} (resp. Proposition \ref{p-typeD-intersec}). 
Thus (2) holds.

%


Let us show (3). 
Since the proofs of (c) and (d) are the same, 
we only show (c). 
Assume that $\tau$ is of type $E_1$. 
Fix a curve $\zeta \subset D$ such that $\tau(\zeta)$ is a point. 
 By $D \sim -\alpha K_{Y^+} -\beta E^+$, we obtain 
$-K_{Y^+} + D \sim (\alpha+1) (-K_{Y^+}) -\beta E^+$, which implies 
\[
0 = (-K_{Y^+} + D) \cdot \zeta = (\alpha+1)  -\beta E^+ \cdot \zeta = \alpha + 1 -\beta\gamma \qquad \text{for}\qquad \gamma := E^+ \cdot \zeta. 
\]
Hence $\alpha +1  =\beta \gamma \in \beta \Z$. 

It suffices to show $\beta  =r_W$. 
We can write $-K_W \sim r_WH_W$ for some ample Cartier divisor $H_W$ on $W$. 
Then 
\[
r_W \tau^*H_W \sim -\tau^*K_W \sim  (-K_{Y^+} + D) \sim \beta \gamma (-K_{Y^+}) -\beta E^+. 
\]
By $\Pic\,Y^+ = \Z K_{Y^+} \oplus \Z E^+$, 
we can uniquely write $\tau^* H_W \sim \delta K_{Y^+} - \beta' E^+$ for some $\delta, \beta' \in \Z$. 
Then we get $r_W \beta'  = \beta$ and 
\[
r_W \tau^*H_W \sim  r_W\beta'( \gamma (-K_{Y^+}) - E^+). 
\]
Taking the pushforward $\tau_*$, we obtain  
\[
r_W H_W \sim r_W \beta' \tau_*( \gamma (-K_{Y^+}) - E^+). 
\]
Since $H_W$ is a primitive element of $\Pic\,W \simeq \Z$ and $\beta' >0$, 
we get $\beta' =1$ by $\beta'>0$. 
Therefore, $\beta = r_W$. 
Thus (c) holds. 
\qedhere

\end{proof}

\begin{lem}\label{l-E1-deg-genus}
We use Notation \ref{n-flop}. 
Assume that $\tau$ is of type $E_1$. 
Set $D := \Ex(\tau)$ and $B := \tau(D)$. 
Recall that we can write $D \sim - \alpha K_{Y^+} - \beta E^+$ 
for some $\alpha, \beta\in \Z_{>0}$ (Lemma \ref{l-beta=r}). 
Then the following hold. 
\begin{enumerate}
\item $(-K_W)^3 = (\alpha+1)^2 (-K_{Y})^3 - 2(\alpha+1)\beta (-K_Y)^2 \cdot E  + \beta^2(-K_Y) \cdot E^2$. 
\item $-K_W \cdot B = \alpha(\alpha+1)(-K_Y)^3 
- (2\alpha +1)\beta (-K_Y)^2 \cdot E 
+\beta^2 (-K_Y) \cdot E^2.$
\item $2g(B) -2 =\alpha^2 (-K_Y)^3 - 2 \alpha \beta(-K_Y)^2 \cdot E + \beta^2(-K_Y) \cdot E^2$. 
\end{enumerate}
\end{lem}

\begin{proof}
Recall that we have $K_{Y^+} = \tau^*K_W +D$. 
Then the assertion (1) holds by the following: 
\begin{align*}
(-K_W)^3 
&= (-\tau^*K_W)^3 \\
&= (-\tau^*K_W)^2 \cdot (-K_{Y^+} +D) \\
&= (-\tau^*K_W)^2 \cdot (-K_{Y^+})\\
&= (-K_{Y^+} +D)^2 \cdot (-K_{Y^+}) \\
&=  (-(\alpha+1)K_{Y^+} -\beta E^+)^2 \cdot (-K_{Y^+})\\
&= (\alpha+1)^2 (-K_{Y})^3 - 2(\alpha+1)\beta (-K_Y)^2 \cdot E  + \beta^2(-K_Y) \cdot E^2,
\end{align*}
where the last equality follows from Lemma \ref{l-blowup-formula2}. 
The assertion (2) holds by 
\begin{align*}
(-K_W) \cdot B &= (\tau^*(-K_W)|_D) \cdot (-K_{Y^+}|_D)\\
&= \tau^*(-K_W) \cdot (-K_{Y^+}) \cdot D \\
&=(-K_{Y^+}) \cdot D  \cdot (-K_{Y^+}+D)\\
&= (-K_{Y^+}) \cdot (\alpha(-K_{Y^+}) - \beta E^+)  \cdot ( (\alpha +1)(-K_{Y^+})-\beta E^+)\\
&= \alpha(\alpha+1)(-K_Y)^3 
- (2\alpha +1)\beta (-K_Y)^2 \cdot E 
+\beta^2 (-K_Y) \cdot E^2, 
\end{align*}
where the last equality follows from Lemma \ref{l-blowup-formula2}. 
The following holds: 
\begin{align*}
D^2 \cdot (-K_{Y^+}) 
&= ( \alpha (-K_{Y^+}) - \beta E^+)^2  \cdot (-K_{Y^+}) \\
&= \alpha^2 (-K_Y)^3 - 2 \alpha \beta (-K_Y)^2 \cdot E + \beta^2(-K_Y) \cdot E^2, 
\end{align*}
where the last equality follows from Lemma \ref{l-blowup-formula2}. 
Thus (3) holds 
{\cred by 
 $2g(B) -2  = D^2 \cdot (-K_{Y^+})$ (Lemma \ref{l-blowup-formula}(2))}.  
\end{proof}

\begin{rem}\label{r-E1-deg-genus}
We use the same notation as in the statement of Lemma \ref{l-E1-deg-genus}. 
By  Lemma \ref{l-beta=r}(3), we can write $\alpha +1  = \beta \gamma$ for some $\gamma \in \Z_{>0}$. 
By $\beta = r_W$, we have $-K_W \sim \beta H_W$ for an ample generator $H_W$ of $\Pic\,W$. 
It follows from  Lemma \ref{l-E1-deg-genus}(1) that 
\[
w:=\frac{1}{2} \beta H_W^3 = \frac{1}{2}(-K_Y)^3 \gamma^2 - (-K_Y)^2 \cdot E \gamma + \frac{1}{2} (-K_Y) \cdot E^2. 
\]
Here we have $w \in \Z_{>0}$ as follows.
\begin{enumerate}
\item Assume $\beta=1$. 
Then $H_W^3 = (-K_W)^3 = 2g_W-2$ and $w = g_W-1$, where 
$g_W$ is the genus of $W$. In particular, $w \in \Z_{>0}$ and $(-K_W)^3 = 2w$. 
Moreover, if $|-K_W|$ is not very ample, then $1 \leq w \leq 2$ \cite[Theorem 1.1]{TanI}. 
\item If $\beta=2$, then we get 
$1 \leq H_W^3 \leq 5$ (Remark \ref{r-index2-pic1}), which implies  $w \in \Z$ and $1 \leq w \leq 5$. 
\item If $\beta =3$, then $(\beta, H_W^3, w)= (3, 2, 3)$. 
\item If $\beta =4$, then $(\beta, H_W^3, w)= (4, 1, 2)$. 
\end{enumerate}
\end{rem}

\begin{rem}\label{r-index2-pic1}
Let $X$ be a Fano threefold with $\rho(X)=1$ and $r_X=2$, 
where $r_X$ denotes the index. 
Let $H$ be a Cartier divisor on $X$ such that $-K_X \sim 2H$. 
Then $1 \leq H \leq 5$. 

Indeed, by \cite[Theorem 2.23]{TanI}, 
$H^3 \geq 6$ implies that $X$ is isomorphic to $\P^1 \times \P^1 \times \P^1$, the blowup of $\P^3$ at a point, or a member $|\MO_{\P^2 \times \P^2}(1, 1)|$. 
Clearly, the first two cases satisfy $\rho(X) \geq 2$. 
As for the last case, i.e., $X \in |\MO_{\P^2 \times \P^2}(1, 1)|$, 
we have a surjection $X \to \P^2$ induced by the first projection $\P^2 \times \P^2$, and hence we get  $\rho(X) \geq 2$. 
\end{rem}

\begin{lem}\label{l-bound--K-E}
Let $X \subset \P^{g+1}$ be an anti-canonically embedded Fano threefold with $\Pic\,X = \Z K_X$. 
Let $\Gamma$ be either a point or a smooth rational curve on $X$. 
Let $\sigma : Y \to X$ be the blowup along $\Gamma$ and set $E := \Ex(\sigma)$. 
Assume that $-K_Y$ is nef and big. 
Then the following hold: 
\begin{eqnarray*}
h^0(Y, -K_Y) &\geq& \chi(Y, -K_Y)= \frac{(-K_Y)^3}{2} +3\\
h^0(E, -K_Y|_E) &=& \chi(E, -K_Y|_E)
=  1 + (-K_Y)^2 \cdot E - \frac{1}{2} (-K_Y) \cdot E^2\\
h^0(Y, -K_Y -E) &\geq& 
h^0(Y, -K_Y) - h^0(E, -K_Y|_E)\\
&\geq& 
\frac{1}{2}(-K_Y)^3 - (-K_Y)^2 \cdot E +\frac{1}{2}(-K_Y) \cdot E^2 +2. 
\end{eqnarray*}
\end{lem}

\begin{proof}
It holds that 
\[
h^0(Y, -K_Y) \geq h^0(Y, -K_Y) - h^1(Y, -K_Y) \overset{{\rm (i)}}{=} \chi(Y, -K_Y) 
\overset{{\rm (ii)}}{=}  \frac{(-K_Y)^3}{2} +3, 
\]
where (i) follows from Proposition \ref{p-kawakami} and (ii) holds by the Riemann--Roch theorem (cf. \cite[(2.5.1), (2.5.2), the proof of Corollary 2.6]{TanI}). 
Since $-K_Y|_E$ is nef and $E$ is toric, the following holds by 
the Riemann--Roch theorem and a toric Kawamata-Viehweg vanishing theorem \cite[Corollary 1.7]{Fuj07} (note that any Cartier divisor on a toric variety is linearly equivalent to a torus-invariant Cartier divisor \cite[Theorem 4.2.1]{CLS11}): 
\[
h^0(E, -K_Y|_E) = \chi(E, -K_Y|_E)=  1 + \frac{1}{2} (-K_Y|_E) \cdot (-K_Y|_E -K_E).  
\]
The last term can be computed as follows: 
\[
(-K_Y|_E) \cdot (-K_Y|_E -K_E) = (-K_Y) \cdot (-2K_Y-E) \cdot E  =
2K_Y^2 \cdot E +K_Y \cdot E^2. 
\]
By the exact sequence 
\[
0 \to \MO_Y(-K_Y-E) \to \MO_Y(-K_Y) \to \MO_Y(-K_Y)|_E \to 0, 
\]
we obtain 
\begin{eqnarray*}
h^0(Y, -K_Y-E) 
&\geq& h^0(Y, -K_Y) - h^0(E, -K_Y|_E)\\
&\geq& \left(\frac{(-K_Y)^3}{2} +3\right) - 
\left( 1 + \frac{1}{2}(2K_Y^2 \cdot E +K_Y \cdot E^2)\right)\\
&=& \frac{1}{2}(-K_Y)^3 - (-K_Y)^2 \cdot E +\frac{1}{2}(-K_Y) \cdot E^2 +2. 
\end{eqnarray*}
\end{proof}

The following vanishing theorem is essentially due to Kawakami \cite{Kaw21}. 

\begin{prop}\label{p-kawakami}
Let $X$ be a Fano threefold. 
Let $\sigma: Y \to X$ be a birational morphism from a smooth projective threefold $Y$. 
\begin{enumerate}
    \item Let $D$ be a nef  Cartier divisor on $Y$ such that $\kappa(Y, D) \geq 0$ and 
$\nu(Y, D) \geq 2$. 
Then $H^1(Y, \MO_Y(-D))=0$. 
\item If $-K_Y$ is nef and big, then $H^2(Y, -K_Y) = H^3(Y, -K_Y)=0$. 
\end{enumerate}
\end{prop}

\begin{proof}
It is enough to show (1), because (1) implies (2) by Serre duality. 
By \cite[Lemma 2.5]{Kaw21}, 
it suffices to prove that $H^0(Y, \Omega_Y^1 \otimes \MO_Y(-p^eD))=0$ for every $e \in \Z_{>0}$. 
For $D_X := f_*D$, we have 
$\kappa(X, D_X) \geq \kappa(Y, D) \geq 0$ and 
the following inclusion \cite[Lemma 2.3]{Kaw21}: 
\[
H^0(Y, \Omega_Y^1 \otimes \MO_Y(-p^eD)) \subset H^0(X, \Omega_X^1 \otimes \MO_X(-p^eD_X)). 
\]
Finally, $H^0(X, \Omega_X^1 \otimes \MO_X(-p^eD_X)) =0$ by \cite[Theorem 3.5]{Kaw21}. 
Thus (1) holds. 
\end{proof}

\subsection{Existence of  flopping curves}\label{ss-flopping}

As explained in Subsection \ref{ss-2ray-intro}, 
it is important 
to compute $D^3$ in order to show the existence of flopping curves. 
The following result summarises the list of the values for $D^3$. 
Note that $-K_W \cdot \tau(D)$ appearing below can be computed by  Lemma \ref{l-E1-deg-genus}(2).

\begin{lem}\label{l-D^3}
We use Notation \ref{n-flop}. 
Then 
the following hold. 
\begin{enumerate}
\item 
If $\tau$ is of type $C$ or $D$, 
then $D^3 =0$ holds for $D :=\tau^*D_W$ and 
a Cartier divisor $D_W$ on $W$. 
\item 
If $\tau$ is of type $E_1$ and $D = \Ex(\tau)$, then 
$K_Y^3 = K_W^3 +3 (-K_W) \cdot \tau(D) +D^3.$
\item 
If $\tau$ is of type $E_2$ and $D = \Ex(\tau)$, then $D^3 =1$. 
\item 
If $\tau$ is of type $E_3$ and $D = \Ex(\tau)$, then $D^3 =2$. 
\item 
If $\tau$ is of type $E_4$ and $D = \Ex(\tau)$, then $D^3 =2$. 
\item 
If $\tau$ is of type $E_5$ and $D = \Ex(\tau)$, then $D^3 =4$. 
\end{enumerate}
\end{lem}

\begin{proof}
The assertion (1) is clear. 
The assertions (3)--(6) follow from Proposition \ref{p-typeE-intersec}. 

Let us show (2). 
We have 
\[
K_{Y^+} = \tau^*K_W + D, 
\]
which implies 
\[
K_Y^3 = (\tau^*K_W + D)^3 = K_W^3 + 3 \tau^*K_W \cdot D^2 +D^3 
=K_W^3 +3 (-K_W) \cdot \tau(D) +D^3. 
\]
Thus (2) holds. 
\end{proof}

\section{Two-ray game for lines I}\label{s-line-I}

Let $X \subset \P^{g+1}$ be an anti-canonically embedded Fano threefold with $\Pic\,X = \Z K_X$. 
The main results of this section and its continuation (Section \ref{s-line-II}) are as follows: 
\begin{enumerate}
\item Let $\Gamma$ be a line of $X$ and let $\sigma: Y \to X$ be the blowup along $\Gamma$. 
Then $-K_Y$ is not ample (Corollary \ref{c-line-not-ample}). 
\item If $g=5$ or $g \geq 7$, then 
$X$ is not covered by the lines on $X$ 
(Proposition \ref{p-line-Hilb2}). 
\end{enumerate}
The strategy of the proof of (1) was explained in Subsection \ref{ss-2ray-intro}, i.e., 
we apply the two-ray game starting with the blowup along the line $\Gamma$. 
The proof of (1) will be postponed to Section \ref{s-line-II}, 
because Diophantine equations, involving $g$, will be simpler after establishing the inequality $g \leq 12$ given in Section \ref{s-conic}. 
On the other hand, we prove (2) in this section, as (2) will be used in Section \ref{s-pt}.

We now overview how to show (2). 
Let $\Hilb_X^{\line}$ be the Hilbert scheme of $X$ parametrising the lines on $X$. 
Let $\Univ_X^{\line}$ be its universal family. 
Then it is enough to show that 
 the induced morphism $\lambda : \Univ_X^{\line} \to X$ is not surjective. 
It is not so hard to check that  
$\Hilb_X^{\line} = \emptyset, \dim \Hilb_X^{\line} =1$, or $\dim \Hilb_X^{\line} =2$ (Proposition \ref{p-line-Hilb1}). 
Assume that  $\dim \Hilb_X^{\line} =2$ and  $\lambda : \Univ_X^{\line} \to X$ is surjective. 
It suffices to derive a contradiction. 
For simplicity, we further assume that $\Hilb_X^{\line}$ is integral, i.e., a projective surface. 
Then $\Univ_X^{\line}$ is a projective threefold. 
Then $\lambda : \Univ_X^{\line} \to X$ is a generically finite morphism. 
If $\lambda : \Univ_X^{\line} \to X$ is not a finite morphism, then we can find a line $\Gamma$ on $X$ 
intersecting infinitely many lines on $X$ (Step \ref{s2-line-Hilb2} of Proposition \ref{p-line-Hilb2}). 
Then, for the blowup $\sigma : Y \to X$ along $\Gamma$, we can check that 
the Stein factorisation $\psi : Y \to Z$ of $\varphi_{|-K_Y|}$ satisfies $\Ex(\psi)=2$, because it contracts the proper transforms of infinitely many lines. 
However, we can show that $\Ex(\psi)=2$ occurs only when $g  \in \{4, 6\}$ (Proposition \ref{p-line-div-cont}), which contradicts our assumption: $g =5$ or $g \geq 7$. 
Thus $\lambda : \Univ_X^{\line} \to X$ is a finite surjective morphism. 
If $\lambda : \Univ_X^{\line} \to X$ is not purely inseparable, then we can again find a line $\Gamma$ on $X$ 
intersecting infinitely many lines on $X$ (Step \ref{s3-line-Hilb2} of Proposition \ref{p-line-Hilb2}), which leads to a contradiction as above. 
Then $\lambda : \Univ_X^{\line} \to X$ is a finite purely inseparable morphism. 
In this case, there exists $e \in \Z_{>0}$ such that the $e$-th iterated Frobenius morphism $F^e :X \to X$ factors through $\lambda$: 
\[
F^e :X \to \Univ_X^{\line} \xrightarrow{\lambda} X. 
\]
In particular, we get the following composite surjection $X \to \Univ_X^{\line} \to \Hilb_X^{\line}$ to a surface $\Hilb_X^{\line}$, which contradicts $\rho(X)=1$.

\subsection{Basic properties}

\begin{prop}\label{p-line-basic}
Let $X \subset \P^{g+1}$ be an anti-canonically embedded Fano threefold with $\Pic\,X = \Z K_X$. 
Let $\Gamma$ be a line on $X$. 
Let $\sigma : Y \to X$ be the blowup along $\Gamma$ and set $E := \Ex(\sigma)$. 
Then the following hold. 
\begin{enumerate}
    \item 
    \begin{enumerate}
    \item $(-K_Y)^3 = 2g -6$. 
    \item $(-K_Y)^2 \cdot E = 3$. 
    \item $(-K_Y) \cdot E^2 = -2$. 
    \item $E^3 = 1$. 
\end{enumerate}
    \item $|-K_Y|$ is base point free. Moreover, if $g \geq 4$, then $-K_Y$ is big. 
    \item 
        $h^0(E, -K_Y|_E) = 5$, 
    $h^0(Y, -K_Y) \geq g$, and 
    $h^0(Y, -K_Y -E) \geq g-5$.
    \item One of the following holds: 
\[
N_{\Gamma/X} \simeq \MO_{\Gamma} \oplus \MO_{\Gamma}(-1), \quad 
N_{\Gamma/X} \simeq \MO_{\Gamma}(1) \oplus \MO_{\Gamma}(-2), \quad 
\]
\item 
The following hold. 
\begin{enumerate}
\item 
If $N_{\Gamma/X} \simeq \MO_{\Gamma} \oplus \MO_{\Gamma}(-1)$, then $E \simeq \F_1$ and 
$-K_Y|_E$ is ample. 
\item 
If $N_{\Gamma/X} \simeq \MO_{\Gamma}(1) \oplus \MO_{\Gamma}(-2)$, 
then $E \simeq \F_3$ and $-K_Y|_E$ is not ample but nef and big. 
In particular, for a curve $B$ on $E$, 
$(-K_Y|_E) \cdot B=0$ if and only if $B$ is the section of the $\P^1$-bundle 
$E \to \Gamma$ whose self-intersection number is negative. 
\end{enumerate}
\end{enumerate}
\end{prop}

\begin{proof}
The assertion (1)  follows from Lemma \ref{l-blowup-formula}. 

Let us show (2).  
Let $\Lambda \subset |-K_X|$ be the linear system consisting of all the hyperplane sections containing the line $\Gamma$. 
Then the blowup $\sigma : Y \to X$ along $\Gamma$ coincides with the resolution of the indeterminacies of $\varphi_{\Lambda}$. 
By the scheme-theoretic equalities $\sigma^{-1}(\Gamma)=E$ and  
$\Bs\,\Lambda = \Gamma$, 
the complete linear system 
$|-K_Y| = |-\sigma^*K_X -E|$ is base point free. 
If $g \geq 4$, then we have $(-K_Y)^3 = 2g - 6 >0$ by (1), and hence $-K_Y$ is big. 
Thus (2) holds. 
The assertion (3) follows from (1) and Lemma \ref{l-bound--K-E}. 

Let us show (4) and (5). 
By $\deg N_{\Gamma/X} = 2g(\Gamma) -2 +\deg \Gamma=-1$ (Lemma \ref{l-blowup-formula}), 
we have $N_{\Gamma/X} \simeq \MO_{\Gamma}(m) \oplus \MO_{\Gamma}(-m-1)$ for some $m \in \Z_{\geq 0}$, and hence $E \simeq \F_{2m+1}$. 
Fix a fibre $f$ of the $\P^1$-bundle $\sigma|_E : E \to \Gamma$ and 
let $s$ be the section of $\sigma|_E$ with $s^2 =-(2m+1)$. 
We can write $-K_Y|_E \sim a s + b f$ for some $a, b \in \Z$. 
By $-K_Y|_E \cdot f =1$, we have $a=1$, and hence $-K_Y|_E \sim s + b f$. 
Since $-K_Y|_E$ is nef by (2), we get $0 \leq (-K_Y|_E) \cdot s = s^2 +b = -(2m+1) +b$, which implies 
$2m+1 \leq b$.   
We obtain 
\[
3 \overset{{\rm (1)}}{=} (-K_Y)^2 \cdot E =( -K_Y|_E)^2 = (s+bf)^2 = s^2 +2b 
\]
\[
= -(2m+1) + 2b 
\geq -(2m+1) +2 \cdot (2m+1) = 2m+1. 
\]
Therefore, it holds that $0 \leq m \leq 1$, which implies (4). 
The equality $-(2m+1) + 2b = 3$ deduces the following. 
\begin{itemize}
\item If $E \simeq \F_1$, i.e., $N_{\Gamma/X} \simeq \MO_{\Gamma} \oplus \MO_{\Gamma}(-1)$, then $-K_Y|_E \sim s + 2f$. 
\item If $E \simeq \F_3$, i.e., 
$N_{\Gamma/X} \simeq \MO_{\Gamma}(1) \oplus \MO_{\Gamma}(-2)$, then $-K_Y|_E \sim s + 3f$. 
\end{itemize}
Then we can check, by Kleiman's criterion for ampleness, that (5) holds. 
\qedhere

\end{proof}

\subsection{Divisorial contractions}

\begin{prop}\label{p-line-div-cont}
Let $X \subset \P^{g+1}$ be an anti-canonically embedded Fano threefold with $\Pic\,X = \Z K_X$ and $g \geq 4$. 
Let $\Gamma$ be a line on $X$. 
Let $\sigma : Y \to X$ be the blowup along $\Gamma$ and set $E := \Ex(\sigma)$. 
Let 
\[
\psi : Y \to Z 
\]
be the Stein factorisation of $\varphi_{|-K_Y|}$, 
which is a birational morphism (cf. Proposition \ref{p-line-basic}(2)). 
Assume that $\dim \Ex(\psi) = 2$. Set $D := \Ex(\psi)$. 
Recall that $D$ is a prime divisor and we have 
\[
D \sim -\alpha K_Y -\beta E
\]
for some $\alpha, \beta \in \Z_{>0}$ (Proposition \ref{p-div-cont1}). 
Then $(g, \alpha, \beta)=(4, 3, 2)$ or $(g, \alpha, \beta)=(6, 1, 2)$. 
\end{prop}

 \begin{proof}
Recall that 
we have $(2g-6)\alpha = 3\beta$ and $\beta \in \{1, 2\}$ (Proposition \ref{p-div-cont1}(5)(8)). 
If $\beta =1$, then 
\[
3   = (2g-6)\alpha \in 2\Z, 
\]
which is a contradiction. Hence $\beta =2$. 
Then $6 = (2g-6)\alpha$, i.e., 
\[
3 = (g-3)\alpha. 
\]
Hence we get $(g, \alpha, \beta) = (4, 3, 2)$ or $(g, \alpha, \beta) = (6, 1, 2)$. 
\end{proof}


\subsection{Hilbert schemes of lines}

Given an anti-canonically embedded Fano threefold $X \subset \P^{g+1}$, 
$\Hilb^{\line}_{X}$ and $\Univ^{\line}_X$ denote the Hilbert scheme of lines and its universal family, respectively \cite[Chapter 5]{FGI05}.

\begin{prop}\label{p-line-Hilb1}
Let $X \subset \P^{g+1}$ be an anti-canonically embedded Fano threefold with $\Pic\,X = \Z K_X$. 
Let $\Gamma$ be a line on $X$. 
Recall that 
$N_{\Gamma/X} \simeq \MO_{\Gamma} \oplus \MO_{\Gamma}(-1)$ or 
$N_{\Gamma/X} \simeq \MO_{\Gamma}(1) \oplus \MO_{\Gamma}(-2)$ (Proposition \ref{p-line-basic}). 
Then the following hold. 
\begin{enumerate}
\item 
If $N_{\Gamma/X} \simeq \MO_{\Gamma} \oplus \MO_{\Gamma}(-1)$, 
then  $\Hilb_X^{\line}$ is smooth at $[\Gamma]$ and $\dim \MO_{\Hilb_X^{\line}, [\Gamma]}=1$, 
where $\MO_{\Hilb_X^{\line}, [\Gamma]}$ denotes the local ring at $[\Gamma]$. 
\item If $N_{\Gamma/X} \simeq \MO_{\Gamma}(1) \oplus \MO_{\Gamma}(-2)$, 
then  
\begin{itemize}
    \item $1 \leq \dim \MO_{\Hilb_X^{\line}, [\Gamma]} \leq 2$, and 
    \item $\dim \MO_{\Hilb_X^{\line}, [\Gamma]}=2$ if and only if $\Hilb_X^{\line}$ is smooth at $[\Gamma]$. 
\end{itemize}
\end{enumerate}
In particular, either $\Hilb^{\line}_X = \emptyset$ or $1 \leq \dim \Hilb^{\line}_X \leq 2$. 
\end{prop}
 
\begin{proof}
Recall that 
\[
h^0(\Gamma, N_{\Gamma/X}) - h^1(\Gamma, N_{\Gamma/X})
\leq \dim \MO_{\Hilb^{\line}_X, [\Gamma]} \leq h^0(\Gamma, N_{\Gamma/X}) = \dim T_{[\Gamma]}\Hilb^{\line}_X, 
\]
where $T_{[\Gamma]}\Hilb^{\line}_X$ denotes the 
tangent space of $\Hilb^{\line}_X$ at $[\Gamma]$. 
Then (1) holds by 
$h^0(\Gamma, N_{\Gamma/X}) =1$ and $h^1(\Gamma, N_{\Gamma/X}) =0$ 
when $N_{\Gamma/X} \simeq \MO_{\Gamma} \oplus \MO_{\Gamma}(-1)$. 

Let us show (2). 
Assume $N_{\Gamma/X} \simeq \MO_{\Gamma}(1) \oplus \MO_{\Gamma}(-2)$. 
Then we have $h^0(\Gamma, N_{\Gamma/X}) = 2$ and 
\[
h^0(\Gamma, N_{\Gamma/X}) - h^1(\Gamma, N_{\Gamma/X}) = \chi(\Gamma, N_{\Gamma/X})
= \chi(\Gamma, \MO_{\Gamma}(-2)) + \chi(\Gamma, \MO_{\Gamma}(1)) = -1 + 2 = 1. 
\]
Thus (2) holds. 
\end{proof}

\begin{prop}\label{p-line-Hilb2}
Let $X \subset \P^{g+1}$ be an anti-canonically embedded Fano threefold with $\Pic\,X = \Z K_X$. 
Assume $g =5$ or $g \geq 7$. 
Then either $\Hilb^{\line}_X = \emptyset$ or 
every irreducible component of $\Hilb^{\line}_X$ is one-dimensional. 
In particular, there exists a closed point $x$ of $X$ such that no line on $X$ passes through $x$. 
\end{prop}

\begin{proof}
Recall that $1 \leq \dim \MO_{\Hilb^{\line}_X, [\Gamma]} \leq 2$ 
for any line $\Gamma$ on $X$ (Proposition \ref{p-line-Hilb1}). 
Suppose that there exists a line $\Gamma$ satisfying 
$\dim \MO_{\Hilb^{\line}_X, [\Gamma]} = 2$. 
It suffices to derive a contradiction.

Fix a two-dimensional irreducible component $H$ of $\Hilb^{\line}_X$, 
which we equip with the reduced scheme structure.  
Set $U := \Univ^{\line}_X \times_{\Hilb^{\line}_X} H$. 
We have the following diagram consisting of the induced morphisms: 
\[
\begin{tikzcd}[row sep=huge, column sep=huge, text height=1.5ex, text depth=0.25ex]
U \arrow[d, "\pi"] \arrow[r, hook, "j_U"] \arrow[rr, "\lambda", bend left=25] & \Univ^{\line}_X \arrow[d, "\Pi"] \arrow[r, "\Lambda"] & X\\
H  \arrow[r, hook, "j_H"] & \Hilb^{\line}_X.  
\end{tikzcd}
\]
Then $U$ is a projective threefold, because $H$ is a projective surface, $\pi$ is a flat projective morphism, 
and the geometric generic fibre of $\pi$ is an integral scheme. 

\setcounter{step}{0}

\begin{step}\label{s1-line-Hilb2}
Let $\Gamma$ be a line on $X$. Then there are at most finitely many lines on $X$ 
which intersect $\Gamma$. 
In particular, for every closed point $P \in X$, 
there are at most finitely many lines on $X$ 
passing through $P$. 
\end{step}

\begin{proof}[Proof of Step \ref{s1-line-Hilb2}]
Suppose that $\Gamma$ is a line on $X$ such that 
there are infinitely many lines $L_1, L_2, ...$ each of which intersects $\Gamma$. 
Let $\sigma: Y \to X$ be the blowup along $\Gamma$. 
Since $-K_Y$ is nef (Proposition \ref{p-line-basic}),
each proper transform $\sigma_*^{-1}L_i$ satisfies 
\[
0 \geq K_Y \cdot \sigma_*^{-1}L_i =(\sigma^*K_X +E) \cdot  \sigma_*^{-1}L_i 
= -1 + E \cdot  \sigma_*^{-1}L_i \geq 0. 
\]
Let 
\[
\psi : Y \to Z 
\]
be the Stein factorisation of $\varphi_{|-K_Y|}$, 
which is a birational morphism (cf. Proposition \ref{p-line-basic}(2)). 
Since $\Ex(\psi)$ contains infinitely many curves $\sigma_*^{-1}L_1, \sigma_*^{-1}L_2, ...$, 
$\psi$ is a divisorial contraction, i.e., $\dim \Ex(\psi) =2$. 
By Proposition \ref{p-line-div-cont}, we get $g=4$ or $g =6$, 
which contradicts our assumption. 
This completes the proof of Step \ref{s1-line-Hilb2}. 
\end{proof}

\begin{step}\label{s2-line-Hilb2}
$\lambda : U \to X$ is a finite surjective morphism. 
\end{step}

\begin{proof}[Proof of Step \ref{s2-line-Hilb2}]
Suppose that $\lambda : U \to X$ is not surjective. 
Fix a closed point $P \in \lambda(U) \subset X$. 
By $\dim U = \dim X (=3)$, we have $\dim \lambda^{-1}(P) \geq 1$. 
This implies that there exist infinitely many lines passing through $P$, which contradicts Step \ref{s1-line-Hilb2}. 
Hence $\lambda$ is surjective. 
If $\lambda$ is not finite, then we can find a closed point $Q \in X$ such that $\dim \lambda^{-1}(P) \geq 1$, which leads to a contradiction again. 
This completes the proof of Step \ref{s2-line-Hilb2}. 
\end{proof}


\begin{step}\label{s3-line-Hilb2}
$\lambda : U \to X$ is a  finite purely inseparable surjective morphism, 
i.e., the induced field extension $K(X) \subset K(U)$ is purely inseparable. 
\end{step}

\begin{proof}[Proof of Step \ref{s3-line-Hilb2}]
Suppose that the separable degree $d$ of the field extension $K(X) \subset K(U)$ satisfies $d \geq 2$, 
i.e., 
for a general closed point $P \in X$, its inverse image $\lambda^{-1}(P)$ consists of distinct $d$ points. 
Fix a general closed point $[\Gamma] \in H$. 
Since $\lambda : U \to X$ is surjective (Step \ref{s2-line-Hilb2}), 
the line $\Gamma$ corresponding to $[\Gamma]$ intersects a general locus of $X$ (because $\Gamma = \lambda( \pi^{-1}([\Gamma]))$). 
For a general closed point $P$ of $\Gamma$, 
we have $\# \lambda^{-1}(P) = d \geq 2$, and hence we can find another line $\Gamma_P$ passing through $P$. 
Note that $\Gamma \cap \Gamma_P = P$. 
Hence we obtain infinitely many lines $\{ \Gamma_P\}_P$ each of which intersects $\Gamma$. 
This contradicts Step \ref{s1-line-Hilb2}. 
This completes the proof of Step \ref{s3-line-Hilb2}. 
\end{proof}

Since $\lambda : U \to X$ is a finite purely inseparable surjective morphism, 
we can find $e \in \Z_{>0}$ such that 
the $e$-th iterated absolute Frobenius morphism $F^e: X \to X$ factors through $\lambda$: 
\[
F^e : X \to U \xrightarrow{\lambda} X. 
\]
In particular, the induced composite morphism 
$X \to U \to H$ 
is a surjective morphism to a surface $H$. 
However, this contradicts $\rho(X)=1$. 
This completes the proof of Proposition \ref{p-line-Hilb2}. 
\end{proof}

\section{Two-ray game for points}\label{s-pt}


Let $X \subset \P^{g+1}$ be an anti-canonically embedded Fano threefold with $\Pic\,X = \Z K_X$. 
The main results of this section are as follows: 
\begin{enumerate}
\item Let $P$ be a point of $X$ and let $\sigma: Y \to X$ be the blowup at $P$. 
Then $-K_Y$ is not ample (Corollary \ref{c-pt-not-ample}). 
\item If $g \geq 8$, then there exists a conic on $X$ (Corollary \ref{t-pt-to-conic}). 
\end{enumerate}
The strategy of the proof of (1) was explained in Subsection \ref{ss-2ray-intro}, i.e., 
we apply the two-ray game starting with the blowup at $P$. 
We now overview how to show (2). 
By Proposition \ref{p-line-Hilb2}, there exists a point $P$ on $X$ 
such that  no line on $X$ passes through $P$. 
Let $\sigma: Y \to X$ be the blowup at $P$. 
Since $-K_Y$ can be checked to be semi-ample (Proposition \ref{p-pt-basic}), 
(1) enables us to find a curve $\zeta$ on $Y$ satisfying $K_Y \cdot \zeta =0$. 
If $E \cdot \zeta =1$ for $E:=\Ex(\sigma)$, then 
the equality $K_Y =\sigma^*K_X +2E$ implies $-K_X \cdot \sigma(\zeta) =2$. 
Then  $\sigma(\zeta)$ is a conic, which is what we want to find. 
Recall that $E \cdot \zeta >0$ (Lemma \ref{l-beta=r}). Hence, by assuming $E \cdot \zeta \geq 2$, 
it suffices to derive a contradiction. 
A typical case is when $E \cap \zeta$ consists of two points. 
In this case, $E \to \psi(E)$ is not an isomorphism, as these two points are mapping to a single point on $\psi(E)$. 
Then we shall derive a contradiction by showing the surjectivity of $H^0(Y, -K_Y) \to H^0(E, -K_Y|_E)$ 
(Proposition \ref{p-max-rank-pt}).

\subsection{Basic properties}

\begin{prop}\label{p-pt-basic}
Let $X \subset \P^{g+1}$ be an anti-canonically embedded Fano threefold with $\Pic\,X = \Z K_X$ and $g \geq 6$. 
Let $P$ be a point on $X$ such that no line on $X$ passes through $P$. 
Let $\sigma : Y \to X$ be the blowup at $P$ and set $E := \Ex(\sigma)$. 
Then the following hold. 
\begin{enumerate}
    \item 
    \begin{enumerate}
    \item $(-K_Y)^3 = 2g -10$. 
    \item $(-K_Y)^2 \cdot E = 4$. 
    \item $(-K_Y) \cdot E^2 = -2$. 
    \item $E^3 = 1$. 
\end{enumerate}
    \item $-K_Y$ is semi-ample and big. 
\item 
    $h^0(E, -K_Y|_E) = 6$, 
    $h^0(Y, -K_Y) \geq g-2$, and 
    $h^0(Y, -K_Y -E) \geq g-8$. 
    \item Let $C \subset Y$ be a curve such that $-K_Y \cdot C =0$. 
Then one of the following holds. 
\begin{enumerate}
\item $C \cdot E = 1$, $C \not\subset E$, and $\sigma(C)$ is a conic on $X$. 
\item $C \cdot E \geq 2$ and $C \not\subset E$. 
\end{enumerate}
\end{enumerate}
\end{prop}

\begin{proof}
The assertion (1) follows from Lemma \ref{l-blowup-formula}(1). 
Assuming (2), the assertion (3) holds by (1) and Lemma \ref{l-bound--K-E}. 
Since $-E$ is $\sigma$-ample and $E$ is $\psi$-ample 
({\cred Proposition \ref{p-div-cont1},} Lemma \ref{l-beta=r}), we get $C \not\subset E$, 
{\cred where $\psi : Y \to Z$  denotes the birational morphism obtained by the Stein factorisation of $\varphi_{|-mK_Y|}$ for some (every) integer $m>0$ such that $|-mK_Y|$ is base point free.} 
Hence {\cred (2) implies} (4). 
It suffices to show (2). 


Recall that $X$ is an intersection of quadrics (Theorem \ref{t-inter-quads}). 
We can write $X = \bigcap_{X \subset Q}Q$, where $Q$ runs over all the quadric hypersurfaces in $\P^{g+1}$ containing $X$. 
We have 
$T_PX =\P^3$, i.e., $T_PX$ is a three-dimensional linear subvariety of $\P^{g+1}$. 
We  obtain 
\[
X \cap T_PX = \bigcap_{X \subset Q} (Q \cap T_PX) \qquad \text{and}\qquad  
T_PX = \bigcap_{X \subset Q} T_PQ.
\]

\setcounter{step}{0}

\begin{step}\label{s1-pt-basic}
    Let $Q$ be a quadric hypersurface $Q \subset \P^{g+1}$ such that $X \subset Q$ and $T_PX \not\subset Q$. 
Then $Q \cap T_PX$ is an effective Cartier divisor on $T_PX = \P^3$ such that $\deg (Q \cap T_PX)=2$ and $Q \cap T_PX$ is singular at $P$. 
\end{step}

\begin{proof}[Proof of Step \ref{s1-pt-basic}]
By $T_PX \not\subset Q$, $Q \cap T_PX$ is an effective Cartier divisor on $T_PX = \P^3$ such that $\deg (Q \cap T_PX)=2$. 
It suffices to show that $P$ is a singular point of $Q \cap T_PX$. 
 Suppose that $Q \cap T_PX$ is smooth at $P$. 
It holds that 
$ T_P(Q \cap T_PX) = T_PQ \cap T_P(T_PX) = T_PQ \cap T_PX$, 
and hence 
\[
T_PX= \bigcap_{X \subset Q} T_PQ \subset T_PQ \cap T_PX = T_P(Q \cap T_PX).  
\]
However, we have  $\dim T_PX =3$ and $\dim T_P(Q \cap T_PX) =2$, 
because $X$ and $Q \cap T_PX$ are smooth at $P$. 
This is absurd. 
This completes the proof of Step \ref{s1-pt-basic}. 
\end{proof}

\begin{step}\label{s2-pt-basic}
    $(X \cap T_PX)_{\red} =  P$. 
\end{step}

\begin{proof}[Proof of Step \ref{s2-pt-basic}]
    Suppose  $(X \cap T_PX)_{\red} \neq  P$. 
Then we can find a point $P' \in X \cap T_PX$ with $P \neq P'$. 
In order to derive a contradiction, it suffices to show that $\overline{PP'} \subset X$, 
where $\overline{PP'}$ denotes the line on $\P^{g+1}$ passing through $P$ and $P'$. 
We fix a quadric hypersurface $Q \subset \P^{g+1}$ with $X \subset Q$. 
By $X = \bigcap_{X \subset Q} Q$, it is enough to prove that $\overline{PP'} \subset Q$. 
If $T_PX \subset Q$, then this 
holds by $\overline{PP'} \subset T_PX \subset Q$. 
We may assume $T_PX \not\subset Q$. 
Then Step \ref{s1-pt-basic} implies that $Q \cap T_PX$ is an effective Cartier divisor on $T_PX =\P^3$ which is of degree two and singular at $P$. 
Hence either $Q$ is a singular quadric surface or a sum of two planes $S$ and $S'$ whose intersection contains $P$ (in the latter case, $S = S'$ might occur, i.e., $Q \cap T_PX = 2S$). 
In any case, $Q \cap T_PX$ can be written (set-theoretically) as a union of lines passing through $P$: 
$Q \cap T_PX= \bigcup_{R \in Q \cap T_PX} \overline{PR}$. 
By $P' \in X \cap T_PX \subset Q \cap T_PX$, 
we obtain  $\overline{PP'} \subset Q \cap T_PX \subset Q$. 
This completes the proof of Step \ref{s2-pt-basic}. 
\end{proof}

\begin{step}\label{s3-pt-basic}
$\Bs\,|-K_Y| \subset E$. 
\end{step}

\begin{proof}[Proof of Step \ref{s3-pt-basic}]
Fix a closed point $P'_Y \in Y$ satisfying $P'_Y \not\in E$. 
Set $P':= \sigma(P'_Y)$. 
It is enough to find an effective divisor $F$ on $Y$ such that $-K_Y \sim F$ and $P'_Y \not\in \Supp\,F$. 
Note that $T_PX$ is an intersection of hyperplanes in $\P^{g+1}$. 
Then, by $P \neq P'$ and Step \ref{s2-pt-basic},  there exists a hyperplane 
$H \subset \P^{g+1}$ such that 
$T_PX \subset H$ and $P' \not\in H$. 
For $H_X := H \cap X$, we have that 
$\sigma^*H_X = \sigma_*^{-1}H_X + m E$ for some integer  $m \geq 2$, 
where $\sigma_*^{-1}H_X$ denotes the proper transform of $H_X$ on $Y$. 
By $-K_X \sim H_X$,  we get 
\[
-K_Y \sim -\sigma^*K_X -2E
\sim \sigma^*H_X -2E= \sigma_*^{-1}H_X + (m-2) E. 
\]
Then the last term $F := \sigma_*^{-1}H_X + (m-2) E$ is 
an effective divisor such that $-K_Y \sim F$ and 
$P'_Y \not\in \Supp\,F$. 
This completes the proof of Step \ref{s3-pt-basic}. 
\end{proof}

\begin{step}\label{s4-pt-basic}
(2) holds, i.e., $-K_Y$ is semi-ample and big. 
\end{step}

\begin{proof}[Proof of Step \ref{s4-pt-basic}]
By $(-K_Y)^3 = 2g -10 >0$, it is enough to prove that $-K_Y$ is semi-ample. 
It suffices to show that $-K_Y$ is nef \cite[Theorem H]{BMPSTWW}. 
Fix a curve $C$ on $Y$. 
If $C \not\subset E$, then we get $-K_Y \cdot C \geq 0$ 
(Step \ref{s3-pt-basic}). 
Hence it is enough to show $-K_Y \cdot C \geq 0$ when $C \subset E$. 
This follows from  
\[
-K_Y \cdot C = -(\sigma^*K_X + 2E) \cdot C =-2E \cdot C>0. 
\]
This completes the proof Step \ref{s4-pt-basic}. 
\end{proof}
This completes the proof of Proposition \ref{p-pt-basic}. 
\qedhere
\end{proof}






\subsection{Divisorial contractions}

\begin{prop}\label{p-pt-div-cont}
Let $X \subset \P^{g+1}$ be an anti-canonically embedded Fano threefold with $\Pic\,X = \Z K_X$ and $g \geq 6$. 
Let $P$ be a point on $X$ such that 
there exists no line on $X$ passing through $P$. 
Let $\sigma : Y \to X$ be the blowup at $P$ and set $E := \Ex(\sigma)$. 
Let 
\[
\psi : Y \to Z 
\]
be the Stein factorisation of $\varphi_{|-mK_Y|}$ for some (every) positive integer $m$ such that $|-mK_Y|$ is base point free (Proposition \ref{p-pt-basic}). 
Assume that $\dim \Ex(\psi) = 2$. Set $D := \Ex(\psi)$. 
Recall that $D$ is a prime divisor and we have 
\[
D \sim -\alpha K_Y -\beta E
\]
for some $\alpha, \beta \in \Z_{>0}$ (Proposition \ref{p-div-cont1}). 
Then the following hold. 
\begin{enumerate}
\item 
One of the following holds. 
\begin{enumerate}
\item $(g, \alpha, \beta)= (6, 2, 1)$. 
\item $(g, \alpha, \beta)= (7, 1, 1)$. 
\item $(g, \alpha, \beta)= (6, 4, 2)$. 
\item $(g, \alpha, \beta)= (7, 2, 2)$. 
\item $(g, \alpha, \beta)= (9, 1, 2)$. 
\end{enumerate}
\item 
There exists a conic on $X$ passing through $P$. 
\end{enumerate}
\end{prop}

\begin{proof}
Let us show (1). 
By Proposition \ref{p-div-cont1}(5)(8) and Proposition \ref{p-pt-basic}(1), we have $\beta \in \{1, 2\}$ and 
$(2g-10) \alpha = 4\beta$. Then 
\[
(g-5)\alpha = 2\beta. 
\]
If $\beta =1$, then $(g, \alpha, \beta) = (6, 2, 1), (7, 1, 1)$. 
If $\beta =2$, then $(g, \alpha, \beta) = (6, 4, 2), (7, 2, 2), (9, 1, 2)$. 
Thus (1) holds. 

Let us show (2). 
If there exists a curve $\zeta$ on $Y$ such that $K_Y \cdot \zeta =0$ and $E \cdot \zeta =1$, 
then 
$\sigma(\zeta)$ is a required conic (Proposition \ref{p-pt-basic}(4)). 
Suppose that $E \cdot \zeta \neq 1$ 
for every curve $\zeta$ on $Y$ satisfying $K_Y \cdot \zeta =0$. 
It suffices to derive a contradiction. 
We use the same notation as in Proposition \ref{p-div-cont2}. 
By Proposition \ref{p-div-cont2}(1){\cred (5)}, we have $\beta =1$ and $\omega_D^2 \in 4\Z \subset 2\Z$. 
Recall that 
the following hold (Proposition \ref{p-pt-basic}(1)):
\[
E^3=1\qquad \text{and}\qquad 
K_Y^3 \equiv K_Y^2 \cdot E \equiv K_Y \cdot E^2 \equiv 0 \mod 2.  
\]
By $D \sim -\alpha K_Y -E$, we get  
\[
\omega_D^2 = (K_Y+D)^2 \cdot D = ( (1-\alpha)K_Y-E)^2 \cdot (-\alpha K_Y-E) 
\equiv (-E)^3 \equiv 1  \mod 2, 
\]
which contradicts $\omega_D^2 \in 2\Z$. 
Thus (2) holds. 
\end{proof}

\subsection{Flops}

\begin{thm}\label{t-pt-flop}
Let $X \subset \P^{g+1}$ be an anti-canonically embedded Fano threefold with $\Pic\,X = \Z K_X$ and $g \geq 6$. 
Let $P$ be a point on $X$ such that there is no line on $X$ passing through $P$. 
Let $\sigma : Y \to X$ be the blowup at $P$ and set $E := \Ex(\sigma)$. 
Since $-K_Y$ is semi-ample and big (Proposition \ref{p-pt-basic}), we have the birational morphism 
\[
\psi : Y \to Z 
\]
which is the Stein factorisation of 
$\varphi_{|-mK_Y|}$ for some (every) positive integer $m$ such that $|-mK_Y|$ is base point free. 
Assume that $\dim \Ex(\psi) = 1$ or $\psi$ is an isomorphism. 
\begin{itemize}
    \item If $\dim \Ex(\psi) = 1$, then let $\psi^+ : Y^+ \to Z$ be the flop of $\psi$, set $E^+$ to be the proper transform of $E$ on $Y^+$, and 
let $\tau : Y^+ \to W$ be the contraction of the $K_{Y^+}$-negative extremal ray. 
    \item If $\psi$ is an isomorphism, then we set $Z :=Y$, $Y^+ :=Y$, $\psi := {\rm id}$, 
$\psi^+ := {\rm id}, E^+ :=E$, and 
let $\tau : Y^+ \to W$ be the contraction of the $K_{Y^+}$-negative extremal ray not corresponding to $\sigma$. 
\end{itemize}
\[
\begin{tikzcd}
Y \arrow[d, "\sigma"'] \arrow[rd, "\psi"]& & Y^+ \arrow[ld, "\psi^+"'] \arrow[d, "\tau"]\\
X & Z & W
\end{tikzcd}
\]
Then $g$, the type of $\tau$, $D$, and $D^3$ satisfy one of the possibilities in 
Table \ref{table-2ray-pt}, where we use the following notation. 
\begin{itemize}
\item 
If $\tau$ is of type $E$, then $D := \Ex(\tau)$. 
If $\tau$ is of type $C$ or $D$, then $D$ is the pullback of an ample Cartier divisor which is a generator of $\Pic\,W$. 
The divisors below $D$ means the linear equivalence, e.g., 
if $g=11$, then $D \sim -K_{Y^+} -E^+$. 
\item If $\tau$ is of type $C$, then $\Delta$ denotes the discriminant bundle. 
\item If $\tau$ is of type $E_1$ or $E_2$, 
then let $r_W$ be the index of the Fano threefold $W$. 
\item If $\tau$ is of type $E_1$, then we set 
$B := \tau(D)$ and  $g(B)$ denotes the genus of $B$. 
\end{itemize}
\begin{table}[h]
\caption{Two-ray game for points}\label{table-2ray-pt}
     \centering
{\renewcommand{\arraystretch}{1.35}%
      \begin{tabular}{|c|c|c|c|c|}
      \hline
$g$  & type of $\tau$ & $D$ & $D^3$ & Other properties\\      \hline
$6$  & $E_2$ & $-4K_{Y^+}-E^+$ & $1$ & $(r_W, -K_W^3)=(1, 10)$\\           \hline
$7$  & $E_2$ & $-2K_{Y^+}-E^+$ & $1$ & $(r_W, -K_W^3)=(1, 12)$\\           \hline
$9$  & $E_2$ & $-K_{Y^+}-E^+$ & $1$ & $(r_W, -K_W^3)=(1, 16)$\\           \hline
$10$  & $D_1$ & $-K_{Y^+}-E^+$ & $0$ &$(-K_{Y^+})^2 \cdot D =6$ \\           \hline
$11$  & $C_1$  & $-K_{Y^+}-E^+$ & $0$ &$\deg \Delta =4$ \\           \hline
$6$  & $E_1$ & $-4K_{Y^+}-E^+$ & $0$ & $(r_W, -K_W^3, -K_W\cdot B, g(B))=(1, 8, 2, 0)$\\           \hline
$6$  & $E_1$ & $-5K_{Y^+}-E^+$ & $-22$ & $(r_W, -K_W^3, -K_W\cdot B, g(B))=(1, 22, 14, 5)$\\           \hline
$6$  & $E_1$ & $-9K_{Y^+}-2E^+$ & $-30$ & $(r_W, -K_W^3, -K_W\cdot B, g(B))=(2, 32, 20, 6)$\\           \hline
$7$  & $E_1$ & $-2K_{Y^+}-E^+$ & $0$ & $(r_W, -K_W^3, -K_W\cdot B, g(B))=(1, 10, 2, 0)$\\           \hline
$7$  & $E_1$ & $-5K_{Y^+}-2E^+$ & $-36$ & $(r_W, -K_W^3, -K_W\cdot B, g(B))=(2, 40, 24, 7)$\\           \hline
$8$  & $E_1$ & $-3K_{Y^+}-2E^+$ & $-6$ & $(r_W, -K_W^3, -K_W\cdot B, g(B))=(2, 24, 8, 0)$\\           \hline
$8$  & $E_1$ & $-5K_{Y^+}-3E^+$ & $-42$ & $(r_W, -K_W^3, -K_W\cdot B, g(B))=(3, 54, 30, 7)$\\           \hline
$9$  & $E_1$ & $-K_{Y^+}-E^+$ & $0$ & $(r_W, -K_W^3, -K_W\cdot B, g(B))=(1, 14, 2, 0)$\\           \hline
$12$  & $E_1$ & $-3K_{Y^+}-4E^+$ & $-22$ & $(r_W, -K_W^3, -K_W\cdot B, g(B))=(4, 64, 24, 0)$\\           \hline
$13$  & $E_1$ & $-2K_{Y^+}-3E^+$ & $-16$ &  $(r_W, -K_W^3, -K_W\cdot B, g(B))=(3, 54, 18, 0)$\\           \hline
      \end{tabular}}
    \end{table}
\end{thm}


\begin{proof}
Recall that we have 
\[
D \sim -\alpha K_{Y^+} -\beta E^+, 
\]
where $\alpha, \beta \in \Z_{>0}$ (Lemma \ref{l-beta=r}). 
The following hold (Lemma \ref{l-blowup-formula2}(1)): 
\begin{align*}
D^2 \cdot (-K_{Y^+})
&= (\alpha (-K_{Y^+}) -\beta E^+)^2 \cdot (-K_{Y^+})\\
&= \alpha^2 (-K_{Y^+})^3  -2 \alpha \beta (-K_{Y^+})^2 \cdot E^+ + \beta^2 (-K_{Y^+}) \cdot (E^+)^2\\
&= (2g-10) \alpha^2 -8 \alpha \beta -2 \beta^2\\
D \cdot (-K_{Y^+})^2 & = (\alpha (-K_{Y^+}) -\beta E^+) \cdot (-K_{Y^+})^2\\
&= \alpha (-K_{Y^+})^3  - \beta (-K_{Y^+})^2 \cdot E^+\\
&= (2g-10) \alpha -4\beta. 
\end{align*}


\setcounter{step}{0}

\begin{step}\label{s1-pt-flop}
If $\tau$ is of type $D$, 
then 
$(g, \alpha, \beta, D \cdot (-K_{Y^+})^2)
=(10, 1, 1, 6)$ and $\tau$ is of type $D_1$. 
\end{step}

\begin{proof}[Proof of Step \ref{s1-pt-flop}]
Recall that $\alpha \in \Z_{>0}$ and 
$\beta \in \{1, 2, 3\}$ (Lemma \ref{l-beta=r}). 
The following holds (Proposition \ref{p-typeD-intersec}): 
\[
f(\alpha, \beta) := (g-5) \alpha^2 -4 \alpha \beta -\beta^2 
=0. 
\]
Then the quadratic equation 
\[
(g-5)x^2 -4x -1 =0
\]
has a solution $x = \alpha/\beta \in \Q$. 
Hence the discriminant $D$ satisfies 
\[
\Q \ni D = \sqrt{(-4)^2 +4(g-5) } = 2 \sqrt{g -1}, 
\]
which implies 
\[
g \in \{ n^2 +1 \,|\, n \in \Z_{>0}\} = \{2, 5, 10, {\cred 17}, 26, 37, ...\}. 
\]
By $g \geq 6$, we get $g=10$ or $g \geq {\cred 17}$. 
{\cred The latter case $g \geq 17$ does not occur, 
because the inequality $g \geq 17$ leads to the following contradiction
 (Proposition \ref{p-typeD-intersec}):}
\[
{\cred 
9 \geq  D \cdot (-K_{Y^+})^2 =  (2g-10) \alpha - 4\beta 
\geq  (2 \cdot 17 -10) \cdot 1 -4 \cdot 3  = 24 -12 =12.}  
\]
{\cred Thus it holds that} $g=10$. 
Then we have 
\[
0 = 5\alpha^2 -4\alpha \beta -\beta^2 = (5\alpha +\beta)(\alpha -\beta), 
\]
which implies $\alpha = \beta$. 
In this case, the following holds (Proposition \ref{p-typeD-intersec}):
\[
9 \geq  D \cdot (-K_{Y^+})^2 =  
 (2g-10) \alpha - 4\beta = (2g-14)\alpha =  6\alpha. 
\]
Therefore, $(g, \alpha, \beta, D \cdot (-K_{Y^+})^2) = (10, 1, 1, 6)$. 
This completes the proof of Step \ref{s1-pt-flop}. 
\end{proof}

\begin{step}\label{s2-pt-flop}
If $\tau$ is of type $C$, 
then 
$(g, \alpha, \beta, \deg \Delta) =(11, 1, 1, 4)$ and  $\tau$ is of type $C_1$. 
\end{step}

\begin{proof}[Proof of Step \ref{s2-pt-flop}]
Recall that 
$\alpha \in \Z_{>0}$ and $\beta \in \{1, 2\}$ (Lemma \ref{l-beta=r}). 
The following holds (Proposition \ref{p-typeC-intersec}): 
\[
f(\alpha, \beta) := (g-5) \alpha^2 -4 \alpha \beta -\beta^2 
=1. 
\]
We get 
\[
\alpha = \frac{2\beta \pm \sqrt{4\beta^2 +(g-5)(\beta^2+1)}}{g-5}.
\]
There is no solution when $g=6$ or $g=7$, 
because $\sqrt{4\beta^2 +(g-5)(\beta^2+1)}$ is not a rational number as follows: 
\begin{itemize}
\item $(g, \beta)=(6, 1) \Rightarrow \sqrt{4\beta^2 +(g-5)(\beta^2+1)} = 
\sqrt{4 +1 \cdot (1+1)}  =\sqrt{6} \not\in \Q$. 
\item $(g, \beta)=(6, 2)\Rightarrow 
\sqrt{4\beta^2 +(g-5)(\beta^2+1)} = 
\sqrt{16 +1 \cdot (4+1)}  =\sqrt{21} \not\in \Q$. 
\item $(g, \beta)=(7, 1)\Rightarrow 
\sqrt{4\beta^2 +(g-5)(\beta^2+1)} = 
\sqrt{4 +2 \cdot (1+1)}  =\sqrt{8} \not\in \Q$. 
\item  $(g, \beta)=(7, 2)\Rightarrow 
\sqrt{4\beta^2 +(g-5)(\beta^2+1)} = 
\sqrt{16 +2 \cdot (4+1)}  =\sqrt{26} \not\in \Q$. 
\end{itemize}
In what follows, we assume $g \geq 8$.

Assume that $\beta=1$. 
In this case, we get 
\[
f(\alpha, 1) = (g-5) \alpha^2 -4 \alpha -1 =1. 
\]
By $g \geq 8$, we obtain $f(1, 1) < 
f(2, 1) < f(3, 1) < \cdots$. 
We have 
\begin{itemize}
\item $f(1, 1) = (g-5)  -4 -1 = g-10$.  
\item $f(2, 1) = 4(g-5)  -8  -1 = 4g -29 \geq 4 \cdot 8 -29 =3>1$. 
\end{itemize}
Thus $f(\alpha, 1)=1$ has a unique solution 
$(g, \alpha, \beta)= (11, 1, 1)$. 

Assume that $\beta=2$. 
In this case, we get 
\[
f(\alpha, 2) = (g-5) \alpha^2 -8 \alpha -4 =1. 
\]
By $g \geq 8$, we obtain 
$f(2, 2) < f(3, 2) < \cdots$. 
We have 
\begin{itemize}
\item $f(1, 2) = (g-5)  -8 -4 = g-17$.  
\item $f(2, 2) = 4f(1, 1) = 4g-40$.
\item $f(3, 2) =  9(g-5) -24 -4 = 9g-73$.
\item $f(4, 2) = 4f(2, 1) = 4(4g-29) \geq 4 \cdot 3>1$. 
\end{itemize}
Then $f(\alpha, 2)=1$ has a unique solution 
$(g, \alpha, \beta) = (18, 1, 2)$.



\medskip

To summarise, we obtain 
$(g, \alpha, \beta) \in \{(11, 1, 1), (18, 1, 2)\}$.
Recall that the following holds (Proposition \ref{p-typeC-intersec}):
\[
12 - \deg \Delta = D \cdot (-K_{Y^+})^2 
= (2g-10) \alpha - 4\beta. 
\]
\begin{itemize}
\item If $(g, \alpha, \beta) = (11, 1, 1)$, then 
$12 - \deg \Delta = (22-10)\cdot 1 - 4 \cdot 1 = 8$. 
\item 
If $(g, \alpha, \beta) = (18, 1, 2)$, then  
$12 - \deg \Delta = (36-10)\cdot 1 - 4 \cdot 2 = 18$, 
which  contradicts $\deg \Delta \geq 0$ (Proposition \ref{p-typeC-intersec}(3)).  
\end{itemize}
This completes the proof of Step \ref{s2-pt-flop}. 
\end{proof}

\begin{step}\label{s3-pt-flop}
Assume that $\tau$ is of type $E_2, E_3, E_4$, or $E_5$.  
Then $\tau$ is of type $E_2$ and $(g, \alpha, \beta) = (9, 1, 1), (7, 2, 1), (6, 4, 1)$. 
Furthermore, $r_W=1$ and $-K_W^3 = -K_X^3$. 
\end{step}

\begin{proof}[Proof of Step \ref{s3-pt-flop}]
The following hold (Proposition \ref{p-typeE-intersec}): 
\begin{eqnarray*}
(2g-10) \alpha -4\beta &=& D \cdot (-K_{Y^+})^2 
=: u \in \{1, 2, 4\}\\
(2g-10) \alpha^2 - 8\alpha \beta -2 \beta^2  &=&
D^2 \cdot (-K_{Y^+}) =-2.
\end{eqnarray*}
We then get 
\[
\alpha(4\beta +u) = (2g-10)\alpha^2 = 8\alpha \beta +2\beta^2 -2, 
\]
which implies 
\[
\alpha (-4\beta +u) = 2\beta^2 -2 \geq 0.
\]
By $u \in \{1, 2, 4\}$, we obtain $u=4$ and $\beta =1$. 
In particular, $\tau$ is of type $E_2$  (Proposition \ref{p-typeE-intersec}) and 
$r_W = \beta =1$ (Lemma \ref{l-beta=r}). 
We have 
\[
(g-5) \alpha =\frac{1}{2}(4\beta + u) =4. 
\]
Hence $(g, \alpha, \beta) = (9, 1, 1), (7, 2, 1), (6, 4, 1)$. 
Finally, we obtain $K_X^3 = K_Y^3 -8 = K^3_{Y^+} -8 = K_W^3$. 
This completes the proof of Step \ref{s3-pt-flop}. 
\end{proof}

\begin{step}\label{s4-pt-flop}
Assume that $\tau$ is of type $E_1$. 
Then one of the possibilities appearing in Table \ref{table-2ray-pt} holds. 
\end{step}

\begin{proof}[Proof of Step \ref{s4-pt-flop}]
In this case, the following hold. 
\begin{enumerate}
\renewcommand{\labelenumi}{(\roman{enumi})}
\item $W$ is a Fano threefold with $\rho(W) =1$ and $\beta =r_W$ for the index $r_W$ of $W$ (Lemma \ref{l-beta=r}). 
In particular, $1 \leq \beta \leq 4$. Fix a Cartier divisor $H_W$ on $W$ with $-K_W \sim \beta H_W$. 
\item $\alpha+1 = \beta \gamma$ for some $\gamma \in \Z_{>0}$  (Lemma \ref{l-beta=r}). 
\item $w := \frac{1}{2}\beta H_W^3 = (g-5)\gamma^2 - 4\gamma  -1$ (Remark \ref{r-E1-deg-genus}). 
\item $h^0(Y^+, -K_{Y^+}-E^+) = h^0(Y, -K_Y-E) \geq g-8$ (Proposition \ref{p-pt-basic}). 
\item $2g(B) -2 =(2g-10)\alpha^2 - 8 \alpha \beta  -2 \beta ^2$ (Lemma \ref{l-E1-deg-genus}). 
\item $-K_W \cdot B = (2g-10)\alpha(\alpha+1) 
- 4(2\alpha +1)\beta  -2\beta^2$  (Lemma \ref{l-E1-deg-genus}).
\item $D^3 = -(2g-10)+ (-K_W)^3 -3 (-K_W) \cdot B$ (Lemma \ref{l-D^3}). 
\end{enumerate}

We now prove the following assertion $(*)$. 
\begin{enumerate}
\item[$(*)$] If $\alpha < \beta$ or $g \geq 10$, then $(g, \alpha, \beta, -K_W^3)  \in \{ (12,  3, 4, 64), (13, 2, 3, 54)\}$. 
In particular, $g \leq 13$. 
\end{enumerate}
We first show that $g \geq 10 \Rightarrow \alpha < \beta$. 
Assume  $g \geq 10$. 
We obtain 
\[
h^0(Y^+, -K_{Y^+}-E^+) = h^0(Y, -K_Y-E) \geq g-8 \geq 2. 
\]
If $\alpha \geq \beta$, then we would get the following contradiction: 
\[
1 = h^0(Y^+, D) = h^0(Y^+, -\alpha K_{Y^+} - \beta E^+) \geq  
h^0(Y^+, - \beta K_{Y^+} - \beta E^+)  \geq g-8 \geq 2. 
\]
This completes the proof of the implication $g \geq 10 \Rightarrow \alpha < \beta$. 
In order to prove $(*)$, we may assume that $\alpha +1 \leq  \beta$. 
This, together with $\alpha+1 =\beta \gamma \geq \beta$, implies  
$\alpha +1= \beta$ and $\gamma =1$. 
By (iii), we get $g= 10 + \frac{1}{2}\beta H_W^3$. 
If $\beta=4$, then $\alpha =3$, $(-K_W)^3 = 64$, and $g=12$. 
If $\beta=3$, then $\alpha =2$, $(-K_W)^3 = 54$, and $g=13$. 
Assume $\beta=2$. Then $(\alpha, \beta, \gamma) =(1, 2, 1)$. 
By (iii) and (v), we get 
\[
5 \geq H_W^3 = g -10 \qquad \text{and}\qquad 
-2 \leq 2g(B) -2  = 2g-34. 
\]
This is absurd, 
because these two inequalities 
deduce $15 \geq g \geq 16$. 
This completes the proof of $(*)$. 
Since the proof of Theorem \ref{t-pt-flop} is completed when $g \geq 10$, 
we may use $g_W \leq 13$ when $r_W=1$. 
In particular, the following holds (Remark \ref{r-E1-deg-genus}): 
\[
1 \leq w \leq 12. 
\]

\medskip

In what follows, we assume $g \leq 9$ and $\alpha \geq \beta$. 
The solution of $w = (g-5)\gamma^2 -4 \gamma -1$ is given by 
\[
\gamma = \frac{2 + \sqrt{4 +(g-5)\left(w+1\right)}}{g-5}, 
\]
because we have $g-5 >0$ and $4 +(g-5)\left(w+1\right) >4$.

Assume $g=9$. 
Then $\gamma = \frac{1 + \sqrt{w+2}}{2}$. 
By $3 \leq w+2 \leq 14$, 
we need $w+2 \in \{4, 9\}$ to assure $\sqrt{w+2} \in \Z$. 
It follows from  $\gamma \in \Z$ that $(w, \gamma) = (7, 2)$. 
By (ii), (iii), and Remark \ref{r-E1-deg-genus}, 
we get  $(g, \alpha, \beta,  -K_W^3) = (9, 1, 1, 14)$. 


Assume $g=8$. 
Then $\gamma = \frac{2 + \sqrt{3w+7 }}{3}$. 
By $10 \leq 3w+7 \leq 43$ and $3w+7 \in 3\Z +1$, 
we need $3w+7 \in \{16, 25\}$ to assure $\sqrt{3w+7} \in \Z$. 
It follows from $\gamma \in \Z$ that $(w, \gamma) = (3, 2)$. 
By (ii), (iii), and Remark \ref{r-E1-deg-genus}, 
we get 
\[
(g, \alpha, \beta, -K_W^3) \in \{(8, 5, 3, 54), (8, 3, 2, 24), (8, 1, 1, 6)\}.
\]
The case $(g, \alpha, \beta, -K_W^3)=(8, 1, 1, 6)$ can be excluded by (vi): 
\[
0< -K_W \cdot B = 2 \cdot 6 - 4 \cdot 3 -2 =-2. 
\]

Assume $g=7$. 
Then $\gamma = \frac{2 + \sqrt{2w+6 }}{2}$. 
By $8 \leq 2w+6 \leq 30$ and $2w+6 \in 2\Z$, 
we need $2w+6  =16$ to assure $\sqrt{2w+6} \in \Z$. 
It holds that  $(w, \gamma) = (5, 3)$. 
By (ii), (iii), and Remark \ref{r-E1-deg-genus}, 
we get 
\[
(g, \alpha, \beta, -K_W^3) \in \{(7, 5, 2, 40), (7, 2, 1, 10)\}.
\]

Assume $g=6$. 
Then $\gamma = 2 + \sqrt{w+5 }$. 
By $6 \leq w+5 \leq 17$, 
we need $w+5  \in \{ 9, 16\}$ to assure $\sqrt{w+5} \in \Z$. 
It holds that  $(w, \gamma) \in \{(4, 5), (11, 6)\}$. 
By (ii), (iii), and Remark \ref{r-E1-deg-genus}, 
we get 
\[
(g, \alpha, \beta, -K_W^3) \in \{(6, 9, 2, 32), (6, 5, 1, 22), (6, 4, 1, 8)\}.
\]

\medskip

The remaining invariants $g(B), -K_W \cdot B, D^3$ can be computed by (v), (vi), (vii), respectively. 
For the sake of completeness, we here include the computations for all the cases. \\
If $(g, \alpha, \beta, -K_W^3) = (13, 2, 3, 54)$, then the following hold.  
\begin{enumerate}
\item[(v)] $2g(B) -2 = (26-10) \cdot 2^2 -8 \cdot 2 \cdot 3 -2 \cdot 3^2 = 64 -48 - 18 =-2$. 
\item[(vi)] $-K_W \cdot B = (26-10) \cdot 2 \cdot 3 -4 \cdot 5 \cdot 3 -2 \cdot 3^2 = 96 -60 - 18 =18$. 
\item[(vii)] $D^3 = -(26-10) + 54 -3 \cdot 18 = -16 +54-54=-16$. 
\end{enumerate}
If $(g, \alpha, \beta, -K_W^3) = (12, 3, 4, 64)$, then the following hold. 
\begin{enumerate}
\item[(v)] $2g(B) -2 = (24-10) \cdot 3^2 -8 \cdot 3 \cdot 4 -2 \cdot 4^2 = 126 -96 - 32 =-2$. 
\item[(vi)] $-K_W \cdot B = (24-10) \cdot 3 \cdot 4 -4 \cdot 7 \cdot 4 -2 \cdot 4^2 = 168-112-32 =24$. 
\item[(vii)] $D^3 = -(24-10) + 64 -3 \cdot 24 = -14 +64-72=-22$. 
\end{enumerate}
If $(g, \alpha, \beta, -K_W^3) = (9, 1, 1, 14)$, then the following hold. 
\begin{enumerate}
\item[(v)] $2g(B) -2 = (18-10) \cdot 1^2 -8 \cdot 1 \cdot 1 -2 \cdot 1^2 = 8-8-2 =-2$. 
\item[(vi)] $-K_W \cdot B = (18-10) \cdot 1 \cdot 2 -4 \cdot 3 \cdot 1 -2 \cdot 1^2 = 16-12-2 =2$. 
\item[(vii)] $D^3 = -(18-10) + 14 -3 \cdot 2 = -8+14-6=0$. 
\end{enumerate}
If $(g, \alpha, \beta, -K_W^3) = (8, 5, 3, 54)$, then the following hold. 
\begin{enumerate}
\item[(v)] $2g(B) -2 = (16-10) \cdot 5^2 -8 \cdot 5 \cdot 3 -2 \cdot 3^2 = 150-120-18 =12$. 
\item[(vi)] $-K_W \cdot B = (16-10) \cdot 5 \cdot 6 -4 \cdot 11 \cdot 3 -2 \cdot 3^2 = 180-132-18 =30$. 
\item[(vii)] $D^3 = -(16-10) + 54 -3 \cdot 30 = -6+54-90=-42$. 
\end{enumerate}
If $(g, \alpha, \beta, -K_W^3) = (8, 3, 2, 24)$, then the following hold. 
\begin{enumerate}
\item[(v)] $2g(B) -2 = (16-10) \cdot 3^2 -8 \cdot 3 \cdot 2 -2 \cdot 2^2 = 54-48-8 =-2$. 
\item[(vi)] $-K_W \cdot B = (16-10) \cdot 3 \cdot 4 -4 \cdot 7 \cdot 2 -2 \cdot 2^2 = 72-56-8 =8$. 
\item[(vii)] $D^3 = -(16-10) + 24 -3 \cdot 8 = -6+24-24=-6$. 
\end{enumerate}
If $(g, \alpha, \beta, -K_W^3) = (7, 5, 2, 40)$, then the following hold. 
\begin{enumerate}
\item[(v)] $2g(B) -2 = (14-10) \cdot 5^2 -8 \cdot 5 \cdot 2 -2 \cdot 2^2 = 100-80-8=12$. 
\item[(vi)] $-K_W \cdot B = (14-10) \cdot 5 \cdot 6 -4 \cdot 11 \cdot 2 -2 \cdot 2^2 = 120-88-8 =24$. 
\item[(vii)] $D^3 = -(14-10) + 40 -3 \cdot 24 = -4+40-72=-36$. 
\end{enumerate}
If $(g, \alpha, \beta, -K_W^3) = (7, 2, 1, 10)$, then the following hold. 
\begin{enumerate}
\item[(v)] $2g(B) -2 = (14-10) \cdot 2^2 -8 \cdot 2 \cdot 1 -2 \cdot 1^2 = 16-16 -2=-2$. 
\item[(vi)] $-K_W \cdot B = (14-10) \cdot 2 \cdot 3 -4 \cdot 5 \cdot 1 -2 \cdot 1^2 = 24-20-2 =2$. 
\item[(vii)] $D^3 = -(14-10) + 10 -3 \cdot 2 = -4+10-6=0$. 
\end{enumerate}
If $(g, \alpha, \beta, -K_W^3) = (6, 9, 2, 32)$, then the following hold. 
\begin{enumerate}
\item[(v)] $2g(B) -2 = (12-10) \cdot 9^2 -8 \cdot 9 \cdot 2 -2 \cdot 2^2 = 162 -144 -8 =10$. 
\item[(vi)] $-K_W \cdot B = (12-10) \cdot 9 \cdot 10 -4 \cdot 19 \cdot 2 -2 \cdot 2^2 = 180-152-8 =20$. 
\item[(vii)] $D^3 = -(12-10) + 32 -3 \cdot 20 = -2+32-60=-30$. 
\end{enumerate}
If $(g, \alpha, \beta, -K_W^3) = (6, 5, 1, 22)$, then the following hold. 
\begin{enumerate}
\item[(v)] $2g(B) -2 = (12-10) \cdot 5^2 -8 \cdot 5 \cdot 1 -2 \cdot 1^2 = 50-40-2 =8$. 
\item[(vi)] $-K_W \cdot B = (12-10) \cdot 5 \cdot 6 -4 \cdot 11 \cdot 1 -2 \cdot 1^2 = 60-44-2 =14$. 
\item[(vii)] $D^3 = -(12-10) + 22 -3 \cdot 14 = -2+22-42=-22$. 
\end{enumerate}
If $(g, \alpha, \beta, -K_W^3) = (6, 4, 1, 8)$, then the following hold. 
\begin{enumerate}
\item[(v)] $2g(B) -2 = (12-10) \cdot 4^2 -8 \cdot 4 \cdot 1 -2 \cdot 1^2 = 32-32-2 =-2$. 
\item[(vi)] $-K_W \cdot B = (12-10) \cdot 4 \cdot 5 -4 \cdot 9 \cdot 1 -2 \cdot 1^2 = 40-36 -2 =2$. 
\item[(vii)] $D^3 = -(12-10) + 8 -3 \cdot 2 = -2+8-6=0$. 
\end{enumerate}
This completes the proof of Step \ref{s4-pt-flop}. 
\end{proof}
This completes the proof of Theorem \ref{t-pt-flop}. 
\end{proof}


\subsection{Existence of flopping curves and conics}\label{ss-pt-flopping}

\begin{prop}\label{p-pt-defect}
We use the same notation as in Theorem \ref{t-pt-flop}. 
Then $E^3 \neq (E^+)^3$. 
In particular, there exists a curve $\zeta$ on $Y$ such that  $K_Y \cdot \zeta =0$. 
\end{prop}

\begin{proof}
Suppose $(1=) E^3 = (E^+)^3$. 
Let us derive a contradiction. 
By $D \sim -\alpha K_{Y^+} -\beta E^+$ (Table \ref{table-2ray-pt} in Theorem \ref{t-pt-flop}), we have 
    \[
D^3 = ( -\alpha K_{Y^+} -\beta E^+)^3 \equiv -\beta^3 (E^+)^3 =-\beta^3  \mod \alpha. 
\]
By Table \ref{table-2ray-pt} in Theorem \ref{t-pt-flop}, 
the remaining cases are as follows: 
\begin{enumerate}
\renewcommand{\labelenumi}{(\roman{enumi})}
    \item $(\text{type of }\tau, g, \alpha, \beta, D^3)=(E_2, 7, 2, 1, 1)$.
    \item $(\text{type of }\tau, g, \alpha, \beta, D^3)=(E_1, 8, 5, 3, -42)$.
    \item $(\text{type of }\tau, g, \alpha, \beta, D^3)=(E_1, 12, 3, 4, -22)$.
    \item $\alpha =1$. 
\end{enumerate} 
We have 
\begin{eqnarray*}
D^3 &=& \alpha^3 (-K_Y)^3 - 3 \alpha^2 \beta (-K_Y)^2 \cdot E + 3 \alpha \beta^2 (-K_Y) \cdot E^2 -\beta^3 E^3\\
&=&  (2g-10)\alpha^3 - 12 \alpha^2 \beta -6 \alpha \beta^2 -\beta^3.
\end{eqnarray*}
For (i)-(iii), we get the following contradictions: 
\begin{enumerate}
\renewcommand{\labelenumi}{(\roman{enumi})}
\item $1 = D^3 = 4 \cdot 2^3 - 12 \cdot 2^2 \cdot 1 -6 \cdot 2 \cdot 1^2 -1^3 \equiv -1 \mod 4$. 
\item $-42 = D^3 = 6 \cdot 5^3 - 12 \cdot 5^2 \cdot 3 -6 \cdot 5 \cdot 3^2 -3^3  \equiv -7 \mod 10$. 
\item $-22  =D^3 = 14 \cdot 3^3 -12 \cdot 3^2 \cdot 4 -6 \cdot 3 \cdot 4^2 -4^3 \equiv -1 \mod 9$. 
\end{enumerate} 
Suppose (iv). 
By Table \ref{table-2ray-pt} in Theorem \ref{t-pt-flop}, we have $\beta =1$ and 
$D^3 \in \{0, 1\}$. We get 
\[
\{0, 1\} \ni  D^3 = (2g-10)- 12 -6 -1 = 2g -29,    
\]
which is absurd, because we have $g \leq 13 $ by Table \ref{table-2ray-pt} in Theorem \ref{t-pt-flop}. 
\qedhere





\end{proof}

\begin{cor}\label{c-pt-not-ample}
Let $X \subset \P^{g+1}$ be an anti-canonically embedded Fano threefold with $\Pic\,X = \Z K_X$. 
Let $P$ be a point on $X$ and let $\sigma : Y \to X$ be the blowup at $P$. 
Then $-K_Y$ is not ample. 
\end{cor}

\begin{proof}
Suppose that $-K_Y$ is ample. 
Let us derive a contradiction. 
Since $-K_Y$ is ample, we obtain $2g -10 = (-K_Y)^3 >0$ (Proposition \ref{p-pt-basic}(1)), which implies $g \geq 6$. 
If there exists a line $L$ on $X$ passing through $P$, 
then we would get the following contradiction for the proper transform $L_Y :=\sigma_*^{-1}L$ of $L$ on $Y$:   
\[
K_Y \cdot L_Y = (\sigma^*K_X +2E) \cdot L_Y = K_X \cdot L +2E \cdot L_Y \geq  -1 + 2 >0. 
\]
Hence there is no line on $X$ passing through $P$. 
Then we may apply Proposition \ref{p-pt-defect}. 
There exists a curve $\zeta$ on $Y$ satisfying $K_Y \cdot \zeta =0$, which is absurd. 
\end{proof}

\begin{prop}\label{p-max-rank-pt}
We use the same notation as in Theorem \ref{t-pt-flop}. 
Then the following hold. 
\begin{enumerate}
\item 
If $g \in \{ 6, 7, 8\}$, then $\rho : H^0(Y, -K_Y) \to H^0(E, -K_Y|_E)$ is injective. 
\item 
If $g \geq 8$, then $\rho : H^0(Y, -K_Y) \to H^0(E, -K_Y|_E)$ is surjective. 
\end{enumerate}
\end{prop}

\begin{proof}
Note that $h^0(Y, -K_Y -E) = h^0(Y^+, -K_{Y^+}-E^+)$. 
Let us show the following: 
\begin{enumerate}
\item[(1)'] If $g \in \{ 6,7, 8\}$, then $H^0(Y, -K_Y-E)=0$. 
\item[(2)'] If $g \geq 8$, then $h^0(Y, -K_Y-E)=g-8$. 
\end{enumerate}
We now finish the proof of Proposition \ref{p-max-rank-pt} by assuming (1)' and (2)'. 
Clearly, (1)' implies (1). 
By $h^0(Y, -K_Y) \geq g-2$ and $h^0(E, -K_Y|_E)=6$  (Proposition \ref{p-pt-basic}(1)), 
(2)' implies (2). 
It suffices to show (1)' and (2)'.

Let us show (1)'. 
Assume $g \in \{ 6, 7, 8\}$. 
Then $H^0(Y, -K_Y-E) =0$, as otherwise we would get the following contradiction: 
\[
2 \geq \kappa(Y^+, D) = \kappa(Y^+, -\alpha K_{Y^+} -\beta E^+) 
\geq \kappa(Y^+, -(\alpha - \beta) K_{Y^+}) = 3, 
\]
where the last equality holds by $\alpha >\beta$ (Table \ref{table-2ray-pt} in Theorem \ref{t-pt-flop}). 
Thus (1)' holds. 

Let us show (2)'. 
If $g = 8$, then (1)' implies (2)' 
by $h^0(Y, -K_Y) \geq g-2$ and $h^0(E, -K_Y|_E)=6$  (Proposition \ref{p-pt-basic}(1)). 
If $g =9$, then the linear equivalence $D \sim -K_{Y^+} -E^+$  
(Table \ref{table-2ray-pt}) implies 
$h^0(Y^+, -K_{Y^+}-E^+) = h^0(Y^+, D)=1$, because $\tau$ is of type $E$ and $D = \Ex(\tau)$. 
If $g=10$, then 
$\tau$ is of type $D$ and $D \sim -K_{Y^+}-E^+$ (Table \ref{table-2ray-pt}), and hence 
$h^0(Y^+, -K_{Y^+}-E^+) = h^0(Y^+, D) = h^0(\P^1, \MO_{\P^1}(1)) =2$, 
where we used the fact that $W=\P^1$ and $\MO_{Y^+}(D) \simeq \tau^*\MO_{\P^1}(1)$. 
If $g=11$, then 
$\tau$ is of type $C$ and $D \sim -K_{Y^+}-E^+$ (Table \ref{table-2ray-pt}),  and hence 
$h^0(Y^+, -K_{Y^+}-E^+) = h^0(Y^+, D) = h^0(\P^2, \MO_{\P^2}(1)) =3$, 
where we used the fact that $W=\P^2$ and $\MO_{Y^+}(D) \simeq \tau^*\MO_{\P^2}(1)$. 
Hence we may assume that $g=12$ or $g=13$. 
We only treat the case when $g=13$, as the proofs are very similar. 
In this case, we have $D \sim -2K_{Y^+} -3E^+$ (Table \ref{table-2ray-pt}). 
For a curve $\ell_{\tau}$ on $Y^+$ contracted by $\tau$ (i.e., a fibre of the induced $\P^1$-bundle structure $D \to \tau(D)$), we obtain 
\[
K_{Y^+} \cdot \ell_{\tau} =D \cdot \ell_{\tau} = -1.
\]
Hence $-1 = D\cdot \ell_{\tau} = ( -2K_{Y^+} -3E^+) \cdot \ell_{\tau} = 2 -3E^+ \cdot \ell_{\tau}$, which implies $E^+ \cdot \ell_{\tau}=1$. 
In particular, we get $(-K_{Y^+} -E^+) \cdot \ell_{\tau}=0$, 
and hence we can write $-K_{Y^+} -E^+ \sim \tau^*F$ for some Cartier divisor $F$ on $W$. 
By $D -K_{Y^+} \sim -3K_{Y^+} -3E^+$, we get $-K_W \sim \tau_*(D-K_{Y^+}) \sim \tau_*(-3K_{Y^+} -3E^+) \sim 3F$. 
Hence we obtain $F \sim H_W$ and $-K_{Y^+} -E^+ \sim \tau^*H_W$, where $H_W$ is an ample generator of the quadric hypersurface $W \subset \P^4$ (Table \ref{table-2ray-pt}). 
Therefore, 
\[
h^0(Y, -K_Y-E) =h^0(Y^+, -K_{Y^+}-E^+) = h^0(W, H_W) = h^0(\P^4, \MO_{\P^4}(1)) = 5.  
\]
Thus (2)' holds. 
\end{proof}

\begin{thm}\label{t-pt-to-conic}
Let $X \subset \P^{g+1}$ be an anti-canonically embedded Fano threefold with $\Pic\,X = \Z K_X$ and $g \geq 8$. 
If $P$ is a general point on $X$, then there exists a conic on $X$ passing through $P$. 
\end{thm}

\begin{proof}
By Proposition \ref{p-line-Hilb2}, 
there exists no line on $X$ passing through $P$. 
Let $\sigma : Y \to X$ be the blowup at $P$. 
Then $-K_Y$ is semi-ample and big (Proposition \ref{p-pt-basic}). 
Let $\psi : Y \to Z$ be the Stein factorisation of $\varphi_{|-mK_Y|}$ 
for some $m \in \Z_{>0}$ such that $|-mK_Y|$ is base point free. 
As $-K_Y$ is big, $\psi$ is birational. 
By Corollary \ref{c-pt-not-ample}, $\psi$ is not an isomorphism, i.e., 
$\dim \Ex(\psi)=2$ or $\dim \Ex(\psi)=1$. 
If $\Ex (\psi) =2$, then there exists a conic on $X$ passing through $P$ 
(Proposition \ref{p-pt-div-cont}). 
We may assume that $\Ex (\psi) =1$. 
By Proposition \ref{p-pt-defect}, 
there exists a curve $\zeta$ on $Y$ such that $-K_Y \cdot \zeta =0$. 
It is enough to show that $E \cdot \zeta =1$. 
Indeed, this implies 
\[
0 = K_Y \cdot \zeta = (\sigma^*K_X+2E) \cdot \zeta = K_X \cdot \sigma(\zeta) +2,
\]
and hence $\sigma(\zeta)$ is a conic on $X$.


Let us show  that $E \cdot \zeta =1$. 
Set $z := \psi(\zeta)$ and $E_Z := \psi(E)$. 
Since $E$ is $\psi$-ample (Lemma \ref{l-beta=r}), 
we have $E \cdot \zeta \geq 1$. 
Suppose that  $E \cdot \zeta \geq 2$. 
It suffices to derive a contradiction. 
Since $H^0(Y, -K_Y) \to H^0(E, -K_Y|_E)$ is surjective (Proposition \ref{p-max-rank-pt}), 
it follows from $E \simeq \P^2$ and the very ampleness of $|-K_Y|_E|$ that $\psi|_E : E \to E_Z$ is an isomorphism. 
By $z = \psi(\zeta)$, the composite morphism $E\cap \zeta \hookrightarrow E \xrightarrow{\psi|_E, \simeq} E_Z$ factors through $z$: 
\[
\begin{tikzcd}
E \arrow[d, "\psi|_E"', "\simeq"] & E \cap \zeta \arrow[l, hook'] \arrow[d]\\  
E_Z & z. \arrow[l, hook']
\end{tikzcd}
\]
Then we obtain two closed immersions: 
\[
E_Z \hookleftarrow z \hookleftarrow (\psi|_E)(E \cap \zeta), 
\]
which is a contradiction, because $h^0(z, \MO_z)=1$ and 
\[
h^0((\psi|_E)(E \cap \zeta), \MO_{(\psi|_E)(E \cap \zeta)}) = h^0(E \cap \zeta, \MO_{E \cap \zeta}) = E \cdot \zeta \geq 2.
\]
\end{proof}

\section{Two-ray game for conics}\label{s-conic}

Let $X \subset \P^{g+1}$ be an anti-canonically embedded Fano threefold with $\Pic\,X = \Z K_X$. 
The main results of this section are as follows: 
\begin{enumerate}
\item $3 \leq g \leq 12$ (Theorem \ref{t-g-bound}). 
\item Let $\Gamma$ be a conic of $X$ and let $\sigma: Y \to X$ be the blowup along $\Gamma$. 
Then $-K_Y$ is not ample (Corollary \ref{c-conic-not-ample}). 
\end{enumerate}
Both the strategies were explained in Subsection \ref{ss-2ray-intro}.

\subsection{Basic properties}

For an anti-canonically embedded Fano threefold $X \subset \P^{g+1}$, 
a {\em conic} on $X$ is a curve $\Gamma$ on $X$  such that $-K_X \cdot \Gamma =2$. 
Note that a conic is always assumed to be an integral scheme under our terminology.







\begin{lem}\label{l-conic-basic}
Let $X \subset \P^{g+1}$ be an anti-canonically embedded Fano threefold with 
$\Pic\,X = \Z K_X$ and $g \geq 5$. 
Let $\Gamma$ be a conic on $X$. 
Then the following hold. 
\begin{enumerate}
\item $\langle \Gamma \rangle = \P^2$, i.e., the smallest linear subvariety $\langle \Gamma \rangle$ of $\P^{g+1}$ containing $\Gamma$ is two-dimensional. 
\item $\Gamma$ is isomorphic to $\P^1$. 
\item 
The scheme-theoretic equality $\langle \Gamma \rangle \cap X = \Gamma$ holds. 
\end{enumerate}

\end{lem}

\begin{proof}
Let us show (1). 
Fix distinct closed points $P, Q \in \Gamma$. 
Since $\Gamma$ is not a line, we have $\overline{PQ} \neq \Gamma$, where $\overline{PQ}$ denotes the line 
passing through $P$ and $Q$. 
Pick $R \in \Gamma \setminus \overline{PQ}$. 
Then $\Gamma$ must be contained in the plane $\overline{PQR}$ containing $P, Q, R$, 
as otherwise we would get $\Gamma \not\subset H$ for a general hyperplane $H \subset \P^N$, 
containing $\overline{PQR}$, 
which leads to a contradiction: $2=\Gamma \cdot H \geq \#(\Gamma \cap H) \geq 3$. 
Thus (1) holds. 
Let us show (2). 
By (1), we have a plane $V = \P^2 \subset \P^N$ containing $\Gamma$. 
Then $\Gamma$ is a conic on $V=\P^2$ (i.e., a prime divisor on $\P^2$ of degree two), 
which is isomorphic to $\P^1$. 
Thus (2) holds.

Let us show (3). 
By (1), we have $V :=  \langle \Gamma \rangle = \P^2$.  
Since $X$ is an intersection of quadrics, we have $X=  \bigcap_{X \subset Q} Q$, where $Q$ runs over all the quadric hypersurfaces of $\P^{g+1}$ containing $X$. 
Hence 
\[
\Gamma \subset X \cap V = \left( \bigcap_{X \subset Q} Q\right) \cap V = \bigcap_{X \subset Q} (Q \cap V) 
\subset V.  
\]
Since $X$ does not contain any plane by adjunction formula, we get $X \cap V \subsetneq V$. 
In particular, there exists a quadric hypersurface $Q_1 \subset \P^{g+1}$ such that $X \subset Q_1$ and $Q_1 \cap V\subsetneq V$. 
Then 
\[
\Gamma \subset X \cap V = \bigcap_{X \subset Q} (Q \cap V) \subset Q_1 \cap V \subsetneq V. 
\]
Since both $\Gamma$ and $Q_1 \cap V$ are effective Cartier divisors on $V=\P^2$ of degree two, 
they must coincide. Thus (3) holds. 
\end{proof}

\begin{prop}\label{p-conic-basic}
Let $X \subset \P^{g+1}$ be an anti-canonically embedded Fano threefold with $\Pic\,X = \Z K_X$ and $g \geq 5$. 
Let $\Gamma$ be a conic on $X$. 
Let $\sigma : Y \to X$ be the blowup along $\Gamma$ and set $E := \Ex(\sigma)$. 
Then the following hold. 
\begin{enumerate}
    \item 
    \begin{enumerate}
    \item $(-K_Y)^3 = 2g -8$. 
    \item $(-K_Y)^2 \cdot E = 4$. 
    \item $(-K_Y) \cdot E^2 = -2$. 
    \item $E^3 = 0$. 
\end{enumerate}
    \item $|-K_Y|$ is base point free and $-K_Y$ is big. 
    \item 
    $h^0(E, -K_Y|_E) = 6$, 
    $h^0(Y, -K_Y) \geq g-1$, and 
    $h^0(Y, -K_Y -E) \geq g-7$. 
\item One of the following holds: 
\[
N_{\Gamma/X} \simeq \MO_{\Gamma} \oplus \MO_{\Gamma}, \quad 
N_{\Gamma/X} \simeq \MO_{\Gamma}(1) \oplus \MO_{\Gamma}(-1), \quad 
N_{\Gamma/X} \simeq \MO_{\Gamma}(2) \oplus \MO_{\Gamma}(-2). 
\]
\item 
The following hold. 
\begin{enumerate}
\item 
If $N_{\Gamma/X} \simeq \MO_{\Gamma} \oplus \MO_{\Gamma}$ or 
$N_{\Gamma/X} \simeq \MO_{\Gamma}(1) \oplus \MO_{\Gamma}(-1)$, then $-K_Y|_E$ is ample. 
\item 
If $N_{\Gamma/X} \simeq \MO_{\Gamma}(2) \oplus \MO_{\Gamma}(-2)$, 
then $-K_Y|_E$ is not ample but nef and big. 
In this case, for a curve $B$ on $E$, 
$(-K_Y|_E) \cdot B=0$ if and only if $B$ is the section of the $\P^1$-bundle 
$E \to \Gamma$ whose self-intersection number is negative. 
\end{enumerate}
    \item Let $C \subset Y$ be a curve such that $-K_Y \cdot C =0$. 
Then one of the following holds. 
\begin{enumerate}
\item $\sigma(C)$ is a line such that $C \cdot E = 1$  and $C \not\subset E$. 
\item $C \cdot E \geq 2$ and $C \not\subset E$. 
\item 
$N_{\Gamma/X} \simeq \MO_{\Gamma}(2) \oplus \MO_{\Gamma}(-2)$ and 
$C$ is the section of the $\P^1$-bundle $\sigma|_E : E \to \Gamma$ whose self-intersection number, as a divisor on $E$, is negative. 
\end{enumerate}
\end{enumerate}
\end{prop}

\begin{proof}
The assertion (1) follows from Lemma \ref{l-blowup-formula}(2). 
Let us show (2). 
By the scheme-theoretic equalities $\sigma^{-1}(\Gamma)=E$ and  
$\langle \Gamma \rangle \cap X = \Gamma$  (Lemma \ref{l-conic-basic}), 
the complete linear system 
$|-K_Y| = |-\sigma^*K_X -E|$ is base point free.  
Since $-K_Y$ is nef, it follows from $(-K_Y)^3 = 2g -8>0$ that $-K_Y$ is big. 
Thus (2) holds. 
The assertion (3) holds by (1), (2), and Lemma \ref{l-bound--K-E}.

Let us show (4) and (5). 
Note that we have $N_{\Gamma/X} \simeq \MO_{\Gamma}(m) \oplus \MO_{\Gamma}(-m)$ for some $m \in \Z_{\geq 0}$ 
(Lemma \ref{l-blowup-formula}(2)), and hence $E \simeq \F_{2m}$. 
Fix a fibre $f$ of the $\P^1$-bundle $\sigma|_E : E \to \Gamma$ and 
let $s$ be a section of $\sigma|_E$ with $s^2 =-2m$. 
We can write $-K_Y|_E \sim a s + b f$ for some $a, b \in \Z$. 
By $-K_Y|_E \cdot f =1$, we have $a=1$, and hence $-K_Y|_E \sim s + b f$. 
Since $-K_Y|_E$ is nef, we get $0 \leq (-K_Y|_E) \cdot s = s^2 +b = -2m +b$, which implies 
$2m \leq b$.   
We obtain 
\[
4 \overset{{\rm (1)}}{=} (-K_Y)^2 \cdot E =( -K_Y|_E)^2 = (s+bf)^2 = s^2 +2b = -2m + 2b \geq -2m +2 \cdot (2m). 
\]
Therefore, it holds that $0 \leq m \leq 2$, which implies (4). 
The equality $-2m + 2b = 4$ deduces the following. 
\begin{itemize}
\item If $E \simeq \F_0$, i.e., $N_{\Gamma/X} \simeq \MO_{\Gamma} \oplus \MO_{\Gamma}$, then $-K_Y|_E \sim s + 2f$. 
\item If $E \simeq \F_2$, i.e., 
$N_{\Gamma/X} \simeq \MO_{\Gamma}(1) \oplus \MO_{\Gamma}(-1)$, then $-K_Y|_E \sim s + 3f$. 
\item If $E \simeq \F_4$, i.e., 
$N_{\Gamma/X} \simeq \MO_{\Gamma}(2) \oplus \MO_{\Gamma}(-2)$, then $-K_Y|_E \sim s + 4f$. 
\end{itemize}
Then we can check, by Kleiman's criterion for ampleness, that (5) holds.

Let us show (6). 
If $C$ is contained in $E$, then (c) holds by (5). 
Assume $C\not\subset E$. In particular, $\sigma(C)$ is a curve and $-K_X \cdot \sigma(C) >0$. 
By $K_Y \cdot C =0$ and $K_Y = \sigma^*K_X + E$, we obtain $E \cdot C >0$. 
If $C \cdot E =1$ (resp. $C \cdot E \geq 2$), then (a) (resp. (b)) holds. 
Thus (6) holds. 
\end{proof}





\subsection{Divisorial contractions}

\begin{prop}\label{p-conic-div-cont}
Let $X \subset \P^{g+1}$ be an anti-canonically embedded Fano threefold with $\Pic\,X = \Z K_X$ and $g \geq 5$. 
Let $\Gamma$ be a conic on $X$. 
Let $\sigma : Y \to X$ be the blowup along $\Gamma$ and set $E := \Ex(\sigma)$. 
Since $|-K_Y|$ is base point free and $-K_Y$ is big (Proposition \ref{p-conic-basic}), we have the birational morphism 
\[
\psi : Y \to Z 
\]
which is  the Stein factorisation of $\varphi_{|-K_Y|}$. 
Assume that $\dim \Ex(\psi) = 2$. Set $D := \Ex(\psi)$ and $B := \psi(D)$. 
Recall that $D$ is a prime divisor and we have 
\[
D \sim -\alpha K_Y -\beta E
\]
for some $\alpha, \beta \in \Z_{>0}$ (Proposition \ref{p-div-cont1}). 
Then the following hold. 
\begin{enumerate}
\item 
One of the following holds. 
\begin{enumerate}
\item $(g, \alpha, \beta)= (5, 2, 1)$. 
\item $(g, \alpha, \beta)= (6, 1, 1)$. 
\item $(g, \alpha, \beta)= (5, 4, 2)$. 
\item $(g, \alpha, \beta)= (6, 2, 2)$. 
\item $(g, \alpha, \beta)= (8, 1, 2)$. 
\end{enumerate}
\item 
If $g \geq 6$, then there exists a line $L$ on $X$ such that $\Gamma \cap L \neq \emptyset$. 
\end{enumerate}
\end{prop}

\begin{proof}
Let us show (1). 
By Proposition \ref{p-div-cont1}(5)(8) and Proposition \ref{p-conic-basic}(1), we have $\beta \in \{1, 2\}$ and 
$(2g-8) \alpha = 4\beta$. Then 
\[
 (g-4)\alpha =2\beta. 
\]
\begin{itemize}
    \item If $\beta =1$, then $(g, \alpha, \beta) \in \{(5, 2, 1), (6, 1, 1)\}$. 
    \item If $\beta =2$, then $(g, \alpha, \beta) \in \{(5, 4, 2), (6, 2, 2), (8, 1, 2)\}$. 
\end{itemize}
Thus (1) holds.

Let us show (2). 
If there exists a curve $\zeta$ on $Y$ such that $K_Y \cdot \zeta =0$ and $E \cdot \zeta =1$, 
then we get $-K_X \cdot \sigma(\zeta)=1$ by $K_Y = \sigma^*K_X+E$. 
Suppose that $E \cdot \zeta \neq 1$ 
for every curve $\zeta$ on $Y$ such that $K_Y \cdot \zeta =0$. 
It suffices to derive a contradiction. 
We use the same notation as in Proposition \ref{p-div-cont2}. 
By $g \geq 6$ and $\beta =1$ (Proposition \ref{p-div-cont2}(1)), 
(1) implies  $(g, \alpha, \beta) = (6, 1, 1)$. 
It follows from Proposition \ref{p-div-cont2}(5) that $\omega_D^2 \in 4\Z$. 
On the other hand, recall that $D \sim - K_Y -E$, 
$-K_Y \cdot E^2 =-2$, and $E^3=0$ (Proposition \ref{p-conic-basic}(1)). 
Hence we get  
\[
\omega_D^2 = (K_Y+D)^2 \cdot D = (-E)^2 \cdot (-K_Y -E) 
=-K_Y \cdot E^2 -E^3 = -2 \not\equiv 0 \mod 4, 
\]
which is a contradiction. 
Thus (2) holds. 
\end{proof}

\subsection{Flops}

\begin{thm}\label{t-conic-flop}
Let $X \subset \P^{g+1}$ be an anti-canonically embedded Fano threefold with $\Pic\,X = \Z K_X$ and $g \geq 5$. 
Let $\Gamma$ be a conic on $X$. 
Let $\sigma : Y \to X$ be the blowup along $\Gamma$ and set $E := \Ex(\sigma)$. 
Since $|-K_Y|$ is base point free and $-K_Y$ is big (Proposition \ref{p-conic-basic}), we have the birational morphism 
\[
\psi : Y \to Z 
\]
which is  the Stein factorisation of $\varphi_{|-K_Y|}$. 
Assume that $\dim \Ex(\psi) = 1$ or $\psi$ is an isomorphism. 
\begin{itemize}
    \item If $\dim \Ex(\psi) = 1$, then let $\psi^+ : Y^+ \to Z$ be the flop of $\psi$, set $E^+$ to be the proper transform of $E$ on $Y^+$, and 
let $\tau : Y^+ \to W$ be the contraction of the $K_{Y^+}$-negative extremal ray. 
    \item If $\psi$ is an isomorphism, then we set $Z :=Y$, $Y^+ :=Y$, $\psi := {\rm id}$, 
$\psi^+ := {\rm id}$, $E^+ :=E$, and let $\tau : Y^+ \to W$ be the contraction of the $K_{Y^+}$-negative extremal ray not corresponding to $\sigma$. 
\end{itemize}
\[
\begin{tikzcd}
Y \arrow[d, "\sigma"'] \arrow[rd, "\psi"]& & Y^+ \arrow[ld, "\psi^+"'] \arrow[d, "\tau"]\\
X & Z & W
\end{tikzcd}
\]
Then $g$, the type of $\tau$, $D$, and $D^3$ satisfy one of the possibilities in 
Table \ref{table-2ray-conic}, where we use the following notation. 
\begin{itemize}
\item 
If $\tau$ is of type $E$, then $D := \Ex(\tau)$. 
If $\tau$ is of type $C$ or $D$, then $D$ is the pullback of an ample Cartier divisor which is a generator of $\Pic\,W$. 
The divisors below $D$ means the linear equivalence, e.g., 
if $g=10$, then $D \sim -K_{Y^+} -E^+$. 
\item If $\tau$ is of type $C$, then $\Delta$ denotes the discriminant bundle. 
\item If $\tau$ is of type $E_1$ or $E_2$, 
then let $r_W$ be the index of the Fano threefold $W$. 
\item If $\tau$ is of type $E_1$, then we set 
$B := \tau(D)$ and  $g(B)$ denotes the genus of $B$. 
\end{itemize}
\begin{table}[h]
\caption{Two-ray game for conics}\label{table-2ray-conic}
     \centering
{\renewcommand{\arraystretch}{1.35}%
      \begin{tabular}{|c|c|c|c|c|}
      \hline
$g$  & type of $\tau$ & $D$ & $D^3$ & Other properties\\      \hline
$5$  & $E_2$ & $-4K_{Y^+}-E^+$ & $1$ & $(r_W, -K_W^3)=(1, 10)$\\           \hline
$6$  & $E_2$ & $-2K_{Y^+}-E^+$ & $1$ & $(r_W, -K_W^3)=(1, 12)$\\           \hline
$8$  & $E_2$ & $-K_{Y^+}-E^+$ & $1$ & $(r_W, -K_W^3)=(1, 16)$\\           \hline
$9$  & $D_1$ & $-K_{Y^+}-E^+$ & $0$ &$(-K_{Y^+})^2 \cdot D =6$ \\           \hline
$10$  & $C_1$  & $-K_{Y^+}-E^+$ & $0$ &$\deg \Delta =4$ \\           \hline
$5$  & $E_1$ & $-4K_{Y^+}-E^+$ & $0$ & $(r_W, -K_W^3, -K_W\cdot B, g(B))=(1, 8, 2, 0)$\\           \hline
$5$  & $E_1$ & $-5K_{Y^+}-E^+$ & $-22$ & $(r_W, -K_W^3, -K_W\cdot B, g(B))=(1, 22, 14, 5)$\\           \hline
$5$  & $E_1$ & $-9K_{Y^+}-2E^+$ & $-30$ & $(r_W, -K_W^3, -K_W\cdot B, g(B))=(2, 32, 20, 6)$\\           \hline
$6$  & $E_1$ & $-2K_{Y^+}-E^+$ & $0$ & $(r_W, -K_W^3, -K_W\cdot B, g(B))=(1, 10, 2, 0)$\\           \hline
$6$  & $E_1$ & $-5K_{Y^+}-2E^+$ & $-36$ & $(r_W, -K_W^3, -K_W\cdot B, g(B))=(2, 40, 24, 7)$\\           \hline
$7$  & $E_1$ & $-3K_{Y^+}-2E^+$ & $-6$ & $(r_W, -K_W^3, -K_W\cdot B, g(B))=(2, 24, 8, 0)$\\           \hline
$7$  & $E_1$ & $-5K_{Y^+}-3E^+$ & $-42$ & $(r_W, -K_W^3, -K_W\cdot B, g(B))=(3, 54, 30, 7)$\\           \hline
$8$  & $E_1$ & $-K_{Y^+}-E^+$ & $0$ & $(r_W, -K_W^3, -K_W\cdot B, g(B))=(1, 14, 2, 0)$\\           \hline
$11$  & $E_1$ & $-3K_{Y^+}-4E^+$ & $-22$ & $(r_W, -K_W^3, -K_W\cdot B, g(B))=(4, 64, 24, 0)$\\           \hline
$12$  & $E_1$ & $-2K_{Y^+}-3E^+$ & $-16$ &  $(r_W, -K_W^3, -K_W\cdot B, g(B))=(3, 54, 18, 0)$\\           \hline
      \end{tabular}}
    \end{table}
\end{thm}


\begin{proof}
The proof is identical to 
that of Theorem \ref{t-pt-flop} 
after replacing $g$ by 
$\widetilde g:= g+1$. 
We here explain how to apply the same argument as in Theorem \ref{t-pt-flop}. 
By Lemma \ref{l-blowup-formula2}(2) and Proposition \ref{p-conic-basic}
we have 
\begin{itemize}
    \item $(-K_Y)^3 = (-K_{Y^+})^3 = 2\widetilde{g} -10$,
    \item $(-K_Y)^2 \cdot E = (-K_{Y^+})^2 \cdot E^+ =4$, and 
    \item $(-K_Y) \cdot E^2 = (-K_{Y^+}) \cdot (E^+)^2=-2$. 
\end{itemize}
Note that these equations coincide, after replacing $\widetilde g$ by 
$g$, with the ones in Proposition \ref{p-pt-basic}(1). 
Let $D$ be a Cartier divisor on $Y^+$ such that  the linear equivalence  
\[
D \sim -\alpha K_{Y^+} - \beta E^+ 
\]
holds for some $\alpha, \beta \in \Z$. 
Then the following hold: 
\begin{align*}
D^2 \cdot (-K_{Y^+})
&= (\alpha (-K_{Y^+}) -\beta E^+)^2 \cdot (-K_{Y^+})\\
&= \alpha^2 (-K_{Y^+})^3  -2 \alpha \beta (-K_{Y^+})^2 \cdot E^+ + \beta^2 (-K_{Y^+}) \cdot (E^+)^2\\
&= (2\widetilde{g}-10) \alpha^2 -8 \alpha \beta -2 \beta^2\\
D \cdot (-K_{Y^+})^2 & = (\alpha (-K_{Y^+}) -\beta E^+) \cdot (-K_{Y^+})^2\\
&= \alpha (-K_{Y^+})^3  - \beta (-K_{Y^+})^2 \cdot E^+\\
&= (2\widetilde{g}-10) \alpha -4\beta. 
\end{align*}

In what follows, we only give a sketch of a proof for the case when $\tau :Y^+ \to W$ is of type $E_1$. 
In this case, the following hold. 
\begin{enumerate}
\renewcommand{\labelenumi}{(\roman{enumi})}
\item $W$ is a Fano threefold with $\rho(W) =1$ and $\beta =r_W$ for the index $r_W$ of $W$ (Lemma \ref{l-beta=r}). 
In particular, $1 \leq \beta \leq 4$. Fix a Cartier divisor $H_W$ on $W$ with $-K_W \sim \beta H_W$. 
\item $\alpha+1 = \beta \gamma$ for some $\gamma \in \Z_{>0}$  (Lemma \ref{l-beta=r}). 
\item $w := \frac{1}{2}\beta H_W^3 = (\widetilde g-5)\gamma^2 - 4\gamma  -1$ (Remark \ref{r-E1-deg-genus}). 
\item $h^0(Y^+, -K_{Y^+}-E^+) = h^0(Y, -K_Y-E) \geq \widetilde g-8$ (Proposition \ref{p-conic-basic}). 
\item $2g(B) -2 =(2\widetilde g-10)\alpha^2 - 8 \alpha \beta  -2 \beta ^2$ (Lemma \ref{l-E1-deg-genus}). 
\item $-K_W \cdot B = (2\widetilde g-10)\alpha(\alpha+1) 
- 4(2\alpha +1)\beta  -2\beta^2$  (Lemma \ref{l-E1-deg-genus}).
\item $D^3 = -(2\widetilde g-10)+ (-K_W)^3 -3 (-K_W) \cdot B$ (Lemma \ref{l-D^3}). 
\end{enumerate}
Note that these equations coincide, after replacing $\widetilde g$ by 
$g$,  with the ones in Step \ref{s4-pt-flop} in the proof of Theorem \ref{t-pt-flop}. 
Then the remaining argument is identical to that of Theorem \ref{t-pt-flop}. 
For example, if $(\alpha, \beta) = (3, 4)$, 
then we get $(-K_W)^3 =64$ and $\widetilde{g}=12$ (i.e., $g =11$), 
which imply $(-K_W \cdot B, g(B), D^3) = (24, 0, -22)$ by (v)-(vii). 
\qedhere

\end{proof}

\begin{thm}\label{t-g-bound}
Let $X \subset \P^{g+1}$ be an anti-canonically embedded Fano threefold with $\Pic\,X = \Z K_X$. 
Then $3 \leq g \leq 12$. 
\end{thm}


\begin{proof}
By $g \geq 3$ (Proposition \ref{p-RR}(3)), 
we may assume that $g \geq 8$. 
In this case, we can find a conic $\Gamma$ on $X$ (Theorem \ref{t-pt-to-conic}). 
Take the blowup $\sigma: Y \to X$ along $\Gamma$. 
Then $|-K_Y|$ is base point free and big (Proposition \ref{p-conic-basic}(2)). 
Let $\psi : Y \to Z$ be the Stein factorisation of $\varphi_{|-K_Y|}$. 
Then we get $g \leq 12$ by Proposition \ref{p-conic-div-cont} ($\Ex(\psi) =2$) and Theorem \ref{t-conic-flop} ($\Ex(\psi) \neq 2$). 
\end{proof}

\begin{cor}\label{c-g-bound}
Let $X$ be a Fano threefold  with $\Pic\,X = \Z K_X$. 
Set $g := \frac{1}{2}(-K_X)^3+1$. 
Then $2 \leq g \leq 12$. 
\end{cor}

\begin{proof}
If $|-K_X|$ is not very ample, then we get $(-K_X)^3 \in \{ 2, 4\}$, i.e., $g \in \{ 2, 3\}$ 
\cite[Theorem 1.1]{TanI}. 
If $|-K_X|$ is very ample, then the assertion follows from Theorem \ref{t-g-bound}. 
\end{proof}

\begin{rem}\label{r-w-bound}
Under the same notation as in Remark \ref{r-E1-deg-genus}, 
we obtain $1 \leq w \leq 11$. 
Indeed, if $\beta \geq 2$, then this is explained in Remark \ref{r-E1-deg-genus}. 
If $\beta =1$, then $1 \leq w \leq 11$ holds by Corollary \ref{c-g-bound}. 
\end{rem}

\subsection{Existence of flopping curves}\label{ss-conic-flopping}

\begin{prop}\label{p-conic-defect}
We use the same notation as in Theorem \ref{t-conic-flop}. 
Then $E^3 \neq (E^+)^3$. 
In particular, there exists a curve $\zeta$ on $Y$ such that  $K_Y \cdot \zeta =0$. 
\end{prop}

\begin{proof}
Suppose $(0=) E^3 = (E^+)^3$ (Proposition \ref{p-conic-basic}). 
Let us derive a contradiction. 
By $D \sim -\alpha K_{Y^+} -\beta E^+$ (Table \ref{table-2ray-conic} in Theorem \ref{t-conic-flop}), we have 
    \[
D^3 = ( -\alpha K_{Y^+} -\beta E^+)^3 \equiv -\beta^3 (E^+)^3 =-\beta^3E^3 =0   \mod \alpha. 
\]
By Table \ref{table-2ray-conic} in Theorem \ref{t-conic-flop}, 
the remaining cases are as follows: 
\begin{enumerate}
\renewcommand{\labelenumi}{(\roman{enumi})}
\item  $(\text{type of }\tau, g, \alpha, \beta, D^3)=(E_1, 5, 4, 1, 0)$.
\item  $(\text{type of }\tau, g, \alpha, \beta, D^3)=(E_1, 6, 2, 1, 0)$.
\item  $(\text{type of }\tau, g, \alpha, \beta, D^3)=(E_1, 7, 3, 2, -6)$.
\item  $(\text{type of }\tau, g, \alpha, \beta, D^3)=(E_1, 12, 2, 3, -16)$.
    \item $\alpha =1$. 
\end{enumerate} 
The following holds (Proposition \ref{p-conic-basic}):
\begin{eqnarray*}
D^3 &=& \alpha^3 (-K_Y)^3 - 3 \alpha^2 \beta (-K_Y)^2 \cdot E + 3 \alpha \beta^2 (-K_Y) \cdot E^2 -\beta^3 E^3\\
&=&  (2g-8)\alpha^3 - 12 \alpha^2 \beta -6 \alpha \beta^2.
\end{eqnarray*}
For (i)--(iv), we get the following contradictions: 
\begin{enumerate}
\renewcommand{\labelenumi}{(\roman{enumi})}
\item $0 = D^3 = 2 \cdot 4^3 - 12 \cdot 4^2 \cdot 1 -6 \cdot 4 \cdot 1^2 \not\in 16\Z$. 
\item $0 = D^3 = 4 \cdot 2^3 - 12 \cdot 2^2 \cdot 1 -6 \cdot 2 \cdot 1^2 \not\in 8\Z$. 
\item $-6 = D^3 = 6 \cdot 3^3 - 12 \cdot 3^2 \cdot 2 -6 \cdot 3 \cdot 2^2 \in 18\Z$. 
\item $-16  =D^3 = 16 \cdot 2^3 -12 \cdot 2^2 \cdot 3 -6 \cdot 2 \cdot 3^2 \not\in 8\Z$. 
\end{enumerate} 
Suppose (v). 
By Table \ref{table-2ray-conic}, we have $\beta =1$ and $D^3 \in \{0, 1\}$. We get 
\[
\{0, 1\} \ni D^3 = (2g-8)- 12 -6 = 2g -26, 
\]
which contradicts $g \leq 12$ (Theorem \ref{t-g-bound}).
\qedhere 





\end{proof}

\begin{cor}\label{c-conic-not-ample}
Let $X \subset \P^{g+1}$ be an anti-canonically embedded Fano threefold with $\Pic\,X = \Z K_X$. 
Let $\Gamma$ be a conic on $X$ and let $\sigma : Y \to X$ be the blowup along $\Gamma$. 
Then $-K_Y$ is not ample. 
\end{cor}

\begin{proof}
Suppose that $-K_Y$ is ample. 
Let us derive a contradiction. 
Since $-K_Y$ is ample, we obtain $2g -8 = (-K_Y)^3 >0$ (Proposition \ref{p-conic-basic}(1)), which implies $g \geq 5$. 
Then we may apply Proposition \ref{p-conic-defect}. 
There exists a curve $\zeta$ on $Y$ satisfying $K_Y \cdot \zeta =0$, which is absurd. 
\end{proof}

\subsection{Hilbert schemes of conics}

{\cred 
Given an anti-canonically embedded Fano threefold $X \subset \P^{g+1}$, 
$\Hilb^{\conic}_{X}$  denotes the Hilbert scheme of conics \cite[Chapter 5]{FGI05}.} 
Although the following result  will not be used in this paper, we include it for a future reference.


\begin{prop}\label{p-conic-Hilb1}
Let $X \subset \P^{g+1}$ be an anti-canonically embedded Fano threefold with $\Pic\,X = \Z K_X$. 
Let $\Gamma$ be a conic on $X$. 
Recall that 
$N_{\Gamma/X} \simeq \MO_{\Gamma} \oplus \MO_{\Gamma}, 
N_{\Gamma/X} \simeq \MO_{\Gamma}(1) \oplus \MO_{\Gamma}(-1),$ or 
$N_{\Gamma/X} \simeq \MO_{\Gamma}(2) \oplus \MO_{\Gamma}(-2)$ (Proposition \ref{p-conic-basic}(4)). 
Then the following hold. 
\begin{enumerate}
\item 
If $N_{\Gamma/X} \simeq \MO_{\Gamma} \oplus \MO_{\Gamma}$ 
or $N_{\Gamma/X} \simeq \MO_{\Gamma}(1) \oplus \MO_{\Gamma}(-1)$, 
then  $\Hilb_X^{\conic}$ is smooth at $[\Gamma]$ and 
$\dim \MO_{\Hilb_X^{\conic}, [\Gamma]}=2$. 
\item If $N_{\Gamma/X} \simeq \MO_{\Gamma}(2) \oplus \MO_{\Gamma}(-2)$, 
then then $\dim \MO_{\Hilb^{\conic}_X, [\Gamma]} \in \{2, 3\}$ and 
one of the following holds. 
\begin{itemize}
\item $\Hilb^{\conic}_X$ is smooth at $[\Gamma]$ and $\dim \MO_{\Hilb^{\conic}_X, [\Gamma]} = 3$. 
\item $\Hilb^{\conic}_X$ is not smooth at $[\Gamma]$ and $\dim \MO_{\Hilb^{\conic}_X, [\Gamma]} = 2$. 
\end{itemize}
\end{enumerate}
In particular, either $\Hilb^{\conic}_X = \emptyset$ or 
$2 \leq \dim \Hilb^{\line}_X \leq 3$. 
\end{prop}

\begin{proof}
The proof is identical to that of Proposition \ref{p-line-Hilb1}. 
\qedhere

\end{proof}

\begin{rem}
Let $X \subset \P^{g+1}$ be an anti-canonically embedded Fano threefold with $\Pic\,X = \Z K_X$ and $g \geq 8$. 
By Theorem \ref{t-pt-to-conic}, the induced morphism $\Univ^{\conic}_X \to X$ is surjective, where $\Univ^{\conic}_X$ denotes  the universal family of $\Hilb_X^{\conic}$.  
\end{rem}

\section{Two-ray game for lines II}\label{s-line-II}

This section is the continuation of Section \ref{s-line-I}.

\subsection{Flops}


\begin{thm}\label{t-line-flop}
Let $X \subset \P^{g+1}$ be an anti-canonically embedded Fano threefold with $\Pic\,X = \Z K_X$  
and $g \geq 4$. 
Let $\Gamma$ be a line on $X$. 
Let $\sigma : Y \to X$ be the blowup along $\Gamma$ and set $E := \Ex(\sigma)$. 
Since $|-K_Y|$ is base point free and $-K_Y$ is big (Proposition \ref{p-line-basic}), we have the birational morphism 
\[
\psi : Y \to Z 
\]
which is  the Stein factorisation of $\varphi_{|-K_Y|}$. 
Assume that $\dim \Ex(\psi) = 1$ or $\psi$ is an isomorphism. 
\begin{itemize}
    \item If $\dim \Ex(\psi) = 1$, then let $\psi^+ : Y^+ \to Z$ be the flop of $\psi$, set $E^+$ to be the proper transform of $E$ on $Y^+$, and 
let $\tau : Y^+ \to W$ be the contraction of the $K_{Y^+}$-negative extremal ray. 
    \item If $\psi$ is an isomorphism, then we set $Z :=Y$, $Y^+ :=Y$, $\psi := {\rm id}$, 
$\psi^+ := {\rm id}$, $E^+ :=E$, and let $\tau : Y^+ \to W$ be the contraction of the $K_{Y^+}$-negative extremal ray not corresponding to $\sigma$. 
\end{itemize}
\[
\begin{tikzcd}
Y \arrow[d, "\sigma"'] \arrow[rd, "\psi"]& & Y^+ \arrow[ld, "\psi^+"'] \arrow[d, "\tau"]\\
X & Z & W
\end{tikzcd}
\]
Then $g$, the type of $\tau$, $D$, and $D^3$ satisfy one of the possibilities in 
Table \ref{table-2ray-line}, where we use the following notation. 
\begin{itemize}
\item 
If $\tau$ is of type $E$, then $D := \Ex(\tau)$. 
If $\tau$ is of type $C$ or $D$, then $D$ is the pullback of an ample Cartier divisor which is a generator of $\Pic\,W$. 
The divisors below $D$ means the linear equivalence, e.g., 
if $g=12$, then $D \sim -K_{Y^+} -2E^+$. 
\item If $\tau$ is of type $C$, then $\Delta$ denotes the discriminant bundle. 
\item If $\tau$ is of type $E_1$, 
then let $r_W$ be the index of the Fano threefold $W$. 
\item If $\tau$ is of type $E_1$, then we set 
$B := \tau(D)$ and  $g(B)$ denotes the genus of $B$. 
\end{itemize}
    \begin{table}[h]
\caption{Two-ray game for lines}\label{table-2ray-line}
     \centering
{\renewcommand{\arraystretch}{1.35}%
      \begin{tabular}{|c|c|c|c|c|}
      \hline
$g$  & type of $\tau$ & $D$ & $D^3$ & Other properties\\      \hline
$5$  & $C_1$ & $-2K_{Y^+}-E^+$ & ${\cred 0}$ & $\deg \Delta =7$\\           \hline
$7$  & $D_1$ & $-K_{Y^+}-E^+$ & ${\cred 0}$ & $(-K_{Y^+})^2 \cdot D =5$\\           \hline
$8$  & $C_1$ & $-K_{Y^+}-E^+$ & ${\cred 0}$ & $\deg \Delta =5$\\           \hline
$4$  & $E_1$ & $-4K_{Y^+}-E^+$ & $-17$ &$(r_W, -K_W^3, -K_W\cdot B, g(B))=(1, 18, 11, 4)$ \\   \hline
$4$  & $E_1$ & $-3K_{Y^+}-E^+$ & $1$ &$(r_W, -K_W^3, -K_W\cdot B, g(B))=(1, 6, 1, 0)$ \\   \hline
$4$  & $E_1$ & $-7K_{Y^+}-2E^+$ & $-20$ &$(r_W, -K_W^3, -K_W\cdot B, g(B))=(2, 24, 14, 4)$ \\   \hline
$4$  & $E_1$ & $-11K_{Y^+}-3E^+$ & $-65$ &$(r_W, -K_W^3, -K_W\cdot B, g(B))=(3, 54, 39, 14)$ \\   \hline
$5$  & $E_1$ & $-2K_{Y^+}-E^+$ & $-9$ &$(r_W, -K_W^3, -K_W\cdot B, g(B))=(1, 16, 7, 2)$\\   \hline
$6$  & $E_1$ & $-K_{Y^+}-E^+$ & $1$ &$(r_W, -K_W^3, -K_W\cdot B, g(B))=(1, 10, 1, 0)$ \\   \hline
$9$  & $E_1$ & $-3K_{Y^+}-4E^+$ & $-32$ &$(r_W, -K_W^3, -K_W\cdot B, g(B))=(4, 64, 28, 3)$ \\   \hline
$10$  & $E_1$ & $-2K_{Y^+}-3E^+$ & $-23$ &$(r_W, -K_W^3, -K_W\cdot B, g(B))=(3, 54, 21, {\cred 2})$ \\  \hline    
$12$  & $E_1$ & $-K_{Y^+}-2E^+$ & $-8$ &$(r_W, -K_W^3, -K_W\cdot B, g(B))=(2, 40, 10, 0)$ \\           \hline
      \end{tabular}}
    \end{table}
\end{thm}

\begin{proof}
Recall that we have 
\[
D \sim -\alpha K_{Y^+} -\beta E^+, 
\]
where $\alpha, \beta \in \Z_{>0}$ (Lemma \ref{l-beta=r}). 
The following hold (Lemma \ref{l-blowup-formula2}(2), Proposition \ref{p-line-basic}(1)): 
\begin{align*}
D^2 \cdot (-K_{Y^+})
&= (\alpha (-K_{Y^+}) -\beta E^+)^2 \cdot (-K_{Y^+})\\
&= \alpha^2 (-K_{Y^+})^3  -2 \alpha \beta (-K_{Y^+})^2 \cdot E^+ + \beta^2 (-K_{Y^+}) \cdot (E^+)^2\\
&= (2g-6) \alpha^2 -6 \alpha \beta -2 \beta^2\\
D \cdot (-K_{Y^+})^2 & = (\alpha (-K_{Y^+}) -\beta E^+) \cdot (-K_{Y^+})^2\\
&= \alpha (-K_{Y^+})^3  - \beta (-K_{Y^+})^2 \cdot E^+\\
&= (2g-6) \alpha -3\beta. 
\end{align*}

\setcounter{step}{0}

\begin{step}\label{s1-line-flop}
If $\tau$ is of type $D$, 
then 
$(g, \alpha, \beta, D \cdot (-K_{Y^+})^2)
=(7, 1, 1, 5)$ and $\tau$ is of type $D_1$. 
\end{step}

\begin{proof}[Proof of Step \ref{s1-line-flop}]
Recall that $\alpha \in \Z_{>0}$ and 
$\beta \in \{1, 2, 3\}$ (Lemma \ref{l-beta=r}). 
The following holds (Proposition \ref{p-typeD-intersec}): 
\[
f(\alpha, \beta) := (g-3) \alpha^2 -3 \alpha \beta -\beta^2 =0. 
\]
Then the quadratic equation 
\[
(g-3)x^2 -3x -1 =0
\]
has a solution $x = \alpha/\beta \in \Q$. 
Then the discriminant $D$ satisfies 
\[
\Q \ni D = \sqrt{(-3)^2 +4(g-3) } = \sqrt{4g -3}. 
\]
By $4 \leq g \leq 12$ (Corollary \ref{c-g-bound}), 
we have $13 \leq 4g-3 \leq 45$, 
which implies $g=7$. 
Then we get 
\[
0 =(g-3) \alpha^2 -3 \alpha \beta -\beta^2 = 4\alpha^2 -3\alpha \beta -\beta^2 = (4\alpha +\beta) (\alpha -\beta), 
\]
and hence $\alpha = \beta$. 
The following holds (Proposition \ref{p-typeD-intersec}):
\[
9\geq  D \cdot (-K_{Y^+})^2 =  
 (2g-6) \alpha - 3\beta = (2g-9)\alpha  = 5\alpha.  
\]
Therefore, $(g, \alpha, \beta, D \cdot (-K_{Y^+})^2) = (7, 1, 1, 5)$. 
This completes the proof of Step \ref{s1-pt-flop}. 
\end{proof}

\begin{step}\label{s2-line-flop}
If $\tau$ is of type $C$, 
then 
$(g, \alpha, \beta, \deg \Delta) \in \{(5, 2, 1, 7), (8, 1, 1, 5)\}$ 
and  $\tau$ is of type $C_1$. 
\end{step}

\begin{proof}[Proof of Step \ref{s2-line-flop}]
Recall that 
$\alpha \in \Z_{>0}$ and $\beta \in \{1, 2\}$ (Lemma \ref{l-beta=r}). 
The following holds (Proposition \ref{p-typeC-intersec}): 
\[
f(\alpha, \beta) := (g-3) \alpha^2 -3 \alpha \beta -\beta^2 
=1. 
\]
By $\alpha>0$ and $g\geq 4$, we get 
\[
\alpha = \frac{3\beta + \sqrt{9\beta^2 +4(g-3)(\beta^2+1)}}{2g-6}.
\]

Assume $g=5$. 
We then get $\alpha = \frac{3\beta + \sqrt{17\beta^2 +8}}{4}$. 
By $\beta \in \{1, 2\}$, we obtain $\beta = 1$ and $\alpha =2$, 
i.e., $(g, \alpha, \beta) =(5, {\cred 2, 1})$. 
When $g=6$ or $g=7$, there is no solution, 
because $\sqrt{9\beta^2 +4(g-3)(\beta^2+1)}$ is not a rational number as follows: 
\begin{itemize}
\item $(g, \beta)=(6, 1) \Rightarrow \sqrt{9\beta^2 +4(g-3)(\beta^2+1)} = 
\sqrt{9 +12 \cdot (1+1)}  =\sqrt{33} \not\in \Q$. 
\item $(g, \beta)=(6, 2)\Rightarrow 
 \sqrt{9\beta^2 +4(g-3)(\beta^2+1)} =
2 \sqrt{9 +3 \cdot ({\cred 4} +1)}  =2\sqrt{{\cred 24}} \not\in \Q$. 
\item $(g, \beta)=(7, 1)\Rightarrow 
\sqrt{9\beta^2 +4(g-3)(\beta^2+1)} =
\sqrt{9 +16 \cdot (1+1)}  =\sqrt{41} \not\in \Q$. 
\item  $(g, \beta)=(7, 2)\Rightarrow 
\sqrt{9\beta^2 +4(g-3)(\beta^2+1)} =
2 \sqrt{9 + 4 \cdot (4+1)}  =2\sqrt{29} \not\in \Q$. 
\end{itemize}
In what follows, we assume $g \geq 8$.

Assume that $\beta=1$. 
In this case, we get 
\[
f(\alpha, 1) = (g-3) \alpha^2 -3 \alpha -1 =1. 
\]
By $g \geq 8$, we have  $f(1, 1) < 
f(2, 1) < f(3, 1) < \cdots$. 
The following hold. 
\begin{itemize}
\item $f(1, 1) = (g-3)  -3 -1 = g-7 \geq 1$.  
\item $f(2, 1) > f(1, 1) \geq 1$. 
\end{itemize}
Thus $f(\alpha, 1)=1$ has a unique solution 
$(g, \alpha, \beta)= (8, 1, 1)$. 

Assume that $\beta=2$. 
In this case, we get 
\[
f(\alpha, 2) = (g-3) \alpha^2 -6 \alpha -4 =1. 
\]
By $g \geq 8$, we have 
$f(1, 2) < f(2, 2) < \cdots$. 
We have 
\begin{itemize}
\item $f(1, 2) = (g-3)  -6 -4 = g-{\cred 13}<1$ (Corollary \ref{c-g-bound}), and 
\item $f(2, 2) = 4f(1, 1) = 4(g-7) \geq 4 >1$.
\end{itemize}
Hence $f(\alpha, 2)=1$ has no solution.  

\medskip

To summarise, we obtain 
$(g, \alpha, \beta) \in \{(5, 2, 1), (8, 1, 1)\}$.
Recall that the following holds (Proposition \ref{p-typeC-intersec}):
\[
12 - \deg \Delta = D \cdot (-K_{Y^+})^2 = ( -\alpha K_{Y^+} - \beta E^+)\cdot (-K_{Y^+})^2
= (2g-6) \alpha - 3\beta. 
\]
\begin{itemize}
\item If $(g, \alpha, \beta) = (5, 2, 1)$, then 
$12 - \deg \Delta = (10-6)\cdot 2 - 3 \cdot 1 = 5$. 
\item 
If $(g, \alpha, \beta) = (8, 1, 1)$, then  
$12 - \deg \Delta = (16-6)\cdot 1 - 3 \cdot 1 = 7$. 
\end{itemize}
This completes the proof of Step \ref{s2-pt-flop}. 
\end{proof}

\begin{step}\label{s3-line-flop}
$\tau$ is not of type $E_2, E_3, E_4$, nor $E_5$. 
\end{step}

\begin{proof}[Proof of Step \ref{s3-line-flop}]
The following hold (Proposition \ref{p-typeE-intersec}): 
\begin{eqnarray*}
(2g-6) \alpha -3\beta &=& D \cdot (-K_{Y^+})^2 
=: u \in \{1, 2, 4\}\\
(2g-6) \alpha^2 - 6 \alpha \beta -2 \beta^2  &=&
D^2 \cdot (-K_{Y^+}) =-2.
\end{eqnarray*}
We then get 
\[
\alpha(3\beta +u) = (2g-6)\alpha^2 = 6\alpha \beta +2\beta^2 -2, 
\]
which implies 
\[
\alpha (-3\beta +u) = 2\beta^2 -2 \geq 0.
\]
By $u \in \{1, 2, 4\}$, we obtain $u=4$ and $\beta =1$. 
Then $\alpha= 2\beta^2 -2 = 2\cdot 1^2 -2 =0$, which is absurd. 
This completes the proof of Step \ref{s3-line-flop}. 
\end{proof}

\begin{step}\label{s4-line-flop}
Assume that $\tau$ is of type $E_1$. 
Then one of the possibilities appearing in Table \ref{table-2ray-line} holds. 
\end{step}

\begin{proof}[Proof of Step \ref{s4-line-flop}]
In this case, the following hold. 
\begin{enumerate}
\renewcommand{\labelenumi}{(\roman{enumi})}
\item $W$ is a Fano threefold with $\rho(W) =1$ and $\beta =r_W$ for the index $r_W$ of $W$ (Lemma \ref{l-beta=r}). 
In particular, $1 \leq \beta \leq 4$. Fix a Cartier divisor $H_W$ on $W$ satisfying $-K_W \sim \beta H_W$. 
\item $\alpha+1 = \beta \gamma$ for some $\gamma \in \Z_{>0}$  (Lemma \ref{l-beta=r}). 
\item $w := \frac{1}{2}\beta H_W^3 = (g-3)\gamma^2 - 3\gamma  -1$ (Remark \ref{r-E1-deg-genus}). 
We have $w \in \Z$ and $1 \leq w \leq 11$ (Remark \ref{r-E1-deg-genus}, 
Remark \ref{r-w-bound}). 
\item $h^0(Y^+, -K_{Y^+}-E^+) = h^0(Y, -K_Y-E) \geq g-5$ (Proposition \ref{p-line-basic}). 
\item $2g(B) -2 =(2g-6)\alpha^2 - 6 \alpha \beta  -2 \beta ^2$ (Lemma \ref{l-E1-deg-genus}). 
\item $-K_W \cdot B = (2g-6)\alpha(\alpha+1) 
- 3(2\alpha +1)\beta  -2\beta^2$  (Lemma \ref{l-E1-deg-genus}).
\item $D^3 = -(2g-6)+ (-K_W)^3 -3 (-K_W) \cdot B$ (Lemma \ref{l-D^3}). 
\end{enumerate}

We now prove the following assertion $(*)$. 
\begin{enumerate}
\item[$(*)$] If $\alpha < \beta$ or $g \geq 7$, then 
\[
(g, \alpha, \beta, -K_W^3)  \in \{ (9,  3, 4, 64), (10, 2, 3, 54), (12, 1, 2, 40)\}. 
\]
\end{enumerate}
We first show the implication $g \geq 7 \Rightarrow \alpha < \beta$. 
Assume  $g \geq 7$. 
By (iv), we obtain 
\[
h^0(Y^+, -K_{Y^+}-E^+) = h^0(Y, -K_Y-E) \geq g-5 \geq 2. 
\]
If $\alpha \geq \beta$, then we would get the following contradiction: 
\[
1 = h^0(Y^+, D) = h^0(Y^+, -\alpha K_{Y^+} - \beta E^+) \geq  
h^0(Y^+, - \beta K_{Y^+} - \beta E^+)  \geq g-5 \geq 2. 
\]
This completes the proof of the implication $g \geq 7 \Rightarrow \alpha < \beta$. 
In order to prove $(*)$, we may assume that $\alpha+1 \leq \beta$. 
This, together with $\alpha+1 =\beta \gamma \geq \beta$, implies  
$\alpha +1= \beta$ and $\gamma =1$. 
By (iii), we get $g= 7 + w$. 
If $\beta=4$, then $(\alpha, \beta, \gamma, w) = (3, 4, 1, 2)$, 
$-K_W^3= 64$, and $g=9$. 
If $\beta=3$, then 
$(\alpha, \beta, \gamma, w) = (2, 3, 1, 3)$,  $(-K_W)^3 = 54$, and $g=10$. 
Assume $\beta=2$. Then $(\alpha, \beta, \gamma) =(1, 2, 1)$ and $1 \leq w = H_W^3 \leq 5$. 
We then get 
\[
5 \geq w= g -7 \quad \text{and}\quad 
-2 \leq 2g(B) -2 \overset{{\rm (v)}}{=} 
(2g-6)\alpha^2 - 6 \alpha \beta  -2 \beta ^2 = 2g-26. 
\]
Hence we get a unique solution 
$(g, \alpha, \beta, -K_W^3, g(B))= (12, 1, 2, 40, 0)$. 
This completes the proof of $(*)$. 


\medskip

In what follows, we assume $g \leq 6$ and $\alpha \geq \beta$. 
The solution of $w = (g-3)\gamma^2 -3 \gamma -1$ is given by 
\[
\gamma = \frac{3 + \sqrt{9 +4(g-3)\left(w+1\right)}}{2(g-3)}, 
\]
because $g-3 >0$ and $\sqrt{9 +4(g-3)\left(w+1\right)} >3$.

Assume $g=6$. 
Then $\gamma = \frac{3 + \sqrt{12w + 21}}{6}$. 
By $1 \leq w \leq 11$ (Remark \ref{r-E1-deg-genus}, Remark \ref{r-w-bound}), we have 
\[
12w +21 \in  \{ 33, 45, 57, 69, 81, 93, 105, 117, 129, 141, 153\}. 
\]
Hence we get $(w, 12w +21) = (5, 81)$, which implies $\gamma =  \frac{3 + \sqrt{81}}{6} =2$. 
Then the following holds (Remark \ref{r-E1-deg-genus}): 
$(g, \alpha, \beta,  -K_W^3) \in \{ (6, 3, 2, 40), (6, 1, 1, 10)\}.$ 
The case $(g, \alpha, \beta,  -K_W^3)=(6, 3, 2, 40)$ does not occur, as otherwise 
we would get the following contradiction: 
\[
0 =\kappa(Y^+, D) = \kappa(Y^+, -3K_{Y^+} -2E^+) \geq \kappa(Y^+, -K_{Y^+}) =3. 
\]
where the inequality follows from $h^0(Y^+, -K_{Y^+} -E^+) \geq g-5 =1$ (Proposition \ref{p-line-basic}). 
Hence we get 
\[
(g, \alpha, \beta,  -K_W^3) = (6, 1, 1, 10). 
\]

Assume $g=5$. 
Then $\gamma = \frac{3 + \sqrt{8w + 17}}{4}$. 
By $1 \leq w \leq 11$ (Remark \ref{r-w-bound}), we have 
\[
8w +17 \in  \{ 25, 33, 41, 49, 57, 65, 73, 81, 89, 97, 105\}. 
\]
Hence we need $(w, 8w +17) = (1, 25), (4, 49), (8, {\cred 81})$. 
By $\gamma =  \frac{3 + \sqrt{8w + 17}}{4}$, we get 
$(w, 8w +17, \gamma) \in \{ (1, 25, 2), {\cred (8, 81, 3)}\}$. 
Then the following holds (Remark \ref{r-E1-deg-genus}): 
\[
(g, \alpha, \beta,  -K_W^3) \in \{(5, 3, 2, 8), (5, 1, 1, 2), {\cred (5, 2, 1, 16)}\}.
\]

Assume $g=4$. 
Then $\gamma = \frac{3 + \sqrt{4w + 13}}{2}$. 
By $1 \leq w \leq 11$ (Remark \ref{r-w-bound}), we have 
\[
4w +13 \in  \{ 17, 21, 25, 29, 33, 37, 41, 45, 49, 53, 57\}. 
\]
Hence we need $(w, 4w +13) = (3, 25), (9, 49)$. 
By $\gamma =  \frac{3 + \sqrt{4w + 13}}{2}$, we get 
$(w, 4w +13, \gamma) \in \{ (3, 25, 4), (9, 49, 5)\}$. 
Then the following holds (Remark \ref{r-E1-deg-genus}): 
\[
(g, \alpha, \beta,  -K_W^3) \in \{ 
(4, 11, 3, 54), (4, 7, 2, 24), (4, 3, 1, 6), 
(4, 4, 1, 18)\}. 
\]


The remaining invariants $g(B), -K_W \cdot B, D^3$ can be computed by (v), (vi), (vii), respectively. 
For the sake of completeness, we here include the computations for all the cases.\\ 
If $(g, \alpha, \beta, -K_W^3) = (12, 1, 2, 40)$, then the following hold.  
\begin{enumerate}
\item[(v)] $2g(B) -2 = (24-6) \cdot 1^2 -6 \cdot 1 \cdot 2 -2 \cdot 2^2 = 18-12 -8 =-2$. 
\item[(vi)] $-K_W \cdot B = (24-6) \cdot 1 \cdot 2 -3 \cdot 3 \cdot 2 -2 \cdot 2^2 = 36-18-8 =10$. 
\item[(vii)] $D^3 = -(24-6) + 40 -3 \cdot 10 = -18 +40-30=-8$. 
\end{enumerate}
If $(g, \alpha, \beta, -K_W^3) = (10, 2, 3, 54)$, then the following hold. 
\begin{enumerate}
\item[(v)] $2g(B) -2 = (20-6) \cdot 2^2 -6 \cdot 2 \cdot 3 -2 \cdot 3^2 = 
56-36 -18 = {\cred 2}$. 
\item[(vi)] $-K_W \cdot B = (20-6) \cdot 2 \cdot 3 -3 \cdot 5 \cdot 3 -2 \cdot 3^2 = 84-45-18 =21$. 
\item[(vii)] $D^3 = -(20-6) + 54 -3 \cdot 21 = -14 +54-63=-23$. 
\end{enumerate}
If $(g, \alpha, \beta, -K_W^3) = (9, 3, 4, 64)$, then the following hold. 
\begin{enumerate}
\item[(v)] $2g(B) -2 = (18-6) \cdot 3^2 -6 \cdot 3 \cdot 4 -2 \cdot 4^2 = 108 -72 -32 =4$. 
\item[(vi)] $-K_W \cdot B = (18-6) \cdot 3 \cdot 4 -3 \cdot 7 \cdot 4 -2 \cdot 4^2 = 144 -84 -32 =28$. 
\item[(vii)] $D^3 = -(18-6) + 64 -3 \cdot 28 = -12 +64 - 84=-32$. 
\end{enumerate}
If $(g, \alpha, \beta, -K_W^3) = (6, 1, 1, 10)$, then the following hold. 
\begin{enumerate}
\item[(v)] $2g(B) -2 = (12-6) \cdot 1^2 -6 \cdot 1 \cdot 1 -2 \cdot 1^2 = 6 - 6 - 2 =-2$. 
\item[(vi)] $-K_W \cdot B = (12-6) \cdot 1 \cdot 2 -3 \cdot 3 \cdot 1 -2 \cdot 1^2 = 12 - 9 - 2 =1$. 
\item[(vii)] $D^3 = -(12-6) + 10 -3 \cdot 1 = -6 +10 - 3=1$. 
\end{enumerate}
If $(g, \alpha, \beta, -K_W^3) = (5, 3, 2, 8)$, then we get the following contradiction: 
\begin{enumerate}
\item[(v)] $-2 \leq 2g(B) -2 = (10-6) \cdot 3^2 -6 \cdot 3 \cdot 2 -2 \cdot 2^2 = 36 - 36 - 8 =-8$. 
\end{enumerate}
If $(g, \alpha, \beta, -K_W^3) = (5, 1, 1, 2)$, then we get the following contradiction: 
\begin{enumerate}
\item[(v)] $-2 \leq 2g(B) -2 = (10-6) \cdot 1^2 -6 \cdot 1 \cdot 1 -2 \cdot 1^2 = 4 - 6 - 2 =-4$. 
\end{enumerate}
{\cred 
If $(g, \alpha, \beta, -K_W^3) = (5, 2, 1, 16)$, then the following hold. 
\begin{enumerate}
\item[(v)] $2g(B) -2 = (10-6) \cdot 2^2 -6 \cdot 2 \cdot 1 -2 \cdot 1^2 = 
16 - 12 - 2 =2$. 
\item[(vi)] $-K_W \cdot B = (10-6) \cdot 2 \cdot 3 -3 \cdot 5 \cdot 1 -2 \cdot 1^2 = 24 -15 -2 =7$. 
\item[(vii)] $D^3 = -(10-6) + 16 -3 \cdot 7 = -4 +16 - 21=-9$. 
\end{enumerate}
}
If $(g, \alpha, \beta, -K_W^3) = (4, 11, 3, 54)$, then the following hold. 
\begin{enumerate}
\item[(v)] $2g(B) -2 = (8-6) \cdot 11^2 -6 \cdot 11 \cdot 3 -2 \cdot 3^2 = 242-198-18 =26$. 
\item[(vi)] $-K_W \cdot B = (8-6) \cdot 11 \cdot 12 -3 \cdot 23 \cdot 3 -2 \cdot 3^2 = 264 - 207 -18 =39$. 
\item[(vii)] $D^3 = -(8-6) + 54 -3 \cdot 39 = -2 +54-117=-65$. 
\end{enumerate}
If $(g, \alpha, \beta, -K_W^3) = (4, 7, 2, 24)$, then the following hold. 
\begin{enumerate}
\item[(v)] $2g(B) -2 = (8-6) \cdot 7^2 -6 \cdot 7 \cdot 2 -2 \cdot 2^2 = 98 - 84 - 8=6$. 
\item[(vi)] $-K_W \cdot B = (8-6) \cdot 7 \cdot 8 -3 \cdot 15 \cdot 2 -2 \cdot 2^2 = 112 - 90 - 8 = 14$. 
\item[(vii)] $D^3 = -(8-6) + 24 -3 \cdot 14 = -2 + 24 -42 =-20$. 
\end{enumerate}
If $(g, \alpha, \beta, -K_W^3) = (4, 3, 1, 6)$, then the following hold. 
\begin{enumerate}
\item[(v)] $2g(B) -2 = (8-6) \cdot 3^2 -6 \cdot 3 \cdot 1 -2 \cdot 1^2 = 18 - 18 -2 =-2$. 
\item[(vi)] $-K_W \cdot B = (8-6) \cdot 3 \cdot 4 -3 \cdot 7 \cdot 1 -2 \cdot 1^2 = 24 - 21 -2  =1$. 
\item[(vii)] $D^3 = -(8-6) + 6 -3 \cdot 1 = -2 +6 - 3=1$. 
\end{enumerate}
If $(g, \alpha, \beta, -K_W^3) = (4, 4, 1, 18)$, then the following hold. 
\begin{enumerate}
\item[(v)] $2g(B) -2 = (8-6) \cdot 4^2 -6 \cdot 4 \cdot 1 -2 \cdot 1^2 = 32 - 24 -2 =6$. 
\item[(vi)] $-K_W \cdot B = (8-6) \cdot 4 \cdot 5 -3 \cdot 9 \cdot 1 -2 \cdot 1^2 = 40-27-2  =11$. 
\item[(vii)] $D^3 = -(8-6) + 18 -3 \cdot 11 = -2 +18 - 33=-17$. 
\end{enumerate}
This completes the proof of Step \ref{s4-line-flop}.     
\end{proof}
This completes the proof of Theorem \ref{t-line-flop}. 
\qedhere

\end{proof}




\subsection{Existence of flopping curves}

\begin{prop}\label{p-line-defect}
We use the same notation as in Theorem \ref{t-line-flop}. 
Then $E^3 \neq (E^+)^3$. 
In particular, there exists a curve $\zeta$ on $Y$ such that  $K_Y \cdot \zeta =0$. 
\end{prop}

\begin{proof}
Suppose $(1=)E^3 = (E^+)^3$ (Proposition \ref{p-line-basic}). 
Let us derive a contradiction. 
By $D \sim -\alpha K_{Y^+} -\beta E^+$ 
(Table \ref{table-2ray-line} in Theorem \ref{t-line-flop}), we have 
    \[
D^3 = ( -\alpha K_{Y^+} -\beta E^+)^3 \equiv -\beta^3 (E^+)^3 =-\beta^3E^3 =-\beta^3   \mod \alpha. 
\]
By Table \ref{table-2ray-line} in Theorem \ref{t-line-flop}, 
the remaining cases are as follows: 
\begin{enumerate}
\renewcommand{\labelenumi}{(\roman{enumi})}
\item  
{\cred $(\text{type of }\tau, g, \alpha, \beta, D^3)=(E_1, 4, 4, 1, -17)$.}
\item  
{\cred $(\text{type of }\tau, g, \alpha, \beta, D^3)=(E_1, 5, 2, 1, -9)$.}
\item  $(\text{type of }\tau, g, \alpha, \beta, D^3)=(E_1, 10, 2, 3, -23)$.
\item $(\text{type of }\tau, g, \alpha, \beta, D^3)=(E_1, 12, 1, 2, -8)$.
    \item {\cred $\alpha =\beta =1$, $D^3 \in \{0, 1\}$, and $g\neq 11$.} 
\end{enumerate} 
The following holds (Proposition \ref{p-line-basic}): 
\begin{eqnarray*}
D^3 &=& \alpha^3 (-K_Y)^3 - 3 \alpha^2 \beta (-K_Y)^2 \cdot E + 3 \alpha \beta^2 (-K_Y) \cdot E^2 -\beta^3 E^3\\
&=&  (2g-6)\alpha^3 - 9 \alpha^2 \beta -6 \alpha \beta^2-\beta^3.
\end{eqnarray*}
For (i)--(iv), we get the following contradictions: 
\begin{enumerate}
\renewcommand{\labelenumi}{(\roman{enumi})}
\item 
{\cred $-17 = D^3 = 2 \cdot 4^3 - 9 \cdot 4^2 \cdot 1 -6 \cdot 4 \cdot 1^2 - 1^3 \equiv 0 + 0 + 8 -1=7 \mod 16$.} 
\item 
{\cred 
$-9 = D^3 = 4 \cdot 2^3 - 9 \cdot 2^2 \cdot 1 -6 \cdot 2 \cdot 1^2 - 1^3 \equiv 
2 + 0 + 0 -1 = 1  \mod 6$.} 
\item $-23 = D^3 = 14 \cdot 2^3 - 9 \cdot 2^2 \cdot 3 -6 \cdot 2 \cdot 3^2 - 3^3 = 112 - 108 - 108 -27 <-23$. 
\item $-8 = D^3 = 18 \cdot 1^3 - 9 \cdot 1^2 \cdot 2 -6 \cdot 1 \cdot 2^2 - 2^3 =18 -18-24-8 \neq -8$. 
\end{enumerate}
{\cred 
Suppose (v). 
We then get 
\[
\{0, 1\} \ni  D^3 = (2g-6)- 9 -6-1 =2g -22, 
\]
which contradicts $g\neq 11$. 
}
\end{proof}

\begin{cor}\label{c-line-not-ample}
Let $X \subset \P^{g+1}$ be an anti-canonically embedded Fano threefold with $\Pic\,X = \Z K_X$. 
Let $\Gamma$ be a line on $X$ and let $\sigma : Y \to X$ be the blowup along $\Gamma$. 
Then $-K_Y$ is not ample. 
\end{cor}

\begin{proof}
Suppose that $-K_Y$ is ample. 
Let us derive a contradiction. 
Since $-K_Y$ is ample, we obtain $2g -6 = (-K_Y)^3 >0$ (Proposition \ref{p-line-basic}(1)), 
which implies $g \geq 4$. 
Then we may apply Proposition \ref{p-line-defect}. 
There exists a curve $\zeta$ on $Y$ satisfying $K_Y \cdot \zeta =0$, which is absurd. 
\end{proof}

\section{Non-existence of Fano blowups}\label{s-nonFano}

Let  $X \subset \P^{g+1}$ be an anti-canonically embedded Fano threefold with $\Pic\,X =\Z K_X$. 
Take a smooth curve $\Gamma$ on $X$. 
The purpose of this section is to prove that 
$Y$ is not Fano for the blowup $\sigma: Y \to X$ along $\Gamma$ (Theorem \ref{t-blowup-main}). 
Recall that we have $3 \leq g \leq 12$ (Theorem \ref{t-g-bound}). 
Since this result is well known for the case when $3 \leq g \leq 5$  (cf. Subsection \ref{ss-nonFano-g<6}), 
it is enough to treat the the case when $6 \leq g \leq 12$. 
In characteristic zero, this theorem was proven by showing the existence of lines on $X$. 
We here use a different strategy. 
The main idea is similar to those of former sections, i.e., 
by assuming that $-K_Y$ is ample, 
consider the Diophantine equations obtained from the contraction $\tau : Y \to W$ of the extremal ray $R$ of $\NE(Y)$ not corresponding to $\sigma : Y \to X$. 
We shall apply  case study depending on the type of $R$. 
The case when $\tau$ is not of type $E_1$ (resp. is of type $E_1$) is treated in Subsection \ref{ss-nonFano-notE1} (resp. Subsection \ref{ss-nonFano-E1}). 
For an overview of some ideas, we refer to Subsection \ref{ss-overview}. 




\begin{nota}\label{n-678}
Let $X \subset \P^{g+1}$ be an anti-canonically embedded Fano threefold with $\Pic\,X =\Z K_X$ and
 $6 \leq g \leq 12$. 
Let $\Gamma$ be a smooth curve on $X$. 
Set $d := -K_X \cdot \Gamma$ and $h := h^1(\Gamma, \MO_{\Gamma})$. 
Let $\sigma : Y \to X$ be the blowup along $\Gamma$. 
Assume that $-K_Y$ is ample. 
Let $\tau : Y \to W$ be the contraction of the extremal ray that is not corresponding to $\sigma : Y \to X$. 
\begin{itemize}
    \item If $\tau$ is of type $C$, then set $D := \tau^*L$ for a line $L$ on $W =\P^2$. 
    \item If $\tau$ is of type $D$, then set $D := \tau^*P$ for a point $P$ on $W =\P^1$. 
    \item If $\tau$ is of type $E$, then set $D := \Ex(\tau)$, which is a prime divisor on $Y$. 
    \end{itemize}   
We can write 
\[
D \sim -\alpha K_Y -\beta E. 
\]
for some $\alpha, \beta \in \Z_{>0}$ (Lemma \ref{l-beta=r}). 
We have $\alpha \Z + \beta \Z =\Z$ when $\tau$ is of type $C, D$, or $E_1$ (Lemma \ref{l-beta=r}). 
\[
\begin{tikzcd}
    & X \arrow[ld, "\sigma"'] \arrow[rd, "\tau"]\\
Y &  & W
\end{tikzcd}
\]
\end{nota}


\begin{nota}\label{n-678-E1}
We use Notation \ref{n-678}. Moreover, assume that $\tau$ is of type $E_1$. 
In this case, the following hold (Lemma \ref{l-beta=r}).
\begin{itemize}
    \item $W$ is a Fano threefold with $\rho(W)=1$. 
    \item $\beta = r_W$ for the index $r_W$ of $W$. 
    \item $\alpha +1 = \beta \gamma$ for some $\gamma \in \Z_{>0}$. 
\end{itemize}
We fix a Cartier divisor $H_W$ on $W$ such that $-K_W \sim \beta H_W$. Set $B :=\tau(D)$. 
\end{nota}

\begin{lem}\label{l-678-h0-d3}
We use Notation \ref{n-678}. 
Then the following hold. 
\begin{enumerate}
\item $d \geq 3$. 
\item 
If $h \geq 1$, then $d \geq 4$. 
\end{enumerate}
\end{lem}

\begin{proof}
Let us show (1). 
Suppose that $d \in \{1, 2\}$. 
Then $\Gamma$ is either a line or a conic. 
This is absurd (Corollary \ref{c-conic-not-ample} and Corollary \ref{c-line-not-ample}). 
Thus (1) holds. 

Let us show (2). 
Assume $h \geq 1$. 
By (1), it is enough to exclude the case $d =3$. 
Suppose $d=3$. 
Let us derive a contradiction. 
Recall that $X$ is an intersection of quadrics, i.e., we have $X = \bigcap_{X \subset Q} Q$, 
where $Q$ runs over all the quadric hypersurfaces in $\P^{g+1}$ containing $X$. 
Set $V :=\langle \Gamma \rangle$, which is the smallest linear subvariety of $\P^{g+1}$ containing $\Gamma$. 

We first treat the case when $ V= \P^2$. 
In this case, $\Gamma$ is an elliptic curve on $V=\P^2$. 
We get 
\[
\Gamma \subset V \cap X = \P^2 \cap \left(\bigcap_{X \subset Q} Q\right)
= \bigcap_{X \subset Q} (\P^2 \cap Q). 
\]
Note that $\P^2 \cap Q$ is either $\P^2$ or an effective divisor on $\P^2$ of degree two. 
By $V \not\subset X$, there exists $Q$ such that $\P^2\cap Q$ is an effective divisor on $\P^2$ of degree two. This is a contradiction, because $\Gamma$ is of degree three. 
Thus $V \neq \P^2$.

Pick general four points $P_1, P_2, P_3, P_4 \in \Gamma$. 
By $\langle \Gamma \rangle \neq \P^2$, the smallest linear subvariety $V'$ of $\P^{g+1}$ containing $P_1, P_2, P_3, P_4$ is  $\P^3$. 
By $\deg \Gamma = -K_X \cdot \Gamma =d=3$, $\Gamma$ must be contained in $V'=\P^3$, 
i.e., $V=V'=\P^3$. 
Since {\cred $X \neq V=\P^3$ and} $X$ is an intersection of quadrics, 
we can find a quadric hypersurface $Q \subset \P^{g+1}$ such that 
\begin{itemize}
\item $\Gamma \subset Q \cap V$ and 
\item $Q \cap V$ is a (possibly singular) quadric surfaces on 
$V =\P^3$. 
\end{itemize}
{\cred This contradicts Lemma \ref{l-cubic-on-qurd-surf} below by 
$\deg \Gamma = d=3$ and $h^1(\Gamma, \MO_{\Gamma})=h \geq 1$.} 
\end{proof}

\begin{lem}\label{l-cubic-on-qurd-surf}
Let $S$ be a (possibly singular) quadric surface on $\P^3$. 
Let $\Gamma$ be a  curve on $S$ of degree three, i.e., $\deg (\MO_{\P^3}(1)|_{\Gamma})=3$. 
Then $\Gamma \simeq \P^1$. 
\end{lem}

\begin{proof}
If $S$ is smooth, then $S \simeq \P^1 \times \P^1$. 
In this case, we obtain $C \in |\MO_{\P^1 \times \P^1}(1, 2)|$ or $C \in |\MO_{\P^1 \times \P^1}(2, 1)|$. 
We then get $(K_S+C) \cdot C <0$, which implies $C \simeq \P^1$. 
Thus we may assume that $S$ is singular. 
In this case, $S$ is obtained by contracting the $(-2)$-curve on $\P_{\P^1}(\MO_{\P^1} \oplus \MO_{\P^1}(2))$. 
In particular, $S$ is $\Q$-factorial and ${\rm Cl}(S) =\Z \ell$, where 
$\ell$ is a line on $S$ and ${\rm Cl}(S)$ denotes the divisor class group of $S$. 
By $\deg \Gamma=3$ and $\deg \ell =1$, we obtain $\Gamma \equiv 3 \ell$. 
In particular, 
\[
(K_S+\Gamma) \cdot \Gamma  = (-4\ell +3\ell) \cdot ({\cred 3}\ell) <0. 
\]
We then get $\Gamma \simeq \P^1$ \cite[Theorem 3.19(1)]{Tan14}. 
\end{proof}

\begin{lem}\label{l-678-easy-bound}
We use Notation \ref{n-678}. 
Then the following hold. 
\begin{enumerate}
\item 
\begin{enumerate}
\item $(-K_Y)^3 = 2g-4-2d +2h$. 
\item $(-K_Y)^2 \cdot E =d -2h+2$. 
\item $(-K_Y) \cdot E^2 = 2h-2$.
\item $-E^3 = d + 2h-2$. 
\end{enumerate}
    \item $h \leq g-2 \leq 10$. 
    \item $2h-1 \leq d \leq h+g-3 \leq 2g -5 \leq 19$. 
\item \begin{enumerate}
\item $D^3 = (2g-4 -2 d +2h)\alpha^3 +3(2h-2-d)\alpha^2\beta 
+3(2h-2)\alpha \beta^2 +(d+2h-2)\beta^3$.
\item $D^2 \cdot (-K_Y) 
    = (2g-4-2d+2h) \alpha^2 -2(d-2h+2) \alpha\beta +(2h-2)\beta^2$. 
\item $D \cdot (-K_Y)^2  = (2g-4-2d+2h) \alpha -(d-2h+2)\beta$. 
\end{enumerate}
\item It holds that 
\begin{eqnarray*}
-D^2 \cdot (-K_Y) \alpha +D^3 
&=& (-d+2h-2)\alpha^2\beta +(4h-4)\alpha\beta^2 + (d+2h-2)\beta^3\\
&=& d\beta (-\alpha+\beta)(\alpha+\beta) + (2h-2)\beta (\alpha+\beta)^2. 
\end{eqnarray*}
In particular, $-D^2 \cdot (-K_Y) \alpha +D^3  \in \beta(\alpha+\beta)\Z$. 
\end{enumerate}
\end{lem}


\begin{proof}
The assertion (1) follows from Lemma \ref{l-blowup-formula}(2). 
Let us show (2) and (3). 
Since $-K_Y$ is ample, 
it follows from 
the Nakai-Moishezon criterion that $-K_Y^3 >0$ and $(-K_Y)^2 \cdot E >0$. 
{\cred By (1),} 
we get $g-d +h -2 >0$ and $d -2h+2>0$, 
which imply 
\[
2h -2 <  d < h +g-2. 
\]
In particular, it holds that $h <g$. 
If $h =g-1$, these inequalities would become 
\[
2g-4<  d < 2g-3, 
\]
which is a contradiction. 
Thus (2) and (3) hold. We obtain (4)  by the following: 
\begin{eqnarray*}
D^3 &=& (-\alpha K_Y -\beta E)^3\\
&=& (-K_Y)^3 \alpha^3 -3 (-K_Y)^2 \cdot E \alpha^2\beta 
+3 (-K_Y) \cdot E^2 \alpha \beta^2 + (-E)^3\beta^3\\
&\overset{{\rm (1)}}{=}& (2g-4 -2 d +2h)\alpha^3 +3(2h-2-d)\alpha^2\beta 
+3(2h-2)\alpha \beta^2 +(d+2h-2)\beta^3. 
\end{eqnarray*}
\begin{eqnarray*}
    D^2 \cdot (-K_Y) &=& (\alpha^2 K^2_Y +2\alpha\beta K_Y \cdot E +\beta^2E^2) \cdot (-K_Y)\\
    &=& (-K_Y)^3 \alpha^2 - 2 (-K_Y)^2 \cdot E \alpha \beta + (-K_Y) \cdot E^2 \beta^2\\
    &\overset{{\rm (1)}}{=}& (2g-4-2d+2h) \alpha^2 -2(d-2h+2) \alpha\beta +(2h-2)\beta^2
\end{eqnarray*}
\begin{eqnarray*}
    D \cdot (-K_Y)^2 &=& (-\alpha K_Y -\beta E) \cdot (-K_Y)^2\\
    &=& (-K_Y)^3 \alpha -  (-K_Y)^2 \cdot E \beta\\
    &\overset{{\rm (1)}}{=}& (2g-4-2d+2h) \alpha -(d-2h+2)\beta.
\end{eqnarray*}
The assertion (5) holds by the following: 
\begin{eqnarray*}
-D^2 \cdot (-K_Y) \alpha +D^3 
&\overset{{\rm (4)}}{=}& (-d+2h-2)\alpha^2\beta +(4h-4)\alpha\beta^2 + (d+2h-2)\beta^3\\
&=& d(-\alpha^2 \beta +\beta^3) + (2h-2) (\alpha^2\beta +2\alpha\beta^2 + \beta^3)\\
&=& d\beta (-\alpha+\beta)(\alpha+\beta) + (2h-2)\beta (\alpha+\beta)^2 \in \beta(\alpha+\beta)\Z. 
\end{eqnarray*}

\end{proof}


\subsection{Not type $E_1$ ($g \geq 6$)}\label{ss-nonFano-notE1}

\begin{lem}\label{l-678-nonE1-E}
We use Notation \ref{n-678}. 
Then $\tau$ is not of type $E_2, E_3, E_4$, nor $E_5$. 
\end{lem}

\begin{proof}
Suppose that $\tau$ is of type $E_2, E_3, E_4$, or $E_5$. 
Let us derive a contradiction. 
We have $v := D^3 \in \{1, 2, 4\}$, $D^2 \cdot (-K_Y) =-2$, and $D \cdot (-K_Y)^2 = 4/v$ (Proposition \ref{p-typeE-intersec}). 

\begin{claim*}
One of {\rm (I)} and {\rm (II)} holds. 
\begin{enumerate}
    \item[(I)] $(v, \beta) = (4, 2)$. 
    \item[(II)] $(v, \beta) = (2, 1)$ and   $(d - 2h+2) \alpha = d+2h-4$. 
\end{enumerate}
\end{claim*}

\begin{proof}[Proof of Claim] 
Lemma \ref{l-678-easy-bound}(5) implies 
\[
2\alpha +v = 
-D^2\cdot (-K_Y) \alpha + D^3 
= d\beta(-\alpha + \beta)(\alpha+\beta) +(2h-2)\beta(\alpha+\beta)^2 \in \beta(\alpha +\beta)\Z. 
\]
In particular, we can write $2\alpha +v = \beta(\alpha +\beta)m$ for some $m \in \Z_{>0}$.

Assume $v=1$. 
Then $2\alpha +1 = \beta(\alpha +\beta)m$. 
If $\beta \geq 2$ or $m\geq 2$, then we get the following contradiction: 
$2\alpha + 1 = \beta(\alpha +\beta)m \geq 2(\alpha +\beta) \geq 2\alpha +2$. 
Hence $\beta = m =1$. 
Then $2\alpha +1 = \beta(\alpha +\beta)m = \alpha +1$, i.e., $\alpha =0$, which is absurd. 
Thus $v \neq 1$.

Assume $v=2$. 
Then $2\alpha +2 = \beta(\alpha +\beta)m$. 
If $\beta \geq 2$, then we get the following contradiction: 
$2\alpha +2 = \beta(\alpha +\beta)m \geq 2(\alpha+2) =2\alpha +4$. 
Hence we get $\beta =1$ and 
\[
2\alpha +2= d(-\alpha + 1)(\alpha+1) +(2h-2)(\alpha+1)^2.  
\]
Dividing by $\alpha+1$, we obtain $2 = d(-\alpha + 1) +(2h-2)(\alpha+1)$, i.e., 
\[
(d - 2h+2) \alpha = d+2h-4.
\]
Thus (II) holds.

Assume $v=4$. 
In order to prove (I), it suffices to  show $\beta =2$. 
If $\beta \geq 3$, 
then we would get 
the following contradiction: 
$2\alpha +4 = \beta(\alpha +\beta)m \geq 3(\alpha +3)=3\alpha +9$. 
Suppose $\beta =1$. 
We get $2\alpha +4 = \beta(\alpha +\beta)m = (\alpha+1)m \in (\alpha+1)\Z$. 
Since $\alpha +1$ and $\alpha+2$ are coprime, we get $2 \in (\alpha+1)\Z$. 
Then $\alpha =\beta =1$, which leads to the following contradiction:  
\[
6= 2\alpha +v = d\beta(-\alpha + \beta)(\alpha+\beta) +(2h-2)\beta(\alpha+\beta)^2 = 0 +4(2h-2) \in 4\Z.
\]
This completes the proof of Claim. 
\end{proof}


Recall that the following hold (Lemma \ref{l-678-easy-bound}(4)): 
\begin{eqnarray*}
    -2=D^2 \cdot (-K_Y) &=&  (2g-4-2d+2h) \alpha^2 -2(d-2h+2) \alpha\beta +(2h-2)\beta^2\\
    4/v = D \cdot (-K_Y)^2 &=& (2g-4-2d+2h) \alpha -(d-2h+2)\beta
\end{eqnarray*}
Assume (I).  In this case, we obtain the following contradiction: 
\[
1 =4/v = D \cdot (-K_Y)^2 = (2g-4-2d+2h) \alpha -2(d-2h+2) \in 2\Z. 
\]
Thus (II) holds. Then the above equations become
\begin{eqnarray*}
    -2&=&  (2g-4-2d+2h) \alpha^2 -2(d-2h+2) \alpha +2(h-1)\\
    2  &=& (2g-4-2d+2h) \alpha -(d-2h+2).
\end{eqnarray*}
We get 
\[
2(d-2h+2) \alpha -2(h-1) -2 = (2g-4-2d+2h) \alpha^2 = \alpha (d-2h+4), 
\]
which implies  
\[
(d -2h)\alpha = 2h. 
\]
This, together with $(d - 2h+2) \alpha = d+2h-4$ (assured by (II)), 
implies  $2\alpha = d-4$. 
Eliminating $d$ from $(d -2h)\alpha = 2h$ by using $2\alpha = d-4$, 
we get $\alpha ( 2\alpha -2h +4) =2h$, i.e., $\alpha^2 +(-h+2)\alpha -h=0$. 
The discriminant $\Delta$ of this quadratic equation is given by 
$\Delta = (-h+2)^2 +4h = h^2 +4$. 
Since $\Delta$ is a square,  we get $(h, \Delta) = (0, 4)$ 
(because  $(h+1)^2 > h^2+4 =\Delta >h^2$ for $h \geq 2$). 
However, $(h, \Delta) = (0, 4)$ implies $\alpha^2 +2\alpha=0$, which contradicts $\alpha \in \Z_{>0}$. 
\qedhere




\end{proof}

\begin{lem}\label{l-678-nonE1-C}
We use Notation \ref{n-678}. 
Then $\tau$ is not of type $C$. 
\end{lem}

\begin{proof}
Suppose that $\tau$ is of type $C$. 
Let us derive a contradiction. 
In this case, we have $D^3 =0$ and $D^2 \cdot (-K_Y)=2$ (Proposition \ref{p-typeC-intersec}). 
By Lemma \ref{l-678-easy-bound}(5), we obtain  
\begin{equation}\label{e1-678-nonE1-C}
-2\alpha = -D^2 \cdot (-K_Y) \alpha + D^3= (-d+2h-2)\alpha^2\beta +(4h-4)\alpha\beta^2 + (d+2h-2)\beta^3. 
\end{equation}
Then $-2\alpha \in \beta \Z$ and $(d+2h-2)\beta^3 \in \alpha \Z$. 
Since $\alpha$ and $\beta$ are coprime, we get $\beta \in \{1, 2\}$ 
and $d+2h-2 \in \alpha \Z$. 


\medskip

We now treat the case when  $\beta =1$. 
Then Lemma \ref{l-678-easy-bound}(5) implies 
\[
-2\alpha = -D^2 \cdot (-K_Y) \alpha + D^3
= (1-\alpha^2)d + 2 (\alpha +1)^2  (h-1) \in (\alpha +1) \Z. 
\]
Since $\alpha$ and $\alpha+1$ are corpime, we obtain $-2 \in (\alpha+1) \Z$. 
Then $\alpha =1$ and $-2 =-2\alpha = 8(h-1) \in 8 \Z$, which is a contradiction. 

\medskip

We then get $\beta =2$. 
Then 
the following holds: 
\begin{eqnarray*}
-\alpha \overset{{\rm (i)}}{=} 
\frac{1}{2}(-D^2 \cdot (-K_Y) \alpha + D^3) \overset{{\rm (ii)}}{=} 
\frac{1}{\beta}(-D^2 \cdot (-K_Y) \alpha + D^3) \overset{{\rm (iii)}}{\in}  (\alpha +\beta) \Z \overset{{\rm (iv)}}{=} (\alpha +2)\Z, 
\end{eqnarray*} 
where (i) follows from (\ref{e1-678-nonE1-C}), 
(ii) and (iv) hold by $\beta =2$, and we have (iii) by Lemma \ref{l-678-easy-bound}(5). 
However, this is impossible by $1 \leq \alpha <\alpha +2$. 
\qedhere

\end{proof}

\begin{lem}\label{l-678-nonE1-D}
We use Notation \ref{n-678}. 
Then $\tau$ is not of type $D$. 
\end{lem}

\begin{proof}
Suppose that $\tau$ is of type $D$. 
Let us derive a contradiction. 
In this case, we have $D^3 =0$ and  $D^2 \cdot (-K_Y)=0$ 
(Proposition \ref{p-typeD-intersec}). 
By Lemma \ref{l-678-easy-bound}(4)(5), we obtain  
\begin{equation}\label{e1-678-nonE1-D}
(g-2-d+h) \alpha^2 -(d-2h+2) \alpha\beta +(h-1)\beta^2=0
\end{equation}
\begin{equation}\label{e2-678-nonE1-D}
d (-\alpha+\beta) + (2h-2) (\alpha+\beta) =0 
\end{equation}
We treat the following three cases separately. 
\begin{enumerate}
    \item[(0)] $h =0$. 
    \item[(1)] $h=1$.
    \item[(2)] $h \geq 2$. 
\end{enumerate}

(0) 
Assume $h=0$. Then (\ref{e1-678-nonE1-D}) implies $\beta^2 \in \alpha \Z$. 
Since $\alpha$ and $\beta$ are coprime, we get $\alpha =1$. 
Again by (\ref{e1-678-nonE1-D}), we obtain 
\begin{eqnarray*}
\varphi(\beta) := (g-2-d) -(d+2) \beta -\beta^2 = 0. 
\end{eqnarray*}
Then the quadratic function $\varphi(x)$ is decreasing on $x \geq 0$. 
We have $\varphi(0) = g-2 -d = \frac{1}{2}(-K_Y)^3 >0$ (Lemma \ref{l-678-easy-bound}(1)) and 
$\varphi(1) =  (g-2-d) -(d+2) -1 = g -5-2d \leq g-5 -2 \cdot 3 = g-11$ 
(Lemma \ref{l-678-h0-d3}(1)). 
Hence $g = 11$ or $g=12$. 
If $g =11$,  then we obtain $(g, d, h, \alpha, \beta) =(11, 3, 0, 1, 1)$, 
which contradicts (\ref{e2-678-nonE1-D}) as follows: 
\[
0 = d (-\alpha+\beta) + (2h-2) (\alpha+\beta) = 3 \cdot (-1+1) + (-2)(1+1) \neq 0.
\]
Thus $g =12$. 
We then get $\varphi(1) = 7-2d \neq 0$, and hence $\varphi(1) = 7-2d >0$. 
By $d \geq 3$ (Lemma \ref{l-678-h0-d3}(1)), we get $d=3$. 
Then $\varphi(\beta) = 7 -5\beta  - \beta^2 =0$, which contradicts $\beta \in \Z$.

\medskip

(1) Assume $h=1$. 
By (\ref{e2-678-nonE1-D}), we get $\alpha = \beta$. 
Since $\alpha$ and $\beta$ are coprime, we obtain $\alpha = \beta  =1$, i.e., 
\[
\tau^*\MO_{\P^1}(1) \sim D \sim -K_Y -E, 
\]
where 
$\tau : Y \to W = \P^1$. 
Then (\ref{e1-678-nonE1-D}) implies 
$g-1-2d =0$. 
We shall prove (i) and (ii) below. 
\begin{enumerate}
    \item[(i)]  $h^0(Y, -K_Y)  \geq g-d+2$. 
    \item[(ii)] $h^0(E, -K_Y|_E)=d$. 
\end{enumerate}   
Assuming (i) and (ii), we now finish the proof. 
By the exact sequence 
\[
0 \to H^0(Y, -K_Y-E) \to H^0(Y, -K_Y) \to H^0(E, -K_Y|_E), 
\]
(i) and (ii), together with $g-1-2d=0$, imply 
\[
h^0(Y, -K_Y -E)  \geq 
h^0(Y, -K_Y) -h^0(E, -K_Y|_E)\geq  
(g-d+2) -d = (g-1-2d) +3 = 3.
\]
This contradicts $h^0(Y, -K_Y -E) = h^0(Y, \tau^*\MO_{\P^1}(1)) = h^0(\P^1, \MO_{\P^1}(1)) =2$.



It is enough to show (i) and (ii). 
By Lemma \ref{l-678-easy-bound}(1), the following hold: 
\begin{itemize}
\item $(-K_Y)^3 = 2g-2d +2h -4 = 2g-2d-2$. 
\item $(-K_Y)^2 \cdot E =d -2h+2 = d$.
\item $(-K_Y) \cdot E^2 = 2h-2 =0$.
\end{itemize}
Then the assertion (i) holds as follows: 
\[
h^0(Y, -K_Y) \geq h^0(Y, -K_Y) - h^1(Y, -K_Y) \overset{(\star)}{=} \chi(Y, -K_Y) 
= \frac{(-K_Y)^3}{2} +3 = g-d+2, 
\]
where 
$(\star)$ is assured by \cite[Corollary 2.6]{TanI}. 
Let us show (ii). 
By the Riemann--Roch theorem, we obtain 
\[
\chi(E, -K_Y|_E) = \chi(E, \MO_E) + \frac{1}{2}(-K_Y|_E) \cdot ((-K_Y|_E)-K_E) 
\]
\[
= \frac{1}{2}E \cdot (-K_Y) \cdot (-2K_Y-E)= (-K_Y)^2 \cdot E =d.
\]
Hence it suffices to prove $H^i(E, -K_Y|_E)=0$ for every $i>0$. 
Since $E$ satisfies Kodaira vanishing \cite[Theorem 3(a)]{Muk13}, 
it is enough to show that the divisor $-K_Y|_E -K_E =(-2K_Y-E)|_E$ is ample. 
This follows from the fact that $-K_Y$ is ample and $-K_Y -E \sim \tau^*\MO_{\P^1}(1)$ is nef.





\medskip

(2) 
Assume $h \geq 2$. 
By $h \geq 2$ and $(\ref{e2-678-nonE1-D})$, we get $\alpha >\beta \geq 1$. 
We have $-(h-1)\beta^2 \in \alpha \Z$ (\ref{e1-678-nonE1-D}), 
and hence $h-1 \in \alpha \Z$. 
By (\ref{e2-678-nonE1-D}), we get  $d (-\alpha +\beta) \in \alpha \Z$, 
which implies $d \in \alpha \Z$. 
By $h-1 \in \alpha \Z$, $d \in \alpha \Z$, and (\ref{e1-678-nonE1-D}), we obtain $(h-1)\beta^2 \in \alpha^2 \Z$, and hence $h-1 \in \alpha^2 \Z$. 
By  $\alpha \geq 2$, $h-1 \in \alpha^2 \Z$, and $2 \leq h \leq 10$ (Lemma \ref{l-678-easy-bound}(2)), we get $\alpha \in \{2, 3\}$. 
By $\alpha > \beta \geq 1$, it holds that $(\alpha, \beta) \in \{ (2, 1), (3, 1), (3, 2)\}$.

Assume $\beta =1$. 
Then (\ref{e2-678-nonE1-D}) implies 
\[
 (2h-2) (\alpha +1) = d(\alpha -1).
\]
If $(\alpha, \beta) =(2, 1)$, then  
$d = 6(h-1) \in 6\alpha^2 \Z = 24 \Z$, which contradicts $1 \leq d \leq 19$ (Lemma \ref{l-678-easy-bound}). 
If $(\alpha, \beta) =(3, 1)$, then $d = 4(h-1) \in 4\alpha^2  \Z = 36 \Z$, which contradicts  $1 \leq d  \leq 19$ (Lemma \ref{l-678-easy-bound}). 

Thus we get $\beta \neq 1$, i.e., $(\alpha, \beta) =(3, 2)$. 
Then (\ref{e2-678-nonE1-D}) implies 
\[
d = d(\alpha -\beta)  =  (2h-2) (\alpha +\beta) = 10(h-1) \in 10 \alpha^2 \Z =90 \Z. 
\]
This contradicts  $1 \leq d  \leq 19$ (Lemma \ref{l-678-easy-bound}).  
\qedhere

\end{proof}

\begin{prop}\label{p-678-nonE1-main}
We use Notation \ref{n-678}. 
Then $\tau$ is of type $E_1$. 
\end{prop}

\begin{proof}
The assertion follows from Lemma \ref{l-678-nonE1-E}, Lemma \ref{l-678-nonE1-C}, 
and  Lemma \ref{l-678-nonE1-D}. 
\end{proof}

\subsection{Type $E_1$ ($g \geq 6$)}\label{ss-nonFano-E1}

\begin{lem}\label{l-678-E1-ABC}
We use Notation \ref{n-678-E1}. 
Then the following hold. 
\begin{enumerate}
    \item[(A)] $\frac{1}{2} \beta H_W^3 = (g-d +h -2)
\gamma^2 - (d -2h+2)\gamma   + (h-1)$. 
\item[(B)] $H_W^3 = (2g-2d +2h -4)\gamma^3 -3(d -2h+2)\gamma^2+ 3(2h-2)\gamma +(d+2h-2)$. 
In particular, $H_W^3 \equiv d+2h -2 \mod \gamma$. 
\item[(C)] $\alpha H_W^3 = (\beta \gamma -1)H_W^3 = (d -2h+2)\gamma^2 -2(2h-2)\gamma  -(d+2h-2)$. 
\item[(D)] 
\begin{itemize}
    \item $8h-8 \in \alpha \Z$ if $\beta=1$. 
    \item $3d +18h-18 \in \alpha \Z$ if $\beta=2$. 
    \item $8d +32h-32 \in \alpha \Z$ if $\beta=3$. 
    \item $15d +50h-50 \in \alpha \Z$ if $\beta=4$. 
\end{itemize}
\end{enumerate}
\end{lem}


\begin{proof}
Recall that the following hold (Lemma \ref{l-E1-deg-genus}(1)(2), Lemma \ref{l-D^3}(2)). 
\begin{enumerate}
\item $(-K_W)^3 = (\alpha+1)^2 (-K_{Y})^3 - 2(\alpha+1)\beta (-K_Y)^2 \cdot E  + \beta^2(-K_Y) \cdot E^2$. 
    \item $D^3 = \alpha^3 (-K_Y)^3 - 3 \alpha^2 \beta (-K_Y)^2 \cdot E + 3 \alpha \beta^2 (-K_Y) \cdot E^2 -\beta^3 E^3$. 
    \item $D^3 = (-K_W)^3 - (-K_Y)^3 -3(-K_W) \cdot B$. 
    \item $-K_W \cdot B = \alpha (\alpha+1)(-K_Y)^3 
- (2\alpha +1) \beta(-K_Y)^2 \cdot E 
+\beta^2 (-K_Y) \cdot E^2.$
\end{enumerate}

Let us show  (A). 
By (1), $(-K_W)^3 = \beta^3 H_W^3$ (Notation \ref{n-678-E1}), and $\alpha + 1 = \beta \gamma$, we obtain 
\[
\beta H_W^3 = 
\gamma^2 (-K_{Y})^3 - 2\gamma (-K_Y)^2 \cdot E  + (-K_Y) \cdot E^2.
\]
Thus (A) holds by Lemma \ref{l-678-easy-bound}(1). 


Let us show (B). 
Let us erase $-K_W \cdot B$ and $D^3$ from the above equations (2)-(4): 
\begin{eqnarray*}
&& \alpha^3 (-K_Y)^3 - 3 \alpha^2 \beta (-K_Y)^2 \cdot E + 3 \alpha \beta^2 (-K_Y) \cdot E^2 -\beta^3 E^3\\
&\overset{(2)}{=}& D^3\\
&\overset{(3)}{=}&(-K_W)^3 - (-K_Y)^3 -3(-K_W) \cdot B\\
&\overset{(4)}{=}& (-K_W)^3 - (-K_Y)^3 -3(\alpha(\alpha+1)(-K_Y)^3 
- (2\alpha +1) \beta(-K_Y)^2 \cdot E 
+\beta^2(-K_Y) \cdot E^2)\\
&=& (-K_W)^3 - (-K_Y)^3 -3\alpha(\alpha+1)(-K_Y)^3 
+ 3(2\alpha +1)\beta (-K_Y)^2 \cdot E 
-3\beta^2 (-K_Y) \cdot E^2.
\end{eqnarray*}
Then 
\begin{eqnarray*}
(-K_W)^3 
&=& (\alpha^3 +3\alpha^2 +3\alpha+1)(-K_Y)^3
+ (-3\alpha^2\beta -3(2\alpha+1)\beta)(-K_Y)^2\cdot E\\
&& + (3\alpha\beta^2 +3\beta^2)(-K_Y) \cdot E^2 -\beta^3 E^3\\
&=& (\alpha+1)^3(-K_Y)^3 -3(\alpha+1)^2\beta(-K_Y)^2\cdot E
+ 3(\alpha+1)\beta^2(-K_Y) \cdot E^2 +\beta^3 (-E^3).
\end{eqnarray*}
By $(-K_W)^3 = \beta^3 H_W^3$ and $\alpha +1 = \beta \gamma$, we get 
\[
H_W^3 =  (-K_Y)^3\gamma^3 -3(-K_Y)^2\cdot E \gamma^2
+ 3(-K_Y) \cdot E^2 \gamma +(-E^3).
\]
Then (B) holds  by Lemma \ref{l-678-easy-bound}(1). 
The equation  (C) is obtained by taking ${\rm (A)} \times 2\gamma - {\rm (B)}$.

Let us show (D). We first compute $-\beta^3 E^3$ modulo $\alpha$: 
\begin{eqnarray*}
-\beta^3 E^3  
&\overset{{\rm (2)}}{\equiv} & D^3\\
&\overset{{\rm (3)}}{=}& (-K_W)^3 - (-K_Y)^3 -3(-K_W) \cdot B\\
&\overset{{\rm (1)}}{\equiv}&- 2\beta (-K_Y)^2 \cdot E  + \beta^2(-K_Y) \cdot E^2
-3(-K_W) \cdot B\\
&\overset{{\rm (4)}}{\equiv}&- 2\beta (-K_Y)^2 \cdot E  + \beta^2(-K_Y) \cdot E^2
-3 (- \beta (-K_Y)^2 \cdot E +\beta^2 (-K_Y) \cdot E^2)\\
&=& \beta (-K_Y)^2 \cdot E  - 2\beta^2(-K_Y) \cdot E^2 \mod \alpha. 
\end{eqnarray*}
Since $\alpha$ and $\beta$ are coprime (Notation \ref{n-678-E1}), 
Lemma \ref{l-678-easy-bound}(1) implies 
\[
\beta^2(2h -2 +d) \equiv (d-2h+2) -2\beta  (2h-2) \mod \alpha, 
\]
i.e., $(\beta^2 -1)d + (\beta +1)^2 (2h-2) \in \alpha \Z$. 
Then all the equations in (D) hold by substituting $\beta = 1, 2, 3, 4$. 
\qedhere

\end{proof}

\subsubsection{Cases $h=0, 1$}

\begin{lem}\label{l-678-E1-h=0}
We use Notation \ref{n-678-E1}. 
Then $h \neq 0$. 
\end{lem}

\begin{proof}
Suppose $h=0$. 
Let us derive a contradiction. 
By Lemma \ref{l-678-E1-ABC}(C), 
we get 
\[
(\beta \gamma -1) H_W^3 = (d+2)\gamma^2 +4\gamma  -(d-2),  
\]
i.e., 
\[
q(\gamma) := (d+2) \gamma^2 +(4-\beta H_W^3)\gamma + (-d+2+H_W^3) =0. 
\]
It is enough to show that $q(x) >0$ for every $x \in \Z_{>0}$. 
To this end, we shall prove (1)-(4) below.  
\begin{enumerate}
\item If $\beta \geq 2$, then $q(1) < q(2) < q(3) < \cdots$. 
\item If $\beta =1$, then $q(2) < q(3) < q(4) < \cdots$. 
\item $q(1)>0$. 
\item If $\beta =1$, then $q(2) >0$. 
\end{enumerate}
We now finish the proof by assuming (1)-(4). 
If $\beta \geq 2$, then 
(1) and (3) imply $0< q(1) < q(2)< q(3) < \cdots$. 
Hence we may assume $\beta =1$. By (3), it is enough to show $q(x) >0$ for every integer $x \geq 2$. This follows from (2) and (4). 
Therefore, what is remaining is  to prove  (1)-(4). 

\medskip

Let us show (1) and (2). 
By $d \geq 3$ (Lemma \ref{l-678-h0-d3}(1)) and $\beta H_W^3 \leq 22$ (Remark \ref{r-w-bound}), the symmetry axis $x=s$ of the quadratic function $q(x)$ satisfies 
\[
s:=\frac{\beta H_W^3-4}{2(d+2)} \leq \frac{22-4}{2(3+2)}<2, 
\]
and hence $q(2) < q(3) < q(4) < \cdots$. 
Thus (2) holds. 
If $\beta \geq 2$, then the inequality $\beta H_W^3 \leq 10$ (Remark \ref{r-E1-deg-genus}) implies 
\[
s=\frac{\beta H_W^3-4}{2(d+2)} \leq \frac{10-4}{2(3+2)}<1. 
\]
Hence (1) holds. 

Let us show (3). 
We have  
\[
q(1) = (d+2) \cdot 1^2  +(4-\beta H_W^3) \cdot 1 + (-d+2+H_W^3) = 8 -(\beta -1)H_W^3. 
\]
If $\beta=1$, then $q(1) =8 >0$. 
If $\beta =2$, then $q(1) >0$ holds by $H_W^3 \leq 5$. 
If $\beta = 3$ (resp. $\beta=4$), then $q(1) >0$ follows from
$H_W^2 =2$ (resp. $H_W^3=1$). 
Thus (3) holds. 

Let us show (4). Assume $\beta =1$. 
By $d \geq 3$ (Lemma \ref{l-678-h0-d3}(1)) and $H_W^3 =\beta H_W^3  \leq 22$ (Remark \ref{r-w-bound}),
we obtain 
  \[
q(2) = (d+2) \cdot 2^2  +(4-H_W^3) \cdot 2 + (-d+2+H_W^3) 
\]
\[
= 3d +18-H_W^3
\geq 3 \cdot 3 + 18 -22>0. 
\]  
Thus (4) holds. 
\qedhere






\end{proof}


\begin{lem}\label{l-678-E1-h=1}
We use Notation \ref{n-678-E1}. 
Then $h \neq 1$. 
\end{lem}

\begin{proof}
Suppose $h=1$. 
Let us derive a contradiction. 
By Lemma \ref{l-678-E1-ABC}(C), 
we get 
$(\beta \gamma -1) H_W^3 = d\gamma^2 -d$, 
i.e., 
\[
q(\gamma) := d \gamma^2 -  \beta H_W^3 \gamma + (-d+H_W^3) =0. 
\]
The symmetry axis $x=s$ of the quadratic function $q(x)$ is given by  
\[
s := \frac{\beta H_W^3}{2d}.
\]
By $h=1$, we have $d \geq 4$ (Lemma \ref{l-678-h0-d3}(2)). 
The following hold. 
\begin{enumerate}
\renewcommand{\labelenumi}{(\roman{enumi})}
    \item $q(1) = d - \beta H_W^3  +(-d+H_W^3) = -(\beta-1)H_W^3$. 
    \item $q(2) = 4d -2\beta H_W^3  +(-d+H_W^3) =  3d -(2\beta -1)H_W^3$. 
    \item $q(3) = 9d -3\beta H_W^3  +(-d+H_W^3) =  8d -(3\beta -1)H_W^3\geq 32 -(3\beta -1)H_W^3$. 
\end{enumerate}

\setcounter{step}{0}
\begin{step}\label{s1-678-E1-h=1}
It holds that $\beta \neq 1$. 
\end{step}

\begin{proof}[Proof of Step \ref{s1-678-E1-h=1}]
Suppose $\beta =1$. 
In this case, $\gamma =\beta\gamma= 1+\alpha \geq 2$ and 
\[
0 = q(\gamma) = d\gamma^2 - H_W^3 \gamma +(-d+H_W^3) = 
 d(\gamma^2-1) -H_W^3(\gamma -1), 
\]
which implies $d(\gamma+1) = H_W^3$. 
By Lemma \ref{l-678-E1-ABC}(A), we obtain $\frac{1}{2}H_W^3 = (g-d-1)\gamma^2 -d\gamma$, and hence 
\begin{equation}\label{e1-678-E1-h=1}
d(\gamma+1) = H_W^3 = 2\gamma ( (g-d-1)\gamma -d). 
\end{equation}

We now show 
\begin{enumerate}
    \item[(I)] $d \in \gamma \Z$ and 
    \item[(II)] $2(g-1) \in (\gamma +1)\Z$. 
\end{enumerate}
The assertion (I) is obtained by $d(\gamma+1) = 2\gamma ( (g-d-1)\gamma -d) \in \gamma \Z$. 
Let us show (II). 
We have  $2\gamma ( (g-d-1)\gamma -d) = d(\gamma+1) \in (\gamma+1)\Z$, 
and hence 
\[
0 \equiv 2 \gamma ((g-d-1)\gamma -d) \equiv -2((g-d-1)(-1) -d) = 2(g-1) \mod \gamma +1.
\]
Thus (II) holds. 

By Lemma \ref{l-678-h0-d3}(2) and  Lemma \ref{l-678-easy-bound}(3), we obtain 
\begin{equation}\label{e2-678-E1-h=1}
4 \leq d \leq h+ g-3 = g-2. 
\end{equation}

Assume $g=6$. 
In this case, we have $4 \leq d \leq g-2 = 4$, i.e., $d=4$. 
By (I) and $\gamma \geq 2$, we get $\gamma \in \{ 2, 4\}$. 
It follows from (II) that $\gamma =4$. Then (\ref{e1-678-E1-h=1}) leads to the following contradiction: 
\[
20 = 4 (4+1) = d(\gamma+1) = H_W^3 = 2 \cdot 4 \cdot ( (g-d-1)\gamma -d) \in 8\Z. 
\]

Assume $g =7$. 
By (II), $12 \in (\gamma+1)\Z$. 
This, together with $\gamma \geq 2$, implies $\gamma \in \{2, 3, 5, 11\}$. 
It follows from (\ref{e2-678-E1-h=1}) that $4 \leq d \leq  5$. 
By (I) and $\gamma \in \{2, 3, 5, 11\}$, we obtain 
$(d, \gamma) \in \{(4, 2), (5, 5)\}$. 
For each case, (\ref{e1-678-E1-h=1}) leads to the following contradiction. 
\begin{itemize}
\item If $(d, \gamma)=(4, 2)$, then  
$0< H_W^3 = 2\gamma ( (g-d-1)\gamma -d) = 2 \gamma ( (7-4-1) \cdot 2 -4)=0$. 
\item 
If $(d, \gamma)=(5, 5)$, then  
$0< H_W^3 = 2\gamma ( (g-d-1)\gamma -d) = 2 \gamma ( (7-5-1) \cdot 5 -5)=0$. 
\end{itemize}


Assume $g=8$. 
By (II), $14 \in (\gamma+1)\Z$. 
This, together with $\gamma \geq 2$, implies $\gamma \in \{6, 13\}$. 
It follows from (\ref{e2-678-E1-h=1}) that $4 \leq d \leq 6$. 
By (I) and $\gamma \in \{6, 13\}$,  we obtain $d=\gamma =6$. 
However, we get $22 \geq H_W^3 = d(\gamma +1) =6(6+1)$  ((\ref{e1-678-E1-h=1}), Corollary \ref{c-g-bound}), which is a contradiction.

Assume $g=9$. 
By (II), $16 \in (\gamma+1)\Z$. 
This, together with $\gamma \geq 2$, implies $\gamma \in \{3, 7, 15\}$. 
It follows from (\ref{e2-678-E1-h=1}) that $4 \leq d \leq 7$. 
By (I) and $\gamma \in \{3, 7, 15\}$, this implies 
$(\gamma, d) \in \{ (3, 6), (7, 7)\}$. 
However, we get $22 \geq H_W^3 = d(\gamma +1) \geq 6(3+1)$ ((\ref{e1-678-E1-h=1}), Corollary \ref{c-g-bound}), which is a contradiction.

Assume $g=10$. 
By (II), $18 \in (\gamma+1)\Z$. 
This, together with $\gamma \geq 2$, implies $\gamma \in \{2, 5, 8, 17\}$. 
It follows from (\ref{e2-678-E1-h=1}) that $4 \leq d \leq 8$. 
By (I) and $\gamma \in  \{2, 5, 8, 17\}$, this implies 
$(\gamma, d) \in \{ (2, 4), (2, 6),  (2, 8), (5, 5), (8, 8)\}$. 
By $22 \geq H_W^3 = d(\gamma +1)$  ((\ref{e1-678-E1-h=1}), Corollary \ref{c-g-bound}), 
we get $(\gamma, d, H_W^3) \in \{ (2, 4, 12), (2, 6, 18)\}$. 
By Lemma \ref{l-678-E1-ABC}(A), we obtain the following contradiction: 
\[
4\Z \not\ni \frac{1}{2}H_W^3 =
\frac{1}{2} \beta H_W^3 = (g-d +h -2)
\gamma^2 - (d -2h+2)\gamma   + (h-1) \in 4\Z. 
\]

Assume $g=11$. 
By (II), $20 \in (\gamma+1)\Z$. 
This, together with $\gamma \geq 2$, implies $\gamma \in \{3, 4, 9, 19\}$. 
It follows from (\ref{e2-678-E1-h=1}) that $4 \leq d \leq  9$. 
By (I) and $\gamma \in \{3, 4, 9, 19\}$, this implies 
$(\gamma, d) \in \{ (3, 6), (3, 9), (4, 4), (4, 8), (9, 9)\}$. 
By $22 \geq H_W^3 = d(\gamma +1)$ ((\ref{e1-678-E1-h=1}), Corollary \ref{c-g-bound}), 
we get $(\gamma, d)= (4, 4)$. 
However, (\ref{e1-678-E1-h=1}) leads to the following contradiction: 
\[
20 = d(\gamma+1) = H_W^3 = 2\gamma ( (g-d-1)\gamma -d) \in 2 \gamma \Z = 8\Z. 
\]

Assume $g=12$. 
By (II), $22 \in (\gamma+1)\Z$. 
This, together with $\gamma \geq 2$, implies $\gamma \in \{10, 21\}$. 
It follows from (\ref{e2-678-E1-h=1}) that $4 \leq d \leq 10$. 
By (I) and $\gamma \in \{10, 21\}$, this implies 
$(\gamma, d) =(10, 10)$. 
However, 
we get the following contradiction: $22 \geq H_W^3 = d(\gamma +1) = 10 \cdot 11$  ((\ref{e1-678-E1-h=1}), Corollary \ref{c-g-bound}). 
This completes the proof of Step \ref{s1-678-E1-h=1}. 
\end{proof}

\begin{step}\label{s2-678-E1-h=1}
If $\beta \geq 2$, then 
$(g, d, h, \beta, \gamma, H_W^3) = (8, 4, 1, 2, 2, 4)$. 
\end{step}

\begin{proof}[Proof of Step \ref{s2-678-E1-h=1}]
Assume $\beta \geq 2$. 
In this case, the following hold (Remark \ref{r-E1-deg-genus},  Lemma \ref{l-678-h0-d3}(2)): 
\[
s = \frac{\beta H_W^3}{2d} \leq \frac{10}{2 \cdot 4} <2.
\]
Therefore, $q(2)<q(3)< \cdots$. 
By (iii), we can check that $q(3) \geq 32 - (3\beta -1)H_W^3>0$ 
(if $\beta =2$ (resp. $\beta=3$, resp. $\beta=4$), then use $H_W^3 \leq 5$ (resp. $H_W^3=2$, resp. $H_W^3=1$)). 
It follows from (i) that $q(1) = -(\beta-1)H_W^3 <0$. 
By $q(\gamma)=0$, we obtain $\gamma =2$. 
Then (ii) implies  
\[
0 =q(2) = 3d -(2\beta -1)H_W^3. 
\]
If $\beta \in \{3, 4\}$, then we get the following contradiction: 
\[
0 = 3d -(2\beta -1)H_W^3 \geq 12 - (2\beta -1)H_W^3 \geq \min \{ 12 - (2\cdot 3-1)\cdot 2, 12-(2 \cdot 4-1)\cdot 1\} >0.
\]
Thus we get $\beta=2$. 
In this case, we obtain 
\[
0 = 3d -(2\beta -1)H_W^3 = 3(d-H_W^3), 
\]
i.e., $H_W^3 = d$. 
Note that $4 \leq d =H_W^3 \leq 5$. 
By Lemma \ref{l-678-E1-ABC}(A), we have 
\[
d = \frac{1}{2}\cdot 2 \cdot d = \frac{1}{2}\beta H_W^3 
= (g-d-1)\gamma^2 -d\gamma = 4( g-d-1)-2d=4g-6d-4. 
\]
In particular, $7d =4g-4 \in 4\Z$. 
By $4 \leq d \leq 5$, we obtain $d =4$, and hence $g = 8$. 
This completes the proof of Step \ref{s2-678-E1-h=1}. 
\qedhere



\end{proof}

\begin{step}\label{s3-678-E1-h=1}
It holds that $(g, d, h, \beta, \gamma, H_W^3) \neq (8, 4, 1, 2, 2, 4)$. 
\end{step}

\begin{proof}[Proof of Step \ref{s3-678-E1-h=1}]
 Suppose    
$(g, d, h, \beta, \gamma, H_W^3) = (8, 4, 1, 2, 2, 4)$.  
We obtain $\alpha = \beta\gamma -1 = 3$ and 
\[
D \sim -3K_Y -2E. 
\]
By Lemma \ref{l-678-easy-bound}(1), the following hold: 
\begin{itemize}
\item $(-K_Y)^3 = 2g-2d +2h -4 = 16 -8 +2-4=6$. 
\item $(-K_Y)^2 \cdot E =d -2h+2 = 4$.
\item $(-K_Y) \cdot E^2 = 2h-2 =0$.
\item $-E^3 = d + 2h-2 = 4$. 
\end{itemize}
We shall prove the following. 
\begin{enumerate}
    \item  $h^0(Y, -K_Y)  \geq 6$. 
    \item $h^0(E, -K_Y|_E)=4$. 
\end{enumerate}   
Assuming (1) and (2), we now finish the proof of Step \ref{s3-678-E1-h=1}. 
By the exact sequence 
\[
0 \to H^0(Y, -K_Y-E) \to H^0(Y, -K_Y) \to H^0(E, -K_Y|_E), 
\]
(1) and (2) imply $h^0(Y, -K_Y -E) >0$. 
This leads to the following contradiction: 
\[
0 = \kappa(Y, D) = \kappa (Y, -3K_Y -2E) 
\geq 
\kappa(Y, -K_Y)=3, 
\]

It is enough to show (1) and (2). The assertion (1) holds by the following \cite[Corollary 2.6]{TanI}: 
\[
h^0(Y, -K_Y) \geq h^0(Y, -K_Y) - h^1(Y, -K_Y) =\chi(Y, -K_Y) 
= \frac{(-K_Y)^3}{2} +3 = 6. 
\]
Let us show (2). 
By the Riemann--Roch theorem, we obtain 
\[
\chi(E, -K_Y|_E) = \chi(E, \MO_E) + \frac{1}{2}(-K_Y|_E) \cdot ((-K_Y|_E)-K_E) 
\]
\[
= \frac{1}{2}E \cdot (-K_Y) \cdot (-2K_Y-E)= (-K_Y)^2 \cdot E =4.
\]
Hence it suffices to prove $H^i(E, -K_Y|_E)=0$ for every $i>0$. 
Since $E$ satisfies Ramanujam vanishing \cite[Theorem 3(a)]{Muk13}, 
it is enough to show that $-K_Y|_E -K_E$ is nef and big. 
We can write $-K_Y|_E -K_E = (-2K_Y-E)|_E$. 
On the other hand, we have $K_Y =\tau^*K_W +D$, and hence 
\[
-\tau^*K_W = -K_Y+D \sim -4K_Y -2E. 
\]
Therefore, 
\[
-K_Y|_E -K_E = (-2K_Y-E)|_E \equiv \frac{1}{2}\tau^*(-K_W)|_E. 
\]
Since $E \to \tau(E)$ is birational, $\tau^*(-K_W)|_E$ is nef and big, as required. 
This completes the proof of Step \ref{s3-678-E1-h=1}. 
\end{proof}
Step \ref{s1-678-E1-h=1}, Step \ref{s2-678-E1-h=1}, Step \ref{s3-678-E1-h=1} complete the proof of Lemma \ref{l-678-E1-h=1}. 
\qedhere







\end{proof}

\subsubsection{Cases $h=2, 3, 4$}

\begin{nota}\label{n-f(x)}
We use Notation \ref{n-678-E1}. 
Set 
\[
f(x) := (g-d +h -2)
x^2 - (d -2h+2)x   + (h-1). 
\]
Let $x=s$ be the symmetry axis of the quadratic function $f(x)$, i.e.,  
\[
s :=
\frac{1}{2} \cdot  
\frac{d-2h+2}{g-d +h -2} \in \R. 
\]
Note that $d-2h+2 =(-K_Y)^2 \cdot E>0$ and $g-d+h-2 = \frac{(-K_Y)^3}{2}>0$ (Lemma \ref{l-678-easy-bound}(1)). In particular, $s>0$. 
By $f(\gamma) =\frac{1}{2}\beta H_W^3$ (Lemma \ref{l-678-E1-ABC}(A)), we obtain $1 \leq f(\gamma) \leq 11$ 
(Remark \ref{r-w-bound}). 
\end{nota}

\begin{lem}\label{l-s-2gamma+3}
We use Notation \ref{n-678-E1} and Notation \ref{n-f(x)}.  
Then the following hold. 
\begin{enumerate}
\item $f(-3) \geq 12$. 
\item $\gamma <2s+3$. 
\item $\gamma <-h+g+2$. 
\end{enumerate}
\end{lem}
\begin{proof}
Let us show (1). 
Since $-K_Y$ is ample, we have $g-d +h -2 = \frac{1}{2}(-K_Y)^3 >0$ and 
$d -2h+2 = (-K_Y)^2 \cdot E>0$ (Lemma \ref{l-678-easy-bound}(1)). 
These, together with $h \geq 2$ (Lemma \ref{l-678-E1-h=0}, Lemma \ref{l-678-E1-h=1}), imply 
\[
f(-3) = 9(g-d +h -2) + 3(d -2h+2)   + (h-1) \geq 9 \cdot 1 +3 \cdot 1 +1 = 12. 
\]
Thus (1) holds. 
Let us show (2). 
Since we have $f(s + r) = f(s-r)$ for every $r \in \R$, it holds that 
\[
12 \overset{{\rm (1)}}{\leq} f(-3) = f( s -(s+3)) = f(s+(s+3)) = f(2s+3). 
\]
By $1 \leq f(\gamma ) \leq 11$ {\cred and $f(0) = h-1\leq 9 \leq 11$ (Lemma \ref{l-678-easy-bound}(2))}, we obtain $\gamma <2s+3$. Thus (2) holds. 
Let us show (3). 
By $d\leq h+g-3$ (Lemma \ref{l-678-easy-bound}(3)), the following holds:  
\[
s 
=\frac{1}{2} \cdot  \frac{d-2h+2}{g-d +h -2} \leq \frac{1}{2} \cdot \frac{(h+g-3)-2h+2}{g-(h+g-3) +h -2} 
= \frac{1}{2}(-h+g-1). 
\]
This and (2) imply 
\[
\gamma < 2s+3 \leq 2 \cdot \frac{1}{2}(-h+g-1) +3 = -h+g+2. 
\]
Thus (3) holds. 
\end{proof}

\begin{lem}\label{l-678-gamma>1}
We use Notation \ref{n-678-E1}. 
Then $2 \leq \gamma \leq g-1 \leq 11$. 
\end{lem}

\begin{proof}
By $g \leq 12$, we get $g-1 \leq 11$. 
The inequality $\gamma \leq g-1$ follows from 
\[
\gamma < -h+g+2 \leq -2 +g+2 =g, 
\]
where the first (resp. second) inequality holds by 
Lemma \ref{l-s-2gamma+3}(3) 
(resp. Lemma \ref{l-678-E1-h=0} and Lemma \ref{l-678-E1-h=1}). 

It is enough to show $\gamma \neq 1$. 
Suppose $\gamma =1$. 
By $h \geq 2$  (Lemma \ref{l-678-E1-h=0}, Lemma \ref{l-678-E1-h=1}) and Lemma \ref{l-678-E1-ABC}(C), we get the following contradiction: 
\[
0< \alpha H_W^3   =(d -2h+2)\gamma^2 -2(2h-2)\gamma  -(d+2h-2) 
\]
\[
= (d -2h+2) -2(2h-2)  -(d+2h-2)  = -8h +8 <0. 
\]
Thus $\gamma \neq 1$.     
\end{proof}

\begin{lem}\label{l-678-E1-h=2}
We use Notation \ref{n-678-E1}. 
Then $h \neq 2$. 
\end{lem}

\begin{proof}
Suppose $h=2$. Let us derive a contradiction. 
We have $2 \leq \gamma \leq 11$ (Lemma \ref{l-678-gamma>1}). 
The following holds (Lemma \ref{l-678-easy-bound}(3)): 
\[
3= 2h-1 \leq d \leq h+g-3 = g-1 \leq 11. 
\]
By $h=2$, (A) and (C) of Lemma \ref{l-678-E1-ABC} imply the following: 
\begin{enumerate}
    \item[(A)] $\frac{1}{2} \beta H_W^3 = (g-d)
\gamma^2 - (d -2)\gamma   + 1$. 
\item[(C)] $\alpha H_W^3 = (\beta \gamma -1)H_W^3 = (d -2)\gamma^2 -4\gamma  -(d+2)
= (\gamma^2-1) d -2\gamma^2 -4\gamma -2$. 
By $\gamma^2 -1 \neq 0$, we obtain 
\[
d = \frac{2\gamma^2 +4\gamma +2 +(\beta \gamma -1)H_W^3}{\gamma^2-1} = 2 + \frac{4\gamma +4 +(\beta \gamma -1)H_W^3}{\gamma^2-1}. 
\]
\end{enumerate}

\medskip

Assume $\beta=4$. 
By (C) and $H_W^3 =1$, 
\[
d = 2 + \frac{8\gamma +3}{\gamma^2-1}. 
\]
We get $8 \gamma +3 \in (\gamma^2-1)\Z \subset (\gamma -1)\Z$. 
Taking modulo $\gamma -1$, we obtain $11 \equiv 8 \gamma +3 \equiv 0 \mod \gamma -1$. 
Hence $\gamma \in \{2, 12\}$. 
This, together with $2 \leq \gamma \leq 11$, implies $\gamma =2$. 
We then get a contradiction: $\Z \ni  \frac{8\gamma +3}{\gamma^2-1} = \frac{19}{3}\not\in \Z$.


\medskip

Assume $\beta =3$. 
By (C) and $H_W^3 =2$, 
\[
d = 2 + \frac{10\gamma +2}{\gamma^2-1}. 
\]
By $10 \gamma +2 \in (\gamma^2-1)\Z \subset (\gamma +1) \Z$, we obtain $-8 \equiv 10 \gamma +2 \equiv 0 \mod \gamma +1$, 
which implies that $\gamma +1$ is a divisor of $8$. 
This, together with $2 \leq \gamma \leq 11$, implies $\gamma \in \{ 3, 7\}$. 
Since $\gamma =7$ leads to a contradiction: $\Z \ni \frac{10\gamma +2}{\gamma^2-1} = \frac{72}{48}\not\in \Z$, 
we obtain $(\gamma, d) = (3, 6)$. 
Then  (A) leads to the following contradiction:
\[
3 = \frac{1}{2}\beta H_W^3 = (g-d)\gamma^2 -(d-2)\gamma +1 \equiv 1 \mod 3. 
\]


\medskip

Assume $\beta =2$. 
We first treat the case when $\gamma =2$. 
By (C), $d = 2 + \frac{4\gamma +4 +(\beta \gamma -1)H_W^3}{\gamma^2-1}
= 2 + \frac{8+4 + 3H_W^3}{3} = H_W^3 +6$. 
By (A), 
\[
H_W^3 =\frac{1}{2} \beta H_W^3 = (g-d)
\gamma^2 - (d -2)\gamma   + 1 = 4(g -H_W^3 -6)-2(H_W^3 +4)+1 = 4g -6H_W^3 -31. 
\]
Then $4g = 7H_W^3 +31$, and hence $H_W^3 \equiv 3\mod 4$. 
By $1 \leq H_W^3 \leq 5$, we get $(H_W^3, g)= (3, 13)$, which contradicts $6 \leq g \leq 12$. 
In what follows, we assume $\gamma \geq 3$ until the end 
of this paragraph. 
By (C), 
 \[
 d = 2 + \frac{4\gamma +4 +(2\gamma -1)H_W^3}{\gamma^2-1}. 
 \]
Assume $H_W^3=1$. Then $d = 2 + \frac{6\gamma +3}{\gamma^2-1}$. 
If $\gamma \geq 7$, then we get a contradiction: $0< 6\gamma +3 < \gamma^2-1$. 
We then get  $3 \leq \gamma \leq 6$. 
We have $6\gamma +3 \in (\gamma^2-1)\Z \subset (\gamma -1)\Z$. 
Taking modulo $\gamma -1$, we obtain $9 \equiv 6\gamma +3 \equiv 0 \mod \gamma -1$. 
By $3 \leq \gamma \leq 6$, we get $\gamma =4$, 
which leads to a contradiction $\Z \ni \frac{6\gamma +3}{\gamma^2-1}=\frac{27}{15} \not\in \Z$. 
Assume $H_W^3 =2$. 
Then $d = 2 + \frac{8\gamma +2}{\gamma^2-1}$. 
We have $8\gamma +2 \in (\gamma^2-1)\Z \subset (\gamma-1)\Z$. 
Taking modulo $\gamma -1$, we obtain $10 \equiv 8\gamma +2 \equiv 0 \mod \gamma-1$. 
By $3 \leq \gamma \leq 11$, $\gamma \in \{ 3, 6, 11\}$. 
However, none of these  satisfies $\frac{8\gamma +2}{\gamma^2-1} \in \Z$. 
Assume $H_W^3 =3$. 
Then $d = 2 + \frac{10\gamma +1}{\gamma^2-1}$. 
We have $10\gamma +1 \in (\gamma^2-1)\Z \subset (\gamma-1)\Z$. 
Taking modulo $\gamma -1$, we obtain $11 \equiv 10\gamma+1 \equiv 0 \mod \gamma -1$, which implies $\gamma \in \{2, 12\}$. 
By $3 \leq \gamma \leq 11$, this is absurd. 
Assume $H_W^3 =4$. 
Then $d = 2 + \frac{12\gamma}{\gamma^2-1}$. 
Since $\gamma$ and $\gamma^2-1$ are coprime, we obtain $12 \in (\gamma^2-1)\Z$. 
By $3 \leq \gamma \leq 11$, this is absurd.  
Assume $H_W^3 =5$. 
Then $d = 2 + \frac{14\gamma-1}{\gamma^2-1}$. 
We have $14\gamma -1 \in (\gamma^2-1)\Z \subset (\gamma-1)\Z$. 
Taking modulo $\gamma -1$, we see that $13 \equiv 14\gamma -1 \equiv 0 \mod \gamma -1$. 
Then $\gamma \in \{2, 14\}$, which contradicts 
$3 \leq \gamma \leq 11$. 
This completes the proof for the case when $\beta =2$.

\medskip

Assume $\beta =1$. In this case, $H_W^3 =(-K_W)^3 \in 2\Z$. 
By (C), 
\[
d 
= 2 + \frac{4\gamma -4 +(\gamma -1)H_W^3 +8}{\gamma^2-1}. 
\]
In particular, $8 \in (\gamma -1)\Z$, and hence $\gamma \in \{2, 3, 5, 9\}$. 
Moreover, 
\[
d = 2 + \frac{4 +H_W^3 +\frac{8}{\gamma-1}}{\gamma+1}. 
\]

Assume $\gamma =9$. 
Then $d = 2 + \frac{H_W^3 +5}{10}$. 
This is absurd, because $H_W^3 \in 2\Z$ implies $H_W^3+5\not\in 2\Z$.

Assume $\gamma =5$. 
Then $d = 2 + \frac{H_W^3 +6}{6}$, which implies $H_W^3 \in 6\Z$. 
By $2 \leq H_W^3 \leq 22$, we get $(H_W^3, d) = (6, 4), (12, 5), (18, 6)$. 
By (A), we have 
\[
\frac{1}{2}H_W^3 = (g-d)\gamma^2 -(d-2)\gamma +1 =25(g-d) -5(d-2)+1 = 25g -30d +11.  
\]
If $(H_W^3, d)= (6, 4)$, then $3 =\frac{1}{2} H_W^3 =25g -30\cdot 4 +11 \equiv 16 \mod 25$, which is absurd. 
If $(H_W^3, d)= (12, 5)$, then $6=\frac{1}{2} H_W^3 = 25g -30 \cdot 5 +11 \equiv 11 \mod 25$, which is a contradiction. 
If $(H_W^3, d)= (18, 6)$, then  $9=\frac{1}{2} H_W^3 = 25g -30 \cdot 6 +11 \equiv 1 \not\equiv 9 \mod 5$, which is absurd.

Assume $\gamma =3$. 
Then $d = 2 + \frac{H_W^3 +8}{4} =\frac{H_W^3}{4} +4$. By $2 \leq H_W^3 \leq 22$, we get $5 \leq d\leq 9$. 
It follows from (A) and 
$d =\frac{H_W^3}{4} +4$ that 
\[
2d-8  = \frac{1}{2}H_W^3 = (g-d)\gamma^2 -(d-2)\gamma +1 =9(g-d) -3(d-2)+1 = 9g -12d +7.  
\]
Hence $9g = 14d -15$. 
Note that $(g, d)=(3, 3)$ is a solution of this equation. 
We then obtain $d \in 3 + 9\Z$ for an arbitrary solution $(g, d) \in \Z^2$, and hence there is no solution with $5 \leq d \leq 9$.

Assume $\gamma =2$. 
Then $d = 2 + \frac{H_W^3 +12}{3} =\frac{H_W^3}{3} +6$. By $2 \leq H_W^3 \leq 22$, we get $7 \leq d\leq 13$. 
It follows from (A) and 
$d =\frac{H_W^3}{3} +6$ that 
\[
\frac{3}{2}d-9  = \frac{1}{2}H_W^3 = (g-d)\gamma^2 -(d-2)\gamma +1 =4(g-d) -2(d-2)+1 = 4g -6d +5.  
\]
Hence $8g = 15d -28$. 
Note that $(g, d)=(4, 4)$ is a solution of this equation. 
We then obtain $d \in 4 + 8\Z$ for an arbitrary solution $(g, d) \in \Z^2$. 
By  $7 \leq d\leq 13$, we get $d=12$, and hence $g = \frac{15 \cdot 12 -28}{8} = \frac{15 \cdot 3 -7}{2} >12$. 
This contradicts $6 \leq g \leq 12$. 
\end{proof}

\begin{lem}\label{l-678-E1-h=3}
We use Notation \ref{n-678-E1}. 
Then $h \neq 3$. 
\end{lem}

\begin{proof}
Suppose $h=3$. Let us derive a contradiction. 
We have $2 \leq \gamma \leq 11$ (Lemma \ref{l-678-gamma>1}). 
The following holds (Lemma \ref{l-678-easy-bound}(3)): 
\begin{equation}\label{e1-678-E1-h=3}
5= 2h-1 \leq d \leq h+g-3 = g \leq 12. 
\end{equation}
By $h=3$, (A) and (C) of Lemma \ref{l-678-E1-ABC} imply the following: 
\begin{enumerate}
    \item[(A)] 
    $\frac{1}{2} \beta H_W^3 = (g-d +1)\gamma^2 - (d -4)\gamma   + 2$. 
\item[(C)] $\alpha H_W^3 = (\beta \gamma -1)H_W^3 = (d -4)\gamma^2 -8\gamma  -(d+4) 
=(\gamma^2 -1)d -4\gamma^2-8\gamma -4$. 
By $\gamma^2 -1 \neq 0$, we obtain 
\[
d = \frac{4\gamma^2+8\gamma +4 + (\beta \gamma -1)H_W^3 }{\gamma^2-1} = 4 + \frac{8\gamma +8 + (\beta \gamma -1)H_W^3 }{\gamma^2-1}.
\]
\end{enumerate}

\medskip

Assume $\beta=4$. 
By (C) and $H_W^3 =1$, 
\[
d = 4 + \frac{12\gamma +7}{\gamma^2-1}. 
\]
We get $12 \gamma +7 \in (\gamma^2-1)\Z \subset (\gamma -1)\Z$. 
Taking modulo $\gamma -1$, we obtain $19 \equiv 12 \gamma +7 \equiv 0  \mod \gamma -1$. 
Hence $\gamma \in \{2, 20\}$. 
By $2 \leq \gamma \leq 11$, we get $\gamma =2$, which leads to a contradiction $\Z \ni \frac{12\gamma +7}{\gamma^2-1} = \frac{12 \cdot 2 +7}{2^2-1} \not\in \Z$. 

\medskip

Assume $\beta =3$. 
By (C) and $H_W^3 =2$, 
\[
d = 4 + \frac{14\gamma +6}{\gamma^2-1}. 
\]
By $14 \gamma +6 \in (\gamma^2-1)\Z \subset (\gamma +1) \Z$, we obtain $-8 \equiv 14 \gamma +6 \equiv 0 \mod \gamma +1$, 
which implies that $\gamma +1$ is a divisor of $8$. 
This, together with $2 \leq \gamma \leq 11$, implies $\gamma \in \{ 3, 7\}$. 
Since $\gamma =7$ leads to a contradiction: $\Z \ni \frac{14\gamma +6}{\gamma^2-1} = \frac{14 \cdot 7 +6}{48}\not\in \Z$, 
we obtain $\gamma=3$. 
Then  (A) leads to the following contradiction:
\[
3 = \frac{1}{2}\beta H_W^3 = (g-d+1)\gamma^2 -(d-4)\gamma +2 \equiv 2 \mod 3. 
\]

\medskip

Assume $\beta =2$. 
We first treat the case when $\gamma =2$. 
By (C), $d = 4 + \frac{16+8 + 3H_W^3}{3} = H_W^3 +12 \geq 13$. 
This contradicts $d\leq 12$ (\ref{e1-678-E1-h=3}). 
In what follows, we assume $3 \leq \gamma \leq 11$ until the end of this paragraph. 
By (C), 
\[
d = 4 + \frac{8\gamma +8 +(2\gamma -1)H_W^3}{\gamma^2-1}. 
\]
Assume $H_W^3=1$. Then 
$d = 4 + \frac{10\gamma +7}{\gamma^2-1}$. 
We have $10\gamma +7 \in (\gamma^2-1)\Z$. 
Taking modulo $\gamma -1$, we obtain $17 \equiv 10 \gamma +7 \equiv 0 \mod \gamma -1$. 
We then get $\gamma \in \{2, 18\}$, which contradicts $3 \leq \gamma \leq 11$. 
Assume $H_W^3 =2$. 
Then $d = 4 + \frac{12\gamma +6}{\gamma^2-1}$. 
We have $12\gamma +6 \in (\gamma^2-1)\Z$. 
Taking modulo $\gamma +1$, we obtain $-6 \equiv 12 \gamma +6 \equiv 0 \mod \gamma+1$. 
By $3 \leq \gamma \leq 11$, we get $\gamma  = 5$, 
which leads to a contradiction: $\Z \ni \frac{12\gamma +6}{\gamma^2-1} = \frac{66}{24} \not\in \Z$. 
Assume $H_W^3 =3$. 
Then $d = 4 + \frac{14\gamma +5}{\gamma^2-1}$. 
We have $14\gamma +5 \in (\gamma^2-1)\Z$. 
Taking modulo $\gamma -1$, 
we obtain $19 \equiv 14 \gamma +5 \equiv 0 \mod \gamma -1$, which implies $\gamma \in \{2, 20\}$. 
By $3 \leq \gamma \leq 11$, this is absurd. 
Assume $H_W^3 =4$. 
Then $d = 4 + \frac{16\gamma+4}{\gamma^2-1}$. 
We have $16\gamma +4 \in (\gamma^2-1)\Z$. 
Taking modulo $\gamma +1$, we obtain $-12 \equiv 16\gamma +4 
\equiv 0 \mod \gamma +1$. 
This, together with $3 \leq \gamma \leq 11$, implies $\gamma \in \{3, 5, 11\}$. 
If $\gamma \in \{3, 5\}$, then $\gamma^2-1 \in 8 \Z$, which contradicts $16\gamma+4 \not\in 8\Z$. 
If $\gamma =11$, then $\Z \ni \frac{16\gamma+4}{\gamma^2-1} = \frac{180}{120}\not\in \Z$, which is absurd. 
Assume $H_W^3 =5$. 
Then $d = 4 + \frac{18\gamma+3}{\gamma^2-1}$. 
We have $18\gamma +3 \in (\gamma^2-1)\Z$. 
Taking modulo  $\gamma +1$, we see that  $-15 \equiv 18 \gamma +3 \equiv 0 \mod \gamma +1$. 
This, together with $3 \leq \gamma \leq 11$, implies $\gamma =4$. 
By (A), we get the following contradiction: 
\[
5 =
\frac{1}{2} \beta H_W^3 = (g-d +1)
\gamma^2 - (d -4)\gamma   + 2 
\in 2\Z. 
\]

\medskip

Assume $\beta =1$. 
By (C), 
\[
d =  4 + \frac{8\gamma +8 +(\gamma -1)H_W^3}{\gamma^2-1} = 
4 + \frac{8\gamma -8 +(\gamma -1)H_W^3 +16}{\gamma^2-1}. 
\]
In particular, $16 \in (\gamma -1)\Z$. 
This, together with $2 \leq \gamma \leq 11$, implies $\gamma \in \{2, 3, 5, 9\}$. 
Moreover, 
\[
d = 4 + \frac{8 +H_W^3 +\frac{16}{\gamma-1}}{\gamma+1}. 
\]

Assume $\gamma =9$. 
Then $d = 4 + \frac{H_W^3 +10}{10} = 5 + \frac{H_W^3}{10}$. 
By $2 \leq H_W^3 \leq 22$, we get $(H_W^3, d) = (10, 6), (20, 7)$. 
Then  (A) leads to the following contradiction: 
\[
\{ 5, 10\} \ni \frac{1}{2}H_W^3 = (g-d+1)\gamma^2 -(d-4)\gamma +2 \equiv 2 \mod 9. 
\]

\medskip 

Assume $\gamma =5$. 
Then $d = 4 + \frac{H_W^3 +12}{6} = 6 + \frac{H_W^3}{6}$. 
By $2 \leq H_W^3 \leq 22$, we get $H_W^3 \in \{ 6, 12, 18\}$. 
Then  (A) leads to the following contradiction: 
\[
\{3, 6, 9 \} \ni \frac{1}{2}H_W^3 = (g-d+1)\gamma^2 -(d-4)\gamma +2 \equiv 2 \mod 5. 
\]

\medskip

Assume $\gamma =3$. 
Then $d = 4 + \frac{H_W^3 +16}{4} = 8 + \frac{H_W^3}{4}$. 
By $2 \leq H_W^3 \leq 22$ and $d \leq 12$ (\ref{e1-678-E1-h=3}), 
we get $(H_W^3, d) \in \{(4, 9), (8, 10), (12, 11), (16, 12)\}$. 
By (A), we have 
\[
\frac{1}{2}H_W^3 = (g-d+1)\gamma^2 -(d-4)\gamma +2 =9(g-d+1) -3(d-4)+2 = 9g -12d +23.  
\]
In particular, $\frac{1}{2}H_W^3 \equiv 2 \mod 3$. 
Hence we obtain $(H_W^3, d) \in \{(4, 9), (16, 12)\}$. 
If $(H_W^3, d) =(4, 9)$, then $2 = \frac{1}{2}H_W^3 =9g -12d +23 
= 9g -85 \equiv 5 \mod 9$, which is a contradiction. 
If $(H_W^3, d) =(16, 12)$,  
then $8 = \frac{1}{2}H_W^3 =9g -12d +23 = 9g -12 \cdot 12 +23 \equiv 5 \mod 9$, which is absurd. 

Assume $\gamma =2$. 
Then $d = 4 + \frac{H_W^3 +24}{3} = 12 + \frac{H_W^3}{3}>12$. 
This contradicts $d \leq 12$ (\ref{e1-678-E1-h=3}). 
\end{proof}

\begin{lem}\label{l-678-E1-h=4}
We use Notation \ref{n-678-E1}. 
Then $h \neq 4$. 
\end{lem}

\begin{proof}
Suppose $h=4$. Let us derive a contradiction. 
We have $2 \leq \gamma \leq 11$ (Lemma \ref{l-678-gamma>1}). 
The following holds (Lemma \ref{l-678-easy-bound}(3)): 
\begin{equation}\label{e1-678-E1-h=4}
7= 2h-1 \leq d \leq h+g-3 = g +1\leq 13. 
\end{equation}
By $h=4$, (A) and (C) of Lemma \ref{l-678-E1-ABC} imply the following: 
\begin{enumerate}
    \item[(A)] $\frac{1}{2} \beta H_W^3 = (g-d +2)
\gamma^2 - (d -6)\gamma   + 3$. 
\item[(C)] $\alpha H_W^3 = (\beta \gamma -1)H_W^3 = (d -6)\gamma^2 -12\gamma  -(d+6) 
=(\gamma^2 -1)d -6\gamma^2-12\gamma -6$. 
By $\gamma^2 -1 \neq 0$, we obtain 
\[
d = \frac{6\gamma^2+12\gamma +6 + (\beta \gamma -1)H_W^3 }{\gamma^2-1} 
= 6 + \frac{12\gamma +12 + (\beta \gamma -1)H_W^3 }{\gamma^2-1}.
\]
\end{enumerate}

\medskip

Assume $\beta=4$. 
By (C) and $H_W^3=1$, 
\[
d = 6 + \frac{16\gamma +11}{\gamma^2-1}. 
\]
We get $16 \gamma +11 \in (\gamma^2-1)\Z \subset (\gamma +1)\Z$. 
Taking modulo $\gamma +1$, we obtain $-5 \equiv 16 \gamma +11 \equiv 0 \mod \gamma +1$. 
This, together with $2 \leq \gamma \leq 11$, implies 
$\gamma =4$. 
By (A), we get 
\[
2=
\frac{1}{2} \beta H_W^3 = (g-d +2)
\gamma^2 - (d -6)\gamma   + 3 \equiv 3 \mod 4, 
\]
which is a contradiction.


\medskip

Assume $\beta=3$. 
By (C) and $H_W^3 =2$, 
\[
d = 6 + \frac{18\gamma +10}{\gamma^2-1}. 
\]
We get $18 \gamma +10 \in (\gamma^2-1)\Z \subset (\gamma +1)\Z$. 
Taking modulo $\gamma +1$, we obtain $-8 \equiv 18 \gamma +10 \equiv 0 \mod \gamma +1$. 
This, together with $2 \leq \gamma \leq 11$, implies $\gamma \in \{3, 7\}$. 
If $\gamma =7$, then we get $\Z \ni \frac{18\gamma +10}{\gamma^2-1} = \frac{136}{48} \not\in \Z$, 
which is absurd. 
Thus $\gamma =3$, which implies $d  = 6 + \frac{18\gamma +10}{\gamma^2-1} = 14$. 
This contradicts $7 \leq d \leq 13$ (\ref{e1-678-E1-h=4}).

\medskip

Assume $\beta=2$. 
By (C), 
\[
d = 6 + \frac{12\gamma +12 +(2\gamma-1)H_W^3}{\gamma^2-1}. 
\]
If $\gamma=2$, then we have $d = 6 + \frac{24 +12 +3H_W^3}{3} = H_W^3 +18 \geq 14$, 
which contradicts $7 \leq d \leq 13$ (\ref{e1-678-E1-h=4}). 
In what follows, we assume $3 \leq \gamma \leq 11$ until the end of  this paragraph. 
Assume $H_W^3=1$. Then $d = 6 + \frac{14\gamma +11}{\gamma^2-1}$. 
By $14\gamma +11 \in (\gamma^2-1)\Z \subset (\gamma+1)\Z$, 
we get $-3 \equiv 14\gamma +11 \equiv 0 \mod \gamma +1$. 
By $3 \leq \gamma \leq 11$, this is absurd. 
Assume $H_W^3=2$. Then $d = 6 + \frac{16\gamma +10}{\gamma^2-1}$. 
By $16\gamma +10 \in (\gamma^2-1)\Z \subset (\gamma+1)\Z$, 
we get $-6 \equiv 16\gamma +10 \equiv 0 \mod \gamma +1$. 
By $3 \leq \gamma \leq 11$, we obtain $\gamma =5$. 
This is a contradiction, because $16 \gamma +10 \in (\gamma^2-1)\Z=24\Z \subset 8\Z$. 
Assume $H_W^3=3$. Then $d = 6 + \frac{18\gamma +9}{\gamma^2-1}$. 
By $18\gamma +9 \in (\gamma^2-1)\Z \subset (\gamma+1)\Z$, 
it holds that $-9 \equiv 18\gamma +9 \equiv 0 \mod \gamma +1$. 
By $3 \leq \gamma \leq 11$, we get $\gamma =8$. 
This leads to a contradiction: $\Z \ni \frac{18\gamma +9}{\gamma^2-1} = \frac{9 \cdot (2 \cdot 8+1)}{9 \cdot 7} =\frac{17}{7}\not\in \Z$. 
Assume $H_W^3=4$. Then $d = 6 + \frac{20\gamma +8}{\gamma^2-1}$. 
By $20\gamma +8 \in (\gamma^2-1)\Z \subset (\gamma -1)\Z$, 
it holds that  $28\equiv 20\gamma +8 \equiv 0 \mod \gamma -1$. 
By $3 \leq \gamma \leq 11$, we get $\gamma  \in \{ 3, 5, 8\}$. 
If $\gamma =3$, then $\Z \ni \frac{20\gamma +8}{\gamma^2-1} = \frac{68}{8} \not\in \Z$, which is absurd. 
If $\gamma =5$, then $\Z \ni \frac{20\gamma +8}{\gamma^2-1} = \frac{108}{24} \not\in \Z$, which is a contradiction. 
If $\gamma =8$, then $\Z \ni \frac{20\gamma +8}{\gamma^2-1} = \frac{8 \cdot 21}{63} =\frac{8}{3} \not\in \Z$, which is absurd. 
Assume $H_W^3=5$. Then $d = 6 + \frac{22\gamma +7}{\gamma^2-1}$. 
By $22\gamma +7 \in (\gamma^2-1)\Z \subset (\gamma -1)\Z$, 
we get $29 \equiv 22\gamma +7 \equiv 0 \mod \gamma -1$, which contradicts $3 \leq \gamma \leq 11$. 

\medskip

Assume $\beta =1$. In this case, $H_W^3 =(-K_W)^3 \in 2\Z$. 
By (C), 
\[
d = 
6 + \frac{12\gamma -12 +(\gamma -1)H_W^3 +24}{\gamma^2-1}. 
\]
In particular, $24 \in (\gamma -1)\Z$. 
This, together with $2 \leq \gamma \leq 11$, implies $\gamma \in \{2, 3, 4, 5, 7, 9\}$. 
Moreover, 
\[
d = 6 + \frac{12 +H_W^3 +\frac{24}{\gamma-1}}{\gamma+1}. 
\]

Assume $\gamma =9$. 
Then $d = 6 + \frac{H_W^3 +15}{10}$. 
This leads to a contradiction: $1+2\Z  \ni H_W^3+15 \in 10\Z$. 

Assume $\gamma =7$. 
Then $d = 6 + \frac{H_W^3 +16}{8} =8+ \frac{H_W^3}{8}$. 
It follows from $2 \leq H_W^3 \leq 22$ that $H_W^3 \in \{ 8, 16\}$. 
Then  (A) leads to the following contradiction:  
\[
\{4, 8\} \ni 
\frac{1}{2} H_W^3 = (g-d +2)
\gamma^2 - (d -6)\gamma   + 3 
\equiv 3 \mod 7. 
\]

Assume $\gamma =5$. 
Then $d = 6 + \frac{H_W^3 +18}{6} =9+ \frac{H_W^3}{6}$. 
It follows from $2 \leq H_W^3 \leq 22$ that $(H_W^3, d) \in \{(6, 10), (12, 11), (18, 12) \}$. 
By (A), 
\[
\frac{1}{2} H_W^3 = (g-d +2)
\gamma^2 - (d -6)\gamma   + 3 = 25(g-d +2) -5(d -6)   + 3\equiv 3 \mod 5. 
\]
Hence $(H_W^3, d) = (6, 10)$. We then get 
\[
3 =\frac{1}{2} H_W^3 = 25(g-10 +2) -5(10 -6)   + 3 \equiv (-5) \cdot (-6) +3 \mod 25, 
\]
which is absurd.

Assume $\gamma =4$. 
Then $d = 6 + \frac{H_W^3 +20}{5} =10+ \frac{H_W^3}{5}$. 
It follows from $2 \leq H_W^3 \leq 22$ and $H_W^3 \in 2\Z$ that 
$(H_W^3, d) \in \{(10, 12), (20, 14)\}$. 
By $7 \leq d \leq 13$ (\ref{e1-678-E1-h=4}), we get $(H_W^3, d)=(10, 12)$. 
By (A), we get
\[
5=\frac{1}{2} H_W^3 = (g-d +2)
\gamma^2 - (d -6)\gamma   + 3 \equiv 3
\not\equiv 5 \mod 4. 
\]
This is absurd.

Assume $\gamma =3$. 
Then $d = 6 + \frac{H_W^3 +24}{4} =12+ \frac{H_W^3}{4}$. 
It follows from $2 \leq H_W^3 \leq 22$ and $7 \leq d \leq 13$ (\ref{e1-678-E1-h=4}) that 
$(H_W^3, d) =(4, 13)$. 
Then  (A) leads to the following contradiction:  
\[
3\Z \not\ni 
2=\frac{1}{2} H_W^3 = (g-d +2)
\gamma^2 - (d -6)\gamma   + 3 \in 3\Z.
\]

Assume $\gamma =2$. 
Then $d = 6 + \frac{H_W^3 +36}{3} 
>13$. 
This contradicts $7 \leq d \leq 13$ (\ref{e1-678-E1-h=4}). 
\end{proof}

\begin{prop}\label{p-678-E1-h>4}
We use Notation \ref{n-678-E1}. 
Then $5 \leq h \leq 10$. 
\end{prop}

\begin{proof}
By Lemma \ref{l-678-easy-bound}(2), we have $h \leq 10$. 
The other inequality $5 \leq h$ holds as follows: 
\begin{itemize}
\item $h \neq 0$: Lemma \ref{l-678-E1-h=0}. 
\item $h\neq 1$: Lemma \ref{l-678-E1-h=1}. 
\item $h \neq 2$: Lemma \ref{l-678-E1-h=2}. 
\item $h \neq 3$: Lemma \ref{l-678-E1-h=3}. 
\item $h \neq 4$: Lemma \ref{l-678-E1-h=4}. 
\end{itemize}
\end{proof}

\subsubsection{Case $h\geq 5$}

\begin{lem}\label{l-d=h+g-3}
We use Notation \ref{n-678-E1}. 
If $d = h+g-3$, then $\gamma \leq 5$.  
\end{lem}

\begin{proof}
Suppose $d = h+g-3$ and $\gamma \geq 6$. 
It suffices to derive a contradiction. 
We use Notation \ref{n-f(x)}. 
By $d = h+g-3$, we obtain  
\[
f(x) = x^2 - (g -h-1)x   + (h-1)
\]
and 
\begin{equation}\label{e1-d=h+g-3}
s 
=\frac{1}{2} \cdot  \frac{d-2h+2}{g-d +h -2} = \frac{1}{2} \cdot \frac{(h+g-3)-2h+2}{g-(h+g-3) +h -2} 
= \frac{1}{2}(-h+g-1). 
\end{equation}
By $g \leq 12$ and $h \geq 5$ (Proposition \ref{p-678-E1-h>4}), 
(\ref{e1-d=h+g-3}) implies 
\begin{equation}\label{e2-d=h+g-3}
s 
=\frac{1}{2}(-h+g-1) \leq \frac{1}{2}(-5+12-1) = 3.
\end{equation}
In particular, 
\begin{equation}\label{e3-d=h+g-3}
f(6)< f(7)<f(8). 
\end{equation}

We now show (i)-(iii) below. 
\begin{enumerate}
\renewcommand{\labelenumi}{(\roman{enumi})}
    \item $6 \leq \gamma \leq 8$. 
    \item $h+6 \leq g \leq 12$. 
    \item $h =5$ or $h=6$. 
\end{enumerate}
Let us show (i). 
It follows from (\ref{e2-d=h+g-3}) and Lemma \ref{l-s-2gamma+3}(2) that $\gamma <2s+3\leq 9$. 
Hence (i) holds. 
Let us show (ii). 
By (i) and $\gamma <2s+3$, 
we obtain $s > \frac{1}{2} (\gamma-3) \geq \frac{3}{2}$. 
This, together with (\ref{e1-d=h+g-3}), implies $g > h +4$. 
In order to prove (ii), it is enough to exclude the case when $g = h+5$. 
Suppose $g =h+5$. In this case, we obtain 
$f(x) = x^2 -4x +(h-1)$ and $f(6) = 36 -24 +(h-1) = h +11 >11$. 
Hence  $11< f(6)<f(7) < f(8)$ (\ref{e3-d=h+g-3}), which contradicts (i) and 
$f(\gamma) \leq 11$ (Notation \ref{n-f(x)}). 
Thus (ii) holds. 
The assertion (iii) follows from (ii) and $h \geq 5$ (Proposition \ref{p-678-E1-h>4}). 
This completes the proof of (i)-(iii). 

\medskip

\begin{claim}\label{cl-d=h+g-3}
The quintuple $(h, g, \gamma, \frac{1}{2}\beta H_W^3, d)$ 
is equal to one of the following. 
\begin{enumerate}
\item $(6, 12, 6, 11, 15)$. 
\item $(5, 12, 6, 4, 14)$. 
\item $(5, 12, 7, 11, 14)$. 
\item $(5, 11, 6, 10, 13)$. 
\end{enumerate}
\end{claim}

\begin{proof}[Proof of Claim \ref{cl-d=h+g-3}] 
By $d = h+g-3$, it is enough to compute the quadruple 
$(h, g, \gamma, \frac{1}{2}\beta H_W^3)$. 
By (ii) and (iii), one of the following holds. 
\begin{enumerate}
\item[(a)] $h=6$ and $g=12$. 
\item[(b)] $h=5$ and $g=12$. 
\item[(c)] $h=5$ and $g=11$. 
\end{enumerate}

(a) Assume $h=6$ and $g=12$. 
We have $f(x)  = x^2 - (g -h-1)x   + (h-1)= x^2 - 5x +5$. 
Hence $f(6)<f(7)<f(8)$ and $f(6) = 36-30 +5 =11$. 
By (i) and $f(\gamma) \leq 11$ (Notation \ref{n-f(x)}), 
we get $\gamma =6$ and $\frac{1}{2}\beta H_W^3 = f(\gamma) = 11$. 
Hence (1) is a unique solution when (a) holds.

(b) 
Assume  $h=5$ and $g=12$. 
We have 
$f(x)  = x^2 - (g -h-1)x   + (h-1)= x^2 - 6x +4$. 
Hence $f(6)<f(7)<f(8)$, $f(6) = 36-36 +4 =4$, and 
$f(7) = 49 - 42 +4 = 11$. 
In this case, we get the solutions (2) and (3) corresponding to 
$f(6) =4$ and $f(7) = 11$, respectively. 

(c) 
Assume $h=5$ and $g=11$. 
We have 
$f(x)  = x^2 - (g -h-1)x   + (h-1)=  x^2 - 5x +4$. 
Hence $f(6)<f(7)<f(8)$, $f(6) = 36-30+4 =10$, and 
$f(7) = 49 - 35 +4 > 11$. 
 In this case, we get the solution (4) corresponding to 
$f(6) = 10$. 
This completes the proof of Claim \ref{cl-d=h+g-3}. 
\qedhere

\end{proof}

By $\frac{1}{2}\beta H_W^3 \in \{ 4, 10, 11\}$ (Claim \ref{cl-d=h+g-3}), 
we obtain $\beta \in \{1, 2\}$ (cf. Remark \ref{r-E1-deg-genus}).





Suppose $\beta =2$. It follows from  $1 \leq \frac{1}{2}\beta H_W^3 = H_W^3 \leq 5$ that Claim \ref{cl-d=h+g-3}(2) holds. 
Then Lemma \ref{l-678-E1-ABC}(C) leads to the following contradiction: 
\[
44=(2 \cdot 6 -1) \cdot 4=(\beta \gamma -1)H_W^3 = (d -2h+2)\gamma^2 -2(2h-2)\gamma  -(d+2h-2)
\]
\[
\equiv 0 +0 -(14+10-2) \equiv 2 \not\equiv 44 \mod 12. 
\]
Thus $\beta \neq 2$.


Suppose $\beta =1$. 
In this case, 
$8(h-1) \in \alpha \Z= (\gamma -1)\Z$  (Lemma \ref{l-678-E1-ABC}(D)). 
We see that none of (2), (3), (4) of Claim \ref{cl-d=h+g-3} satisfies this condition. 
Hence Claim \ref{cl-d=h+g-3}(1) holds. 
By $H_W^3 \equiv d +2h-2 \mod \gamma$ (Lemma \ref{l-678-E1-ABC}(B)), we get the following contradiction: 
\[
22=H_W^3 \equiv d +2h-2 = 15+12-2 =25 \not \equiv 22 \mod 6. 
\]
Thus $\beta \neq 1$. This completes the proof of Lemma \ref{l-d=h+g-3}. 
\end{proof}

\begin{prop}\label{p-678-h-bound}
We use Notation \ref{n-678-E1}. 
Then $2 \leq \gamma \leq 5$. 
\end{prop}

\begin{proof}
By $\gamma \geq 2$ (Lemma \ref{l-678-gamma>1}), 
it suffices to show $\gamma \leq 5$. 
We use Notation \ref{n-f(x)}. 
If $s \leq \frac{3}{2}$ or $d = h+g-3$, then 
we get $\gamma \leq 5$ (Lemma \ref{l-s-2gamma+3}(2), Lemma \ref{l-d=h+g-3}). 
Recall that we have $d \leq h+g-3$ (Lemma \ref{l-678-easy-bound}(3)). 
Then we may assume that $d \leq  h+g-4$. 
It is enough to show $s \leq \frac{3}{2}$. 
This holds by the following: 
\[
s 
\overset{{\rm (i)}}{=}\frac{1}{2} \cdot \frac{d-2h+2}{g-d +h -2} \leq \frac{1}{2} \cdot \frac{(h+g-4)-2h+2}{g-(h+g-4) +h -2} 
\]
\[
= \frac{1}{4}(-h+g-2) 
\overset{{\rm (ii)}}{\leq}
\frac{1}{4}(-5+12-2)  = \frac{5}{4} < \frac{3}{2}. 
\]
where (i) follows from Notation \ref{n-f(x)} and 
(ii) holds by $g \leq 12$ and $h\geq 5$ (Proposition \ref{p-678-E1-h>4}). 
\qedhere

\end{proof}

\begin{lem}\label{l-678-h-gamma}
We use Notation \ref{n-678-E1}.
Then one of {\rm (1)}-{\rm (4)} holds. 
\begin{enumerate}
\item $\gamma =3, h=5, g=12$. 
\item $\gamma = 4, h =5, g \in \{11, 12\}$. 
\item 
$\gamma =5, h=5, g \in \{ 10, 11, 12\}$. 
\item 
$\gamma =5, h=6, g=12$. 
\end{enumerate}
\end{lem}

\begin{proof}
It follows from 
$d \geq 3$ (Lemma \ref{l-678-h0-d3}(1)) 
that $d +2h-2 >0$. 
By Lemma \ref{l-678-E1-ABC}(C), the following holds: 
\begin{eqnarray*}
0&<& \alpha H_W^3  \\
&=& (d -2h+2)\gamma^2 -2(2h-2)\gamma  -(d+2h-2)\\
&<& \gamma ( (d -2h+2)\gamma -2(2h-2)). 
\end{eqnarray*}
This, together with $d \leq h+g-3$ (Lemma \ref{l-678-easy-bound}(3)), implies 
\[
2(2h-2) <  (d -2h+2)\gamma \leq ((h+g-3) -2h+2)\gamma = (g-1-h)\gamma.
\]
Thus $4h-3 \leq (g-1-h) \gamma$, which implies 
\begin{equation}\label{e1-678-h-gamma}
5 \leq h \leq \frac{(g-1)\gamma +3}{\gamma+4} \leq \frac{11\gamma +3}{\gamma +4}, 
\end{equation}
where the first inequality is guaranteed by Proposition \ref{p-678-E1-h>4}. 
By using $2 \leq \gamma \leq 5$ (Proposition \ref{p-678-h-bound}), we shall prove the assertion 
by case study. 

\medskip

Assume $\gamma =2$. Then (\ref{e1-678-h-gamma}) leads to the following contradiction: 
\[
5 \leq h \leq  \frac{11\gamma +3}{\gamma +4}  = \frac{25}{6} <5. 
\]

Assume $\gamma=3$. Then (\ref{e1-678-h-gamma}) implies 
\[
5 \leq h \leq  \frac{(g-1)\gamma +3}{\gamma+4} = \frac{3g}{7} \leq \frac{36}{7}<6, 
\]
and hence $11 < 5 \cdot \frac{7}{3} \leq g \leq 12$. 
Then  $h=5$ and $g =12$. 
Thus (1) holds.

Assume $\gamma =4$. 
Then  (\ref{e1-678-h-gamma}) implies 
\[
5 \leq h \leq   \frac{(g-1)\gamma +3}{\gamma+4} =  \frac{4g -1}{8} \leq \frac{47}{8} <6,
\]
and hence $10 < \frac{5 \cdot 8+1}{4}  \leq g \leq 12$. 
Then $h=5$ and $g\in \{ 11, 12\}$. Thus (2) holds.

Assume $\gamma =5$. 
Then  (\ref{e1-678-h-gamma}) implies 
\[
5 \leq h \leq   \frac{5g -2}{9} \leq  \frac{5 \cdot 12 -2 }{9}= \frac{58}{9} < 7. 
\]
Thus $h \in \{ 5, 6\}$. 
If $h=6$, then $6 =h \leq   \frac{5g -2}{9} \leq \frac{5 \cdot 12 -2 }{9}$ implies $g=12$, 
and hence  (4) holds. 
If $h=5$, then 
$5 \leq   \frac{5g -2}{9} \leq \frac{5 \cdot 12 -2 }{9}$ implies $g \in \{ 10, 11, 12\}$, and hence (3) holds. 
\qedhere

\end{proof}


\begin{prop}\label{p-678-E1-main}
We use Notation \ref{n-678}. 
Then $\tau$ is not of type $E_1$.  
\end{prop}

\begin{proof}
Suppose that $\tau$ is of type $E_1$. 
Let us derive a contradiction. 
We use Notation \ref{n-678-E1}. 
By Lemma \ref{l-678-h-gamma}, 
 one of {\rm (1)}-{\rm (4)} holds. 
\begin{enumerate}
\item $\gamma =3, h=5, g=12$. 
\item $\gamma = 4, h =5, g \in \{11, 12\}$. 
\item 
$\gamma =5, h=5, g \in \{ 10, 11, 12\}$. 
\item 
$\gamma =5, h=6, g=12$. 
\end{enumerate}
In what follows, we treat the above four cases separately. 


\medskip

(4) Suppose $\gamma =5, h=6, g=12$. 
By  Lemma \ref{l-678-E1-ABC}(B), the following holds: 
\begin{eqnarray*}
H_W^3 &=& (2g-2d +2h -4)\gamma^3 -3(d -2h+2)\gamma^2+ 3(2h-2)\gamma +(d+2h-2)\\
&=& 125(24-2d +12 -4) -75(d -12+2)+ 15(12-2) +(d+12-2)\\
&=& 125(-2d +32) -75(d -10)+ 150 +(d+10) \\
&=& (4000 +750 +150+10) - (250 +75-1)d =4910 -324d. 
\end{eqnarray*}
By $1 \leq H_W^3 \leq 22$ (Remark \ref{r-w-bound}), this is a contradiction, because $4910 \equiv 50 \mod 324$. 

\medskip

(1) Suppose $\gamma =3, h=5, g=12$.
By  Lemma \ref{l-678-E1-ABC}(B), the following holds: 
\begin{eqnarray*}
H_W^3 &=& (2g-2d +2h -4)\gamma^3 -3(d -2h+2)\gamma^2+ 3(2h-2)\gamma +(d+2h-2)\\ 
&=& 27(24-2d +10 -4) -27(d -10+2)+ 9(10-2) +(d+10-2)\\ %
&=& 27(-2d +30) -27(d -8)+ 72 +(d+8) \\
&=& (810 + 216 + 72 +8) - (54+27-1) d
= 1106 -80d. 
\end{eqnarray*}
By $1 \leq H_W^3 \leq 22$ (Remark \ref{r-w-bound}), this is a contradiction, because $1106 \equiv 66 \mod 80$.

\medskip

(2) 
Suppose $\gamma = 4, h =5, g \in \{11, 12\}$. 
By  Lemma \ref{l-678-E1-ABC}(B), the following holds: 
\begin{eqnarray*}
H_W^3 &=& (2g-2d +2h -4)\gamma^3 -3(d -2h+2)\gamma^2+ 3(2h-2)\gamma +(d+2h-2)\\
&=& 64(2g-2d +10 -4) -48(d -10+2)+ 12(10-2) +(d+10-2)\\
&=& 64(2g-2d +6) -48(d -8)+ 96 +(d+8) \\
&=& 128g + (384 +384 +96 +8) - (128 +48 -1)d\\
&=& 128 g + 872 -175d 
= 2280 + 128\widetilde{g} -175d, 
\end{eqnarray*}
where $\widetilde{g} := g-11 \in \{ 0, 1\}$. 
By $1 \leq H_W^3 \leq 22$  (Remark \ref{r-w-bound}) and  $2280 \equiv 5 \mod 175$, we obtain 
$H_W^3 =5$ and $\widetilde g=0$, i.e., $g=11$. 
Thus 
$\beta =2$ (Remark \ref{r-E1-deg-genus}). 
By Lemma \ref{l-678-E1-ABC}(A), we get 
\[
5=\frac{1}{2} \beta H_W^3 = (g-d +h -2)
\gamma^2 - (d -2h+2)\gamma   + (h-1) \in 2\Z, 
\]
which is absurd. 

\medskip

(3) 
Suppose $\gamma =5, h=5, g \in \{ 10, 11, 12\}$. 
By  Lemma \ref{l-678-E1-ABC}(B), the following holds: 
\begin{eqnarray*}
H_W^3 &=& (2g-2d +2h -4)\gamma^3 -3(d -2h+2)\gamma^2+ 3(2h-2)\gamma +(d+2h-2)\\
&=& 125(2g-2d +10 -4) -75(d -10+2)+ 15(10-2) +(d+10-2)\\
&=&  125(2g-2d +6) -75(d -8)+ 120 +(d+8)\\
&=&  250g + (750 +600 +120 +8) - (250+75-1)d\\
&=& 250g + 1478 -324d 
= 3978 + 250\widetilde{g} -324d, %
\end{eqnarray*}
where $\widetilde{g} := g-10 \in \{ 0, 1, 2\}$. 
We have 
\begin{itemize}
\item $3978 \equiv 90 \mod 324$, 
\item $3978+250 \equiv 16 \mod 324$, and 
\item $3978+500 \equiv 266 \mod 324$. 
\end{itemize}
By $1 \leq H_W^3 \leq 22$  (Remark \ref{r-w-bound}), we obtain $H_W^3 =16$ and 
$\widetilde g=1$, i.e., $g=11$. 
Thus $\beta =1$ (Remark \ref{r-E1-deg-genus}).  
By Lemma \ref{l-678-E1-ABC}(A), we get 
\[
8 = \frac{1}{2} \beta H_W^3 = (g-d +h -2)
\gamma^2 - (d -2h+2)\gamma   + (h-1)\\
\equiv 4 \not\equiv 8 \mod 5,  
\]
which is absurd. 
\end{proof}

\begin{thm}\label{t-678-main}
Let $X \subset \P^{g+1}$ be an anti-canonically embedded Fano threefold with $\Pic\,X =\Z K_X$ and
 $g \geq 6$. 
Let $\Gamma$ be a smooth curve on $X$. 
Let $\sigma : Y \to X$ be the blowup along $\Gamma$. 
Then $-K_Y$ is not ample.  
\end{thm}

\begin{proof}
Suppose that $-K_Y$ is ample. 
We use Notation \ref{n-678}. 
Then we get a contradiction by 
Proposition \ref{p-678-nonE1-main} and Proposition \ref{p-678-E1-main}. 
\end{proof}

\begin{rem}
If we use Theorem \ref{t g neq 11}, then 
it was not necessary to treat the cases $g=11$ and $g=12$ in Notation \ref{n-678}. 
However, the author decided to treat these cases in this section, 
as the proof of Theorem \ref{t-678-main} is of characteristic free, whilst 
the one of  Theorem \ref{t g neq 11} is not. 
\end{rem}

\begin{rem}\label{r computer}
By using computer, 
we can check that there exists  
a unique solution 
\[
(g, d, h, \beta, \gamma, v) = (8, 4, 1, 2, 2, 4) 
\]
for $(g, h, d, \beta, \gamma, v) \in \Z^6$ satisfying (1)-(8), (A), and (C) 
below (where 
we set $v := H_W^3$ in Notation \ref{n-678-E1}). 
Therefore, if it is allowed to use computer, 
then Proposition \ref{p-678-E1-main} is proven only by 
Step \ref{s3-678-E1-h=1} in the proof of Lemma \ref{l-678-E1-h=1}. 
\begin{enumerate}
\item $6 \leq g \leq 12$ (Notation \ref{n-678-E1}). 
\item $1 \leq d \leq 19$.
\item $0 \leq h \leq 10$. 
\item $1\leq \gamma \leq 13$. 
\item $h \leq g-2$ (Lemma \ref{l-678-easy-bound}(2)). 
\item $2h -1 \leq d \leq h+g -3$ (Lemma \ref{l-678-easy-bound}(3)). 
\item $(\beta, \gamma) \neq (1, 1)$ (Notation \ref{n-678-E1}). 
\item One of the following holds (Remark \ref{r-E1-deg-genus}, Remark \ref{r-w-bound}). 
\begin{itemize}
\item $(\beta, v) = (4, 1)$. 
\item $(\beta, v) = (3, 2)$. 
\item $\beta =2$ and $1 \leq v \leq 5$. 
\item $\beta =1$, $2 \leq v \leq 22$, and $v \in 2\Z$. 
\end{itemize}
\item[(A)] $\frac{1}{2} \beta v = (g-d +h -2)
\gamma^2 - (d -2h+2)\gamma   + (h-1)$ (Lemma \ref{l-678-E1-ABC}). 
\item[(C)]$(\beta \gamma -1)v = (d -2h+2)\gamma^2 -2(2h-2)\gamma  -(d+2h-2)$ (Lemma \ref{l-678-E1-ABC}).  
\end{enumerate}
Note that (3) holds by (1) and (5). 
Then (2) follows from (1), (3), and (6). 
The bound (4) is assured by 
\[
\gamma < -h + g +2 \leq -0 + 12 + 2 =14, 
\]
where the first (resp. second) inequality holds by Lemma \ref{l-s-2gamma+3}(3) 
(resp. (1) and (3)). 
\end{rem}

\subsection{Case $g \leq 5$}\label{ss-nonFano-g<6}



\begin{prop}\label{p-line-nonFano}
Let $X \subset \P^{g+1}$ be an anti-canonically embedded Fano threefold with $\Pic\,X = \Z K_X$. 
Let $\Gamma$ be a smooth curve on $X$ and let $\sigma : Y \to X$ be the blowup along $\Gamma$. 
Assume that there exists a line on $X$. 
Then $-K_Y$ is not ample. 
\end{prop}

\begin{proof}
By $\Hilb^{\line}_X \neq \emptyset$, we have $\dim \Hilb^{\line}_X \geq 1$ (Proposition \ref{p-line-Hilb1}). 
Hence there exists a prime divisor $D$ on $X$ covered by lines. 
By $\rho(X)=1$, we get $D \cap \Gamma \neq \emptyset$. 
Hence we can find a line $L$ on $X$ such that $L \cap \Gamma \neq \emptyset$. 
If $\Gamma = L$, then we can apply Corollary \ref{c-line-not-ample}. 
Hence we may assume $\Gamma \neq L$. 
For $E:=\Ex(\sigma)$ and the proper transform $L_Y$ of $L$ on $Y$, 
it holds that 
\[
K_Y \cdot L_Y =(\sigma^*K_X +E) \cdot L_Y = K_X \cdot L +E \cdot L_Y =-1 +E \cdot L_Y \geq -1 +1 =0. 
\]
Thus $-K_Y$ is not ample. 
\end{proof}

\begin{lem}\label{l-345-line}
Let $X \subset \P^{g+1}$ be an anti-canonically embedded Fano threefold with $\Pic\,X = \Z K_X$ and $g \not\in \{6, 7\}$. 
Assume that one of the following conditions holds. 
\begin{enumerate}
\item There exists a $W(k)$-lift $\wt{X}$ of $X$, 
i.e., $\wt{X}$ is a flat projective $W(k)$-scheme such that 
$\wt{X} \times_{\Spec\,W(k)} \Spec\,(W(k)/pW(k)) \simeq X$. 
\item  $g \leq 5$. 
\end{enumerate}
Then there exists a line on $X$. 
\end{lem}

\begin{proof}
We first prove that (2) implies (1). Assume (2). 
Since $X$ is a complete intersection in $\P^{g+1}$ (Proposition \ref{p-genus-345}), 
we can find,  by taking lifts of defining equations,  
a  lift $\wt{X} \subset \P^{g+1}_R$ of $X \subset \P^{g+1}_k$ 
to the ring $R:=W(k)$ of Witt vectors. 
This completes the proof of the implication (2) $\Rightarrow$ (1).

In what follows, we assume (1). 
For $Q := {\rm Frac}\,R = {\rm Frac}\,W(k)$ and its algebraic closure $\overline Q$, 
$\wt{X} \times_{R} \overline Q$ contains a line 
(Remark \ref{r exist line}). 
Therefore, we can find a finite extension $R \subset R'$ of discrete valuation rings such that 
the generic fibre of the base change $\wt{X} \times_{R} R'$ contains a line $L$ (defined over ${\rm Frac}\,R'$). 
Taking the closure   $\wt{L} \subset \wt{X} \times_R R'$ of $L$, the fibre $\wt{L}_k \subset X$ over the closed point is a line on $X$. 
\end{proof}


\begin{rem}\label{r exist line}
Let $k_0$ be an algebraically closed field of characteristic zero and 
let $X\subset \P^{g+1}$ be an anti-canonically embedded Fano threefold over $k_0$ with $\Pic X =\Z K_X$. 
If $g \not\in \{ 6, 7\}$, then 
it is well known that there exists a line on $X$: 
\begin{itemize}
    \item $g \leq 5$: \cite[Corollary 1.4]{Can21}
    \item $g \geq 8$: \cite[Theorem 0.2]{Tak89}. 
\end{itemize}
On the other hand, if $g =6$ (resp. $g=7$), 
then the reference  \cite{IP99}
seems to contain  a gap in the proof of the existence of 
lines (resp. conics) in 
\cite[Case (IIa) in page 94]{IP99} (resp. \cite[Case (III) in Page 101]{IP99}). 
\qedhere

\end{rem}

\begin{thm}\label{t-blowup-main}
Let $X \subset \P^{g+1}$ be an anti-canonically embedded Fano threefold with $\Pic\,X =\Z K_X$. 
Let $\Gamma$ be a smooth curve on $X$. 
Let $\sigma : Y \to X$ be the blowup along $\Gamma$. 
Then $-K_Y$ is not ample.  
\end{thm}

\begin{proof}
If $g \geq 6$, then the assertion follows from Theorem \ref{t-678-main}. 
Assume $g \leq 5$. 
By Lemma \ref{l-345-line}, there exists a line on $X$. 
Hence we can apply Proposition \ref{p-line-nonFano}. 
\end{proof}

\section{Non-existence of the case $g = 11$}\label{s-g11}

The purpose of this section is to 
prove $g \neq 11$ 
for an anti-canonically embedded Fano threefold $X \subset \P^{g+1}$ with $\Pic X = \Z K_X$. 
The idea is to construct a $W(k)$-lift of $X$ by using 
the two-ray game for conics (Theorem \ref{t g neq 11}). 



\begin{dfn}
Let $X$ be a projective $k$-scheme. 
We say that $\wt{X}$ is a {\em $W(k)$-lift} of $X$ 
if $\wt{X}$ is a flat projective $W(k)$-scheme such that 
$X \times_{\Spec\,W(k)} \Spec\,W(k)/pW(k) \simeq X$. 
\end{dfn}




\begin{dfn}
Let $f : X \to Y$ be a $k$-morphism of projective $k$-schemes. 
We say that $\wt{f} : \wt{X} \to \wt{Y}$ is a {\em $W(k)$-lift} of $f : X \to Y$  if 
$\wt{f} : \wt{X} \to \wt{Y}$ is a $W(k)$-morphism of flat projective $W(k)$-schemes $\wt{X}$ and $\wt{Y}$ such that there exists 
the following diagram of $W(k)$-schemes in which all the squares are cartesian: 
\[
\begin{tikzcd}
\wt{X} \arrow[d, "\wt{f}"] & X \arrow[d, "f"] \arrow[l]\\
\wt{Y} \arrow[d]  & Y\arrow[l] \arrow[d]\\
\Spec\,W(k) & \Spec W(k)/pW(k).\arrow[l]  
\end{tikzcd}
\]
\end{dfn}

\begin{dfn}[cf. \cite{KM98}*{Definition 6.10}]\label{d-D-flop}
Let $U$ and $V$ be noetherian integral normal separated schemes, 
where $U$ is regular and $V$ is Gorenstein. 
\begin{enumerate}
\item 
We say that $f : U \to V$ is a {\em flopping contraction} if 
\begin{enumerate}
\item $f$ is a proper birational morphism with $\codim\,\Ex(f) \geq 2$, and  
\item $K_U$ is $f$-numerically trivial. 
\end{enumerate}
For a $\Q$-divisor $D$ on $U$, we say that $f$ is a {\em $D$-flopping contraction}  
if (a)-(c) hold. 
\begin{enumerate}
\item[(c)] $-D$ is $f$-ample. 
\end{enumerate}
\item 
Given a $\Q$-divisor $D$ on $U$ and a $D$-flopping contraction $f: U \to V$, 
we say that $f^+ : U^+ \to V$ is the $D$-{\em flop} of $f : U \to V$ if 
\begin{enumerate}
\renewcommand{\labelenumii}{(\roman{enumii})}
\item 
$U^+$ is a noetherian integral normal separated scheme, 
\item $f^+$ is a proper birational morphism with $\codim\,\Ex(f^+) \geq 2$, and 
\item the proper transform $D^+$ of $D$ on $U^+$ is $\Q$-Cartier and $f^+$-ample. 
\end{enumerate}
\end{enumerate}
\end{dfn}

\begin{rem}\label{r-flop-fg}
We use the same notation as in 
Definition \ref{d-D-flop}. 
Let  $f : U \to V$ be a $D$-flopping contraction and set $D_V :=f_*D$. 
By the same argument as in \cite[Section 6.1]{KM98}, 
the $D$-flop is unique 
and given by $U^+ = \Proj_V\bigoplus_{d \geq 0} \MO_V(dD_V)$ 
if it exists \cite[Corollary 6.4(3), the proof of Lemma 6.2]{KM98}. 
Moreover, the following are equivalent 
 \cite[Corollary 6.4(2), Proposition 6.6]{KM98}.  
\begin{enumerate}
\item There exists the $D$-flop of $f$. 
\item $\bigoplus_{d \geq 0} \MO_V(dD_V)$ is a 
finitely generated $\MO_V$-algebra. 
\item 
$\bigoplus_{d \geq 0} \MO_{\widehat{V}_P}(dD|_{\widehat{V}_P})$ is a 
finitely generated 
$\MO_{\widehat{V}_P}$-algebra  
for every closed point $P \in V$, 
where we set $\widehat{V}_P := \Spec\,\widehat{\MO}_{V, P}$ for the completion $\widehat{\MO}_{V, P}$ of $\MO_{V, P}$. 
\end{enumerate}
Note that  the $D$-flopping contraction $f$ 
is a $(K_U+D)$-flipping contraction, and hence  
many arguments in \cite[Section 6.1]{KM98} are applicable. 
\end{rem}

\begin{rem}\label{r-flop-pic1}
We use the same notation as in 
Definition \ref{d-D-flop}. 
\begin{enumerate}
\item 
A typical case is when $\rho(U/V)=1$. 
In this case, a $D$-flopping contraction coincides with $D'$-flopping contraction 
whenever $-D$ and $-D'$ are $f$-ample. 
Therefore, $D$-flop is simply called {\em flop} 
when $\rho(U/V)=1$. 
\item 
In our application, 
we start with the situation satisfying   $\rho(U/V)=1$. 
However, this property is not preserved under taking 
the base change $\Spec\,\widehat{\MO}_{V, P} \to V$ 
for the completion $\widehat{\MO}_{V, P}$ of $\MO_{V, P}$. 
Hence we need the generalised setting introduced above. 
\end{enumerate}
\end{rem}





\begin{lem}\label{l flop criterion}
Let $Y$ and $Z$ be noetherian integral normal separated schemes, 
where $Y$ is regular and $Z$ is Gorenstein. 
Let $\psi : Y \to Z$ be a flopping contraction which is not an isomorphism.  
Assume that there exists a finite surjective morphism 
$\theta: Z \to \overline{Z}$ to a noetherian integral regular separated scheme $\overline Z$ with $[K(Z) : K(\overline Z)]=2$. 
Then the following hold. 
\begin{enumerate}
\item The induced field extension $K(Z) / K(\overline Z)$ is a Galois extension. 
\item 
For every $\Q$-divisor $D$ on $Y$ such that $-D$ is $\psi$-ample, 
the composite morphism 
\[
Y \xrightarrow{\psi} Z \xrightarrow{\iota} Z 
\]
is the $D$-flop of $\psi : Y \to Z$ 
for the Galois involution $\iota : Z \xrightarrow{\simeq} Z$ 
$($i.e., $\iota$ is the automorphism induced by the  element of ${\rm Gal}(K(Z) / K(\overline Z))$ which is not the identity$)$.  
\end{enumerate}
\end{lem}


\begin{proof} 
Fix a Cartier divisor $D$ on $Y$ such that $-D$ is $\psi$-ample.

Let us show (1). 
Suppose that $K(Z) / K(\overline Z)$ is not Galois. 
By  $[K(Z) : K(\overline Z)]=2$, 
$K(Z) / K(\overline Z)$ is purely inseparable. 
As $\overline Z$ is regular, 
it follows from \cite[Lemma 2.5]{Tan18b} that $Z$ is $\Q$-factorial. 
Then we have $D = \psi^*\psi_*D$. 
This is absurd, because $-D$ is $\psi$-ample and $\psi$ is not an isomorphism, whilst $\psi^*\psi_*D$ is $\psi$-numerically trivial. 
Thus (1) holds.



Let us show (2). 
Let 
\[
\iota: Z \xrightarrow{\simeq} Z
\]
be the Galois involution. 
Set $D_Z := \psi_*D$ and $D_{\overline Z} := \theta_*D_Z = \theta_*\psi_*D$. 
Both $D_Z$ and $D_{\overline Z}$ are Weil divisors. 
Since $\overline Z$ is regular, 
the Weil divisor $D_{\overline Z}$ is Cartier. 
It is easy to see  
\[
\theta^*D_{\overline Z}= D_{Z} + D'_{Z}
\]
for the Galois conjugate $D'_{Z} := \iota^*D_Z$ of $D_{Z}$. 
For the proper transform $D'$ of $D'_{Z}$ on $Y$, we obtain 
\[
\psi^*\theta^*D_{\overline Z} = D + D'. 
\]
By $D' = \psi^*\theta^*D_{\overline Z} - D$, 
$D'$ is a $\psi$-ample Cartier divisor on $Y$. 
Therefore, 
$Y = \Proj_Z\bigoplus_{d \geq 0} \psi_*\MO_{Y}(dD') 
= \Proj_Z \bigoplus_{d \geq 0} \MO_Z(dD'_Z)$. 
Then the flop $Y^+$ of $\psi$ is given by 
\[
Y \simeq 
\Proj_Z \bigoplus_{d \geq 0} \MO_{Z}(dD'_Z)
\simeq 
\Proj_Z \bigoplus_{d \geq 0} \MO_{Z}(dD_Z) = Y^+, 
\]
where the middle isomorphism is induced by the Galois involution $\iota$ and 
the last equality holds by Remark \ref{r-flop-fg}. 
Therefore, the $D$-flop $Y^+$ of $\psi : Y \to Z$ exists 
(Remark \ref{r-flop-fg}) 
and 
$Y \xrightarrow{\psi} Z \xrightarrow{\iota} Z$ is the $D$-flop of $\psi : Y \to Z$. 
\end{proof}

{\cred 
\begin{lem}\label{l cont lift}
Let $\psi : Y \to Z$ be a morphism of projective $k$-schemes  and  
let $\wt{\psi} : \wt{Y} \to \wt{Z}$ be a $W(k)$-lift of $\psi$. 
Assume that $\psi_*\MO_Y = \MO_Z$ and $R^1\psi_*\MO_Y =0$. 
Then $\wt{\psi}_*\MO_{\wt Y} = \MO_{\wt Z}$ and $R^1\wt{\psi}_*\MO_{\wt Y} =0$. 
\end{lem}

\begin{proof}
The exact sequence 
\[
0 \to \MO_{\wt Y} \xrightarrow{\times p} \MO_{\wt Y} \to \MO_Y \to 0 
\]
induces another one 
\[
R^1\wt{\psi}_*\MO_{\wt Y}  \xrightarrow{\times p} R^1\wt{\psi}_*\MO_{\wt Y} 
\to  R^1\psi_*\MO_{Y} =0. 
\]
By Nakayama's lemma, we get $R^1\wt{\psi}_*\MO_{\wt Y} =0$. 

By $R^1\wt{\psi}_*\MO_{\wt Y} =0$, we 
have the following commutative diagram in which each horizontal sequence is exact: 
\[
\begin{tikzcd}
0 \arrow[r] & \MO_{\wt Z} \arrow[r, "\times p"] \arrow[d, "\alpha"]&  
\MO_{\wt Z} \arrow[r]  \arrow[d, "\alpha"]& 
\MO_{Z} \arrow[r]  \arrow[d, "\simeq"] &  0\\
0 \arrow[r] & \wt{\psi}_*\MO_{\wt Y} \arrow[r, "\times p"] &  
\wt{\psi}_*\MO_{\wt Y} \arrow[r] & 
\psi_*\MO_{Y} \arrow[r]  &  0.
\end{tikzcd}
\]
It follows from the snake lemma that 
\[
\times p: \Ker(\alpha) \xrightarrow{\simeq} \Ker(\alpha), \qquad 
\times p: \Coker(\alpha) \xrightarrow{\simeq} \Coker(\alpha). 
\]
We then get $\Ker(\alpha)  = \Coker(\alpha) =0$ by Nakayama's lemma. 
\end{proof}



}

\begin{prop}\label{p flop lift}
Let $\psi :Y \to Z$ be a flopping contraction with $\rho(Y/Z)=1$, 
where $Y$ is a smooth projective threefold over $k$ and $Z$ is a Gorenstein normal projective threefold over $k$. 
Let $\wt{\psi}: \wt{Y} \to \wt{Z}$ 
be 
a $W(k)$-lift 
of $\psi : Y \to Z$. 
Then the following hold. 
\begin{enumerate}
\item 
$\wt{Y}$ is regular and $\wt{Z}$ is Gorenstein. 
\item $\rho(\wt{Y}/\wt{Z}) =1$ and 
$\wt{\psi} : \wt{Y} \to \wt{Z}$ is a flopping contraction. 
\item There exists the flop $\wt{\psi}^+ : \wt{Y}^+ \to \wt{Z}$ of $\wt{\psi} : \wt{Y} \to \wt{Z}$. Moreover, for every closed point $P \in \wt{Z}$, 
there exists the following  isomorphism 
over $\Spec\,\widehat{\MO}_{\wt{Z}, P}$: 
\[
\wt{Y} \times_{\wt Z} \Spec\,\widehat{\MO}_{\wt{Z}, P} \simeq 
\wt{Y}^+ \times_{\wt Z} \Spec\,\widehat{\MO}_{\wt{Z}, P}, 
\]
where $\widehat{\MO}_{\wt{Z}, P}$ denotes the completion of 
$\MO_{\wt{Z}, P}$. 

\item The flop $\wt{\psi}^+ : \wt{Y}^+ \to \wt{Z}$ of $\wt{\psi}$ is a 
$W(k)$-lift of 
the flop $\psi^+ : Y^+ \to Z$ of $\psi$. 
\end{enumerate}
\end{prop}







\begin{proof}
Let us show (1). 
Suppose that $\wt{Y}$ is not regular. 
Since $\wt{Y}$ is excellent, 
the non-regular locus $\wt{Y}_{\text{non-reg}}$ is a closed subset of $\wt{Y}$. 
By $\wt{Y}_{\text{non-reg}} \neq \emptyset$, its image to $\Spec\,W(k)$ is a 
non-empty closed subset, 
as $\wt{Y} \to \Spec\,W(k)$ is projective. 
Hence it has the closed point $pW(k)$ of $\Spec\,W(k)$. 
Then 
$\wt{Y}$ has a non-regular point around 
the fibre $Y$ over $pW(k)$, which contradicts the smoothness of $Y$. 
Hence $\wt{Y}$ is regular. 
Since $Z$ is Gorenstein, 
$\wt{Z}$ is Gorenstein around the fibre $Z$. 
Then $\wt{Z}$ is Gorenstein  by the same argument as above. 
Thus (1) holds. 



Let us show (2). 
By the same argument as in \cite[the second paragraph of the proof of Lemma 2.6]{CT20}, 
we get $\rho(\wt{Y}/\wt{Z})=1$ and $K_{\wt{Y}}$ is 
$\wt{\psi}$-numerically trivial 
(e.g., if $L$ is an invertible sheaf on $\wt{Y}$ satisfying $L|_Y \equiv 0$, 
then we get $L \equiv 0$, because any curve on the generic fibre of $\wt{Y} \to \Spec\,W(k)$ degenerates to a curve on the central fibre $Y$). 
{\cred 
By $\psi_*\MO_Y = \MO_Z$ and $R^1\psi_*\MO_Y =0$ \cite[Proposition 6.10]{Tan-flop}, 
we get $\wt{\psi}_*\MO_{\wt Y} = \MO_{\wt Z}$ (Lemma \ref{l cont lift}).} 
Hence $\wt{\psi}: \wt{Y} \to \wt{Z}$ is birational. 
Since $\wt{\psi}(\Ex(\wt{\psi})) \cap Z$ is $0$-dimensional, 
$\wt{\psi}(\Ex(\wt{\psi}))$ is at most $1$-dimensional. 
Every fibre of $\wt{\psi}$ is at most one-dimensional by the 
upper semi-continuity of fibre dimensions. 
Therefore, $\dim \Ex(\wt{\psi}) \leq 2$. 
Then 
\[
\codim\,\Ex(\wt{\psi}) \overset{(\star)}{=} 
\dim \wt{Y} -\dim \Ex(\wt{\psi}) \geq 4-2 =2, 
\]
where $(\star)$ holds by \cite[Proposition 2.8]{CT20}. 
Thus (2) holds.


Let us prove the implication (3) $\Rightarrow$ (4). 
Assume (3). Let us show (4). 
Let $\wt{\psi}^+_p : \wt{Y}^+_p \to \wt{Z}_p = Z$ be the morphism 
obtained from $\wt{Y}^+ \to \wt{Z}$ by taking the modulo $p$ reduction. 
By (3), $\wt{\psi}^+_p : \wt{Y}^+_p \to Z$ is a flopping contraction 
from a smooth projective threefold $\wt{Y}^+_p$. 
Take a Cartier divisor $\wt{E}$ on $\wt Y$ such that $-\wt{E}$ is ample. 
Since the proper transform $\wt{E}^+$ of $\wt{E}$ is $\wt{\psi}^+$-ample, 
its reduction $\wt{E}^+_p$ modulo $p$ is $\wt{\psi}^+_p$-ample. 
Therefore, $\wt{Y}^+_p \to Z$ is a $\wt{E}_p$-flop. 
By uniqueness of flops (Remark \ref{r-flop-fg}), 
we get $Y^+ \simeq \wt{Y}^+_p$. 
This  completes the proof of the implication (3) $\Rightarrow$ (4).



It suffices to show (3). 
Fix a Cartier divisor $\wt{D}$ on $\wt{Y}$ such that $-\wt{D}$ is $\wt{\psi}$-ample (e.g., set $\wt{D} := -\wt{H}$ for an ample Cartier divisor $\wt{H}$ on $\wt{Y}$).  
Pick a closed point $P \in Z \subset \wt{Z}$. 
\begin{itemize}
\item 
{\cred Fix an affine open neighbourhood $P \in \Spec\,\wt{A} \subset \wt{Z}$. 
Let $\Spec\,A$ be the inverse image of $\Spec\,\wt{A}$, 
which is an affine open neighbourhood of $P \in Z$.} 
\item Set $B := A_{\m_A} = \MO_{Z, P}$, where $\m_A$ denotes the  the maximal ideal of $A$ corresponding to $P$. 
\item Take the $\m_B$-adic completion 
$C := \varprojlim_n B/\m_B^n = \widehat{\MO}_{Z, P}$, 
where $\m_B := \m_A B$. 
\end{itemize}
{\cred It holds that $\wt{A}$ is $W(k)$-flat  (i.e., $p$-torsion free) 
and $\wt{A}/p\wt{A}=A$.}
Hence we have the following exact sequence consisting of  $\wt{A}$-module homomorphisms: 
\begin{equation}\label{e1 flop lift}
0 \to \wt{A} \xrightarrow{\times p} \wt{A} \to A \to 0. 
\end{equation}
Let $\m_{\wt{A}}$ be the maximal ideal of $\wt{A}$ 
that is the inverse image of $\m_{{A}}$ by the surjection 
$\wt{A} \to \wt{A}/p\wt{A} =A$. 
For $\wt{B} := \wt{A}_{\m_{\wt{A}}}$ and $\m_{\wt{B}} := \m_{\wt A} \wt B$, 
we have 
\[
\wt{B}/p\wt{B} = 
\wt{A}_{\m_{\wt{A}}} / p \wt{A}_{\m_{\wt{A}}} 
= (\wt{A}/p\wt{A})_{\m_{\wt{A}}} =A_{{\m_{\wt{A}}}} = A_{\m_A} =B.
\]
Hence 
we get the following exact sequence by 
applying $(-) \otimes_{\wt{A}} \wt{B}$ to (\ref{e1 flop lift}):
\begin{equation}\label{e2 flop lift}
0 \to \wt{B} \xrightarrow{\times p} \wt{B} \to B \to 0. 
\end{equation}
We now apply the $\m_{\wt B}$-adic completion $(-)^{\wedge}$ to this exact sequence. 
Set $\wt{C} := (\wt{B})^{\wedge} 
= \varprojlim_n \wt{B}/\m_{\wt B}^n$. 
 By $\m_{\wt B} B = \m_{B}$, we have 
 \[
 \wt{C} / p \wt{C} = (\wt{B})^{\wedge} / p(\wt{B})^{\wedge}  = 
 (\wt{B}/p\wt{B})^{\wedge} = B^{\wedge} = 
 \varprojlim_n B/(\m_{\wt B} B)^n = 
 \varprojlim_n B/\m_B^n = C. 
 \]
Applying  the $\m_{\wt B}$-adic completion 
$(-)^{\wedge} = (-) \otimes_{\wt B} \wt C$ to (\ref{e2 flop lift}),  
we get the following exact sequence: 
\begin{equation}\label{e3 flop lift}
0 \to \wt{C} \xrightarrow{\times p} \wt{C} \xrightarrow{\pi} C \to 0. 
\end{equation}
Therefore, $\wt{C}$ is a $W(k)$-lift of $C$ (i.e., $p$-torsion free and $\wt{C}/p\wt{C} = C$). 
{\cred By construction, we get 
\[
\wt{C} = (\wt{B})^{\wedge} =  (\wt{A}_{\m_{\wt{A}}})^{\wedge} 
= \widehat{\MO}_{\wt{Z}, P}. 
\]}

\medskip

Recall that we have a  double cover. 
\[
\varphi : k[[x, y, z]]\hookrightarrow \widehat{\MO}_{Z, P} =C, 
\]
i.e., $\varphi$ is a $k$-algebra homomorphism such that 
\[
C = k[[x, y, z]] \oplus k[[x, y, z]] c
\]
for some $c \in C$ \cite[Lemma 6.12]{Tan-flop}. 
Fix lifts $\wt{x}, \wt{y}, \wt{z} \in \wt{C}$ of $x, y, z \in C$.  
Since $\pi(\wt{x}) = x \in \m_C$, 
we get $\wt{x} \in \pi^{-1}(\m_C) = \m_{\wt C}$, 
where $\m_{\wt C}$ denotes the maximal ideal  of $\wt{C}$. 
Similarly, we obtain $\wt{x}, \wt{y}, \wt{z} \in \m_{\wt C}$. 
Since $\wt{C}$ is $\m_{\wt C}$-adically complete, 
we get a $W(k)$-algebra homomorphism: 
\[
\wt{\varphi} : W(k)[[\wt{x}, \wt{y}, \wt{z}]] \to \wt{C}, \qquad 
\wt{x} \mapsto \wt{x} \quad \wt{y} \mapsto \wt{y}, \quad \wt{z} \mapsto \wt{z}. 
\]
Fix a lift $\wt{c} \in \wt{C}$ of $c \in C$. 
Consider the natural $W(k)[[\wt{x}, \wt{y}, \wt{z}]]$-module homomorphism
\[
\theta: W(k)[[\wt{x}, \wt{y}, \wt{z}]] \oplus 
W(k)[[\wt{x}, \wt{y}, \wt{z}]]\wt{c} \to 
\wt{C}. 
\]
We now finish the proof by assuming that $\theta$ is an isomorphism. 
Since $\theta$ is {\cred an} isomorphism, the induced morphism 
\[
\wt{\varphi}^* :
\Spec\,\widehat{\MO}_{{\cred {\wt Z}}, P} = \Spec\,\wt{C} \to \Spec\,W(k)[[\wt{x}, \wt{y}, \wt{z}]]. 
\]
is a finite surjective morphism 
of degree $2$, and hence we are done by  Lemma \ref{l flop criterion}. 

It is enough to show that $\theta$ is an isomorphism. 
Note that 
$\theta$ will become a surjection after taking modulo $p$ reduction, 
i.e., $\Im(\theta) + p\wt{C} = \wt{C}$. 
By Nakayama's lemma, $\theta$ is surjective. 
Let us show the injectivity of $\theta$. 
Pick $f, g \in W(k)[[\wt{x}, \wt{y}, \wt{z}]]$ such that 
$\theta (f + g\wt{c})=0$. 
Then, taking modulo $p$, 
 $f + g\wt{c}$ is mapping to  zero in 
$k[[x, y, z]] \oplus k[[x, y, z]] c$. 
This implies that $f +  g\wt{c} = pf' + pg'\wt{c}$ for some 
$f', g' \in W(k)[\wt{x}, \wt{y}, \wt{z}]$. 
Since $p\theta(f' + g'\wt{c})=0$ and $\wt{C}$ is $p$-torsion free, 
we get $\theta(f'+ g'\wt{c})=0$. 
By the same argument as above, we obtain $f'+ g'\wt{c}=pf'' +pg''\wt{c}$  for some 
$f'', g'' \in W(k)[\wt{x}, \wt{y}, \wt{z}]$. 
Then $f = p^2f''$ and $g= p^2 g''$. 
Repeating this procedure, we obtain 
\[
f, g \in \bigcap_{n \geq 1} p^n W(k)[[\wt{x}, \wt{y}, \wt{z}]] =0. 
\]
Therefore, $\theta$ is {\cred an} isomorphism. 
Thus (3) holds. 
\qedhere







\end{proof}


\begin{prop}\label{p cont lift}
Let $f: U \to V$ be a $k$-morphism of normal projective varieties over $k$ 
such that $f_*\MO_U = \MO_V$, $R^1f_*\MO_U =0$, and $H^2(U, \MO_U)=0$. 
Let $\wt{U}$ be a $W(k)$-lift of $U$. 
Then there exists a $W(k)$-lift $\wt{f} : \wt{U} \to \wt{V}$ of $f: U \to V$  
such that $\wt{f}_*\MO_{\wt{U}} = \MO_{\wt{V}}$. 
\end{prop}

\begin{proof}    
Let $\mathfrak U$ be the $p$-adic formal scheme of $\wt{U}$.  
By \cite[Thereom 2.15]{BBKW}, 
there exists a lift $\mathfrak f: \mathfrak U \to \mathfrak V$ 
of $f: U \to V$ in the sense of formal schemes. 
By the same argument as in \cite[Proposition 5.2]{BBKW}, this is algebrisable. 
{\cred Then the equality $\wt{f}_*\MO_{\wt{U}} = \MO_{\wt{V}}$ follows from Lemma \ref{l cont lift}.} 
\end{proof}

\begin{thm}\label{t g neq 11}
Let $X \subset \P^{g+1}$ be an anti-canonically embedded Fano threefold 
with $\Pic X = \Z K_X$. 
Then the following hold. 
\begin{enumerate}
\item $g \neq 11$. 
\item If $g =12$, then there exists a lift $\wt{X}$ of $X$ over $W(k)$, 
and hence there exists a line on $X$ (Lemma \ref{l-345-line}). 
\end{enumerate}
\end{thm}

\begin{proof}
We may assume that $g = 11$ or $g=12$.

\begin{claim*}
There exists a diagram 
\begin{equation}\label{e1 g neq 11}
\begin{tikzcd}
Y \arrow[d, "\sigma"'] \arrow[rd, "\psi"]& & Y^+ \arrow[ld, "\psi^+"'] \arrow[d, "\tau"]\\
X & Z & W=\P^3_k
\end{tikzcd}
\end{equation}
such that   the following hold. 
\begin{enumerate}
\renewcommand{\labelenumi}{(\alph{enumi})}
\item If $g=11$ (resp. $g=12$), 
then 
$\sigma : Y \to X$ is a blowup along a conic (resp. at a point). 
\item 
\begin{itemize}
\item $-K_Y$ is semi-ample and big.
\item $\psi : Y \to Z$ is the birational morphism obtained by the Stein factorisation 
of $\varphi_{|-mK_Y|}$ for some (every) $m \in \Z_{>0}$ such that $|-mK_Y|$ is base point free. 
\item $\dim \Ex(\psi) =1$
\end{itemize}
\item $\psi^+ : Y^+ \to Z$ is the flop of $\psi : Y \to Z$. 
\item $\tau$ is a blowup of $W=\P^3_k$ along a smooth rational curve $B$ of degree $6$. 
\end{enumerate}
\end{claim*}

\begin{proof}[Proof of Claim]
We separately treat the cases (i) $g=11$ and (ii) $g=12$. 

(i) Assume $g=11$. 
Fix a conic $\Gamma$ on $X$, whose existence is guaranteed by 
Theorem \ref{t-pt-to-conic}. 
Let $\sigma : Y \to X$ be the blowup along $\Gamma$. Then (a) holds. 
Note that  $|-K_Y|$ is base point free and $-K_Y$ is big 
(Proposition \ref{p-conic-basic}(2)). 
Let $\psi : Y \to Z$ be the contraction that is obtained by taking the Stein factorisation of $\varphi_{|-K_Y|}$. 
Then $\psi$ is a birational morphism. 
By Proposition \ref{p-conic-div-cont} and Corollary \ref{c-conic-not-ample}, we get $\dim \Ex(\psi) =1$. 
Take the flop $\psi^+ : Y^+ \to Z$. 
Then (b) and (c) hold. 
By Theorem \ref{t-conic-flop}, 
we get $\tau : Y^+ \to W$ satisfying (d).

(ii) Assume $g=12$. 
Fix a point $P$ on $X$ such that no line on $X$ passes through $P$, 
whose existence is guaranteed by Proposition \ref{p-line-Hilb2}. 
Let $\sigma : Y \to X$ be the blowup at $P$. Then (a) holds. 
Note that  $-K_Y$ is semi-ample and big  
(Proposition \ref{p-pt-basic}(2)). 
Let $\psi : Y \to Z$ be the contraction that is obtained by taking the Stein factorisation of $\varphi_{|-mK_Y|}$, 
where $m$ is a positive integer such that $|-mK_Y|$ is base point free. 
Then $\psi$ is a birational morphism. 
By Proposition \ref{p-pt-div-cont}(1) and Corollary \ref{c-pt-not-ample}, we get $\dim \Ex(\psi) =1$. 
Take the flop $\psi^+ : Y^+ \to Z$. 
Then (b) and (c) hold. 
By Theorem \ref{t-pt-flop}, 
we get $\tau : Y^+ \to W$ satisfying (d). 
This completes the proof of Claim. 
\end{proof}


We now finish the proof by assuming $(\star)$ below. 
\begin{enumerate}
\item[$(\star)$] There exists a lift $\wt{X}$ of $X$ over $W(k)$. 
\end{enumerate}
If $g=12$, then $(\star)$ implies (2). 
Assume $g=11$. 
Let $X_0$ be  the geometric generic fibre $X_0$ of $\wt{X} \to \Spec\,W(k)$. 
Then $X_0$  is a Fano threefold over an algebraically closed field of characteristic zero such that 
$\rho(X_0) = \rho(X)=1$, $r_{X_0} = r_X =1$, 
and $(-K_{X_0})^3 = (-K_X)^3 = 2g -2 =20$ \cite[Proposition 3.6]{KT23}. 
However, such a Fano threefold does not exist in characteristic zero 
\cite[\S 12.2]{IP99}.



\medskip

It is enough to show $(\star)$. 
To this end,  we shall construct a $W(k)$-lift of the diagram (\ref{e1 g neq 11}) in reverse order.

\underline{Lifting of $\tau : Y^+ \to W$:}  Set $\wt{W} :=\P^3_{W(k)}$, which is clearly a $W(k)$-lift of $\P^3_k$.  
In order to find a $W(k)$-lift of $\tau : Y^+ \to W$, 
it suffices to show that the closed immersion $B \hookrightarrow W$ lifts to $W(k)$. 
It is enough to prove $H^1(B, N_{B/W})=0$ \cite[Theorem 6.2]{Har10}. 
Recall that we have the following exact sequences for $\MO_W(1) := \MO_{\P^3_k}(1)$ \cite[Ch. II, Theorem 8.17 and Example 8.20.1]{Har77}: 
\[
0 \to \MO_W \to \MO_{W}(1)^{\oplus 4} \to T_{W} \to 0
\]
\[
0 \to T_B \to T_W|_B \to N_{B/W} \to 0. 
\]
Therefore, we get 
 surjections
$\MO_W(1)^{\oplus 4}|_B \twoheadrightarrow  T_W|_B \twoheadrightarrow N_{B/W}$, 
which induce another surjection: 
\[
H^1(B, \MO_W(1)|_B)^{\oplus 4} \twoheadrightarrow H^1(B, N_{B/W}). 
\]
By $B \simeq \P^1_k$ and $\deg (\MO_W(1)|_B) = \deg B  =6$, 
we get $H^1(B, \MO_W(1)|_B)=0$, and hence $H^1(B, N_{B/W})=0$, as required. 

\medskip

\underline{Lifting of $\psi^+: Y^+ \to Z$:}  
By Proposition \ref{p cont lift}, it suffices to show that 
$(\psi^+)_*\MO_{Y^+} = \MO_Z$, $R^1(\psi^+)_*\MO_{Y^+} =0$, and $H^2(Y^+, \MO_{Y^+})=0$. 
The first equality is obvious. 
The second one follows from \cite[Proposition 6.10]{Tan-flop}. 
The third one holds by the fact that $\tau : Y^+ \to \P^3$ is a blowup 
 along a smooth curve $B$. 

\medskip

\underline{Lifting of $\psi: Y \to Z$:} 
Recall that $\psi^+ : Y^+ \to Z$ is the flop of $\psi : Y \to Z$. 
Hence $\psi : Y \to Z$ is the flop of $\psi^+ : Y^+ \to Z$. 
Since $\psi^+ : Y^+ \to Z$ has a $W(k)$-lift over $W(k)$, 
so does $\psi : Y \to Z$ (Proposition \ref{p flop lift}(4)). 

\medskip

\underline{Lifting of $\sigma: Y \to X$:} 
By Proposition \ref{p cont lift}, it suffices to show that 
$\sigma_*\MO_{Y} = \MO_X$, $R^1\sigma_*\MO_{Y} =0$, and 
$H^2(Y, \MO_{Y})=0$. 
The first equality is obvious. 
Since $\sigma : Y \to X$ is a blowup along a smooth curve or at a point  (Claim(a)), 
we get $R^i\sigma_*\MO_Y = 0$ for every $i>0$. 
Hence $H^2(Y, \MO_Y) \simeq H^2(X, \MO_X)=0$, 
as required. 
\end{proof}

\begin{rem}
By the above proof, if $B$ is a smooth rational curve or an elliptic curve on $\P^n$, 
then the blowup $X$ along $B$ always has a $W(k)$-lift. 
The author learned this argument from Burt Totaro. 
\end{rem}

\begin{rem}
Let $X \subset \P^{g+1}$ be an anti-canonically embedded  Fano threefold 
with $\Pic\,X = \Z K_X$. 
It is tempting to hope the existence of a $W(k)$-lift of $X$. 
If $g\leq 5$ or $g=12$, 
then there exists a $W(k)$-lift of $X$ 
(Theorem \ref{t g neq 11} and the first paragraph of the proof of 
Lemma \ref{l-345-line}). 
For the remaining cases $6 \leq g \leq 10$, 
the author does not know whether $X$ has a $W(k)$-lift. 
\end{rem}

\bibliographystyle{skalpha}
\bibliography{reference.bib}

\end{document}